\documentclass[11pt]{article}
\usepackage{amsfonts}
\usepackage{comment}
\usepackage{tabularx}
\usepackage{enumerate}
\usepackage{bbm,bm}
 \usepackage[affil-it]{authblk}
\usepackage{tikz}
\usepackage{latexsym}
\usepackage{mathtools}
\usepackage[pdftitle={A bijective path to Liouville quantum gravity},
pdfauthor={Olivier Bernardi, Nina Holden, Xin Sun}
bookmarks=true,bookmarksopen=true,bookmarksopenlevel=3]{hyperref}
\usepackage{wrapfig}
\usepackage{caption}
\usepackage{paralist}
\usepackage{amsmath, amsthm, amssymb, color,  graphicx, hyperref}
\usepackage[margin=1in]{geometry}
\usepackage[normalem]{ulem}

\DeclareGraphicsExtensions{.jpg,.pdf,.eps}
\newtheorem{thm}{Theorem}[section]
\newtheorem{theorem}[thm]{Theorem}
\newtheorem{proposition}[thm]{Proposition}
\newtheorem{prop}[thm]{Proposition}
\newtheorem{lemma}[thm]{Lemma}
\newtheorem{cor}[thm]{Corollary} 
\newtheorem{corollary}[thm]{Corollary} 

\newtheorem{definition}[thm]{Definition}

\newtheorem{claim}[thm]{Claim}

\newtheorem{notation}[thm]{Notation}
\newtheorem{assumption}[thm]{Assumption}

\theoremstyle{definition}
\newtheorem{remark}[thm]{Remark}
\newtheorem{fact}{Fact}[section]

\newcommand{\ds}{\displaystyle}

\newcommand{\fig}[3]{\begin{figure}[!h]\begin{center}\includegraphics[#1]{#2}\end{center}\caption{#3}\label{fig:#2}\end{figure}}

\newcommand{\xin}[1]{{#1}}
\newcommand{\nina}[1]{{#1}}
\newcommand{\referee}[1]{{#1}}
\newcommand{\old}[1]{{}} 

\newcommand{\cut}{\mathrm{Cut}}
\newcommand{\ans}{\mathrm{AnFr}}

\newcommand{\env}{\mathrm{env}}

\newcommand{\C}{\mathbb{C}}

\newcommand{\1}{\mathbf{1}}
\newcommand{\D}{\mathbb{D}}
\newcommand{\E}{\mathbb{E}}
\newcommand{\N}{\mathbb{N}}
\newcommand{\Q}{\mathbb{Q}}
\newcommand{\Z}{\mathbb{Z}}
\newcommand{\R}{\mathbb{R}}
\renewcommand{\P}{\mathbb{P}}
\newcommand{\bbH}{\mathbb{H}}
\newcommand{\ol}{\overline}
\newcommand{\wt}{\widetilde}
\newcommand{\ul}{\underline}

\newcommand{\eps}{\varepsilon}
\def\P{\mathbb{P}}
\def\E{\mathbb{E}}

\DeclareMathOperator{\SLE}{SLE}
\DeclareMathOperator{\CLE}{CLE}
\DeclareMathOperator{\LQG}{LQG}

\def\cP{\mathcal{P}}

\def\cN{\mathcal{N}}
\def\cM{\mathcal{M}}
\def\cL{\mathcal{L}}
\def\cK{\mathcal{K}}

\def\cI{\mathcal{I}}

\def\cF{\mathcal{F}}

\newcommand{\ob}{{\op{b}}}
\newcommand{\ow}{{\op{w}}}
\newcommand{\hQ}{\wh Q}

\newcommand{\etad}{\etab_{\op D}}
\newcommand{\pivd}{\piv_{\op D}}
\newcommand{\etas}{\etab_{\op S}}

\newcommand{\etadh}{\wh{\etab}_{\op D}}
\newcommand{\Zd}{\Zb^{\op D}}
\newcommand{\Ld}{\Lb^{\op D}}
\newcommand{\Rd}{\Rb^{\op D}}
\newcommand{\Zs}{\Zb^{\op S}}
\newcommand{\Ls}{\Lb^{\op S}}
\newcommand{\Rs}{\Rb^{\op S}}
\newcommand{\Gabd}{\Gab_{\op D}}
\newcommand{\Gabs}{\Gab_{\op S}}
\newcommand{\taud}{\taub^*_{\op D}}
\newcommand{\taus}{\taub^*_{\op S}}
\newcommand{\phib}{\bm{\varphi}}
\newcommand{\nud}{\nub^{\op D}}
\newcommand{\nus}{\nub^{\op S}}
\newcommand{\Eb}{\bm{E}}
\newcommand{\dbl}{\mathrm{dbl}}
\newcommand{\gffd}{\gff_{\op D}}
\newcommand{\gffs}{\gff_{\op S}}
\newcommand{\reg}{\mathrm{reg}}
\newcommand{\piv}{\mathrm{piv}}
\newcommand{\sig}{\mathrm{sig}}
\newcommand{\cardy}{\op{Cardy}}

\newcommand{\ub}{{\mathbf{u}}}
\newcommand{\al}{\alpha}
\newcommand{\be}{\beta}
\newcommand{\ga}{\gamma}

\newcommand{\si}{\sigma}

\newcommand{\Ga}{\Gamma}

\newcommand{\RR}{\mathbb{R}}
\newcommand{\QQ}{\mathbb{Q}}
\newcommand{\ZZ}{\mathbb{Z}}
\newcommand{\NN}{\mathbb{N}}

\newcommand{\mA}{\mathcal{A}}
\newcommand{\mB}{\mathcal{B}}
\newcommand{\mC}{\mathcal{C}}

\newcommand{\mK}{\mathcal{K}}

\newcommand{\mP}{\mathcal{P}}

\newcommand{\mT}{\mathcal{T}}

\newcommand{\tw}{\tilde{w}}

\newcommand{\tT}{\tilde{T}}

\newcommand{\ba}{\bar{a}}
\newcommand{\bb}{\bar{b}}

\newcommand{\hpi}{\widehat{\pi}}

\newcommand{\bmK}{\overline{\mK}}

\newcommand{\bmC}{\overline{\mC}}

\newcommand{\bPhi}{\overline{\Phi}}

\newcommand{\bphi}{\overline{\phi}}
\newcommand{\bOm}{\overline{\Om}}

\newcommand{\suf}[1]{\overleftarrow{#1}}
\newcommand{\smK}{\suf{\mK}}

\newcommand{\wttau}{\widetilde{\tau}}
\newcommand{\ttau}{\widetilde{\tau}}
\newcommand{\wta}{\widetilde{a}}
\newcommand{\wtb}{\widetilde{b}}
\newcommand{\wtw}{\widetilde{w}}

\DeclareMathOperator{\hcode}{h-code}
\DeclareMathOperator{\ccode}{cw-code}
\DeclareMathOperator{\cwcode}{cw-code}
\DeclareMathOperator{\ccwcode}{ccw-code}

\newcommand{\mH}{\mathcal{H}}
\newcommand{\Om}{\Omega}
\newcommand{\Perc}{\mathrm{Perc}}
\newcommand{\DFS}{\mathrm{DFS}}
\newcommand{\gao}{\sigma^{\circ}}
\newcommand{\dfs}{\textrm{dfs-tree}}
\newcommand{\cltree}{\textrm{cluster-tree}}
\newcommand{\dfsdual}{\textrm{dfs-dual}}
\newcommand{\matchingpar}{\textrm{p}}
\newcommand{\treepar}{\vec{\matchingpar}}
\newcommand{\sma}{\textrm{anc}}

\newcommand{\bmT}{\overline{\mT}}
\newcommand{\bmTP}{\bmT_P}
\newcommand{\imK}{\mK^\infty}
\newcommand{\imT}{\mT^\infty}
\newcommand{\imTP}{\mT^\infty_P}
\newcommand{\smTP}{\suf{\mT_P}}

\newcommand{\ta}{{\widetilde{a}}}
\newcommand{\tb}{{\widetilde{b}}}
\newcommand{\LR}{\op{Spine}}
\renewcommand{\tw}{\widetilde{w}}
\renewcommand{\tT}{\widetilde{T}}

\def\@rst #1 #2other{#1}
\newcommand\MR[1]{\relax\ifhmode\unskip\spacefactor3000 \space\fi
		\MRhref{\expandafter\@rst #1 other}{#1}}\newcommand{\MRhref}[2]{\href{http://www.ams.org/mathscinet-getitem?mr=#1}{MR#2}}

\def\MR#1{\href{http://www.ams.org/mathscinet-getitem?mr=#1}{MR#1}}

\newcommand{\aryb}{\begin{eqnarray*}}
\newcommand{\arye}{\end{eqnarray*}}
\def\alb#1\ale{\begin{align*}#1\end{align*}}
\newcommand{\eqb}{\begin{equation}}
\newcommand{\eqe}{\end{equation}}
\newcommand{\eqbn}{\begin{equation*}}
\newcommand{\eqen}{\end{equation*}}

\newcommand{\BB}{\mathbbm}
\newcommand{\op}{\operatorname}

\newcommand{\frk}{\mathfrak}
\newcommand{\eqD}{\overset{d}{=}}

\newcommand{\ep}{\epsilon}
\newcommand{\rta}{\rightarrow}

\newcommand{\wh}{\widehat} 
\newcommand{\mcl}{\mathcal}

\newcommand{\bdy}{\partial}
\newcommand{\val}{\op{val}}

\newcommand{\Trev}{T^{\op{rev}}}
\newcommand{\Wrev}{Z^{\op{rev}}}
\newcommand{\Whrev}{\wh Z^{\op{rev}}}

\newcommand{\cLrev}{L^{\op{rev}}}

\newcommand{\cRrev}{R^{\op{rev}}}

\newcommand{\li}{\op{env}}
\newcommand{\ci}{\op{cone}}

\newcommand{\Ab}{\bm{A}}
\newcommand{\Yb}{\bm{Y}}
\newcommand{\Zb}{\bm{Z}}
\newcommand{\Lb}{\bm{L}}
\newcommand{\Rb}{\bm{R}}
\newcommand{\Xb}{\bm{X}}
\newcommand{\ab}{\bm{a}}
\newcommand{\Tb}{\bm{T}}

\newcommand{\gam}{\gamma}
\newcommand{\Gam}{\Gamma}
\newcommand{\taub}{\bm{\tau}}
\newcommand{\nub}{\bm{\nu}}
\newcommand{\mub}{\bm{\mu}}
\newcommand{\pb}{\bm{p}}
\newcommand{\gab}{\bm{\gamma}}
\newcommand{\Gab}{\bm{\Gamma}}
\newcommand{\etab}{\bm{\eta}}
\newcommand{\etae}{\eta_{\op{e}}}
\newcommand{\etav}{\eta_{\op{v}}}
\newcommand{\etavf}{\eta_{\op{vf}}}
\newcommand{\etalr}{\eta_{\op{lr}}}
\newcommand{\lambdav}{\lambda_{\op{v}}}

\newcommand{\gff}{\bm{h}}
\newcommand{\perc}{\sigma}
\newcommand{\ellb}{\bm{\ell}}
\newcommand{\area}{\op{area}}

\newcommand{\fll}{ \op{fl}_{\op{L}} }
\newcommand{\flr}{\op{fl}_{\op{R}}}

\newcommand{\Lo}{{\op{L}}}
\newcommand{\Ro}{{\op{R}}}

\newcommand{\Xo}{{\op{X}}}
\newcommand{\Io}{{\op{I}}}
\newcommand{\nuloc}{ \nub^{\op{}} }

\newcommand{\lra}{$\Leftrightarrow$}
\newcommand{\erase}[1]{}

\begin {document}
\title{Percolation on triangulations: a bijective path to Liouville quantum gravity}

\author[1]{Olivier Bernardi}
\author[2]{Nina Holden}
\author[3]{Xin Sun
}
\affil[1]{\small Brandeis University}
\affil[2]{\small ETH-ITS Z\"urich}
\affil[3]{\small Columbia University}

\date{\today}
\maketitle

\begin{abstract}
We set the foundation for a series of works aimed at proving strong relations between uniform random planar maps and Liouville quantum gravity (LQG). Our method relies on a bijective encoding of site-percolated planar triangulations by certain 2D lattice paths. Our bijection parallels in the discrete setting the \emph{mating-of-trees} framework of LQG and Schramm-Loewner evolutions (SLE) introduced by Duplantier, Miller, and Sheffield. Combining these two correspondences allows us to relate uniform site-percolated triangulations to $\sqrt{8/3}$-LQG and SLE$_6$. In particular, we establish the convergence of several functionals of the percolation model to continuous random objects defined in terms of $\sqrt{8/3}$-LQG and SLE$_6$. For instance, we show that the exploration tree of the percolation converges to a branching SLE$_6$, and that the collection of percolation cycles converges to the conformal loop ensemble CLE$_6$. We also prove convergence of counting measure on the pivotal points of the percolation. Our results play an essential role in several other works, including a program for showing convergence of the conformal structure of uniform triangulations and works which study the behavior of random walk on the uniform infinite planar triangulation.
\end{abstract}



\tableofcontents

\section{Introduction}\label{sec:intro}
We study critical site-percolation on random planar triangulations, and its relation to \emph{Liouville quantum gravity} (LQG) and \emph{Schramm-Loewner evolutions} (SLE). 
Recall that LQG is a random fractal 2D surface \cite{shef-kpz,rhodes-vargas-review} which is  defined by considering the standard Euclidean metric distorted by the \emph{Gaussian free field} (GFF).
LQG was originally introduced by Polyakov in the 1980s as a model for the random surface corresponding to the space-time evolution of a string \cite{polyakov-qg1,polyakov-qg2,poly-2dqg}. LQG and the GFF have since then appeared in several other mathematical physics contexts. 
Recall also that SLE curves are random fractal curves in the plane \cite{schramm0,werner-notes,lawler-book} which arise as the scaling limit of the interfaces in a wide range of statistical physics models on 2D lattices. See Section~\ref{sec:intro2} for further details. 

In this paper we show the convergence in law of several important observables of the percolation model on random planar triangulations. The continuum limits are expressed in terms of $\sqrt{8/3}$-LQG and $\SLE_6$. We show that, in a precise sense, the continuum limit of the percolation model on random planar maps is the random surface $\sqrt{8/3}$-LQG~\cite{shef-kpz,rhodes-vargas-review}, decorated with an independent instance of the \emph{conformal loop ensemble} $\CLE_6$ (an infinite collection of random loops closely related to $\SLE_6$). Precisely, upon choosing a suitable embedding of the random planar triangulations in the plane, the continuous limit of the vertex distribution has the law of the $\sqrt{8/3}$-LQG area measure, while the continuous limit of the percolation interfaces has the law of $\CLE_6$.

Our work is based on a bijective correspondence between site-percolated triangulations and certain 2D lattice walks known as 
\emph{Kreweras walks}. This bijective correspondence stems from a new interpretation and extension of a bijection of the first author~\cite{bernardi-dfs-bijection}.
A lot of information about the percolated map can be read conveniently from the walk through the bijection (see Table~\ref{table:dictionary-bijection}). Convergence in a strong sense of Kreweras walks to 2D Brownian motion is then exploited to deduce our results about percolated maps.
In particular, we obtain a convergence result for the \emph{exploration tree} of the percolated map, which is the spanning tree obtained by a depth-first search (DFS) exploration of the percolation interfaces. We also obtain convergence of the \emph{percolation cycles}, which separate clusters of different colors. Furthermore, we prove convergence of the counting measure on macroscopic \emph{pivotal points}, which are the vertices of the percolated map whose change of color changes the connectivity of percolation clusters on a macroscopic scale.   Although our convergence results are first expressed in terms of the 2D Brownian motion, a fundamental work of Duplantier, Miller and Sheffield~\cite{wedges} allows us to translate these results to results about $\sqrt{8/3}$-LQG and SLE$_6$ curves. 

Before we state one of our main result in Section~\ref{sec:intro3}, we will briefly introduce the relevant objects in Sections~\ref{subsec:bi} and~\ref{sec:intro2}. The results of the current paper play a fundamental role in several other works on planar maps and LQG, and in Section~\ref{sec:intro4} we give a brief overview of these further developments.

\subsection{Site-percolation on triangulations and Kreweras walks}\label{subsec:bi}
Let us first define more precisely our percolation model on maps.
A \emph{planar map} is a decomposition of the 2D sphere into a finite (or countable) number of vertices, edges and faces, considered up to (orientation preserving) homeomorphism. We only consider planar maps in this article, and call them simply \emph{maps}. Also, our maps will all be \emph{rooted} (see Section~\ref{subsec:prelim} for precise definitions). A \emph{near-triangulation} is a map in which all the non-root faces have degree 3. A \emph{site-percolation configuration} on a map $M$ is any coloring of its vertices in black and white.

Building on a work by the first author~\cite{bernardi-dfs-bijection}, we present in Section~\ref{sec:bijection} a bijective encoding of site-percolated near-triangulations by certain 2D lattice walks. We shall call \emph{Kreweras walk} a lattice walk on $\ZZ^2$ made of the three types of steps $a=(1,0)$, $b=(0,1)$, and $c=(-1,-1)$. This appellation is in honor of Germain Kreweras who first enumerated this type of walks confined in the first quadrant $\NN^2$ in relation with plane partitions~\cite{Kreweras:walks}.\footnote{Kreweras actually considered the reverse steps $(-1,0)$, $(0,-1)$, and $(1,1)$, so our walks would usually be called \emph{reverse Kreweras walks}. We mention that several simplifications and alternatives to Kreweras' original counting technique have been proposed over the years~\cite{bernardi-dfs-bijection,MBM:Kreweras,Gessel:Kreweras,Fayolle:walks-quarter-plane,kauers07v,KR-12,BeBoRa-arxiv}.} 

We now describe informally an important special case of the bijection $\bPhi$ obtained in Section~\ref{subsec:bij-chordal}. Let us call \emph{Kreweras $k$-excursion} a Kreweras walk starting at $(0,0)$ ending at $(0,-k)$ and remaining in the quadrant $\{(i,j)~|~i\geq 0,~j\geq-k\}$. The bijection $\bPhi$ induces a bijection between the set of Kreweras $k$-excursions and the set $\mT_P^{(k)}$ of (2-connected) site-percolated near-triangulations with $k+2$ outer vertices: 1 white vertex, and $k+1$ black vertices. This correspondence is illustrated in Figure~\ref{fig:example-bij-intro-disc}. We shall also extend the bijection $\bPhi$ to the infinite volume setting in Section~\ref{sec:bij-inf}. Recall that the \emph{uniform infinite planar triangulation} (UIPT) is the local limit of uniform triangulations as defined by Angel and Schramm~\cite{angel-schramm-uipt} (see Section \ref{sec:bij-inf} for more details). In the infinite setting, the mapping $\bPhi$ gives a measure-preserving correspondence between bi-infinite Kreweras walks and site-percolated UIPT.

\fig{width=.8\linewidth}{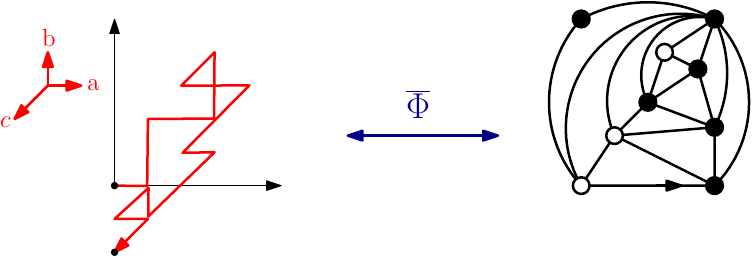}{The bijection $\bPhi$ maps the Kreweras $2$-excursion pictured on the left, to the site-percolated near-triangulation pictured on the right.}

As mentioned earlier, the bijection $\bPhi$ is related to a bijection $\Om$ obtained in~\cite{bernardi-dfs-bijection}.\footnote{There is a rich literature on bijections for planar maps. In particular, a bijection between loopless triangulations and a family of decorated trees was obtained in \cite{Schaeffer:triangulation} (see also \cite{OB-EF:girth} for a far-reaching generalization). However, this bijection is not related to the bijection $\Om$ obtained in \cite{bernardi-dfs-bijection}.} However the bijection $\Om$ is presented as a correspondence between Kreweras walks and pairs made of a near-triangulation $M$ and a \emph{depth-first search} (DFS) tree $\tau^*$ of the dual map $M^*$ (in fact, all the definitions in~\cite{bernardi-dfs-bijection} are in terms of $M^*$). A bijective correspondence between the site-percolation configurations of $M$ and the DFS trees of $M^*$ is established in Section~\ref{sec:exploration-tree-from-map}. This establishes the link between the bijection $\Om$ and (the base case of) our bijection $\bPhi$.

The bijection $\bPhi$ has the nice property that many important quantities about percolated triangulations can be expressed as simple functionals of the corresponding Kreweras walks. For example, we are interested in the DFS tree $\tau^*$, the percolation cycles $\Gamma$, and the counting measure on pivotal points of the percolation.
Using the correspondences between the percolation observables and the Kreweras walk, 
we express the scaling limit of many important quantities about percolated triangulations in terms of 2D Brownian excursions with correlation 1/2. 

This could be the end of a beautiful story, but there is much more. Indeed, in~\cite{wedges}, Duplantier, Miller, and Sheffield establish a measure-preserving correspondence between 2D Brownian excursions with correlation 1/2 and pairs made of an instance of $\sqrt{8/3}$-LQG disk and an instance of a space-filling $\SLE_6$ curve in the disk (see below for more details). Hence we can further express our scaling limit results in terms of such pairs. The motivation for this further expression is that our bijection $\bPhi$ should in fact be considered as the exact discrete analog (defined in terms of random maps) of the correspondence of Duplantier, Miller, and Sheffield (defined in terms of $\sqrt{8/3}$-LQG). By composing the two correspondences, a strong relation between random maps and $\sqrt{8/3}$-LQG is obtained here and in subsequent papers \cite{ghs-metric-peano,hs-quenched}.

\begin{table}
	\begin{center}
		\begin{tabular}{|c|c|}
			\hline
			Percolated triangulations \erase{&} & Kreweras walks\\
			\hline \erase{&}&\\[-4mm]
			$(M,\si)\in\bmT_P$ \erase{&$\bPhi$} &  $w\in \bmK$ \\
			inner triangles \erase{&$\etavf$} & $a$-steps and $b$-steps\\
			triangle incident to an active left/right edge \erase{&\lra} & unmatched $a$-step/$b$-step\\
			white/black in-vertices \erase{&$\etavf$} & $c$-steps of type $a$/$b$\\
			outer in-vertices \erase{&\lra} & unmatched $c$-step\\
			\hline
			spine-looptrees decomposition $\LR(M,\si)$ \erase{&=} & $\LR(w)$\\
			triangles on percolation path  \erase{&$\etavf$} & spine-steps of $w$ (i.e., steps of $\hpi(w)$)\\
			white/black bubbles of $\LR(M,\si)$ of length $k$ \erase{&$\etavf$} & steps $\ba_{k-1},\bb_{k-1}$ of $\hpi(w)$\\
			sub-triangulation inside a bubble \erase{&$\bPhi$} & subwalk enclosed by a spine-step matching\\
			\hline
			exploration tree $\tau^*=\dfs(M,\si)$\erase{&=} &$\dfs(w)$\\
			edges of $\tau^*$\erase{&$\etae^*$} & $a$-steps and $b$-steps\\ 
			in-edges of $M^*/\tau^*$ \erase{&$\etae^*$}& $c$-steps \\ 
			\hline 
			spanning tree $\tau=\dfsdual(M,\si)$ of $M$ \erase{&=} & $\dfsdual(w)$\\
			tree of clusters $\cltree(M,\sigma)$ \erase{&=}& unicolor contraction of  $\dfsdual(w)$\\
			percolation cycles of length $k$ \erase{&$\etae^*$}  & envelope excursions with $k$ spine steps\\
			\hline
		\end{tabular}
		\caption{The correspondences induced by the bijection $\bPhi$ between parameters of site-percolated triangulations  and parameters of Kreweras walks.}\label{table:dictionary-bijection}
	\end{center}
\end{table}

\subsection{$\SLE_6$ on $\sqrt{8/3}$-LQG and mating of trees} \label{sec:intro2}
We now give a very brief and informal introduction to $\SLE_6$ on $\sqrt{8/3}$-LQG and the mating-of-trees construction as needed to state Theorem~\ref{thm2} below. We refer to Section~\ref{sec:dictionary} for a more detailed description.

The \emph{mating-of-trees construction} is a measure-preserving correspondence introduced by Duplantier, Miller and Sheffield in~\cite{wedges}. This construction can be thought  as a continuum analog of the bijection $\bPhi$ defined in Section~\ref{subsec:bi}, as we now explain.
The continuum analog of the planar map is an instance of a $\sqrt{8/3}$-LQG surface. Let $\gff$ be (some variant of) the random distribution known as the Gaussian free field\footnote{For many continuum random variables throughout the paper we use bold letters. In cases where a discrete random variable converges in law to a continuum random variable we will typically use the same symbol for  the discrete and continuum random variable, except that the continuum one is bold. See Table \ref{table-dictionary}.} (GFF)~\cite{shef-gff} in the unit disk $\BB D$. We obtain a \emph{$\sqrt{8/3}$-LQG disk} by considering the random area measure $\mub\equiv\mub_{\gff}$ informally defined by $e^{\sqrt{8/3}\gff}\,dx\,dy$, where $dx\,dy$ denotes Lebesgue area measure in $\BB D$. This definition does not make literal sense since $\gff$ is a distribution and not a function, but $\mub$ may be defined rigorously by considering regularized versions of $\gff$~\cite{shef-kpz,rhodes-vargas-review,berestycki-gmt-elementary}. The area measure $\mub$ almost surely assigns positive measure to any open set. The field $\gff$ also induces a length measure to the boundary $\partial\D$. We assume $\gff$ is conditioned such that $\mub(\D)=1$, and such that the length of $\partial\D$ is 1. 

The continuum analog of the percolation $\perc$ is encoded by a variant of the random curve known as the Schramm-Loewner evolution (SLE). Recall that SLE is a family of random fractal curves which describes the scaling limit of the interfaces of a wide range of statistical physics models on 2D lattices~\cite{lsw-lerw-ust,smirnov-cardy,ss-dgff,chelkak-smirnov-ising,cdhks14,ks-ising,lawler-viklund-16}. SLE curves are indexed by a parameter $\kappa$, where $\kappa=6$ corresponds to the scaling limit of percolation interfaces \cite{smirnov-cardy}. An instance of SLE$_6$ in $\D$ almost surely divides $\D$ into smaller disjoint domains (called ``bubbles'') by hitting its past and $\partial\D$. A \emph{space-filling} variant $\etab$ of SLE$_6$ can be obtained by filling in recursively each bubble by a random space-filling SLE$_6$ loop right after it is enclosed (see Section \ref{sec:sle} and references therein for more details).

In~\cite{wedges} the authors consider a pair $(\gff,\etab)$ as above, where $\gff$ and $\etab$ are independent. The curve $\etab$ is parametrized by $\sqrt{8/3}$-LQG area, meaning that the time-parametrization is chosen so that $\mub(\etab([s,t]) )=t-s$ for any $0\leq s<t\leq 1$. The pair $(\gff,\etab)$ defines a process $\Zb=(\Lb_t,\Rb_t)_{t\in[0,1]}$ which describes how the length of the left and right frontier of $\etab$ evolves in time. It is proved in~\cite{wedges} that $\Zb$ has the law of a two-dimensional Brownian excursion with correlation $1/2$, conditioned to start at $(1,0)$, end at $(0,0)$, and stay in the first quadrant. Furthermore, $(\gff,\etab)$ and $\Zb$ are related by a continuum version of the bijection in Section~\ref{subsec:bi}, since it can be proved that the two objects generate the same $\sigma$-algebra. The construction of $(\gff,\etab)$ using $\Zb$ is known as the \emph{mating-of-trees construction} of LQG\footnote{The name \emph{mating-of-tree} comes from the fact that each coordinate of $\Zb$ encodes an infinite-volume continuum random tree \cite{aldous-crt1,aldous-crt2,aldous-crt3}, and $(\gff,\etab)$ may be viewed as a \emph{gluing} or \emph{mating} of these two trees such that $\etab$ describes the interface between the two trees. In that perspective, the construction in~\cite{wedges} is a ``bijective'' encoding of $(\gff,\etab)$ in terms of a pair of mated trees.}.

The space-filling SLE$_6$ $\etab$ encodes two other random variables of interest: the \emph{conformal loop ensemble} CLE$_6$ $\Gab$ and the \emph{branching} SLE$_6$ $\taub^*=(\wh{\etab}^z)_{z\in\D}$. The three random objects $\etab$, $\Gab$, and $\taub^*$ can be coupled in such a way that they determine each other (in the sense that they generate the same $\sigma$-algebra). The conformal loop ensemble~\cite{shef-cle} $\Gab$ is a random countable collection of non-crossing loops with the same local properties as SLE$_6$ curves. This ensemble is known to describe the scaling limit of percolation interfaces for critical site-percolation on the triangular lattice~\cite{smirnov-cardy,camia-newman-full}. 
The branching SLE$_6$ $\taub^*$ is a tree such that each branch $\wh\etab^z$ for fixed $z\in\D$ has the law of a (non-space-filling) SLE$_6$ from 1 to $z$. 
\emph{In the canonical coupling between $\etab$, $\Gab$, and $\taub^*$,} the space-filling  SLE$_6$ $\etab$ is obtained by a depth-first exploration of the tree $\taub^*$, while the branches $\wh\etab^z$ are defined by an exploration of the loops $\Gab$. 
The field $\gff$ induces a length measure along the CLE$_6$ loops $\gab\in\Gab$ and the branches $\wh\etab^z$ of $\taub^*$, which we may use to parametrize these curves. 

Pivotal points are points where a CLE$_6$ loop hits itself or other loops. The \emph{significance} of a pivotal point is defined in terms of the LQG area enclosed by the relevant loops. For $\eps>0$ the set of pivotal points with significance at least $\eps$ is a random fractal with dimension $3/4$. The field $\gff$ induces an LQG-measure $\nub^\eps$ supported on these points.

\subsection{Scaling limit of percolation observables}
\label{sec:intro3}
In this paper, we prove that a number of percolation observables converge jointly to a continuum limit which may be defined in terms of SLE$_6$ and $\sqrt{8/3}$-LQG. In this introduction we only give an informal statement of our result for finite volume maps with \emph{disk topology} (that is, for finite random triangulations with a single simple boundary). We refer to Section~\ref{sec:conv} for a more formal statement and for the cases of planar maps with whole-plane topology (IUPT)  and sphere topology (finite triangulations without boundary). For $n\in\ZZ^{\geq 2}$, let $M_n$ be a loopless triangulation with a simple boundary of  length $\lceil n^{1/2}\rceil$, having $n$ interior vertices. Let $\perc_n$ be a coloring of the vertices, such that the boundary vertices are white, except a single black vertex. Let the root-edge be the unique boundary edge which is directed in counterclockwise direction towards the black vertex.

Given an \emph{embedding} $\phi_n:V(M_n)\to \D$ of the map $M_n$ into the unit disk $\D$, we study the following observables of the percolation configuration on the embedded map: vertex counting measure, percolation cycles, exploration tree, space-filling exploration path, and pivotal measure. A precise definition of these objects and the topological spaces to which they belong is given in Section~\ref{sec:conv}, but we give a brief description here: 
\begin{compactitem}
	\item The \emph{vertex counting measure} $\mu_n$ is a measure on $\D$ given by counting measure on $V(M_n)$, where each vertex has mass $n^{-1}$.
	\item A \emph{percolation cycle} is a cycle separating two percolation clusters (see Section~\ref{sec:bijection}). We label the percolation cycles $\ga_0^n,\ga_1^n,\dots$, such that the number of vertices enclosed by the cycles is decreasing. We view the cycles as elements in the space of parametrized curves in $\D$, where the parametrization is such that each edge is crossed in $n^{-3/4}$ units of time.
	\item The bijection described in Section~\ref{subsec:bi} defines an ordering of all edges of $M_n$ through a particular depth-first search (DFS), which defines the exploration tree $\tau^*_n$ on the dual map $M^*_n$. For any $z\in\D$ the percolation exploration $\wh\eta^z_n$ from $\frk e_n$ to $z$ is the branch in $\tau^*_n$ from the root-edge to $z$, where the parametrization is such that each edge is traced in $n^{-3/4}$ units of time.
	\item The \emph{space-filling exploration path} $\eta_n$ is a path which visits the edges of $M$ in chronological order as determined by the DFS, starting at the root-edge, such that it takes $n^{-1}$ units of time to go from one edge to the next.
	\item A \emph{pivotal point} $v\in V(M_n)$ of the percolation $\perc_n$ is a vertex such that changing its color makes percolation cycles merge or split. The \emph{significance} of a pivotal point is defined in terms of the number of vertices enclosed by the relevant percolation cycles. We define $\nu_n^\eps$ to be renormalized counting measure in $\D$ on pivotal points of significance at least $\eps$, where each pivotal has mass $n^{-1/4}$. 
\end{compactitem}

The continuum counterpart of the percolated map $(M_n,\perc_n)$ is a pair $(\gff,\etab)$, where $\gff$ is the field associated with a $\sqrt{8/3}$-LQG disk of area 1 and boundary length 1, and $\etab$ is an independent space-filling SLE$_6$. Recall from Section~\ref{sec:intro2} that an instance of $\etab$ determines an instance $\taub^*=(\wh\etab^z)_{z\in\D}$ of a branching SLE$_6$ and an instance $\Gab$ of the conformal loop ensemble CLE$_6$. Also recall that $\gff$ induces an area measure $\mub$ in $\D$ and a measure $\nub^\eps$ supported on the pivotal points of $\Gab$. See Table~\ref{table-dictionary} for an overview of the discrete and continuum correspondences. 
\begin{theorem}\label{thm2}
	For $n\in\N_+$ let $(M_n,\sigma_n)$ be a site-percolated map as above, chosen uniformly at random. There exists a map $\phi_n:V(M_n)\to\D$ such that the following quantities converge jointly in law towards their continuum counterpart $(\gff,\etab)$ as $n\rta\infty$: 
	\begin{compactitem}
		\item The vertex counting measure $\mu_n$ on $\D$ converges weakly to the $\sqrt{8/3}$-LQG area measure $\mub$ associated with $\gff$.
		\item The embedded percolation interfaces $\ga_1^n,\ga_2^n,\dots$ in $\D$ converge to the CLE$_6$ loops $\Gab=(\gab_1,\gab_2,\dots)$ as parametrized curves.
		\item The finite marginals of the DFS tree $\tau^*_n$ converge to the finite marginals of the branching SLE$_6$ $\taub^*$. In particular, for each fixed $t\in(0,1)$ and $z:=\etab(t)$ the percolation exploration $\wh\eta^{z}_n$ converges to the SLE$_6$ $\wh\etab^z$.
		\item The space-filling percolation exploration $\eta_n$ converges uniformly to the space-filling SLE$_6$ $\etab$.
		\item The pivotal measure $\nu_n^\eps$ converges to $\nub^\eps$ for each $\eps>0$.
	\end{compactitem} 
\end{theorem}
We remark that, although the above theorem is stated in terms of a particular embedding $\phi_n$, one can deduce from the theorem that we have joint convergence in law of several interesting functionals of the percolation configuration which are not defined in terms of an embedding. For example, we have joint convergence in law of the lengths, enclosed areas, connectivity properties, and pivotal measure of the macroscopic percolation cycles.
Furthermore, we have convergence of certain crossing events, and, as we will explain in Remark \ref{rmk:lt}, the so-called looptree associated with each percolation cycle converges in the Gromov-Hausdorff topology.


\begin{table}
	\begin{center}
		\begin{tabular}{|c|c|c|c|}
			\hline
			Discrete variable & Notation & Continuum variable & Notation\\
			\hline
			Uniform infinite planar triangulation & $M$ & $\sqrt{8/3}$-LQG cone & $\gff$\\
			Uniform triangulation of disk & $M$ & $\sqrt{8/3}$-LQG disk & $\gff$\\
			Vertex counting measure & $\mu$ & $\sqrt{8/3}$-LQG area measure & $\mub$ \\
			Space-filling percolation exploration & $\eta$ & Space-filling SLE$_6$ & $\etab$\\
			Kreweras walk (bi-infinite) & $w$ & Brownian motion & $(\Zb_t)_{t\in\R}$\\
			Bijections: walk $-$ percolated map & $\Phi,\bPhi,\Phi^\infty$ & Mating-of-trees construction & \\
			Reduced word & $\wh\pi(w^-)$ & L\'evy process & $(\wh \Zb_t)_{t\leq 0}$\\
			Spine-looptrees decomposition & $\LR(w^-)$ & Forested lines of past wedge & $\fll, \flr$
			\\
			Depth-first search tree & $\tau^*$ & Branching SLE$_6$ & $\taub^*$\\ 
			Branch of DFS tree & $\wh\eta^z$ & SLE$_6$ & $\wh\etab^z$\\
			Percolation cycles & $\Gam$	& CLE$_6$ loops & $\Gab$\\
			Envelope interval of perc.\ cycle & & Envelope interval of CLE$_6$ loop & $\env(\ga)$\\
			Pivotal point counting measure & $\nu$ & CLE$_6$ double point measure & $\nub$ \\
			Crossing events (percolation) & $E_\ob,E_\ow$& Crossing events (SLE$_6$) &$\Eb_\ob,\Eb_\ow$\\
			\hline
		\end{tabular}
		\caption{The table illustrates the close correspondence between the discrete and continuum models.}\label{table-dictionary}
	\end{center}
\end{table}

LQG has long been conjectured to be related to the scaling limit of random planar maps, and much more recently some relations have been rigorously proved. 
The conjectures initially appeared in the physics literature, but without precise statements.
The relations which have been rigorously proved depend on various choices of topologies which can be put on the set of maps and on the particular planar map distribution. 
Let us now review briefly these known relations and how they compare to our result. 
Of course, there is a huge literature about both planar maps and LQG surfaces, and we cannot attempt to list all of the relevant references (although we will try to discuss all the references about \emph{concrete relations} between planar maps and LQG surfaces).
The first results about the scaling limits of maps are in terms of their \emph{metric properties} (a planar map defines a metric space, which is obtained by endowing the vertex set with the graph distance). 
Using bijective results do to Schaeffer and others \cite{Schaeffer:these,BDFG:mobiles}, Le Gall and Miermont independently proved that uniformly random quadrangulations 
considered as metric spaces converge in the Gromov-Hausdorff topology to a random metric space known as the \emph{Brownian map} (which is homeomorphic to a sphere and of Hausdorff dimension 4) \cite{LeGall:limitmaps,marc-mokk-tbm,legall-uniqueness,miermont-brownian-map}. This major result was subsequently extended to other classes of planar maps (see e.g. \cite{Albenque:limit-simple-triang,aasw-type2}). 
In another breakthrough,  Miller and Sheffield \cite{lqg-tbm1,lqg-tbm2,lqg-tbm3} proved that the Brownian map is equivalent to $\sqrt{8/3}$-LQG in the sense that the two surfaces can be coupled together so they generate the same $\sigma$-algebra. 
These results are for uniformly random planar maps without statistical physics models, and are not giving any information about embeddings of planar maps in the plane. By contrast the result in \cite{gms-tutte} deals with the embedding of a non-uniform classes of maps which are defined in terms of $\ga$-LQG for $\ga\in(0,2)$: it is shown that a particular random map defined by a coarse-graining of a $\ga$-LQG surface converges to the $\ga$-LQG area measure when embedded into $\C$ using the Tutte embedding. Let us also  mention that the convergence of random planar maps decorated with statistical physics models in the so-called \emph{peanosphere topology} has been established for several universality classes of statistical physics models  \cite{shef-burger,gms-burger-cone,gms-burger-local,gms-burger-finite,gkmw-burger,kmsw-bipolar,ghs-bipolar,lsw-schnyder-wood}. However, the peanosphere topology is defined in terms of the Brownian motion $\Zb$ associated to the pair $(\gff,\etab)$, hence a convergence result in this topology does not imply convergence of the map itself, but rather of the \emph{pair} consisting of the map and an instance of a statistical physics model on the map.

\subsection{Future works and directions}
\label{sec:intro4}
The bijection and convergence results of the current paper already play an essential role in a number of recent works and works in progress:  
\begin{compactitem}
	\item The present paper gives a first convergence result for percolated triangulations to SLE$_6$-decorated $\sqrt{8/3}$-LQG surfaces. Building on this work, the second and third authors of this paper together with collaborators strengthen the notion of convergence in future works, culminating with the proof that the \emph{conformal structure} of the triangulation is converging \cite{hs-quenched}. It is shown there that the embedding used in Theorem \ref{thm2} (which uses the mating-of-trees construction, and depend on a sequence of coupled percolated maps) is not too different from the so-called Cardy embedding (which only depend on the map itself). This proves the convergence of the planar map and several percolation observables for the Cardy embedding. See also Section \ref{sec:flip} for additional details.
	\item In~\cite{ghs-metric-peano} a first such improvement of the notion of convergence is obtained. It will be proved in~\cite{aasw-type2} that loopless triangulations $M$ converge in the Gromov-Hausdorff topology to a limiting metric space known as the Brownian disk. It will be proved in~\cite{ghs-metric-peano} that the convergence in Gromov-Hausdorff topology is \emph{joint} with the convergence established in this paper. This is helpful for studying dynamical percolation on triangulations in \cite{hs-quenched}.
	Indeed, it implies that the limiting LQG surface will stay fixed if we resample the percolation. 
	\item The bijection introduced in this paper allows us to encode properties of the planar map $M$ in terms of the word $w$. Using a strong coupling between $w$ and the mating-of-trees Brownian motion $\Zb$, this allows us to relate properties of the planar map $M$ to properties of $\Zb$ and the associated SLE$_6$-decorated $\sqrt{8/3}$-LQG surface. This approach is used in~\cite{ghs-map-dist} to study distances in planar maps in several universality classes. In particular, the second and the third author of this paper together with Gwynne use known results for the UIPT to study distances in the so-called mated-CRT map for $\kappa=6$. The mated-CRT map is a map defined in terms of $\Zb$ and is studied in for instance \cite{gms-tutte,ghs-dist-exponent,gm-spec-dim,gh-displacement}.
	\item In~\cite{gm-spec-dim,gh-displacement} strong coupling between $w$ and $\Zb$ is used to transfer properties of random walk on the mated-CRT map to properties of random walk on the UIPT. In particular, it is proved in~\cite{gm-spec-dim} that the spectral dimension of the UIPT is 2, and that a random walk on the UIPT typically travels at least $n^{1/4+o(1)}$ units of graph distance in $n$ units of time. In~\cite{gh-displacement} the matching upper bound to the latter result is proved.
\end{compactitem}

The results of this paper also opens several future research directions:
\begin{compactitem}
	\item One may attempt to find bijections between other planar maps models and random walks. As illustrated above, these are powerful tools for analyzing the planar maps and their scaling limits. The \emph{Gessel walks}, which have steps $\{(-1,0),(0,1),(1,1),(-1,-1)\}$, are one potential candidate for such a bijection. This set of walks is exactly solvable~\cite{bkr-gessel-walk,kkz-gessel-walk,mbm-gessel-walk,BeBoRa-arxiv} and is known to have an algebraic generating function. Since the coordinates of the walks have correlation $1/\sqrt{2}=-\cos(3\pi/4)$, it corresponds to the universality class of the Ising model~\cite[Theorem 1.13]{wedges}. One could also try to find bijections for maps decorated with $O(n)$ models.
	\item One could try to extend the bijection of this paper to maps on other surfaces. A bijection for percolated maps on the torus may improve the understanding of the \emph{Brownian torus} (an analogue of the Brownian map with torus topology).
	One goal could be to relate the construction of the Brownian torus given in~\cite{bet-torus1,bet-torus2} to the construction of the $\sqrt{8/3}$-LQG torus given in~\cite{drv-torus}.
\end{compactitem}

\xin{
	\subsection{Global strategy of the paper}
	The rest of the paper can be divided into three parts: discrete (Sections~\ref{sec:bijection}--\ref{sec:appendix}), continuum (Section~\ref{sec:dictionary}), and  convergence (Sections~\ref{sec:conv} and~\ref{app:conv}).

	In Section~\ref{sec:bijection}, we describe the bijection~$\bPhi$ between percolated loopless triangulations and the Kreweras walks, which is mentioned in Section~\ref{subsec:bi}. 
	We present variants of ~$\bPhi$ that correspond to maps on the sphere, disk, and whole plane (infinite triangulations).  
	Then in Sections~\ref{sec:spine}--\ref{sec:dual-DFS}, we describe in detail how the bijection encodes various geometric observables. Namely, we establish the dictionary summarized in  Table~\ref{table:dictionary-bijection}. 
	In Section~\ref{sec:appendix}, we provide the proofs of all these statements except two, which will be proven in the appendix.
	The results proved in the appendix are not used in the rest of the paper and their proofs require inputs from Sections~\ref{sec:dictionary}--\ref{app:conv}.
	
	
	\nina{In Section~\ref{sec:dictionary}, we present mating-of-trees theory for $\gamma=\sqrt{8/3}$ and $\kappa=6$, building on work by Duplantier, Miller, and Sheffield~\cite{wedges}. Based on this we give the continuum analog of the dictionary in Table~\ref{table:dictionary-bijection}, whose objects are as in the right column of Table~\ref{table-dictionary}. Part of this section is a review of~\cite{wedges} and related papers while part of this section is novel.}
	
	In Section~\ref{sec:conv}, we state a refined version of  Theorem~\ref{thm2} as well as its analog in the whole-plane and sphere cases. In particular, we include the convergence of the percolation crossing event, which is defined in Section~\ref{sec:crossing-discrete}. 
	We also state a result (Proposition~\ref{prop23}) concerning the joint convergence of the percolated triangulation  before and after the flip of a macroscopic pivotal point. This result is crucial to the second and the third authors' program on the Cardy embedding.  Section~\ref{app:conv} is devoted to the proofs of the results from Section~\ref{sec:conv}. 
	The strategy of the proofs is to use the discrete and continuum dictionaries along with convergence of random walk to Brownian motion. 
	However, given the intricacy of the geometric observables involved, some arguments in Section~\ref{app:conv} are quite involved.}

\section{A bijection between Kreweras walks and percolated triangulations}\label{sec:bijection}

\subsection{Basic definitions: maps, site-percolation, and Kreweras walks}\label{subsec:prelim}
For a positive integer $n$, we denote by $[n]$ the set $\{1,2,\ldots,n\}$. For $k\in\ZZ$, we define $\ZZ^{\geq k}$ as the set $\{n\in\ZZ~|~n \geq k\}$. The sets  $\ZZ^{\leq k},\ZZ^{>k}$ and $\ZZ^{<k}$ are defined similarly. Lastly, we use the notation $\NN:=\ZZ^{\geq 0}$ and $\NN_+:=\ZZ^{> 0}$.\\

\noindent {\bf Basic definitions about maps.}\\
Unless otherwise specified, our \emph{graphs} are finite and undirected; self-loops and multiple edges are allowed.
The \emph{length} of a path or cycle in a graph is its number of edges. For instance, \emph{self-loops} are cycles of length 1. A \emph{bridge} of a graph is an edge whose deletion would increase the number of connected components.

A \emph{planar map} (or \emph{map} for short) is a proper
embedding of a connected planar graph in the 2-dimensional sphere, considered up to \emph{deformation} (orientation preserving homeomorphism). 
The \emph{faces} of a map are the connected components of the complement of the graph. For a map $M$ we denote by $V(M)$, $E(M)$, and $F(M)$ its set of vertices, edges, and faces respectively.
A map is \emph{rooted} if one of its edge is distinguished as the \emph{root-edge} and oriented. The origin of the root-edge is called \emph{root-vertex}, and the face at the right of the root-edge is called \emph{root-face}.
When drawing maps on the plane (as in Figure~\ref{fig:triangul-with-boundary}), the root-face is usually taken as the infinite face. We call \emph{outer} the vertices and edges incident to the root-face, and \emph{inner} the other. We also call \emph{inner faces} the non-root faces.

\fig{width=.9\linewidth}{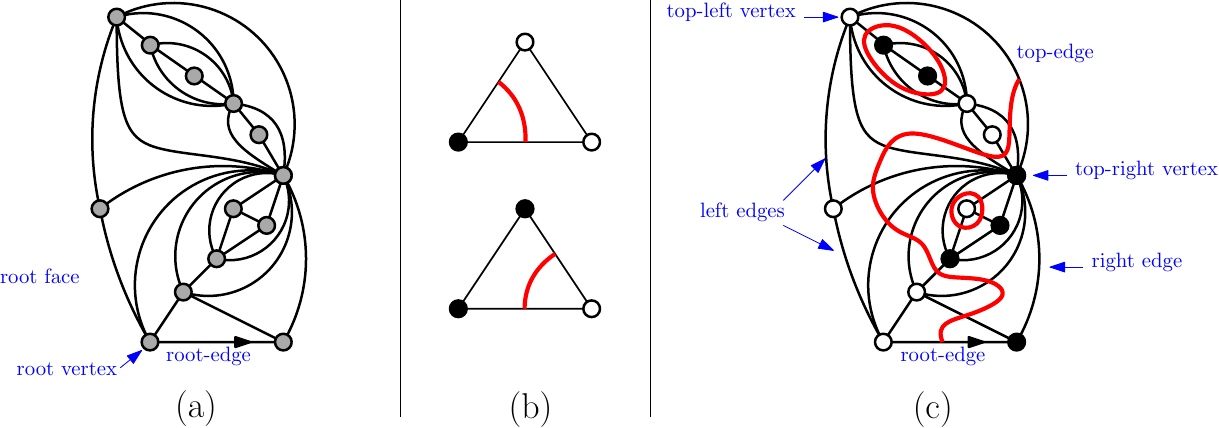}{(a) A rooted triangulation with a simple boundary $M\in\mT$. (b) The rule for drawing the percolation interfaces. (c) A site-percolation configuration $\sigma$ on $M$. The percolation cycles and paths are indicated in bold lines. The configuration satisfies the root-interface condition. In fact, $(M,\sigma)$ is in $\mT_P$ and is equal to $\Phi(w)$ for $w=abbaabbcaccbcaaaabcbbbbacacc$.}

The \emph{degree} of a vertex or face is the number of edges incident to it, counted with multiplicity (an edge is counted twice if it is twice incident). 
A \emph{triangulation} is a map in which every face has degree~3. A \emph{loopless triangulation} is a triangulation without self-loops.
A \emph{near-triangulation} is a map in which every inner face has degree~3. A near-triangulation such that the boundary of root-face is simple (that is, has no cut-point) is called \emph{triangulation with a simple boundary}. 
We denote by $\mT$ the set of rooted loopless triangulations with a simple boundary.\footnote{In the probability literature (e.g. \cite{angel-schramm-uipt}) triangulations are sometimes said to be of \emph{type I, II or III}, depending on whether loops and multiple edges are allowed. 
	In that terminology, our triangulations are of type II (loops are forbidden, but multiple edges are allowed).} 
An example is given in Figure~\ref{fig:triangul-with-boundary}(a).
By convention, the rooted map $M_0$ with a single edge and two vertices is considered to be an element of $\mT$.\\

\noindent {\bf Basic definitions about percolation.}\\
A \emph{site-percolation configuration} on a map $M$ is a coloring of its vertices in two colors: \emph{black} or \emph{white}. A \emph{site-percolated map} (or \emph{percolated map} for short) is a pair $(M,\sigma)$, where $M$ is a map and $\sigma$ is a site-percolation configuration.
An edge of $(M,\sigma)$ is \emph{unicolor} if its endpoints are of the same color, and \emph{bicolor} otherwise. Keeping only unicolor edges of $(M,\sigma)$ gives a disjoint union of maps called \emph{percolation clusters}, which can be either \emph{black} or \emph{white}.

Suppose now that $M$ is a near-triangulation. The inner triangles in $(M,\sigma)$ are either \emph{unicolor}, or \emph{bicolor} in which case they are incident to two bicolor edges. In Figure~\ref{fig:triangul-with-boundary}(b) we have drawn in each bicolor triangle a curve joining the middle of the two incident bicolor edges. The result is a set of disjoint curves which are either cycles, or paths starting and ending on the boundary of $M$. These curves (which can be thought as simple paths and cycles on the dual map $M^*$) are called \emph{percolation interfaces}, and are separating the black and white clusters of $(M,\sigma)$; see Figure~\ref{fig:triangul-with-boundary}(c).
The percolation interfaces which are cycles are called \emph{percolation cycles}, and the ones which are paths starting and ending on the boundary of $M$, are called the \emph{percolation paths} of $(M,\sigma)$. The \emph{length} of a percolation cycle or path is the number of triangles it crosses.\\

\noindent {\bf Basic definitions about Kreweras walks.}\\
We denote by $\{a,b,c\}^*$ the set of finite \emph{words} (sequence of letters) on the \emph{alphabet} $\{a,b,c\}$. 
We identify the words in $\{a,b,c\}^*$ with the finite lattice walks on $\ZZ^2$ starting at $(0,0)$, and made of steps $a=(1,0)$, $b=(0,1)$ and $c=(-1,-1)$. We refer to such lattice walks as \emph{Kreweras walks}. 

We denote by $\mK\subset \{a,b,c\}^*$, the set of words $w$ such that any prefix contains no more $c$'s than $a$'s, and no more $c$'s than $b$'s. Equivalently $w\in\mK$ is a walk staying in the quadrant $\NN^2$.

For a word $w=w_1w_2\cdots w_n\in\{a,b,c\}^*$, we say that $w_i$ is a \emph{$a$-step} if $w_i=a$; and similarly for $b$ and $c$.
We say that an $a$-step $w_i$ and a $c$-step $w_k$ are \emph{matching} if $i<k$, there are as many $a$-steps and $c$-steps in $w_{i}w_{i+1}\cdots w_k$, and for all $j\in\{i,\ldots,k-1\}$ there are strictly more $a$-steps than $c$-steps in $w_{i}w_{i+1}\cdots w_j$ (in terms of lattice walks, the subwalk $w_{i}w_{i+1}\cdots w_k$ is a ``right-excursion'' in the sense that it stays strictly to the ``right'' of the steps $w_i,w_k$). 
The \emph{matching} of $b$-steps and $c$-steps is defined analogously. 
\referee{Figure \ref{fig:matching-letters} illustrates matching steps.}

\fig{width=.5\linewidth}{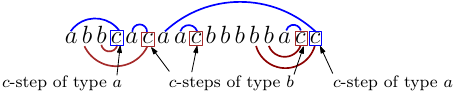}{Matching pairs of steps in a word $w\in\{a,b,c\}^*$. The matching $a$-steps and $c$-steps are indicated by arcs above $w$. The matching $b$-steps and $c$-steps are indicated by arcs below $w$.}

Clearly each $a$-step or $b$-step has at most one matching $c$-step, and we call \emph{unmatched} the ones that do not have a matching $c$-step. Also each $c$-step $w_k$ has at most one matching $a$-step, and at most one matching $b$-step. We call \emph{unmatched} a $c$-step which does not have a matching $a$-step or a matching $b$-step.
Note that $\mK$ is the set of walks such that every $c$-step is \emph{matched}, that is, has both a matching $a$-step and a matching $b$-step. Moreover, if a walk $w\in\mK$ has $x$ unmatched $a$-steps and $y$ unmatched $b$-steps, then it ends at the point $(x,y)$.

\referee{Let $w=w_1w_2\cdots w_n\in\{a,b,c\}^*$ and let  $w_k$ be a $c$-step. If $w_k$ has a matched $b$-step $w_j$, and either no matching $a$-step or a matching $a$-step $w_i$ with $i<j$, then we say that $w_k$ is a $c$-step of \emph{type $a$}. Symmetrically, if  $w_k$ has a matched $a$-step $w_j$ and either no matching $b$-step or a matching $b$-step $w_i$ with $i<j$, then we say that $w_k$ is a $c$-step of \emph{type $b$}. See Figure \ref{fig:matching-letters}.}


\subsection{The bijection: restricted case}\label{subsec:finite}
We will now define a bijection $\Phi$ between site-percolated near-triangulations and Kreweras walks. This bijection is represented in Figure~\ref{fig:example-bij-triang}.

The mapping $\Phi$ is defined on the set $K$ of Kreweras walks staying in the quadrant $\N^2$.
We now define the corresponding set of site-percolated near-triangulations. 
We say that a site-percolated near-triangulation $(M,\sigma)$ satisfies the \emph{root-interface condition} if the root-edge is oriented from a white vertex to a black vertex, and no other outer edge goes from a white vertex to a black vertex in counterclockwise direction around the root-face; see Figure~\ref{fig:triangul-with-boundary}(c). \referee{The root-interface condition is also known as \emph{Dobrushin boundary condition}.}
Observe that this condition ensures that there is a unique percolation path, and that this path starts at the root-edge.
\begin{definition}\label{def:mtp}
	We denote by $\mT_P$ the set of percolated near-triangulations $(M,\sigma)$ such that 
	\begin{compactitem}
		\item[(i)] $M$ is in $\mT$ (that is, $M$ is a rooted loopless triangulation with a simple boundary),
		\item[(ii)] $\sigma$ is a site-percolation configuration satisfying the root-interface condition,
		\item[(iii)] the unique percolation path goes through all the inner triangles incident to an outer edge.
	\end{compactitem}
\end{definition}

A pair $(M,\sigma)$ in $\mT_P$ is represented in Figure~\ref{fig:triangul-with-boundary}(c). \referee{Let us mention that the technical Condition~(iii) in Definition~\ref{def:mtp} is necessary for the restricted bijection $\Phi$ defined below (Theorem~\ref{thm:bij-OB}), but can be bypassed by using the more general setting of Section~\ref{subsec:bij-chordal} (Corollary~\ref{cor:chordal-case}).}

Note that, for any pair $(M,\sigma)$ in $\mT_P$, the percolation path connects the root-edge to another bicolor outer edge. We call this bicolor outer edge the \emph{top-edge} of $(M,\sigma)$, and we call its white and black endpoints the \emph{top-left} and \emph{top-right} vertices respectively. We call \emph{left vertices} (resp. \emph{right vertices}) the white (resp. black) outer vertices.   
We call \emph{left edges} (resp. \emph{right edges}) the unicolor white (resp. black) outer edges. See Figure~\ref{fig:triangul-with-boundary}(c).

\fig{width=\linewidth}{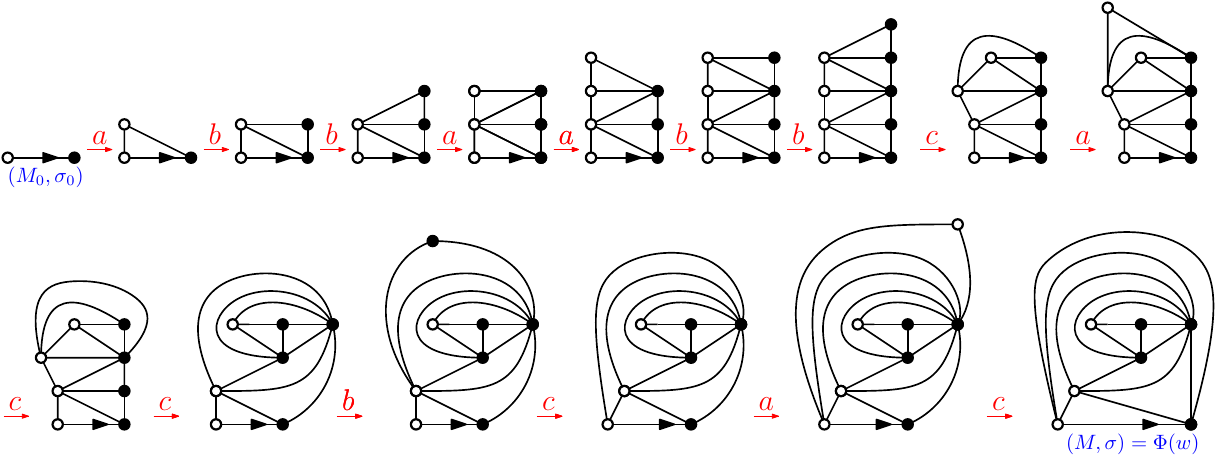}{Bijection $\Phi$ for the word $w=abbaabbcaccbcac$.}

\begin{definition}\label{def:Phi}
	For a walk $w\in\mK$, the pair $(M,\sigma)=\Phi(w)$ is constructed step by step, by reading the steps of $w$ and following the rules illustrated in Figure~\ref{fig:bij-rules-triang}. 
	Formally, the rules for the steps $a,b,c$ are given by mappings $\phi_a,\phi_b,\phi_c$ from $\mT_P$ to $\mT_P$ defined below, and for $w=w_1\cdots w_n$ we define $\Phi(w):=\phi_{w_n}\circ\cdots\circ\phi_{w_1}(M_0,\sigma_0)$, where $(M_0,\sigma_0)\in \mT_P$ is the percolated map with a single root-edge going from a white vertex to a black vertex.
	\begin{compactitem}
		\item For $(M,\sigma)\in\mT_P$, the map $\phi_a(M,\sigma)$ is obtained by gluing a triangle with two white vertices and one black vertex to the top-edge of $(M,\sigma)$.
		\item For $(M,\sigma)\in\mT_P$, the map $\phi_b(M,\sigma)$ is obtained by gluing a triangle with two black vertices and one white vertex to the top-edge of $(M,\sigma)$.
		\item For $(M,\sigma)\in\mT_P$ having both a left edge and a right edge, we define $\phi_c(M,\sigma)$ as follows (see Figure~\ref{fig:bij-rules-triang}). Let $e_\ell$ be the left edge incident to the top-left vertex $v_\ell$, and let $e_r$ be the right edge incident to the top-right vertex $v_r$. Let $P$ be the percolation path of $(M,\sigma)$, and consider $P$ as \emph{starting} at the root-edge and \emph{ending} at the top-edge $e$.
		By definition of $\mT_P$, the inner triangles $t_\ell$ and $t_r$ incident to $e_\ell$ and $e_r$ respectively are on $P$; and one of them $t\in\{t_\ell,t_r\}$ is the last triangle on $P$ incident to a left or right edge. If $t=t_r$, then we recolor the vertex $v_r$ in white, and we glue the edges $e$ and $e_\ell$ together (so that $v_\ell$ becomes an inner white vertex). Symmetrically, if $t=t_\ell$, then we recolor the vertex $v_\ell$ in black, and we glue the edges $e$ and $e_r$ together (so that $v_r$ becomes an inner black vertex). 
	\end{compactitem}
\end{definition}

\fig{width=.9\linewidth}{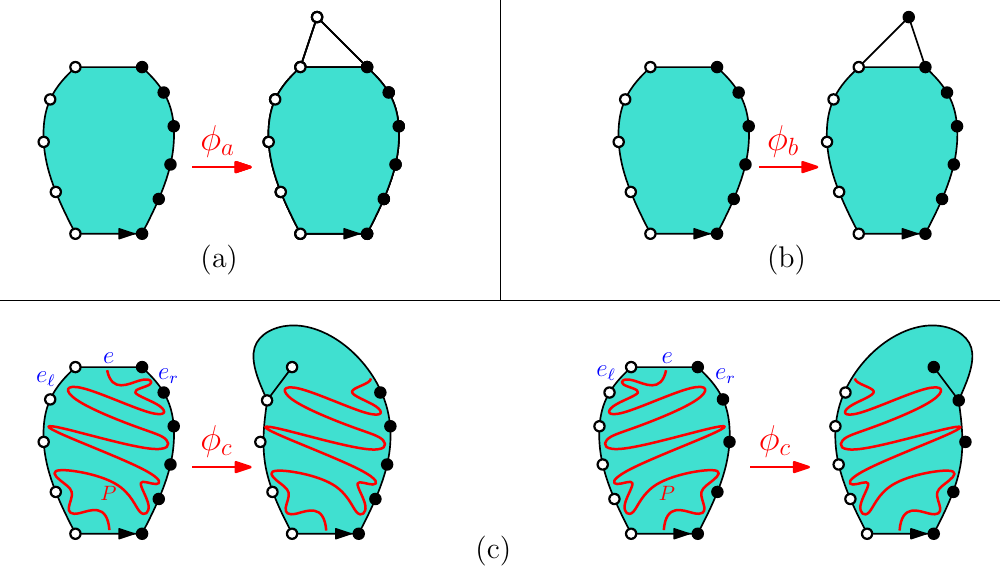}{The mappings $\phi_a$, $\phi_b$, and $\phi_c$. For the mapping $\phi_c$ there are two cases depending on the order in which the percolation path $P$ crosses the triangles $t_\ell$ and $t_r$ incident to the edges $e_\ell$ and $e_r$. On the left (resp. right) we represented the case where $P$ crosses $t_r$ (resp. $t_\ell$) last.}

\begin{remark}\label{rk:phi-well-def}
	Observe that the coordinates of the lattice walk $w=w_1w_2\cdots w_n\in\mK$ (viewed as a lattice walk on $\ZZ^2$) describe the evolution of the number of left and right edges of the map as we apply the mappings $\phi_{w_1},\phi_{w_2},\ldots,\phi_{w_n}$ successively. In particular, if $w$ ends at $(i,j)$, then $\Phi(w)$ has $i$ left edges and $j$ right edges. 
	Indeed, this is clear by induction on the length of $w$, since applying $\phi_a$ (resp. $\phi_b$) increases the number of left (resp. right) edge by 1, while applying $\phi_c$ decreases the number of left and right edges by 1.
	
	In particular, for any word $w=w_1\cdots w_n\in \mK$, the composition $\Phi(w)=\phi_{w_n}\circ\cdots\circ\phi_{w_1}(M_0,\sigma_0)$ is well-defined, in the sense that the mapping $\phi_c$ is only applied to pairs $(M,\sigma)\in\mT_p$ having both a left edge and a right edge. 
\end{remark}

It will be convenient to introduce some explicit correspondences between the steps of $w\in\mK$ and the vertices, faces, and edges of $\Phi(w)$. These correspondences defined below naturally follow from the step-by-step construction of $\Phi(w)$. 

\begin{definition}\label{def:eta}
	Let $w=w_1\cdots w_n\in \mK$ and let $(M,\sigma)=\Phi(w)=\phi_{w_n}\circ\cdots\circ\phi_{w_1}(M_0,\sigma_0)$. 
	\begin{compactitem}
		\item We call \emph{in-edges} of $(M,\sigma)$ the edges which are neither a left edge, a right edge or the top-edge. Let $E$ be the set of in-edges of $(M,\sigma)$. We define the mapping $\etae$ from $[n]$ to $E$ as follows. For each step $w_i$, applying $\phi_{w_i}$ makes the top-edge become an in-edge $e$ of $(M,\sigma)$, and we set $\etae(i)=e$.
		\item Let $V$ and $F$ be the sets of inner vertices and inner triangles of $M$ respectively. 
		We define the mapping $\etavf$ from $[n]$ to $V\cup F$ as follows. If $w_i$ is an $a$-step or a $b$-step, then applying $\phi_{w_i}$ adds one inner triangle $f$ to $(M,\sigma)$ and we set $\etavf(i)=f$. If $w_i$ is a $c$-step, then applying $\phi_{w_i}$ adds one inner vertex $v$ to $(M,\sigma)$ (more precisely, either the top-left or top-right vertex becomes an inner vertex $v$) and we let $\etavf(i)=v$.
		\item For an unmatched $a$-step (resp. $b$-step) $w_i$, it is easy to see that the triangle $\etavf(i)$ is incident to a left (resp. right) edge $e$ of $(M,\sigma)$, and we set $\etalr(i)=e$. 
	\end{compactitem}
\end{definition}

\begin{thm}[Bijection: restricted case]\label{thm:bij-OB}
	The mapping $\Phi$ is a bijection between $\mK$ and $\mT_P$. Moreover, for a walk $w\in\mK$ and its image $\Phi(w)=(M,\sigma)$, we have the following correspondences. 
	\begin{compactitem}
		\item[(i)] The mapping $\etae$ gives a one-to-one correspondence between the steps of $w$ and the in-edges of $(M,\sigma)$.
		\item[(ii)] The mapping $\etavf$ gives a one-to-one correspondence between the $a$-steps and $b$-steps of $w$ and the inner triangles $M$. The mapping $\etavf$ also gives a one-to-one correspondence between the $c$-steps of $w$ of type $a$ (resp. $b$) with the white (resp. black) inner vertices of $(M,\sigma)$.
		\item[(iii)] The mapping $\etalr$ gives a one-to-one correspondence between the unmatched $a$-steps (resp. $b$-steps) of $w$ and the left (resp. right) edges of $(M,\sigma)$.
	\end{compactitem}
\end{thm}

Theorem~\ref{thm:bij-OB} will be proved in Section~\ref{sec:appendix-bijection}. The proof is based on a non-trivial reinterpretation of a bijection of the first author \cite{bernardi-dfs-bijection}.

Note that Theorem~\ref{thm:bij-OB} implies in particular that the subset $\mK^{(0,0)}$ of non-empty walks in $\mK$ ending at $(0,0)$ (i.e.\ non-empty words with no unmatched steps) is in bijection with the set $\mT_P^{(0,0)}$ of percolated near-triangulations $(M,\sigma)$ in $\mT_P$ having two outer edges. Upon removing the top-edge of $(M,\sigma)\in\mT_P^{(0,0)}$, we obtain a percolated triangulation (without boundary). 
Clearly this operation allows to identify $\mT_P^{(0,0)}$ with the set of percolated loopless rooted triangulations such that the root-edge is oriented from a white vertex to a black vertex. We summarize:

\begin{cor} \label{cor:bij-path-returning-0}
	The mapping $\Phi$ induces a bijection $\Phi_0$ between the set $\mK^{(0,0)}$ of non-empty walks in $\mK$ ending at $(0,0)$ and the set of site-percolated loopless rooted triangulations (without boundary) such that the root-edge is oriented from a white vertex to a black vertex. Moreover, if $w\in\mK^{(0,0)}$ and $\Phi(w)=(M,\sigma)$, then the mapping $\etae$ gives a one-to-one correspondence between the steps of $w$ and the edges of $M$.
\end{cor}

\subsection{Bijection: general case}\label{subsec:bij-chordal}
In this section, we generalize the bijection $\Phi$ to a larger class of Kreweras walks $\bmK\supset \mK$. As a special case, we obtain a bijective encoding of all site-percolated triangulations with a simple boundary satisfying the root-interface condition. 

Recall that $\mK$ is the set of words in $\{a,b,c\}^*$ such that every $c$-step has both a matching $a$-step and a matching $b$-step.
We denote by $\bmK$ the set of words in $\{a,b,c\}^*$ such that every $c$-step has at least one matching step. We now define the set $\bmTP$ of percolated near-triangulations which is in bijection with $\bmK$. 

\begin{definition}\label{def:bmTP}
	We denote by $\bmTP$ the set of percolated near-triangulations $(M,\sigma)$ such that 
	\begin{compactitem}
		\item $M$ is in $\mT$ (i.e.\ a rooted loopless triangulation with a simple boundary), and its outer edges are marked as either \emph{active} or \emph{inactive},
		\item the marking of the outer edges is such that the top-edge $e$ is active, the root-edge $e'$ is inactive (except in the case $e'=e$), and the active edges are consecutive around the root-face of $M$, 
		\item $\sigma$ is a site-percolation configuration satisfying the root-interface condition,
		\item the percolation path goes through all of the inner triangles incident to an active outer edge.
	\end{compactitem}
\end{definition}

\fig{width=.8\linewidth}{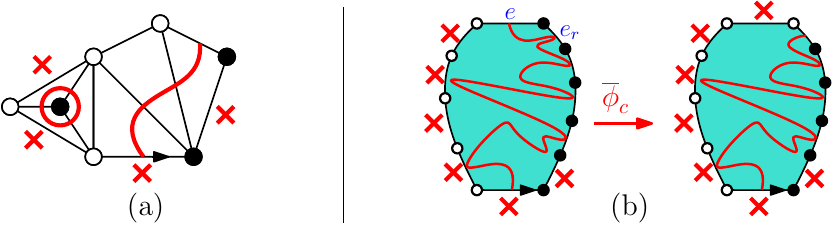}{(a) A pair $(M,\sigma)$ in $\bmTP$, and its percolation interfaces. The inactive edges are indicated by crosses. (b) The mapping $\bphi_c$.}

A pair $(M,\sigma)$ in $\bmTP$ is represented in Figure~\ref{fig:extension-bijection}(a). Note that $\mT_P$ can be thought as the subset of pairs $(M,\sigma)$ in $\bmTP$ such that all the non-root outer edges are active. \referee{Note also that the last condition of Definition \ref{def:bmTP} becomes void when every outer edge other than the top-edge are inactive, and this will be exploited in Corollary \ref{cor:chordal-case}.}
Next, we define the mapping $\bPhi$, which extends $\Phi$ from $\mK$ to $\bmK$.

\begin{definition}
	An example is given in Figure~\ref{fig:example-bPhi}. Let $\phi_a$, $\phi_b$, $\phi_c$, and $(M_0,\sigma_0)$ be as in Definition~\ref{def:Phi}. For $w=w_1\cdots w_n$ in $\bmK$, we define $\bPhi(w)=\bphi_{w_n}\circ\cdots\circ\bphi_{w_1}(M_0,\sigma_0)$, where $\bphi_a$, $\bphi_b$, and $\bphi_c$ are extensions of $\phi_a$, $\phi_b$, $\phi_c$ defined as follows.
	\begin{compactitem}
		\item The mappings $\bphi_a$ and $\bphi_b$ are defined exactly as $\phi_a$ and $\phi_b$: see Figure~\ref{fig:bij-rules-triang}(a,b). The outer edges created by applying $\bphi_a$ and $\bphi_b$ are marked as active.
		\item The mapping $\bphi_c$ is defined on the set of pairs $(M,\sigma)$ in $\bmTP$ with at least one active left or right edge. For pairs $(M,\sigma)$ such that there is both an active left edge and an active right edge, the mapping $\bphi_c$ does the same as $\phi_c$: see Figure~\ref{fig:bij-rules-triang}(c) (and the outer edges keep their active/inactive status). 
		For pairs $(M,\sigma)$ without an active left edge, we consider the top-edge~$e$, the top-right vertex $v_r$ and the incident right edge $e_r$. The pair $\bphi_c(M,\sigma)$ is obtained by recoloring $v_r$ in white (so $e$ becomes a left edge and $e_r$ becomes the top-edge) and setting the edge $e$ to be inactive (while $e_r$ remains active). This is illustrated in Figure~\ref{fig:extension-bijection}(b). For pairs $(M,\sigma)$ without an active right edge, $\bphi_c(M,\sigma)$ is defined symmetrically. 
	\end{compactitem}
\end{definition}

\fig{width=\linewidth}{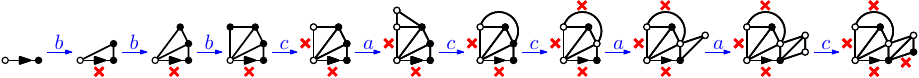}{The construction of $\bPhi(w)$ for $w=bbbcaccaac$. Inactive edges are indicated by crosses.}

\begin{remark}
	It is easy to check (by induction on the length of the walk $w$) that the numbers of active left and right edges of $\bPhi(w)$ are equal to the numbers of unmatched $a$-steps and $b$-steps respectively. Moreover, for $w=w_1\cdots w_n\in \bmK$,  each $c$-step $w_i$ of $w\in\bmK$ has at least one matching $a$-step or $b$-step, so the prefix $w'=w_1\cdots w_{i-1}$ has at least one unmatched $a$-step or $b$-step.
	Thus for $w=w_1\cdots w_n\in \bmK$, the composition $\bPhi(w)=\bphi_{w_n}\circ\cdots\circ\bphi_{w_1}(M_0,\sigma_0)$ is well-defined, in the sense that the mapping $\bphi_c$ is only applied to pairs $(M,\sigma)\in\mT_p$ having at least one active left or right edge.
\end{remark}

Before stating the generalization of Theorem~\ref{thm:bij-OB}, we need some additional definitions.
\begin{definition}\label{def:eta-extend}
	Let $w=w_1\cdots w_n\in \bmK$ and let $(M,\sigma)=\bPhi(w)$. Recall that the $c$-steps are called \emph{matched} if they have both a matching $a$-step and a matching $b$-step, and \emph{unmatched} otherwise.
	We now extend Definition~\ref{def:eta} of the mapping $\etae$, $\etavf$, and $\etalr$. 
	\begin{compactitem}
		\item We call \emph{in-edges} of $(M,\sigma)$ the edges which are not active outer edges, and we let $E$ be the set of in-edges. The mapping $\etae:[n]\to E$ is as in Definition~\ref{def:eta}: for each step $w_i$, the application of $\bPhi_{w_i}$ turns the top-edge into an in-edge $e$ of $(M,\sigma)$ and we denote $\etae(i)=e$.
		\item We call \emph{in-vertices} the vertices not incident to an active outer edge. Let $V$ be the set of in-vertices and let $F$ be the set of inner triangles. The mapping $\etavf:[n]\to V\cup F$ is as in Definition~\ref{def:eta}, except that if $w_i$ is an unmatched $c$-step, then we define $\etae(i)$ to be the top-left or top-right vertex which becomes an in-vertex by application of $\bPhi_{w_i}$. 
		\item The matching $\etalr$ is defined as before (on the set of indices of unmatched $a$-steps and $b$-steps).
	\end{compactitem}
\end{definition}

\begin{thm}[Bijection in the general case]\label{thm:bij-extend} 
	The mapping $\bPhi$ is a bijection between $\bmK$ and $\bmTP$. Moreover, if $w\in\bmK$, and $\bPhi(w)=(M,\sigma)$ then we have the following correspondences.
	\begin{compactitem}
		\item[(i)] The mapping $\etae$ gives a one-to-one correspondence between the steps of $w$ and the in-edges of $(M,\sigma)$. Moreover, the unmatched $c$-steps of type $a$ (resp. type $b$) correspond through $\etae$ to the left (resp. right) inactive outer edges of $M$.
		\item[(ii)] The mapping $\etavf$ gives a one-to-one correspondence between the $a$-steps and $b$-steps of $w$ and the inner triangles $M$. Moreover, the unmatched $a$-steps (resp. $b$-step) correspond through $\etavf$ to the triangles incident to a white (resp. black) active outer-edge.
		The mapping $\etavf$ also gives a one-to-one correspondence between the matched $c$-steps of $w$ of type $a$ (resp. $b$) and the white (resp. black) in-vertices of $(M,\sigma)$. Moreover, the unmatched $c$-steps of type $a$ (resp. $b$) correspond to the white (resp. black) outer in-vertices.
		\item[(iii)] The mapping $\etalr$ gives a one-to-one correspondence between the unmatched $a$-steps (resp. $b$-steps) of $w$ and the left (resp. right) active edges of $(M,\sigma)$.
	\end{compactitem}
\end{thm}

Theorem~\ref{thm:bij-extend} is proved in Section~\ref{sec:appendix-bijection-extended}. Note that the restriction of Theorem~\ref{thm:bij-extend} to the set $\mK\subset \bmK$ of walks is exactly Theorem~\ref{thm:bij-OB}. We now obtain another important specialization by considering the subset $\smK\subset \bmK$ of walks without any unmatched $a$-step nor unmatched $b$-step. 
\begin{definition}\label{def:smTp}
	We denote by $\smTP$ the set of percolated near-triangulations $(M,\sigma)$ such that 
	\begin{compactitem}
		\item $M$ is in $\mT$ (i.e.\ a rooted loopless triangulation with a simple boundary),
		\item $\sigma$ is a site-percolation configuration satisfying the root-interface condition.
	\end{compactitem}
\end{definition}

Clearly, the set $\smTP$ identifies with the set of elements in $\bmTP$ without active left edges nor active right edges. Thus, Theorem~\ref{thm:bij-extend} immediately implies the following correspondence.

\begin{corollary}\label{cor:chordal-case}
	The mapping $\bPhi$ induces a bijection between the set $\smK$ of walks (walks such that all the $a$-steps and $b$-steps are matched, and the all the $c$-steps have at least one matching step), and the set $\smTP$ of percolated near-triangulations defined in Definition~\ref{def:smTp}.
	Moreover, for a percolated near-triangulation $(M,\sigma)\in\smTP$ and the corresponding walk $w\in\smK$, we have the following correspondences.
	\begin{compactitem}
		\item[(i)] The $a$-steps and $b$-steps of $w$ are in one-to-one correspondence with the inner triangles $M$. 
		\item[(ii)] The matched $c$-steps of type $a$ (resp. $b$) are in one-to-one correspondence with the white (resp. black) inner vertices of $(M,\sigma)$.
		\item[(iii)] The unmatched $c$-steps of type $a$ (resp. $b$) are in one-to-one correspondence with the white (resp. black) outer vertices not incident to the top-edge.
	\end{compactitem}
\end{corollary}

\begin{remark}\label{rmk8}
	An important special case of Corollary~\ref{cor:chordal-case} mentioned in the introduction, corresponds to percolated near-triangulations having all the outer vertices of the same color, except one. Precisely, $\bPhi$ induces a bijection between percolated near-triangulations in $\smTP$ having $n$ interior vertices, and $k+2$ outer vertices of which only one is white (hence having $3n+2k+1$ edges by the Euler formula), and the set of all walks $w$ having $3n+2k$ steps in $\{a,b,c\}$, which start at $(0,0)$, end at $(0,-k)$ and stay in the quadrant $\{(x,y)\,|\,x\geq 0,y\geq -k\}$. This correspondence is illustrated in Figure~\ref{fig:example-bij-intro-disc}. 
\end{remark}

\subsection{Future/past decomposition and an alternative description of the bijection}\label{subsec:future-past}
In this subsection we give an alternative presentation of the bijection $\bPhi$. We then define a canonical decomposition of a site-percolated triangulation with a marked edge into two near-triangulations separated by a cycle: a ``past'' near-triangulation and a ``future'' near-triangulation. \\

\noindent\textbf{Alternative description of the bijection $\bPhi$.} 
Consider the percolated maps $T_a$, $T_b$, $T_\ta$, $T_\tb$ represented in Figure~\ref{fig:bricks}, and that we shall call \emph{bricks}. The bricks have both a \emph{root-edge} indicated by a single arrow, and a \emph{top-edge} which is an oriented edge indicated by a double arrow. 
\fig{width=.8\linewidth}{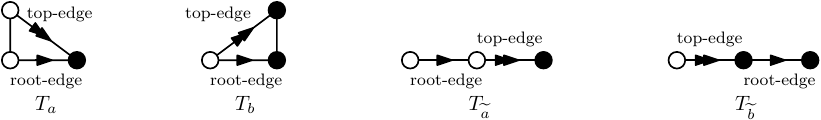}{The bricks $T_a$, $T_b$,$T_\ta$, and $T_\tb$.}

Given a walk $w$ in $\bmK$, we consider the word $\wt{w}=\tw_1\cdots\tw_n\in\{a,b,\ta,\tb\}$ obtained from $w$ by replacing every $c$-step of type $a$ by $\ta$ and every $c$-step of type $b$ by $\tb$. We then construct a percolated near-triangulation $T_{w}$ from the bricks $T_{\tw_1},\ldots ,T_{\tw_n}$ by gluing them by the following process (which is illustrated in Figure~\ref{fig:example-bij-alternative}): 
\begin{compactitem}
	\item For all $i\in[n-1]$ (in this order), we glue the top-edge of $T_{\tw_i}$ to the root-edge of $T_{\tw_{i+1}}$ (respecting their orientation) and if the colors of the glued vertices differ, we keep the color of the vertices of $T_{\tw_{i+1}}$. We denote by $\tT_w$ the percolated near-triangulation obtained at this stage.
	\item The map $\tT_w$ has a root-edge (the root-edge of $T_{\tw_1}$) and a top-edge (the top-edge of $T_{\tw_n}$). We call \emph{left-sides} (resp. \emph{right-sides}) of $\tT_w$ the edges-sides followed when turning around the root-face of $T_w$ in clockwise (resp. counterclockwise) direction from the root-edge to the top-edge. The left-sides are said to be \emph{opening} if they correspond to edges of bricks of the form $T_a$ or $T_b$ and \emph{closing} if they correspond to edges of bricks of the form $T_{\ta}$ or $T_{\tb}$. Note that the sequence of opening and closing left-sides of $\tT_w$ (in clockwise order around the root-face) is the same as the sequence of $a$-steps and $c$-steps in $w$. We say that an opening and a closing left-side are \emph{matching} if the corresponding $a$ and $c$ steps are matching. It is easy to see that one can glue all the pairs of matching left-sides together in a planar manner. Symmetrically, the opening and closing right-sides correspond to the $b$-steps and $c$-steps of $w$, and the pairs of matching right sides can all be glued together in a planar manner. We denote by $T_w$ the percolated \referee{near-triangulation} obtained at the end of this gluing process.
	
	The color of the vertices of $T_w$ are determined as follows: if a vertex $v$ of $T_w$ comes from gluing the vertices $v_1,\ldots,v_k$ of $\tT_w$ together, we consider the last brick $T$ of $T_w$ incident to one of these vertices (the bricks are considered ordered from $T_{\tw_1}$ to $T_{\tw_n}$), and $v$ takes the color of the (unique) vertex $v_k$ incident to the last brick $T$. 
	
	Lastly, we mark as \emph{active} (resp. \emph{inactive}) the unmatched opening (resp. closing) sides of $T_w$. 
\end{compactitem}

\fig{width=.9\linewidth}{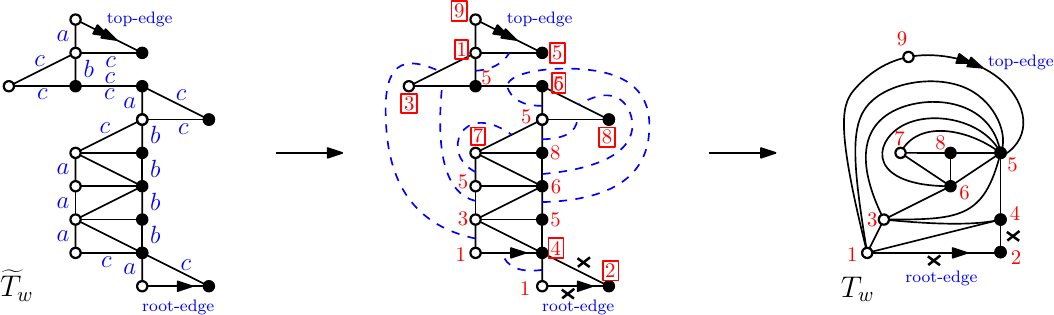}{Alternative description of the bijection $\bPhi$ for the walk $w=acabbaabbcaccbca$. Here $\tw=a\tb abbaabb\ta a\tb\tb b\ta a$. Left: the triangulation with boundary $\widetilde{T}_w$, with the left and right sides labeled by the steps of $w$ they correspond to (so that closing sides are labeled by the letter $c$). Middle: the matching of the opening and closing sides (indicated by dashed lines), and the identification of vertices that gluing these sides in pair will imply. The unmatched closing sides are marked by a cross. Right: the map $T_w$ obtained by gluing the matching sides. The triangulation $T_w$ has 9 vertices, that we labeled (arbitrarily) from 1 to 9; these are the label given in the middle figure. Moreover, for all $v\in[9]$, the ``last vertex'' $v'$ of $\tT_w$ corresponding to $v$ is boxed in the middle figure (as defined above, the color of this last vertex $v'$ in $\tT_w$ determines the color of $v$ in $T_w$).}

The following result is proved in Section~\ref{sec:appendix-future-past}.

\begin{prop}\label{prop:alternative-description} 
	For any walk $w\in \bmK$, the percolated near-triangulations $\bPhi(w)$ and $T_w$ are equal.
\end{prop}

\begin{remark}\label{rk:burger-link}
	An interesting special case of the bijection $\bPhi$ is for walks $w$ in $\bmK$ having no consecutive $c$-step. In this case, the walk $w$ belongs to the set of words $\{a,b,ac,bc\}^*$, and instead of considering the bricks $T_a$, $T_b$, $T_\ta$, and $T_\tb$ as above, one can consider the ``blocks'' $T_a$, $T_b$, $\tT_{ac}$, and $\tT_{bc}$ represented in Figure~\ref{fig:bricks2}. 
	This special case of $\bPhi$ can easily be identified with a bijection described by Sheffield in \cite{shef-burger} (which itself is related to a construction of Mullin \cite{mullin-maps}).
\end{remark}

\fig{width=.7\linewidth}{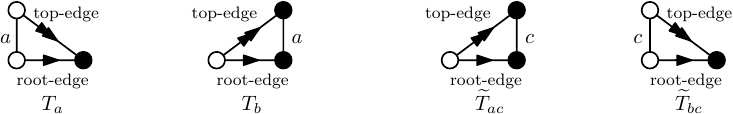}{The ``blocks'' $T_a$, $T_b$, $\tT_{ac}$, $\tT_{bc}$ relevant when $\Phi$ is applied to $w\in\{a,b,ac,bc\}^*\subset \mK$.}

\noindent\textbf{Future/past decomposition.} 
We will now define a canonical decomposition of a percolated triangulation with a marked edge into two near-triangulations separated by a cycle: a ``past'' near-triangulation and a ``future'' near-triangulation. 

Recall that  $\mK^{(0,0)}$ denotes the set of walks in $\mK$ ending at $(0,0)$.
Let $w\in \mK^{(0,0)}$, and let $(M,\sigma)=\Phi(w)$ be the corresponding percolated near-triangulation. Recall that each step of $w$ corresponds to an edge of $M$ via the mapping $\etae$. Given a decomposition $w=uv$ of $w$ into a prefix $u$ and a suffix $v$, we will now describe the corresponding decomposition of $(M,\sigma)$ in two parts: one part called ``past'' made of the edges of $M$ corresponding to the steps in $u$, and another part called ``future'' made of the edges of $M$ corresponding to the steps in $v$.

Observe that all the $c$-steps of $u$ are matched (hence $u$ is in $\mK$), and that all the $a$-steps and $b$-steps of $v$ are matched.
Now consider the $c$-steps of $v$ having neither a matching $a$-step nor a matching $b$-step: these $c$-steps separate (possibly empty) subwords of $v$ in $\smK$. Let us write $v=v_1cv_2c\cdots cv_k$ with $v_i\in \smK$ this decomposition. 
Let $(M,\sigma)=\Phi(w)$, let $(P,\al)=\Phi(u)$, and let $(Q_i,\be_i)=\bPhi(v_i)$. We now describe the process for obtaining $(M,\sigma)$ from the pairs $(P,\al)$ and $(Q_i,\be_i)$. This process is illustrated in Figure~\ref{fig:future-past}.

\fig{width=.75\linewidth}{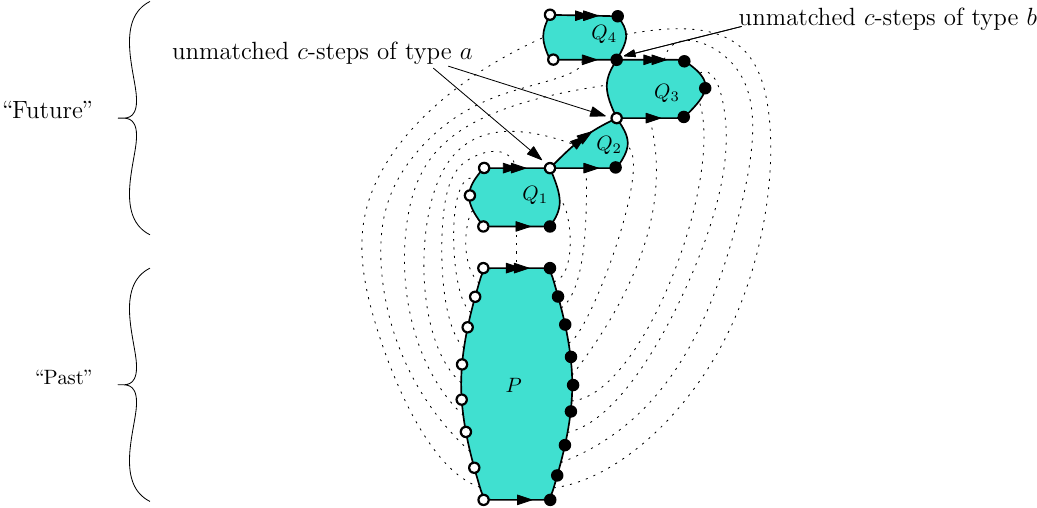}{The future/past decomposition of $(M,\sigma)=\Phi(uv)$. 
	Here $(P,\al)=\Phi(u)$ corresponds to the "past'', while the percolated triangulation with boundary $(Q,\be)$ obtained by gluing $(Q_1,\be_1)=\bPhi(v_1),\ldots,(Q_k,\be_k)=\bPhi(v_k)$ together corresponds to the ``future''. 
}

We denote by $(Q,\be)$ the percolated near-triangulation obtained by gluing together the maps $(Q_1,\be_1),\ldots,(Q_k,\be_k)$ according to the following \emph{future rule}.
\smallskip

\noindent \emph{Future rule:} For all $i\in[k-1]$, if the $c$-step between the subwords $v_i$ and $v_{i+1}$ is of type $a$ (resp. type $b$) in $w$, then the top-right (resp. top-left) vertex of $(Q_i,\beta_i)$ is recolored in white (resp. black) and glued to the white (resp. black) endpoint of the root-edge of $(Q_{i+1},\beta_{i+1})$. 

\smallskip
We take the root-edge of $(Q,\be)$ to be the root-edge of $(Q_1,\be_1)$ and we call \emph{top-edge} of $(Q,\be)$ the top-edge of $(Q_k,\be_k)$. We call \emph{left-sides} (resp. \emph{right-sides}) of $(Q,\be)$ the edge-sides followed while walking around the root-face of $(Q,\be)$ in clockwise (resp.counterclockwise) direction from the root-edge to the top-edge.
We claim that the number $n_\ell$ of left edges of $(P,\al)$ is equal to the number $n_\ell'$ of left-sides of $(Q,\be)$.
Indeed $n_\ell$ is the number of unmatched $a$-steps of $u$, while $n_\ell'$ is the number of $c$-step of $v$ without a matching $a$-step. Similarly, the number $n_r$ of right edges of $(P,\al)$ is equal to the number $n_r'$ of right-sides of $(Q,\be)$.
We can now state our result, which is proved in Section~\ref{sec:appendix-future-past}.

\begin{prop}\label{prop:future-past} Let $w=uv$ be a walk in $\mK^{(0,0)}$.
	With the above notation, the percolated near-triangulation $(M,\sigma)=\Phi(w)$ is obtained by first gluing the top-edge of $(P,\al)$ to the root-edge of $(Q,\be)$, and then gluing the left and right edges of $(P,\alpha)$ to the non-top outer edges of $Q$ in the unique planar manner, while keeping the color of the outer vertices of $(Q,\be)$ (irrespective of the color of the vertices of $(P,\alpha)$ to which they are glued).
\end{prop}

\subsection{Bijection: infinite case}
\label{sec:bij-inf}
\textbf{Infinite planar maps and the UIPT.}
In this section we consider \emph{infinite planar maps}. We call \emph{infinite graph}, a graph with an infinite but countable number of vertices and edges, which is \emph{locally finite} (that is, every vertex has finite degree).  An infinite graph $G$ is said to be \emph{one-ended} if for any finite subgraph $H$, $G\setminus H$ contains exactly one infinite connected component.
We call \emph{infinite planar map} an embedding of an infinite planar graph in the sphere without edge crossing, considered up to orientation-preserving homeomorphism. An infinite planar map is \emph{rooted} if an edge is marked and oriented.
Note that a one-ended infinite planar map $M$ can be drawn in the sphere with a single accumulation point; equivalently $M$ can be drawn in the plane without accumulation point. An \emph{infinite triangulation} is a infinite planar map such that every face has degree~3.
\begin{definition}
	We denote by $\imT$ the set of rooted one-ended infinite loopless triangulations. We denote by $\imTP$ the set of pairs $(M,\si)$ such that $M$ is in $\imT$, and $\si$ is a coloring of its vertices in black and white.
\end{definition}
Recall that in \cite{angel-schramm-uipt}, Angel and Schramm defined the \emph{uniform infinite planar triangulation} (UIPT) of types II, and III. The UIPT of type II, which is also called the \emph{loopless UIPT}, is a probability measure on the set of infinite rooted loopless planar triangulations, which corresponds to the \emph{local limit} of finite loopless triangulations. By this, we mean that for any $k\in \NN$ the probability distribution $P_k$ of the ball of radius $k$ (for the graph distance) around the root-vertex of the UIPT is equal to the limit in law of the probability distribution $P_{k,n}$ of the ball of radius $k$ around the root-vertex of a uniformly random rooted loopless triangulation with $n$ vertices. In \cite{angel-peeling,angel-scaling-limit}, Angel further considered the \emph{critical site-percolated loopless UIPT}, in which the vertices of the UIPT (of type II) are independently colored in black or white with probability 1/2. The UIPT (of any type) is almost surely one-ended, hence the critical site-percolated loopless UIPT is almost surely in $\imTP$. Therefore we may view the law of the UIPT as a probability distribution on $\imT$.
\begin{definition}\label{def:percolated-UIPT}
	The \emph{UIPT distribution} on $\imT$ is defined to be the law of the UIPT viewed as a probability distribution on $\imT$. We say that an element $(M,\sigma)$ of $\imTP$ has been sampled from the \emph{distribution of the percolated UIPT} if $M$ has the law of the UIPT on $\imT$ and $\si$ is a uniform and independent coloring of the vertices of $M$.
\end{definition}
A \emph{near-triangulation with infinite boundary}, is a one-ended infinite map with every face of degree~3, except one face of infinite degree. It is a \emph{near-triangulation with infinite simple boundary} if the boundary of the infinite face is a simple (bi-infinite) path. 
Such maps are sometimes said to have \emph{half-plane topology}, since they can be drawn in the plane without accumulation points with the infinite face being the half plane $\RR\times \RR^+$.\\

\noindent \textbf{Infinite Kreweras walks and normality.}
We denote by $\{a,b,c\}^\ZZ$ the set of \emph{bi-infinite} words on the alphabet $\{a,b,c\}$ (that is, words whose steps are indexed by $\ZZ$). For $w\in\{a,b,c\}^\ZZ$, we denote by $w^-=\dots w_{-2}w_{-1}$ the infinite prefix of $w$ made of all steps indexed by negative integers, and we denote by $w^+=w_0w_1\dots$ the infinite suffix of $w$ made of all steps indexed by non-negative integers. As before we identify words in $\{a,b,c\}^\ZZ$ with infinite lattice walks in $\ZZ^2$ with steps $a=(1,0)$, $b=(0,1)$, and $c=(-1,-1)$ such that the steps are indexed by $\ZZ$ and the position of the walk just before step $w_0$ is $(0,0)$.

We call a word in $\{a,b,c\}^\ZZ$ \emph{fully-matched} if every $a$-step and every $b$-step has a matching $c$-step, and every $c$-step has both a matching $a$-step and a matching $b$-step. Let $w=w^-w^+\in\{a,b,c\}^\ZZ$ be a fully matched word. We call \emph{cut-time} of $w^+$ a $c$-step of $w^+$ such that its matching $a$-step and $b$-step are both in $w^-$. We call \emph{split-time} of $w^-$ a step $w_{-k}$ such that the subword $w_{-k}w_{-k+1}\cdots w_{-1}$ is in $\mK$. 
We call $w\in\{a,b,c\}^\ZZ$ \emph{normal} if it is fully-matched, and moreover $w^+$ has infinitely many cut-times, and $w^-$ has infinitely many split-times. We denote by $\imK$ the set of normal words $w\in \{a,b,c\}^\ZZ$. We will prove the following result in Section~\ref{app:inf}.
\begin{lemma}\label{prop:normal}
	Let $w\in\{a,b,c\}^\ZZ$ be a random word such that the steps are sampled uniformly and independently at random. Then $w\in\imK$ (i.e., $w$ is normal) with probability 1.
\end{lemma}

We will typically sample words from the distribution of the above lemma, which motivates the following definition.
\begin{definition}
	The \emph{uniform distribution} on $\{a,b,c\}^\ZZ$ is defined to be the probability distribution where the steps are sampled uniformly and independently at random. By Lemma \ref{prop:normal}, we may view the uniform distribution as a distribution on $\imK$, and we will call it the  \emph{uniform distribution} on $\imK$.
\end{definition}

\noindent \textbf{Bijection in the infinite setting.}
Let $w=w^-w^+\in\imK$ be an infinite normal Kreweras walk. 
Consider the decomposition $w^+=w^+_1cw^+_2cw^+_3c\ldots$ with the subwords $w^+_i$ separated by cut-times of $w^+$ (and containing no cut-times), and the decomposition $w^-=\ldots w^-_{3}w^-_{2}w^-_{1}$ with the subwords $w^-_i$ starting with a split-time of $w^-$ (and containing no other split-times). Since $w$ is normal, all the subwords $w^-_i$ and $w^+_i$ are finite. Note moreover that for all $i$, $w^+_i$ is in $\smK$ and $w^-_i$ is in $\mK$. We let $\Phi(w^-_i)=(P_i,\al_i)\in\mT_P$ and $\bPhi(w^+_i)=(Q_i,\be_i)\in\bmTP$.
To $w$ we associate two site-percolated near-triangulations:
\begin{compactitem} 
	\item The \emph{past site-percolated near-triangulation} $(M^-,\si^-)$ (corresponding to $w^-$) is obtained by gluing the root-edge of $(P_i,\al_i)$ to the top-edge of $(P_{i+1},\al_{i+1})$ for all $i\geq 1$. 
	\item The \emph{future site-percolated near-triangulation} $(M^+,\si^+)$ (corresponding to $w^+$) is obtained by applying the \emph{future rule} described in Section~\ref{subsec:future-past} to the sequence of pairs $(Q_i,\be_i)$.
\end{compactitem} 
By definition, $(M^-,\si^-)$ and $(M^+,\si^+)$ are near-triangulations with infinite boundaries. 

The past and future site-percolated near-triangulations can be glued as in Figure~\ref{fig:future-past}: the top-edge of $M^-$ (i.e.\ the top-edge of $P_1$) is glued to the root-edge of $M^+$ (i.e.\ the root-edge of $Q_1$), and then the other outer edges of $M^-$ are glued to the other outer edges of $M^+$, with the colors of the vertices in $M^+$ determining the color of the glued vertices. This creates an infinite site-percolated loopless triangulation with a marked edge (the root-edge of $M^+$) that we denote by $\Phi^\infty(w)$. It is easy to see that $\Phi^\infty(w)$ is in $\imTP$.
Given $w\in\imK$, and $(M,\si)=\Phi^\infty(w)$, we define the bijections $\etae:\Z\to E$ and $\etavf:\Z\to V\cup F$ as in Definition~\ref{def:eta}, where $V,E,F$ are the set of vertices, edges and faces of $M$ respectively.

The following result is proved in Section~\ref{app:inf}.
\begin{theorem}\label{thm:UIPT}
	Let $w\in\{a,b,c\}^{\ZZ}$, and assume the steps of $w$ are sampled uniformly and independently at random. Recall that $w$ is normal (equivalently, $w\in\imK$) with probability 1, so the infinite site-percolated triangulation $\Phi^\infty(w)\in\imTP$ is well-defined almost surely. We consider the random infinite rooted triangulation $(M,\sigma)$ obtained from $\Phi^\infty(w)$ by taking the marked edge as the root-edge and orienting it with probability 1/2 in either direction. Then $(M,\sigma)$ is distributed like the percolated UIPT (Definition \ref{def:percolated-UIPT}). Furthermore, the walk $w$ is almost surely determined by $(M,\sigma)$.
\end{theorem}

\begin{remark}
	Observe from Theorem \ref{thm:UIPT}, that although the mapping $\Phi^\infty$ is not a bijection between $\mK^\infty$ and $\imTP$ (for instance, percolated triangulations with every vertex white are not in the image), it is a (measure preserving) bijection between two subsets of measure 1.
	Observe also that for $(M,\sigma)\in\imTP$ chosen according to the percolated UIPT distribution, the past and future  site-percolated near-triangulations $(M^-,\si^-)$ and  $(M^+,\si^+)$ are almost surely determined by $(M,\sigma)$ (because the walk $w$ is almost surely determined by $(M,\sigma)$).
\end{remark}

\section{Discrete dictionary I: Spine-looptrees decomposition}\label{sec:spine}

In this section, we explore the geometry of the percolation path of a percolated triangulation $(M,\sigma)$ in $\bmT_P$ or $\imTP$. In Section~\ref{sec:LR-decomposition-map}, starting from a (finite) percolated triangulation $(M,\sigma)\in\bmT_P$ we define  a submap $\LR(M,\sigma)$ encoding the geometry of percolation path. In Section~\ref{sec:LR-decomposition-walk} we explain how this submap can be obtained directly as a function of the walk $w=\bPhi^{-1}(M,\sigma)$. The infinite volume case is treated in Section~\ref{sec:LR-decomposition-infinite}.

\subsection{Spine-looptrees decomposition of a site-percolated triangulation}\label{sec:LR-decomposition-map}
Let us first define precisely the geometric information we will capture about the percolation path.

\begin{definition} \label{def:LR-M}
	Let $(M,\sigma)\in\bmT_P$, and let $T$ be the set of inner triangles on the percolation path. We denote by $\LR(M,\sigma)$ the site-percolated map obtained from $(M,\sigma)$ by 
	\begin{compactenum}
		\item deleting all the inner edges not incident to one of the triangles in $T$,
		\item replacing each unicolor edge incident to two triangles in $T$ by a double edge,
		\item replacing each unicolor inactive outer edge incident to a triangle in $T$ by a double edge.
	\end{compactenum}
	The outer edges of $\LR(M,\sigma)$ can be identified with the outer edges of $M$, and are still marked as either active or inactive.
\end{definition}

\fig{width=.6\linewidth}{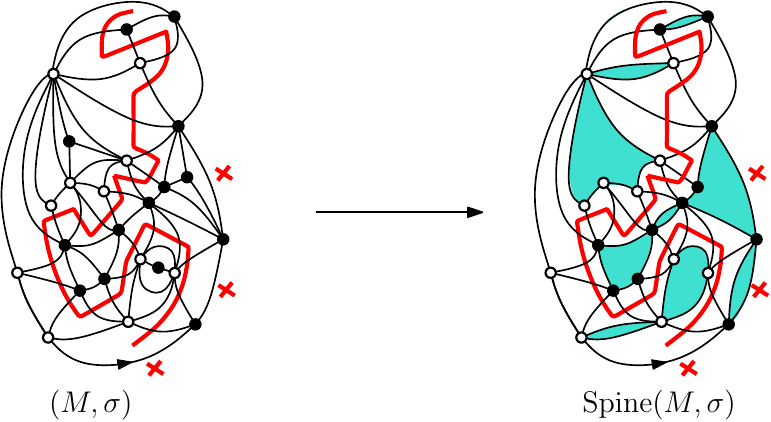}{A percolated near-triangulation $(M,\sigma)\in\bmT_P$, and the associated percolated map $\LR(M,\sigma)$. The percolation path is indicated by a bold line, and the inactive outer edges are indicated by crosses.}

The definition of $\LR(M,\sigma)$ is illustrated in Figure~\ref{fig:LR-map2}. The percolated map $\LR(M,\sigma)$ has three parts: the white cluster, the black cluster, and the set of bicolor edges. 
To describe each cluster we will use the concept of \emph{looptree} which is a collections of cycles glued along a tree structure.
Then it is natural to think of $\LR(M,\sigma)$ as a ``shuffle of two looptrees'' as represented in Figure~\ref{fig:generic-spine}. 
\referee{The concept of looptree is due to  Curien and Kortchemski \cite{curien-kortchemski-looptree-def} who used it to study percolation clusters of the site-percolated UIPT \cite{curien-kortchemski-looptree-perc}. Recently, Richier \cite{Richier:looptrees} also used shuffles of two looptrees to study percolation clusters of the site-percolated UIHPT.}  

The main goal of this section is to show that the shuffle of two looptrees $\LR(M,\sigma)$ can be recovered from the Kreweras walk $w=\bPhi^{-1}(M,\sigma)$ in a very natural manner. In a nutshell, the shuffle of two looptrees is bijectively encoded by the the walk in $\ZZ^2$ obtained from $w$ by replacing each ``cone excursion'' by a single step.

\fig{width=.4\linewidth}{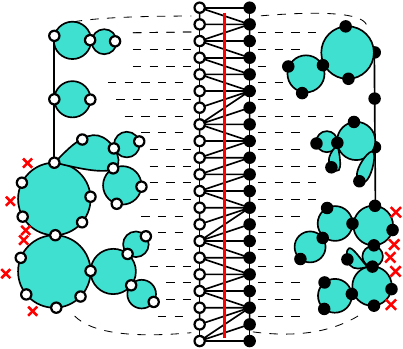}{Generic decomposition of the percolated map $\LR(M,\sigma)$ as a ``shuffle'' of two discrete looptrees. Note that the active left (resp. right) edges of $M$ are bridges of the white (resp. black) cluster of $\LR(M,\sigma)$, but that there is no other bridges in these clusters.}

We start with definitions related to looptrees. A \emph{discrete looptree} is a rooted map such that every edge is incident both to the root-face and to an inner face. A discrete looptree is represented in Figure~\ref{fig:looptree-code}. Clearly, each inner face of a discrete looptree is bounded by a simple cycle, that we call a \emph{bubble}. Discrete looptrees can therefore be thought as obtained by gluing simple cycles along vertices in a tree-like fashion. 

To a discrete looptree $\frk L$ with $k$ edges, we associate the \emph{clockwise code} $\cwcode(\frk L)$, which is a lattice path on $\ZZ$ defined as follows. The path $\cwcode(\frk L)$ starts at $0$ and has $k$ steps $s_1,\ldots, s_k$ corresponding to the $k$ edges $e_1,\ldots,e_k$ followed when turning around $\frk L$ in clockwise direction starting just after the root-edge (so that the root-edge is $e_k$): if $e_i$ is the last edge of a bubble with $d$ edges then $s_i=-d+1$, and otherwise $s_i=1$. The clockwise code is given in Figure~\ref{fig:looptree-code}. 
The \emph{counterclockwise code} $\ccwcode(\frk L)$ is defined similarly excepts one turns around $\frk L$ in counterclockwise direction (the root-edge being still visited last). 
We omit the proof of the following easy lemma.
\begin{lemma} \label{lem:cwcode}
	The mapping $\cwcode$ (resp. $\ccwcode$) is a bijection between discrete looptrees and the set of lattice paths on $\ZZ$ starting and ending at 0, remaining non-negative, and having steps at most $1$. 
\end{lemma}

\fig{width=.9\linewidth}{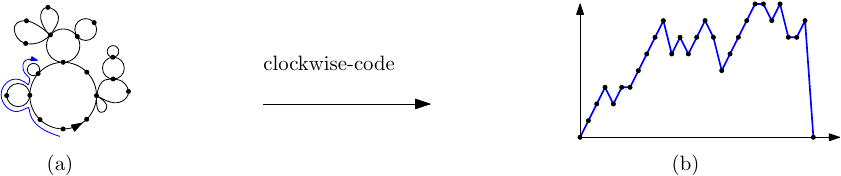}{(a) A discrete looptree $\frk L$. (b) The clockwise-code $\cwcode(\frk L)$.}

Note that for any percolated near-triangulation $(M,\sigma)\in\bmT_P$, the white (resp. black) cluster of $\LR(M,\sigma)$ is made of a sequence of discrete looptrees attached to the path made of the active left (resp. right) outer edges of $M$. The two clusters are attached along the sides of the triangles on the percolation path, so that $\LR(M,\sigma)$ can be thought as a shuffle of two sequences of discrete looptrees. Figure~\ref{fig:generic-spine} represents the generic situation.

\subsection{Spine-looptrees decomposition as a function of the Kreweras walk}\label{sec:LR-decomposition-walk}
We now describe a way to construct a site-percolated map $\LR(w)$ from a walk $w$ in $\bmK$. We will then show that $\LR(M,\sigma)=\LR(w)$ for $(M,\sigma)=\bPhi(w)$. Let $w=w_1\cdots w_n \in\bmK$. If $w_k$ is a $c$-step and $w_i$ is a matching $a$-step or $b$-step, we say that the steps $w_{i+1}\cdots w_{k-1}$ are \emph{enclosed by the matching $w_i,w_k$}. 
Let $w_k$ be a $c$-step and let $w_i$ and $w_j$ be its matching steps with $i<j<k$. We say that $w_i,w_k$ form a \emph{far-matching} (and we say that $w_i$ is \emph{far-matched}), while $w_j,w_k$ form a \emph{close-matching} (and $w_j$ is \emph{close-matched}). If the $c$-step $w_k$ has a unique matching step $w_j$, then the matching $w_j,w_k$ is considered a \emph{close-matching}.

Let us describe the geometry of the subwalk $w'=w_j,w_{j+1},\ldots,w_k$ bounded by a close-matching $w_j,w_k$.
Note first that all the $a$-steps and $b$-steps of $w'$ have a matching $c$-step in $w'$ (because both matching steps of $w_k$ appear before $w_{j+1}$). Moreover, if $w_j$ is an $a$-step (resp. $b$-step), then all the $c$-steps in $w'$ have a matching $a$-step (resp. $b$-step) in $w'$ (because $w_j,w_k$ are matching). Thus $w'$ is in $\smK$. Moreover considering $w'$ as a lattice walk, all the steps $w_{j+1},\ldots,w_{k-1}$ are strictly above and to the right of the step $w_k$, but $w_j$ is not. In fact, these properties exactly characterize subwalks bounded by a close-matching; see Figure~\ref{fig:cone-excursion}. We call $w'$ the \emph{cone excursion starting at} $w_j$. 
The \emph{size} of the cone excursion $w'$ is its number of steps $k-j+1$ (which is $\geq 2$), and its \emph{height} is the distance between its starting point and ending point (which is $\geq 1$).

\fig{width=.6\linewidth}{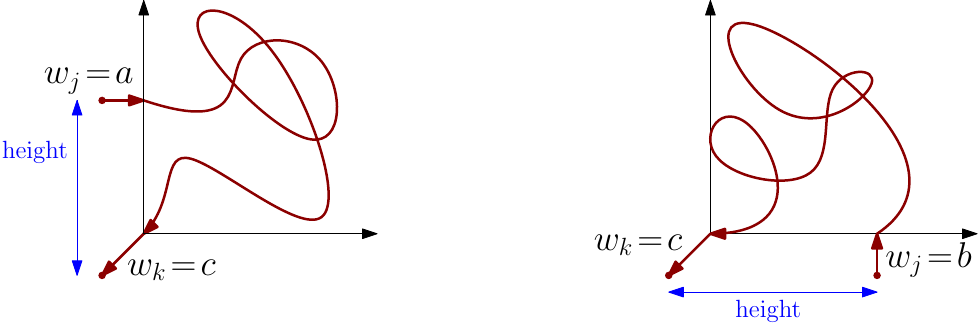}{Left: cone excursion $w_j,w_{j+1},\ldots,w_k$ starting with an $a$-step. Right: cone excursion $w_j,w_{j+1},\ldots,w_k$ starting with a $b$-step. \referee{The steps $w_{j+1}\ldots w_{k-1}$ enclosed by the close matching $w_j,w_k$ stay inside the indicated quadrant.}}

\begin{definition}\label{def:pi-w}
	Let $w\in \bmK$. We call \emph{spine step} of $w$ an $a$-step or $b$-step which is not enclosed by any close-matching. 
	\begin{compactitem}
		\item We denote by $\pi(w)$ the subword of $w$ made of its spine steps.
		\item We denote by $T(0),T(1),\dots,T(s)\in\N$ the indices of the spine steps of $w$, so that $\pi(w)=w_{T(0)}w_{T(1)}w_{T(2)}\cdots w_{T(s)}$.
		\item We denote by $\hpi(w)$ the word on the infinite alphabet $\{a,b\}\cup\{\ba_k,\bb_k,~k\geq 1\}$ obtained from $\pi(w)$ by replacing each close-matched $a$-step (resp. $b$-step) $w_i$ by $\bb_{k}$ (resp. $\ba_{k}$), where $k$ is the height of the cone excursion starting at $w_i$.
	\end{compactitem}
\end{definition}

\fig{width=.8\linewidth}{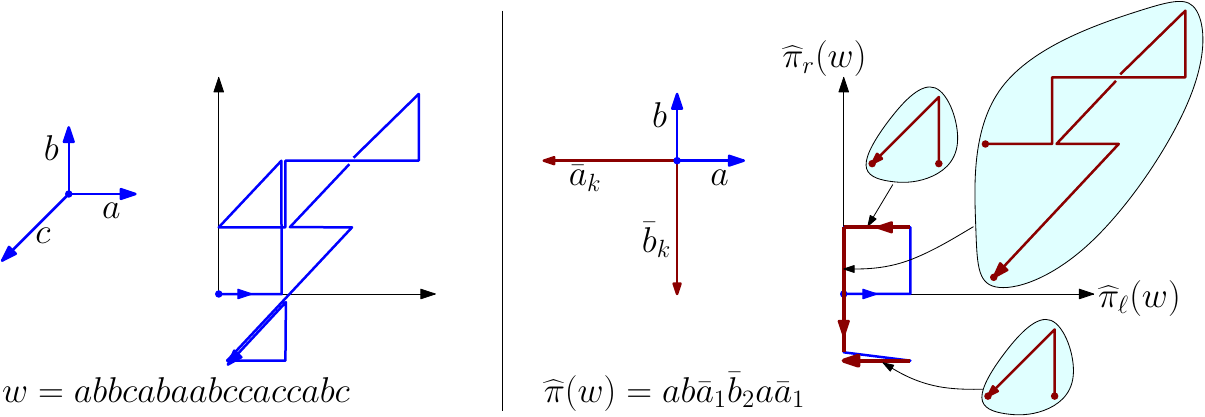}{The walk $w=abbcabaabccaccabc\in\bmK$, and the corresponding walk $\hpi(w)$. \nina{In the example shown, we have $\pi(w)=abbaab$.}}

We identify $\hpi(w)$ with a walk in $\ZZ^2$ starting at $(0,0)$ and with steps $a=(1,0)$, $b=(0,1)$, $\ba_k=(-k,0)$, and $\bb_k=(0,-k)$. This walk is obtained from the walk $w$ by replacing each (maximal) cone excursions by single steps with the same endpoints. This is represented in Figure~\ref{fig:example-path2}. 
Next, we define $\LR(w)$ in terms of $\hpi(w)$.

\begin{definition}\label{def:hpi}
	Let $w\in \bmK$ having $i$ unmatched $a$-steps, $j$ unmatched $b$-steps, $i'$ unmatched $c$-steps of type $a$, and $j'$ unmatched $c$-steps of type $b$. 
	\begin{compactitem}
		\item We denote by $\hpi_\ell(w)$ the word on $\ZZ$ obtained from $\hpi(w)$ by deleting the steps in $\{b\}\cup\{\bb_k,~ k\geq 1\}$, replacing the steps $a$ by $1$ and replacing the steps $\ba_k$ by $-k$ for all $k$, and finally adding $i'$ steps 1 at the beginning and one step $-i$ at the end. 
		\item We denote by $\hpi_r(w)$ the word on $\ZZ$ obtained from $\hpi(w)$ by deleting the steps in $\{a\}\cup\{\ba_k,~ k\geq 1\}$, replacing the steps $b$ by $1$ and replacing the steps $\bb_k$ by $-k$ for all $k$, and finally adding $j'$ steps 1 at the beginning and one step $-j$ at the end. 
		\item We denote by $\hpi_s(w)$ the word on $\{a,b\}$ obtained from $\hpi(w)$ by replacing every step in $\{a\}\cup\{\ba_k,~ k\geq 1\}$ by $a$ and every step in $\{b\}\cup\{\bb_k,~ k\geq 1\}$ by $b$.
	\end{compactitem}
\end{definition}
Definition~\ref{def:hpi} is illustrated in Figure~\ref{fig:LR-theorem2} (left). Roughly speaking, $\hpi_\ell(w)$ and $\hpi_r(w)$ encode the projection of $\hpi(w)$ on the $x$-axis and $y$-axis respectively (modulo time change because we do not record the 0 steps), and $\hpi_s(w)$ encode the way these two walks are shuffled. Note that $\hpi_\ell(w)$ and $\hpi_r(w)$ are walks on $\ZZ$ starting and ending at 0, and remaining non-negative (moreover their steps are at most 1 and only the last step can be 0).

\begin{definition}\label{def:LR-w}
	Let $w\in\bmK$. 
	We denote by $\LR(w)$ the site-percolated map obtained as follows.
	\begin{compactitem}
		\item Let $\frk L_\ell(w)$ and $\frk L_r(w)$ be the discrete looptrees such that $\ccwcode(\frk L_\ell(w))=\hpi_\ell(w)$, and $\cwcode(\frk L_r(w))=\hpi_r(w)$. 
		We color in white the vertices of $\frk L_\ell(w)$ and in black the vertices of $\frk L_r(w)$.
		We call \emph{inactive} the $i'+1$ (resp. $j'+1$) first outer edges of $\frk L_\ell(w)$ (resp. $\frk L_r(w)$) in counterclockwise (resp. clockwise) order around its root-face, starting with the root-edge. The other edges are called \emph{active}.
		\item Let $P_s(w)$ be the site-percolated triangulation $\Phi(\hpi_s(w))$. We denote by $\LR(w)$ the map obtained by gluing the left edges of $P_s(w)$ to the active edges of $\frk L_\ell(w)$, gluing the right edges of $P_s(w)$ to the active edges of $\frk L_r(w)$, and then deleting the root-edges of $\frk L_\ell$ and $\frk L_r$. 
	\end{compactitem}
	The map $P_s(w)$ viewed as a submap of $\LR(w)$ is called the \emph{spine} of $\LR(w)$. 
\end{definition}

\fig{width=\linewidth}{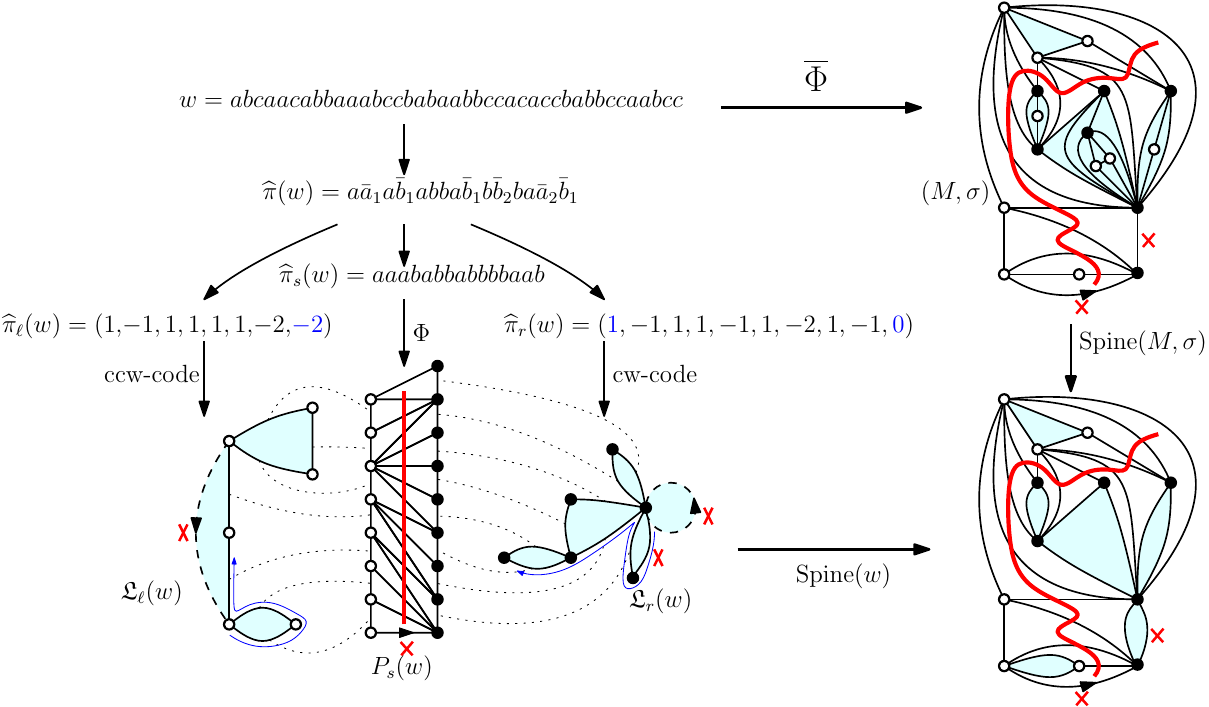}{This figure illustrates Theorem~\ref{thm:LR}. The top row represents a walk $w\in\bmK$ and the corresponding percolated triangulation $(M,\sigma)=\bPhi(w)$. The walk $w$ has $i=2$ (resp. $j=0$) unmatched $a$-steps (resp. $b$-steps) and $i'=0$ (resp. $j'=1$) unmatched $c$-steps of type $a$ (resp. $b$).
	The right column shows how to obtain $\LR(M,\sigma)$. As usual, inactive edges are indicated by crosses. 
	The left column shows how to obtain $\LR(w)$: we find the words $\hpi_\ell(w)$, $\hpi_r(w)$ and $\hpi_s(w)$, and then we compute the maps $\frk L_\ell(w)$, $\frk L_r(w)$, and $P_s(w)$, before gluing them together (bottom row).}

Definition~\ref{def:LR-w} is illustrated in Figure~\ref{fig:LR-theorem2} (left column and bottom row). 
We now state the main result of this section.

\begin{thm}\label{thm:LR}
	Let $w$ be a walk in $\bmK$ and let $(M,\sigma)=\bPhi(w)$ be the associated percolated triangulation. Then the percolated maps $\LR(M,\sigma)$ (obtained by Definition~\ref{def:LR-M}) and $\LR(w)$ (obtained by Definition~\ref{def:LR-w}) are equal. 
	Moreover, we have the following correspondences.
	\begin{compactitem}
		\item[(i)] The spine-steps of $w$, are in one-to-one correspondence, via the mapping $\eta_{\op{vf}}$, with the triangles on the percolation path of $(M,\sigma)$. Equivalently, the steps of $\hpi(w)$ are in one-to-one correspondence with the bicolor triangles of $\LR(M,\sigma)$.
		\item[(ii)] The steps $\ba_k$ (respectively, $\bb_k$) of $\hpi(w)$, or equivalently the spine $b$-steps (respectively, $a$-steps) of $w$ starting a cone excursion of $w$ of height $k$, are in one-to-one correspondence with the white (resp. black) inner faces of $\LR(M,\sigma)$ of degree $k+1$.
		\item[(iii)] By (i) and (ii), each step $x$ of $\hpi(w)$ of the form $\ba_k$ or $\bb_k$ corresponds to both a bicolor triangle $t$ and to a unicolor face $f$ of $\LR(M,\sigma)$. Furthermore, the site-percolated near-triangulation $(M',\sigma')$ of $(M,\sigma)$ formed of $t$ and all the vertices and edges of $(M,\sigma)$ inside and on the boundary of $f$ is equal to $\bPhi(w')$, where $w'$ is the cone excursion of $w$ corresponding to $x$.
	\end{compactitem}
\end{thm}

Theorem~\ref{thm:LR} is illustrated in Figure~\ref{fig:LR-theorem2}. It is proved in Section~\ref{sec:appendix-spine}.
As Theorem~\ref{thm:LR} shows, the percolated map $\LR(M,\sigma)$ is very naturally encoded by the walk $\hpi(w)$ obtained from by $w=\bPhi^{-1}(M,\sigma)$ by replacing each maximal cone excursion by a single step. 

\begin{remark}
	We can reformulate Property (iii) in Theorem~\ref{thm:LR} as follows. Let $f$ be the unicolor face of $\LR(M,\sigma)$ associated with a step $x$ on the form $\ba_k$ or $\bb_k$, let $w'$ be the cone excursion associated with $x$, and let $w''$ be obtained from $w'$ by removing the first letter (which is $a$ or $b$) and the last letter (which is $c$). If $w''$ is the empty word, then $f$ has boundary length 2 and there are no vertices of $M$ which are enclosed by $f$ (recall steps 2 and 3 in Definition~\ref{def:LR-M} for this convention.). 
	If $w''$ is not the empty word, then consider the submap $S$ formed of all the vertices and edges of $(M,\sigma)$ inside and on the boundary of $f$. We choose the root-edge of $S$ to be the edge $e_0$ of $S$ incident to $t$, and we change the color of one of the endpoints of $e_0$ so that this edge goes from a white to a black vertex. The resulting site-percolated near-triangulation $(M'',\sigma'')$ is in $\bmTP$, and Property (iii)
	implies that $(M'',\sigma'')=\bPhi(w'')$. See also Figure \ref{fig:proof-spine2}.
	
	We can now use this property to define a total order on the set $E(M)$ of edges of $M$. This ordering is defined recursively as follows. In the setting above, if $e_1,e_2$ are two edges on the percolation path, then $e_1$ is ordered before $e_2$ if and only if $e_1$ is crossed before $e_2$ by the percolation path (recall that our convention is that the percolation path goes from the root-edge to the top-edge). For any triangle $t$ on the percolation path as above, exactly two edges of $t$ are crossed by the percolation path. Let $e_1$ (resp. $e_2$) be the edges on $t$ which is crossed first (resp.\ last). In the case where $w''$ is empty, we order the third edge of $t$ after $e_1$ and before $e_2$. In the case where $w''$ is non-empty, all the edges of $(M'',\sigma'')$ are ordered after $e_1$ and before $e_2$, and the relative order of the edges of $(M'',\sigma'')$ is determined iteratively by considering $\LR(M'',\sigma'')$. It follows from Property (iii) that the ordering on $E(M)$ defined in this recursive manner equals the one determined by $\etae$.
\end{remark}


\noindent \textbf{Relation with other constructions in the literature.} Let us mention a connection between our bijection $\bPhi$, and the \emph{peeling process} which has been used extensively to study the percolation model on random maps. The peeling process is a way of constructing a percolation path of a random site-percolated triangulation ``one triangle at a time''. If one endows the set of finite site-percolated  triangulations with a Boltzmann distribution,
then the law of the peeling process takes a simple form. See Section~\ref{subsec:cycles} for the definition of the critical Boltzmann triangulation. Here we have shown that $\LR(w)$ encodes the law of the percolation path, and hence is closely related to the peeling process: in fact, the law of the peeling process is given by the law of the path $\hpi(w)$ (via the correspondences (i-iii) of Theorem~\ref{thm:LR}). As explained above, $\bPhi(w)$ can be obtained from $\LR(w)$ by filling each monochromatic face of $\LR(w)$ by the site-percolated triangulations encoded (via the bijection $\bPhi$) by the cone excursions corresponding to the negative steps of $\hpi(w)$. This shows that the bijection $\bPhi$ can be thought as a ``recursive'', ``space-filling'' version of the peeling process. We use this in our proof of Lemma~\ref{prop22}, and refer to that proof for further details.

Let us also mention a connection between our bijection $\Phi$, and Mullin's bijection \cite{mullin-maps} which underlies the influential ``inventory'' construction of Sheffield \cite{shef-burger}. Mullin's bijection can be interpreted as a bijection between the set $\mA$ of walks on $\NN^2$ starting and ending at $(0,0)$ and made of steps in $\{(1,0),(0,1),(-1,0),(0,-1)\}$ and the set $\mB$ of rooted maps with a marked spanning tree. Now consider the subset $\mH$ of $\mK$ made of walks in $w$ ending at $(0,0)$ and having no consecutive $c$-steps. 
For $w\in\mH$, $\hpi(w)$ is simply obtained from $w$ by replacing the subwords $ac$ and $bc$ by $\bb_1=(0,-1)$ and $\ba_1=(-1,0)$ respectively. 
In fact, $\hpi$ induces a bijection between $\mH$ and $\mA$.
Moreover for $w\in\mH$ the unicolor faces of $\LR(w)$ all have degree 2, and $\Phi(w)$ is obtained from $\LR(w)$ by collapsing each unicolor face into a single edge. Thus for $w$ in $\mH$, the map $\Phi(w)$ has one white cluster which is a tree, one black cluster which is a tree, and a path of bicolor triangles ``snaking'' between the two trees. As observed in Remark~\ref{rk:burger-link}, this is essentially the version of Mullin's bijection described in \cite{shef-burger}. In order to get an element of $\mB$ from $\Phi(w)$ (and recover Mullin's bijection) one needs to keep the white vertices, the unicolor white edges (which form the marked spanning tree), and also add the \emph{dual} of every unicolor black edge: every black edge is the diagonal of a square made of two bicolor triangles and needs to be replaced by the other diagonal of that square. Conversely, one can also view the bijection $\Phi$ as a special case of Mullin's bijection by replacing each step $c$ by a sequence of two steps $(-1,0),(0,-1)$ in a suitable order. See \cite[Section 3.2]{ghs-dist-exponent} for details. Using this relation greatly simplifies the analysis of distances in the so-called mated-CRT map for $\gamma=\sqrt{8/3}$ in \cite{ghs-dist-exponent}.

\subsection{Spine-looptrees decomposition in the infinite volume setting}\label{sec:LR-decomposition-infinite}
In this section, we give an analogue of Theorem~\ref{thm:LR} in the infinite volume setting. To be more precise, we will consider the percolation path of the ``past'' near-triangulation associated to $w\in\{a,b,c \}^\ZZ$.

Let $w\in\imK$. We write $w=w^-w^+$ as in Section~\ref{sec:bij-inf}. Let $(M,\si)=\Phi^\infty(w)\in\imTP$, and let $(M^-,\si^-)$ be the past site-percolated near-triangulation associated to $w^-$ (see Section~\ref{sec:bij-inf}). 
Recall that $M^-$ is a near-triangulation with an infinite simple boundary. Moreover the marked \emph{top-edge} separates the boundary into two semi-infinite paths: the \emph{left-boundary} with black vertices, and the \emph{right-boundary} with white vertices. Hence, there is a single infinite \emph{percolation path} of $(M^-,\si^-)$ ending at the top-edge. 
We now define $\LR(M^-,\sigma^-)$ as in Definition~\ref{def:LR-M}.

Recall that $(M^-,\sigma^-)$ is obtained by gluing together the site-percolated near-triangulations $(P_i,\al_i)\in \mT_P$ corresponding to the finite walks $w_i^-\in\mK$ in the decomposition $w^-=\ldots w_3^-w_2^-w_1^-$ at split-times. Hence the percolation path $\gamma$ of $(M^-,\si^-)$ is the concatenation of the percolation paths of the finite maps $(P_i,\al_i)$, and $\gamma$ goes through all the triangles incident to the outer edges of $M^-$. 
Also, $\LR(M^-,\sigma^-)$ is made by concatenating all the percolated maps $\LR(P_i,\al_i)$.

We now explain how to obtain $\LR(M^-,\sigma^-)$ directly from $w^-$.
\begin{definition}\label{def:pi-w-inf}
	Let $w^-=\dots w_{-2}w_{-1} \in \{a,b,c \}^{\ZZ^{< 0}}$ be a word whose letters are indexed by negative integers. We call \emph{spine step} of $w$ an $a$-step or a $b$-step which is not enclosed by any close-matching. 
	\begin{compactitem}
		\item We denote by $\pi(w^-)$ the subword of $w^-$ made of its spine steps.
		\item We denote by $\dots,T(-2),T(-1),T(0)\in\ZZ^{< 0}$ the times associated with the spine steps of $w^-$, so $\pi(w^-)=\dots w_{T(-2)}w_{T(-1)}w_{T(0)}$.
		\item We denote by $\hpi(w^-)$ the word on the infinite alphabet $\{a,b\}\cup\{\ba_k,\bb_k,~k\geq 1\}$ obtained from $\pi(w^-)$ by replacing each close-matched $a$-step (resp. $b$-step) $w_i$ by $\bb_{k}$ (resp. $\ba_{k}$), where $k$ is the height of the cone excursion starting at $w_i$.
		\item We denote by $\wh\pi_s(w^-)$ the word obtained from $\hpi(w^-)$ by replacing the letters in $\{a\}\cup\{\ba_k,~k\geq 1\}$ by $a$ and the letters in $\{b\}\cup\{\bb_k,~k\geq 1\}$ by $b$. We denote by $\wh\pi_\ell(w^-)$ (resp. $\wh\pi_r(w^-)$ ) the word obtained from $\hpi(w^-)$ by replacing the letters $a$ (resp. $b$) by 1, replacing the letters $\ba_k$ (resp. $\bb_k$) by $-k$, and deleting the letters in $\{b\}\cup\{\bb_k,~k\geq 1\}$ (resp. $\{a\}\cup\{\ba_k,~k\geq 1\}$).
	\end{compactitem}
\end{definition}

\begin{figure}
	\centering
	\includegraphics[scale=1]{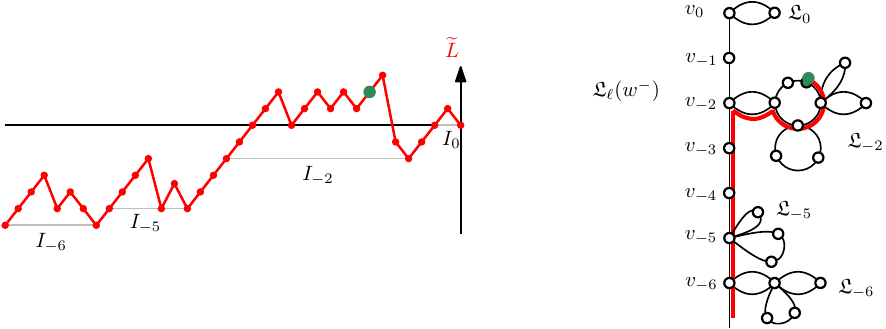}
	\caption{
		Left: The walk $\wt L$ associated to $\wh\pi_\ell(w^-)$. Right: the corresponding forested line $\mathfrak{L}_\ell(w^-)$. Each looptree of $\mathfrak{L}_\ell(w^-)$ corresponds to an excursion of $\wt L$ via the clockwise code. We have also drawn a canonical path $P_v$ (in red) from a vertex $v$ of $\mathfrak{L}_\ell(w^-)$ to $-\infty$. Note that the height of walk $\wt L$ describes how the length the canonical path $P_v$ varies as one trace the contour of $\mathfrak{L}_\ell(w^-)$.}\label{fig:LR-inf}
\end{figure}

A \emph{discrete forested line} is a collection of discrete looptrees, each associated to a different integer $n\in \ZZ$ called its \emph{root}. 
We now associate discrete forested lines to $\wh\pi_\ell(w^-)$ and $\wh\pi_r(w^-)$.
Let $\wt L$ be the lattice walk on $\ZZ$ ending at 0 having steps given by $\hpi_\ell(w^-)$. We denote by $(\wt L_m)_{m\in\ZZ^{\leq 0}}$ the successive values of $\wt L$ (so that, $\wt L_0=0$). 
An \emph{excursion interval} for $\wt L$ is an interval $I=\{m_1,m_1+1,\ldots,m_2\}$ of $\ZZ^{\leq 0}$ such that $m_1\leq m_2$ and
\eqbn
\wt L_{m_1-1}+1=\wt L_{m_1}=\wt L_{m_2}=\inf\{ \wt L_m~|~m_1\leq m\leq m_2\}<\inf\{ \wt L_m~|~m_2< m \leq 0\}.
\eqen
We call \emph{level} of the excursion $I$ the value $\wt L_{m_1}=\wt L_{m_2}$.
Because $w$ is normal, the walk $\wt L$ has one excursion interval $I_n$ of each level $n\in\ZZ^{\leq 0}$.
We now associate to $\wt L$ a discrete forested line made of the discrete looptrees $\mathfrak{L}_n=\ccwcode^{-1}(\wt L_{|I_n})$ with root $n$ for all $n$ in $\ZZ^{\leq 0}$, where $\wt L_{|I_n}$ is the subwalk of $\wt L$ corresponding to the interval $I_n$. We finally create a map $\mathfrak{L}_\ell(w^-)$ as follow: we consider a semi-infinite path $P$ with vertices $v_0,v_{-1},v_{-2},\ldots$ and we attach the root of the looptree $\mathfrak{L}_n$ to $v_n$, on the right of $P$ (with $P$ oriented toward $v_0$).  This is represented in Figure~\ref{fig:LR-inf}. 

Similarly, we can consider the walk $\wt R$ ending at 0 having steps given by $\hpi_r(w^-)$. This walk has one excursion $I'_n$ of each level $n\in\ZZ^{\leq 0}$. We then consider the discrete forested line made of the discrete looptrees $\mathfrak{L}_n'=\cwcode^{-1}(\wt R_{|I'_n})$ of root $n$ for all $n$ in $\ZZ^{\leq 0}$.  We finally create a map $\mathfrak{L}_r(w^-)$ by creating a semi-infinite path $P'$ with vertices $v_0',v_{-1}',v_{-2}',\ldots$ and attaching the root of $\mathfrak{L'}_n$ to $v'_n$, on the left of the path $P'$ (with $P'$ oriented toward $v_0'$).

We now define $\LR(w^-)$ similarly as in the finite case given by Definition~\ref{def:LR-w}; see Figure~\ref{fig:LR-theorem2} for the finite volume case. More precisely, we denote $P_s(w^-)$ the site-percolated triangulation corresponding to the walk $\wh\pi_s(w^-)$. We glue the vertex $v_0$ of $\mathfrak{L}_\ell(w^-)$ to the top-left vertex of $P_s(w^-)$ and we glue the vertex $v_0'$ of $\mathfrak{L}_r(w^-)$ to the top-right vertex of $P_s(w^-)$. We then glue the left edges of $P_s(w^-)$ to the edges on $\mathfrak{L}_\ell(w^-)$ in a planar manner (with all vertices of $\mathfrak{L}_\ell(w^-)$ colored white) and we glue the right edges of the $P_s(w^-)$ to the edges on $\mathfrak{L}_\ell(w^-)$ in a planar manner (with all vertices of $\mathfrak{L}_r(w^-)$ colored black). The resulting map $\LR(w^-)$ is a site-percolated map with half-plane topology (that is, the root-face has an infinite simple boundary).

The following result is the infinite volume analogue of Theorem~\ref{thm:LR}. 
\begin{theorem}\label{thm:spinelt-inf}
	Let $w^-\in \{a,b,c \}^{\ZZ^{<0}}$ be such that the steps are sampled uniformly and independently at random. Almost surely, $w^-$ has infinitely many split-times, hence the associated past percolated near-triangulation $(M^-,\sigma^-)$ is well-defined.
	Moreover the infinite site-percolated maps $\LR(M^-,\sigma^-)$ and $\LR(w^-)$ are the same almost surely. Furthermore, the correspondences (i-iii) of Theorem~\ref{thm:LR} still hold (with $(M,\si)$ replaced by $(M^-,\si^-)$ and $w$ replaced by $w^-$).
\end{theorem}

Theorem~\ref{thm:spinelt-inf} follows immediately from Theorem~\ref{thm:LR} by decomposing $w^-=\ldots w^-_{3}w^-_{2}w^-_{1}$ at split-times, and applying Theorem~\ref{thm:LR} to each finite walk $w^-_i$.

\begin{remark}\label{rk:pivotal-from-spine}
	By definition, a downward jump of $\wt L$ corresponds to the completion of a loop of a looptree of $\mathfrak{L}_\ell(w^-)$. Hence, downward jumps occur when one visits a cut-vertex of $\mathfrak{L}_\ell(w^-)$ for the second (or more) time when tracing the contour of $\mathfrak{L}_\ell(w^-)$; see Figure~\ref{fig:LR-inf}. By Theorem~\ref{thm:spinelt-inf}, the cut-vertices of $\mathfrak{L}_\ell(w^-)$ correspond to \emph{pivotal point} of the outer white cluster of $(M^-,\si^-)$ (see Section~\ref{subsec:pivot} for definitions of pivotal points). Thus, downward jumps of $\wt L$ and $\wt R$ correspond to pivotal points of outer clusters of $(M^-,\si^-)$. This property will be important to prove the convergence of the counting measure on pivotal points.
\end{remark}

\section{Discrete dictionary II: Exploration tree}\label{sec:DFS}
In this section, we describe how to associate to a percolated near-triangulation $(M,\sigma)\in\bmT_P$ a spanning tree $\tau^*$ of the dual map $M^*$. The spanning tree $\tau^*$, called \emph{exploration tree}, is related to an exploration of the map $M^*$ which ``tries to follow the percolation interfaces whenever possible''. The mapping $\Delta_M$ which associates the exploration tree $\tau^*$ to the percolation configuration $\sigma$ is a bijection between \emph{inner colorings} of $M$ and \emph{DFS-trees} of $M^*$ (see Theorem~\ref{thm:DFS-to-perco} for a precise statement).\footnote{The mapping $\Delta_M$ is not entirely new: it was used in several special cases in order to associate a tree to a percolation configuration (see for example \cite{shef-cle}). But its bijective nature seems to have been overlooked (for instance, there is no claim of injectivity nor a characterization of the image). On the other hand, a combinatorial argument showing that the DFS trees of $M^*$ are equinumerous to the inner colorings of $M$ was given in \cite{bernardi-dfs-bijection}, but no bijection was given between the two sets.} 
We will first establish important properties of the exploration tree $\tau^*$ and of the mapping $\Delta_M$ in Section~\ref{sec:exploration-tree-from-map}. Then, in Section~\ref{sec:exploration-tree-from-walk}, we explain how to directly obtain $\tau^*$ as a functional of the walk $w=\bPhi^{-1}(M,\sigma)$.

\subsection{Exploration tree associated to a site-percolation configuration}\label{sec:exploration-tree-from-map}
\noindent \textbf{Basic definitions about trees and planar duality.}
A \emph{tree} is a connected acyclic graph. 
A \emph{rooted tree} is a tree with a vertex distinguished as the \emph{root-vertex}. For a rooted tree, we adopt the usual vocabulary about \emph{parents}, \emph{children}, etc. For instance, an \emph{ancestor} of a vertex $v$ is any vertex on the path between the root-vertex and $v$. By convention, a vertex is considered both an ancestor and a descendant of itself. 
We call \emph{parent-edge} of a non-root vertex $v$, the edge joining $v$ to its parent. A \emph{spanning tree} of a connected graph $G$ is a subgraph containing every vertex of $G$, and which is a tree. We often identify a spanning tree with its edge set.

The \emph{dual} of a map $M$, denoted by $M^*$, is the map obtained by placing a vertex of $M^*$ in each face of $M$ and an edge of $M^*$ across each edge of $M$. Duality is represented in Figure~\ref{fig:duality}. Each vertex, edge, or face $x$ of $M$ naturally corresponds to a face, edge, or vertex of $M^*$ that we call \emph{dual of $x$} and usually denote by $x^*$. 
If $M$ is rooted with root-edge $e_0$, then $M^*$ is rooted with root-edge $e_0^*$ oriented from the right of $e_0$ to the left of $e_0$; see Figure~\ref{fig:duality}. With this convention, if the root-face of $M$ is $f_0$, then the root-vertex of $M^*$ is $f_0^*$.\\

If $T$ is a spanning tree of $M$, then the \emph{dual tree} is the spanning tree $T^*$ of $M^*$ made of the dual of the edges of $M$ not in $T$. This is represented in Figure~\ref{fig:duality}. Note that $T^*$ is indeed a spanning tree of $M^*$ (indeed, $T^*$ is acyclic because $T$ is connected, and $T^*$ is connected because $T$ is acyclic). 

\fig{width=.7\linewidth}{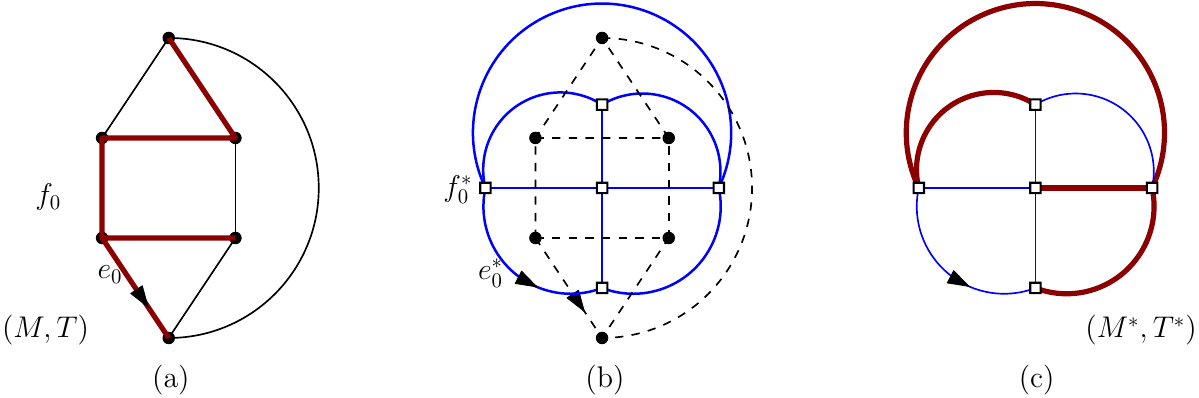}{(a) A rooted map $M$ and a spanning tree $T$ (bold edges). (b) The dual map $M^*$. (c) The dual map $M^*$ and the dual spanning tree $T^*$ (bold edges).}

We call \emph{cubic} a map such that every vertex has degree 3, and \emph{near-cubic} a map such that every non-root vertex has degree 3. Note that duality transforms triangulations (resp. near-triangulations) into cubic maps (resp. near-cubic maps), and transforms loopless maps into bridgeless maps. 
Let $M$ be a near-triangulation, and let $\sigma$ be a site-percolation configuration. We can obviously identify $\sigma$ with a coloring of the faces of $M^*$ in black and white. Note now that the \emph{percolation interface} of a percolated near-triangulation $(M,\sigma)\in\bmT_P$ as defined in Section~\ref{subsec:prelim} can be interpreted as the set of edges of $M^*$ which separate two faces of different colors, and are not incident to the root-vertex of $M^*$.\\

\noindent {\bf From percolation configurations to DFS trees: the mapping $\Delta_M$.}\\
We now establish a bijective correspondence between site-percolation configurations of a near-triangulation $M\in\mT$ and certain types of spanning trees of $M^*$ related to the \emph{depth-first search algorithm}. Let us first recall the definition of the \emph{depth-first search algorithm}. 
\begin{definition}\label{def:DFS-algo}
	Let $G$ be a connected graph and let $v_0$ be a vertex. A \emph{depth-first search} (or \emph{DFS} for short) of $G$ starting at $v_0$ is a visit of its vertices by a ``chip'' according to the following rule. At the beginning of the process, the chip is placed at $v_0$. The vertex $v_0$ is considered \emph{visited}, whereas all the other vertices are considered \emph{unvisited}. After that, we repeat the following step, where $u$ denotes the vertex where the chip is placed:
	\begin{compactitem}
		\item Case (a): there exists some edge between $u$ and an unvisited vertex. In this case, we choose such an edge $e=\{u,v\}$ and move the chip from $u$ to $v$. Then, we mark $v$ as \emph{visited}, call $u$ the \emph{parent} of $v$, and call $e$ the \emph{parent-edge} of $v$.
		\item Case (b): there is no edge between $u$ and an unvisited vertex. In this case, if $u\neq v_0$, then the chips moves to the parent of $u$, while if $u=v_0$, then the DFS stops.
	\end{compactitem}
	The \emph{tree associated to a DFS} of $G$ is the spanning tree of $G$ made of the set $T$ of all the parent-edges.
\end{definition}

We now recall a well-known characterization of the spanning trees which can be obtained by a DFS of $G$ starting at $v_0$ (see for example \cite[Section 23.3]{Cormen:introduction-algorithms}). For a spanning tree $T$ of $G$ rooted at $v_0$, we say that two vertices are \emph{$T$-comparable} if one is the ancestor of the other in the tree $T$.

\begin{claim}[Folklore] \label{claim:characterization-DFStree}
	A spanning tree $T$ of $G$ can be obtained by a DFS of $G$ starting at $v_0$ if and only if any two adjacent vertices of $G$ are \emph{$T$-comparable}. We call such a tree a \emph{$v_0$-DFS tree} of $G$.
\end{claim}

Next, we define a set of DFS trees and a set of percolation configurations, and the bijection between them.

\begin{definition}\label{def:Perc-and-DFS-sets}
	Let $M$ be a near-triangulation in $\mT$, let $M^*$ be the dual map, and let $v_0$ be the root-vertex of $M^*$.
	We call \emph{inner-coloring} of $M$ a coloring of the inner vertices of $M$ in black or white.
	We denote by $\Perc_M$ the set of inner-colorings of $M$.
	We denote by $\DFS_{M^*}$ the set of $v_0$-DFS trees $T$ of $M^*$ such that the root-edge $e_0^*$ of $M^*$ is in $T$. 
\end{definition}
\fig{width=\linewidth}{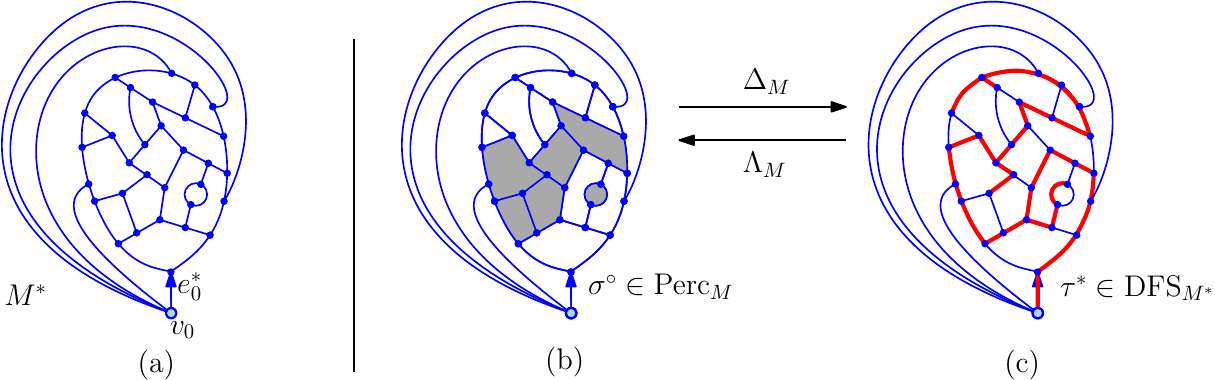}{(a) A near-cubic map $M^*$ with root-vertex $v_0$. The map $M^*$ is the dual of a map $M\in\mT$. (b) An inner-coloring $\gao\in \Perc_M$ represented as a coloring of the faces of $M^*$ not incident to the root-vertex $v_0$. (c) The $v_0$-DFS tree $\tau^*=\Delta_M(\gao)$.}

We point out that, for all $M\in\mT$, the trees in $\DFS_{M^*}$ contain no edge incident to $v_0$ apart from the root-edge $e_0^*$ (otherwise $v_0$ would have several children and the characterization of $v_0$-DFS trees given by Claim~\ref{claim:characterization-DFStree} would be violated somewhere). Observe that this implies that during a DFS of $M^*$ starting at $v_0$, the chip does not visit $v_0$ except at the first and last step. 
We now introduce some additional vocabulary. During a DFS of $M^*$, for any vertex $v\neq v_0$ which has already been visited, we consider the incident edges $e_1,e_2,e_3$ in clockwise order around $v$, with $e_1$ being the parent-edge of $v$. We call $e_2$ the \emph{left forward edge} of $v$, and $e_3$ the \emph{right forward edge} of $v$. We call the face $f$ containing the corner between $e_2$ and $e_3$, the \emph{forward face} of $v$.
We can now define the bijection $\Delta_M$ between $\Perc_M$ and $\DFS_{M^*}$.

\begin{definition}\label{def:Delta}
	Let $M$, $M^*$ and $v_0$ be as in Definition~\ref{def:Perc-and-DFS-sets}.
	Given an inner-coloring $\gao\in \Perc_M$, we consider the corresponding coloring of the faces of $M^*$ \emph{with the convention that the faces of $M^*$ dual to the outer vertices of $M$ are colored white}. We define $\Delta_M(\gao)$ as the spanning tree of $M^*$ obtained by the DFS of $M^*$ defined as follows:
	\begin{compactitem} 
		\item[(i)] The chip starts at the root-vertex $v_0$ of $M^*$, and first moves along the root-edge $e_0^*$ of $M^*$.
		\item[(ii)] Subsequently, each time the algorithm is in Case (a) of Definition~\ref{def:DFS-algo} and several edges are possible to move the chip, the choice is made according to the following rule. Let $u$ be the current position of the chip, and let $f$ be its forward face. If $f$ is black, then the chip moves along the forward left edge. If $f$ is white, then the chip moves along the forward right edge. This rule is represented in Figure~\ref{fig:DFS-rule-2}(a). 
	\end{compactitem}
\end{definition}

The mapping $\Delta_M$ is illustrated in Figure~\ref{fig:near-cubic-map}.

\fig{width=\linewidth}{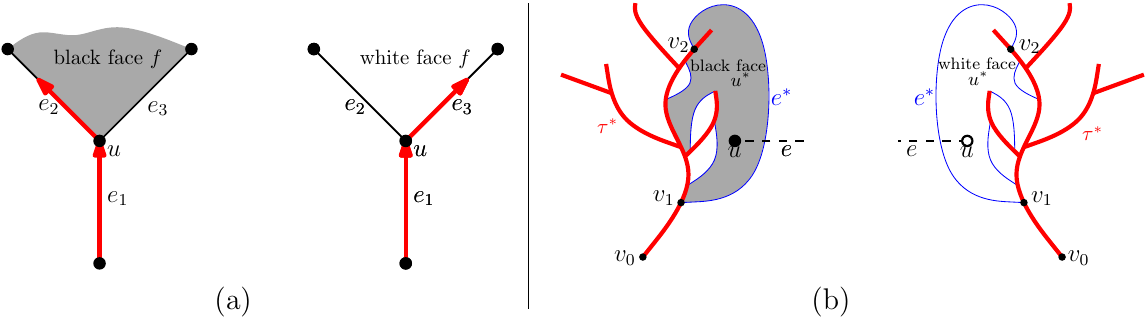}{(a) Rule of the mapping $\Delta_M$ for choosing the next edge $e=\{u,v\}$ of the DFS (Case (a) of Definition~\ref{def:DFS-algo}): the chip ``turns left'' when encountering a black face of $M^*$, and ``turns right'' when encountering a white face of $M^*$. Here $f$ is the forward face and $e_2,e_3$ are the left and right forward edges respectively.
	(b) Rule of the mapping $\Lambda_M$ for choosing the color of the inner vertex $u$ of $M$. The edge $e^*$ of $M^*$ is the unique edge of $M^*$ incident to the face $u^*$ such that the cycle inside $\tau^*\cup \{e^*\}$ separate $u^*$ from the root-face. Let $v_1,v_2$ be the endpoints of $e^*$ with $v_1$ the ancestor of $v_2$. 
	If $e^*$ is on the left (resp. right) of $\tau^*$ at $v_1$, then $u$ is colored white (resp. black).}

We now state a technical fact highlighting that certain choices made during the DFS algorithm do not affect the resulting DFS tree (even if they affect the order in which vertices are visited).

\begin{claim} \label{claim:color-doesnt-matters} 
	Let $M$, $M^*$ and $v_0$ be as in Definition~\ref{def:Perc-and-DFS-sets} and let $\gao\in \Perc_M$. 
	Let $\tau^*$ be a spanning tree of $M^*$ obtained from a DFS of $M^*$ satisfying the rules of Definition~\ref{def:Delta} except that Rule (ii) is loosened as follows:
	\begin{compactitem} 
		\item[({ii}$\,^{\star}$)] Subsequently, each time the algorithm is in Case (a) of Definition~\ref{def:DFS-algo} and several edges are possible to move the chip, the choice is made according to the following rule. Let $u$ be the current position of the chip, and let $f$ be its forward face. If no vertex $v\neq u$ incident to $f$ has been visited and $f$ is black (resp. white), then the chip moves along the forward left (resp. right) edge. If a vertex $v\neq u$ incident to $f$ has already been visited, then  we may choose the left forward edge or the right forward edge to move the chip, independently of the color of $f$.
	\end{compactitem}
	In this case, the resulting \referee{DFS tree} $\tau^*$ is equal to $\Delta_M(\gao)$.
\end{claim}


Claim~\ref{claim:color-doesnt-matters} will be proved in Section~\ref{sec:appendix-DFS-Perco}. We will define later (Definition~\ref{def:space-filling}) a particular DFS algorithm of $M^*$, called \emph{space-filling-DFS}, which satisfies the hypothesis of Claim~\ref{claim:color-doesnt-matters} (hence give the same DFS tree $\Delta_M(\gao)$), and is such that the order in which the vertices of $M^*$ are visited corresponds to the order in which they are created during the bijection $\bPhi$.

We also mention the following consequence of Claim~\ref{claim:color-doesnt-matters}.

\begin{remark}\label{rk:outer-coloring-does-not-matter} In Definition~\ref{def:Delta} of $\Delta_M(\gao)$ we adopted the convention that the outer vertices of $M$ were all colored white. But, by Claim~\ref{claim:color-doesnt-matters}, any convention for the color of the outer vertices of $M$ would result in the same spanning tree $\Delta_M(\gao)$. This is because the faces of $M^*$ dual to the outer vertices of $M$ are incident to the vertex $v_0$ which is the first vertex visited by the DFS.
\end{remark}


\noindent {\bf From DFS trees to percolation configurations: the mapping $\Lambda_M$.}
\begin{definition}\label{def:Gamma}
	Let $M$, $M^*$ and $v_0$ be as in Definition~\ref{def:Perc-and-DFS-sets}.
	Consider a $v_0$-DFS tree $\tau^*\in \DFS_{M^*}$, and the dual spanning tree $\tau$ of $M$. 
	We define the inner-coloring $\Lambda_M(\tau^*)$ of $M$ by deciding the color of any inner vertex $u$ of $M$ according to the following rule:
	
	Let $e$ be the parent-edge of $u$ in the spanning tree $\tau$ of $M$ dual to the spanning tree $\tau^*$ of $M^*$. By Claim~\ref{claim:characterization-DFStree}, the edge $e^*\in M^*\setminus \tau^*$ joins a vertex $v_1$ to one of its descendants $v_2$ in $\tau^*$. Note that $v_0,v_1,v_2$ are distinct vertices. Then we color $u$ white, if the edge $e^*$ is on the left of the path of $\tau^*$ from $v_0$ to $v_2$ at the vertex $v_1$, and we color $u$ black otherwise.
\end{definition}

In summary, the inner-coloring $\Lambda_M(\tau^*)$ considered as a coloring of the faces of $M^*$ is obtained as follows: the faces of $M^*$ ``on the left of $\tau^*$'' are colored white, while the faces of $M^*$``on the right of $\tau^*$'' are colored black. This rule is represented in Figure~\ref{fig:DFS-rule-2}(b). We now state the main result of this subsection.

\begin{thm}\label{thm:DFS-to-perco}
	Let $M$, $M^*$ and $v_0$ be as in Definition~\ref{def:Perc-and-DFS-sets}.
	The mapping $\Delta_M$ is a bijection from $\Perc_M$ to $\DFS_{M^*}$, and $\Lambda_M$ is the inverse bijection. 
	
	Moreover, for any inner-coloring $\gao\in \Perc_M$ and for any site-percolation configuration $\sigma$ of $M$ which extends $\gao$ (by attributing a color to the outer vertices) and which satisfies the root-interface condition, the tree $\tau^*=\Delta_M(\gao)$ satisfies the following properties.
	\begin{compactitem}
		\item[(i)] The percolation path of $(M,\sigma)$ is contained in $\tau^*$. 
		\item[(ii)] For any percolation cycle $C$ of $(M,\sigma)$, every edge of $C$ except one is in $\tau^*$.
		\item[(iii)] Consider the coloring of the faces of $M^*$ corresponding to the configuration $\sigma$. Then any edge of $\tau^*$ separating a black face and a white face of $M^*$ has the white face on its left when oriented from parent to child.
	\end{compactitem}
\end{thm}

Theorem~\ref{thm:DFS-to-perco} is illustrated in Figure~\ref{fig:near-cubic-map}. It will be proved in Section~\ref{sec:appendix-DFS-Perco}. We now state an immediate corollary. 
\begin{definition}\label{def:Perc-and-DFS-sets-e} 
	Let $M$, and $M^*$ and $v_0$ be as in Definition~\ref{def:Perc-and-DFS-sets}. 
	\begin{compactitem}
		\item For a non-root outer edge $e$ of $M$, we denote by $\Perc_M^e$ the set of site-percolation configurations of $M$ satisfying the root-interface condition, such that $e$ is bicolor, and such that the percolation path visits every inner triangle of $M$ incident to an outer edge.
		\item For a non-root outer edge $e^*$ of $M^*$ incident to $v_0$, we denote by $\DFS_{M^*}^{e^*}$ the set of trees $\tau^*\in \DFS_{M^*}$ such that the non-root vertex $v_1$ incident to $e^*$ is the descendant in $\tau^*$ of every vertex of $M^*$ adjacent to $v_0$. 
	\end{compactitem}
	For $\tau^*\in\DFS_{M^*}^{e^*}$, we denote by $\Lambda_M^{e^*}(\tau^*)$ the unique site-percolation configuration of $M$ extending the inner-coloring $\gao=\Lambda_M(\tau^*)$ and satisfying the root-interface condition with $e$ bicolor.
\end{definition}
Note that a pair $(M,\sigma)$ is in $\mT_P$ if and only if $M\in\mT$ and $\sigma\in \Perc_M^e$ for some non-root outer edge $e$ of $M$.

\begin{corollary}\label{cor:DFS-site-perco}
	Let $M$, and $M^*$ and $v_0$ be as in Definition~\ref{def:Perc-and-DFS-sets}.
	Let $e$ be a non-root outer edge of $M$, and let $e^*$ be the dual edge. 
	The mapping $\Lambda_M^{e^*}$ is a bijection between $\DFS_{M^*}^{e^*}$ and $\Perc_M^e$. 
\end{corollary}

\begin{proof}
	Let $\tau^*\in\DFS_{M^*}^{e^*}$. Let $\sigma=\Lambda_M^{e^*}(\tau^*)$. By Theorem~\ref{thm:DFS-to-perco}, $\tau^*=\Delta_M(\sigma^0)$, where $\sigma^0$ is the inner-coloring of $M$ induced by $\sigma$.
	By Property (i) of Theorem~\ref{thm:DFS-to-perco}, the percolation path $P$ of $(M,\sigma)$ is the path of the tree $\tau^*$ going from the root-vertex $v_0$ to the other endpoint $v_1$ of $e^*$. Thus $P$ visits every inner triangle of $M$ incident to an outer edge if and only if $v_1$ is the descendant in $\tau^*$ of every vertex of $M^*$ adjacent to $v_0$. Hence Corollary~\ref{cor:DFS-site-perco} follows immediately from Theorem~\ref{thm:DFS-to-perco}. 
\end{proof}

\subsection{Exploration tree as a function of the Kreweras walk}\label{sec:exploration-tree-from-walk}
For $(M,\sigma)\in \bmT_P$, we denote by $\dfs(M,\sigma):=\Delta_{M}(\gao)$ the DFS tree of $M^*$ corresponding to the inner-coloring $\gao$ of $M$ induced by $\sigma$. In this subsection, we will describe how to obtain the tree $\dfs(M,\sigma)$ directly from the walk $w=\bPhi^{-1}(M,\sigma)$.

We first recall the definition of the \emph{height-code} of a tree. A \emph{rooted plane tree} is a rooted map whose underlying graph is a tree. The \emph{contour} of a rooted plane tree $T$ is the walk around $T$ in clockwise direction (that is, keeping $T$ on the right of the walker) starting at the root-vertex $v_0$, just before the root-edge. The contour of a tree is represented in Figure~\ref{fig:height-code}(a). The \emph{prefix order} of the vertices of $T$ is the order in which vertices appear along the contour of $T$ (the root-vertex is first in this order). 
The \emph{height} of a vertex $v$ in $T$ is the number of edges on the path from $v$ to $v_0$ (so that the height of $v_0$ is 0).
The \emph{height code} of $T$ is the tuple $\hcode(T):=(h_0,h_1,h_2,\ldots,h_n)$, where $h_i$ is the height of the vertices $v_0,v_1,\ldots,v_n$ in prefix order; see Figure~\ref{fig:height-code}(b).
It is well-known that the mapping $\hcode$ establishes a bijection between rooted plane trees with $n+1$ vertices and $(n+1)$-tuples of integers $(h_0,h_1,\ldots,h_n)$ such that $h_0=0$, and for all $i\in [n]$, $1\leq h_{i}\leq h_{i-1}+1$. 

\fig{width=.6\linewidth}{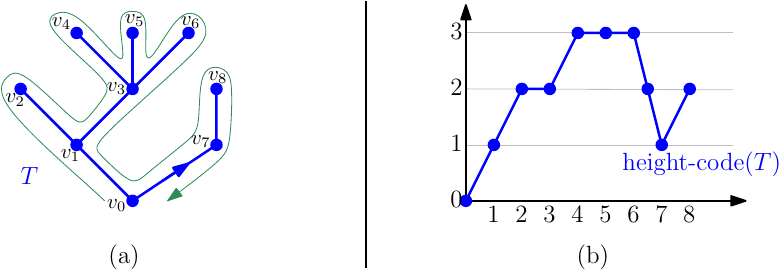}{(a) A rooted plane tree $T$ with root-vertex $v_0$, its contour, and its vertices $v_0,v_1,\ldots,v_8$ labeled in prefix order. (b) The height-code of $T$.} 

\begin{definition} \label{def:DFS-w}
	Let $w\in \bmK$. Let $n$ be the total number of $a$-steps and $b$-steps of $w$.
	For $i\in[n]$, let $w^{(i)}$ be the prefix of $w$ ending just after the $i^{\textrm{th}}$ step of $w$ which is either an $a$-step or a $b$-step.
	
	We denote by $\dfs(w)$ the rooted plane tree having height code $(h_0,\ldots,h_{n})$, where $h_0=0$ and for all $i\in[n]$, $h_i$ is the number of spine steps of the walk $w^{(i)}$ (see Definition~\ref{def:pi-w} of spine steps).
\end{definition}

The tree $\dfs(w)$ is represented in Figure~\ref{fig:height-code-DFStree}.

\fig{width=\linewidth}{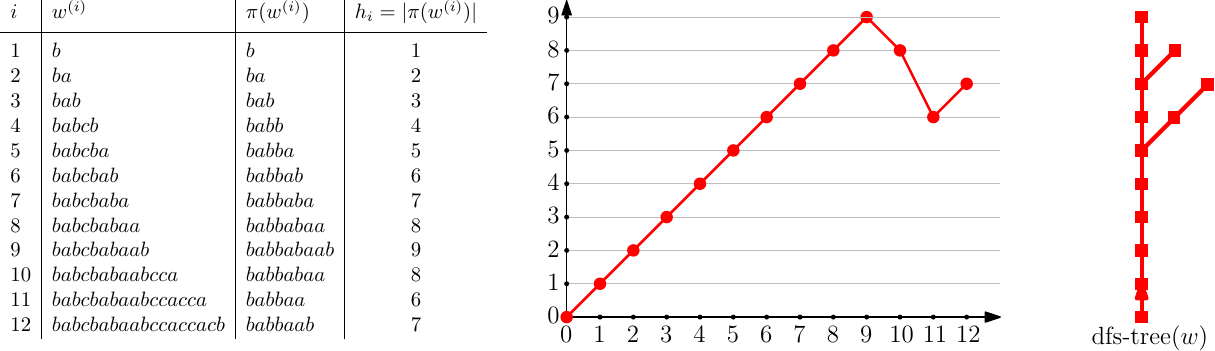}{The tree $\dfs(w)$ for the word $w=babcbabaabccaccacb\in\bmK$.}

\begin{thm}\label{thm:DFS}
	Let $w$ be a walk in $\bmK$, let $(M,\sigma)=\bPhi(w)$ be the associated site-percolated triangulation, and let $\gao$ be the inner-coloring of $M$ induced by $\sigma$. Then, the rooted plane tree $\dfs(w)$ (obtained by Definition~\ref{def:DFS-w}) and the spanning tree $\dfs(M,\sigma):=\Delta_M(\gao)$ (obtained by Definition~\ref{def:Delta}) have the same underlying rooted tree (although their planar embeddings may differ).
	Moreover, we have the following correspondences.
	\begin{compactitem}
		\item[(i)] The mapping $\etae$ gives a one-to-one correspondence, via duality, between the $a$-steps and $b$-steps of $w$ and the edges of $M^*$ in the tree $\dfs(M,\sigma)$. Moreover, the spine steps of $w$ correspond to the edges of $M^*$ on the percolation path. 
		\item[(ii)] The mapping $\etae$ also gives a one-to-one correspondence, via duality, between the $c$-steps of $w$ and the edges of $M^*$ which are not in the tree $\dfs(M,\sigma)$ and are not dual to active outer edges of $M$. 
	\end{compactitem}
\end{thm}

\fig{width=.25\linewidth}{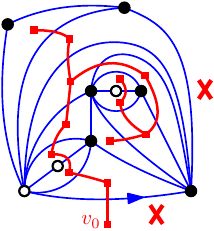}{The percolated near-triangulation $(M,\sigma)\in\bmT_P$ corresponding to the walk $w=babcbabaabccaccacb\in\bmK$. The DFS tree $\dfs(M,\sigma):=\Delta_M(\gao)$ is indicated in bold red lines (where $\gao$ is the inner-coloring of $M$ induced by $\sigma$). The tree $\dfs(w)$ had been computed in Figure~\ref{fig:height-code-DFStree}, and one can check it has the same underlying rooted tree as $\dfs(M,\sigma)$.}

Theorem~\ref{thm:DFS} is illustrated in Figure~\ref{fig:DFS-thm}. It is proved in Section~\ref{sec:appendix-DFS-from-walk}. We mention that Theorem~\ref{thm:DFS} is actually more transparent in terms of the bijection $\Om$ from \cite{bernardi-dfs-bijection} which is recalled in Section~\ref{sec:appendix-bijection}.

In the rest of this subsection we describe, for any pair $(M,\sigma)\in \bmT_P$, a DFS of $M^*$ for which the order of visit of the vertices corresponds to the order of creation of these vertices during the bijection~$\bPhi$.

\begin{definition}\label{def:space-filling} Let $(M,\sigma)\in\bmTP$.
	We call \emph{space-filling exploration of $M^*$} the DFS of $M^*$ defined as follows. 
	\begin{compactitem} 
		\item[(i)] The chip starts at the root-vertex $v_0$ of $M^*$, and first moves along the root-edge $e_0^*$ of $M^*$.
		\item[(ii')]
		Subsequently, each time the algorithm is in Case (a) of Definition~\ref{def:DFS-algo} and several edges are possible to move the chip, the choice is made according to the following rule. Let $u$ be the current position of the chip, and let $f$ be its forward face.
		\begin{compactitem} 
			\item[(1)] If none of the ancestors of $u$ is incident to $f$, and $f$ is black (resp. white), then the chip moves along the left (resp. right) forward edge of $u$.
			\item[(2)] If some ancestor of $u$ is incident to $f$, then we consider the \emph{last} ancestor $v$ of $u$ incident to $f$ (by \emph{last} we mean the one closet to $u$ in the DFS tree). If $v=v_0$ and $f$ is black (resp. white), then the chip moves along the right (resp. left) forward edge of $u$.
			If $v\neq v_0$, we consider the edge $e$ incident to $v$ on the path from $v$ to $u$ in the DFS tree. If $e$ is the left (resp. right) forward edge of $v$, then the chip moves along the right (resp. left) forward edge of $u$.
		\end{compactitem}
	\end{compactitem}
\end{definition}

Observe that by Claim~\ref{claim:color-doesnt-matters}, the DFS tree defined by the space-filling exploration of $M^*$ is $\dfs(M,\sigma)$. The space-filling exploration defines an order of visits of the vertices of $M^*$. We now define an order of \emph{treatment} of the edges of $M^*$. We call \emph{in-edge} of $M^*$ the edges of $M^*$ which are dual to in-edges of $M$, and \emph{active edges} the edges of $M^*$ which are dual to active outer edges of $M$. An edge of the DFS tree $\tau^*$ is said to be \emph{treated} during the space filling exploration of $M^*$, when the chip first moves along it (for the first time) during the space-filling exploration of $M^*$. An in-edge $e=\{u,v\}$ of $M^*$ not in $\tau^*$, with $u$ the ancestor of $v$, is said to be \emph{treated} the first time that the chip is at $v$ and either all the other edges incident to $v$ have been treated, or there is a path of yet untreated edges (not using $e$) between $v$ and an ancestor of $u$. It is easy to see that all the in-edges of $M^*$ will be treated during the space-filling exploration of $M^*$. By convention, the active edges of $M^*$ are considered to be treated after the space-filling exploration is complete (that is, after all the in-edges have been treated). 

\begin{proposition}\label{prop32}
	Let 	$(M,\sigma)\in\bmTP$.
	\begin{compactitem}
		\item[(i)] Let $t_1,\dots,t_k$ be the inner triangles of $M$, and let $t_1^*,\ldots,t_k^*$ be the corresponding vertices of $M^*$ visited in this order during the space-filling exploration. Then, $\etavf^{-1}(t_1),\dots,\etavf^{-1}(t_k)$ is increasing. In other words, the triangles $t_1,\dots,t_k$ are created in this order during the bijection $\bPhi$. 
		\item[(ii)] Let $e_1,\dots,e_n$ be the in-edges of $M$, and let $e_1^*,\ldots,e_n^*$ be the corresponding in-edges of $M^*$ treated in this order during the space-filling exploration. Then, for all $i\in [n]$, $\etae(i)=e_i$.
	\end{compactitem}
\end{proposition}

Proposition~\ref{prop32} will be proved in Section~\ref{sec:appendix-DFS-from-walk}.

\subsection{Exploration tree in the infinite setting}\label{sec:DFS-infinite}
In this subsection we define the infinite volume analogues of the bijections $\Delta_M$ and $\Lambda_M$, and obtain an analogue of~\ref{thm:DFS} for describing the exploration tree as a function of the Kreweras walk. All the proofs are given in Section~\ref{app:inf}. 

Let $G$ be an infinite graph. Given a one-ended spanning tree $T$ of $G$, we say that a vertex $u$ is an \emph{ancestor} of $v$ if $u$ is on the path of $T$ from $v$ to $\infty$ \referee{(in other words, we consider the ``root'' of the tree to be at infinity)}. We say that $u,v$ are \emph{$T$-comparable} if one is an ancestor of the other. We call \emph{DFS tree} of $G$ a one-ended spanning tree $T$ of $G$ such that every edge of $G$ joins $T$-comparable vertices. For $M\in\imT$, we denote by $\DFS_{M^*}$ the set of (one-ended) DFS trees of $M^*$. 

We recall that the dual tree of a one-ended spanning tree on an infinite planar graph is also a one-ended spanning tree (see e.g.\ \cite[Lemma 7.1]{bp93}). Hence for $\tau^*$ in $\DFS_{M^*}$ and $u$ a vertex of $M$, the notion of \emph{parent-edge} of $u$ in the dual spanning tree $\tau$ of $M$ is well defined.
Given $\tau^*$ in $\DFS_{M^*}$, we define the percolation configuration $\Lambda_M(\tau^*)$ of $M$ as in Definition~\ref{def:Gamma}. In short, vertices of $M$ are white if the corresponding face of $M^*$ is ``on the left'' of the tree $\tau^*$, and are black otherwise.
\begin{theorem}\label{thm:dfs-inf}
	Let $(M,\si)\in\imTP$ be chosen according to the percolated UIPT distribution. Then, the following properties hold almost surely.
	\begin{compactitem}
		\item[(i)] There exists a unique spanning tree $\tau^*\in \DFS_{M^*}$ such that $\Lambda_M( \tau^*)=\sigma$. In this case, we write $\tau^*=\Delta_M(\sigma)$. 
		\item[(ii)] The tree $\tau^*$ satisfies Properties $(ii)$ and $(iii)$ of Theorem~\ref{thm:DFS-to-perco}.
		\item[(iii)] The percolation path of  $(M^-,\sigma^-)$ is contained in $\tau^*$: it is the path of $\tau^*$ from $\infty$ to the vertex $v_0$ defined as follows.
		Let $e_0^*$ be the edge of $M^*$ dual to the root-edge of $M$, and let $u,v$ be its endpoints with $u$ the ancestor of $v$ in $\tau^*$. Then $v_0=u$ if $e_0^*$ is in $\tau^*$, and $v_0=v$ otherwise. 
	\end{compactitem}
\end{theorem}

The proof of Theorem~\ref{thm:dfs-inf} will provide an explicit way to obtain $\Delta_M(\sigma)$ from $(M,\sigma)$ by a certain DFS procedure (see the proof of Lemma~\ref{prop9}). The recursive step of this procedure is as described in Definition~\ref{def:Delta}.

Let $w\in\imK$, and let $(M,\si)=\Phi^\infty(w)$. As we now explain, the infinite DFS tree $\Delta_M(\sigma)$ (assuming it is uniquely determined) can be obtained directly from the walk $w$ similarly as in the finite volume case. 
We first need to extend the notion of \emph{height-code} to infinite plane trees.
An \emph{infinite rooted plane tree} is an infinite one-ended planar map with a single face, together with a \emph{marked vertex}. Note that the graph underlying an infinite rooted plane tree is a one-ended tree.
Let $T$ be an infinite rooted plane tree, and let $v_0$ be its marked vertex. 
The \emph{relative height} of a vertex $v$ of $T$ is the ``difference of length'' between the path $P$ from $v$ to $\infty$ in $T$ and the path $P'$ from $v_0$ to $\infty$ in $T$ (this difference makes sense because $P,P'$ coincide except on a finite portion). The \emph{height-code} of $T$ is the bi-infinite sequence $\hcode(T)=(h_i)_{i\in\ZZ}$ of the relative heights of the vertices seen in clockwise order around $T$, with time 0 corresponding to being at the first corner of $v_0$.
We adopt the convention $h_i=-\infty$ if there are less than $i$ vertices appearing after $v_0$ around the tree.
It is not hard to see that $\hcode$ is a bijection between the set of infinite rooted plane trees and the set of sequences $(h_i)_{i\in\ZZ}$ with $h_i\in \ZZ\cup\{-\infty\}$ such that $h_0=0$, $\liminf_{i\to -\infty}(h_i)=\liminf_{i\to +\infty}(h_i)=-\infty$, and for all $i\in \ZZ$, $h_{i+1}\leq h_i+1$.

Let $w\in\imK$. Let $w^{(0)}$ be the prefix of $w$ ending right after the last $a$-step or $b$-step $x_0$ of $w^-$. For all $i\in\ZZ$, let $w^{(i)}$ be the prefix of $w$ ending just after the $i$th step of $w$ which is either an $a$-step or a $b$-step, where we count these steps relative to step $x_0$. 
Let $(T^i(k))_{k\in\ZZ^{\leq 0}}$ be the sequence of indices of the spine steps of $w^{(i)}$, and let $\pi(w^{(i)})=\ldots w_{T^i(-2)} w_{T^i(-1)}w_{T^i(0)}$.
Since $w$ is normal there exists a split time $t$ of $w^{(\max(0,i))}$ preceding $\min(0,i)$, and it is not hard to see that the spine steps of $w^{(0)}$ and $w^{(i)}$ before the split time $t$ are equal. Hence there are numbers $k_i,k_i'\in\ZZ^{\leq 0}$ such that 
\begin{equation}
	\label{eq:branching-times}
	\begin{split}
		&T^i(k_i-j)=T^0(k_i'-j),\quad\forall j\in\N,\\
		&\{T^i(k_i+1),\dots,T^i(0) \}\cap \{T^0(k_i'+1),\dots,T^0(0) \}=\emptyset.
	\end{split}
\end{equation}
We define $h_i:=k_i-k_i'$, and call $(h_i)_{i\in\ZZ}$ the \emph{spine length sequence} of $w$.

\begin{definition}\label{def:DFS-w-inf}
	For $w\in\imK$, we let $\dfs(w)$ be the infinite rooted plane tree having height-code equal to the spine length sequence $(h_i)_{i\in\ZZ}$ of $w$. 
\end{definition}
We now state the infinite volume analogue of Theorem~\ref{thm:DFS}.
\begin{theorem}\label{thm:dfs-inf2}
	Let $w\in \imK$, let $(M,\sigma)=\Phi^\infty(w)$ be the associated percolated triangulation, and assume that the spanning tree $\Delta_M(\sigma)$ (described in Theorem~\ref{thm:dfs-inf}) is uniquely defined.
	Then, the infinite rooted plane tree $\dfs(w)$ (obtained by Definition~\ref{def:DFS-w-inf}) and the spanning tree 
	$\Delta_M(\sigma)$ 
	have the same underlying graphs (the underlying trees are equal although their planar embeddings may differ). 
	Moreover, the marked vertex of $\dfs(w)$ corresponds to the vertex of $\Delta_M(\sigma)$ denoted by $v_0$ in Theorem \ref{thm:dfs-inf}(iii) (that is, $v_0$ is the end of the percolation path of $(M^-,\sigma^-)$).
\end{theorem}

In the rest of this section we state the infinite volume analogue of Proposition~\ref{prop32}. Given $w\in\mK^\infty$ and $(M,\sigma)=\Phi^{\infty}(w)$, we want to describe a DFS of $M^*$ whose order of visit of vertices coincides with the order of creation of these vertices during $\Phi^{\infty}$. We first need to define DFS of infinite graphs.

\begin{definition}\label{def:DFS-algo-inf}
	Let $G$ be an infinite graph. A \emph{depth-first search} of $G$ is a surjective map $s:\Z\to V(G)$ which describes a walk of a ``chip'' on $G$ (so that $s(k),s(k+1)$ are adjacent vertices for all $k\in\ZZ$) which must satisfy the following conditions. First, for any $v\in V(G)$, $m_v=\min(k,~s(k)=v)$ must be finite, and we call $s(m_v-1)$ the \emph{parent of $v$}. Second, for any $k\in\Z$ the value $s(k+1)$ must be compatible with $s_{|\Z^{\leq k}}$ in the following sense. 
	Call $s(k)$ the \emph{position of the chip} at time $k$, and call $\{s(i),~i\leq k\}$ the set of \emph{visited vertices} at time $k$. Then $s(k+1)$ must satisfy the following:
	\begin{compactitem}
		\item Case (a): if there is a neighbor of $s(k)$ which is unvisited at time $k$, then $s(k+1)$ is one of these unvisited vertices.
		\item Case (b): otherwise $s(k+1)$ is the parent of $s(k)$.
	\end{compactitem}
	The \emph{tree associated to a DFS} of $G$ is the spanning tree having one edge between each vertex of $G$ and its parent.
\end{definition}

\referee{Note that not every infinite graph admits a depth-first search. A necessary condition is the existence of a one-ended spanning tree. In fact, one can easily check the following analogue of Claim~\ref{claim:characterization-DFStree}.}

\begin{claim}\label{claim:flokloreDFS-infinite}
	The tree associated to a DFS of an infinite graph $G$ is indeed a spanning tree of $G$. Moreover Claim~\ref{claim:characterization-DFStree} still holds: a one-ended spanning tree $T$ of an infinite graph $G$ is associated to a DFS if and only if every edge of $G$ join $T$-comparable vertices.
\end{claim}


The following analogue of Proposition~\ref{prop32} holds in the infinite volume setting.
\begin{prop}\label{prop:dfs-ordering-inf}
	Let $(M,\sigma)\in\imTP$ be chosen according to the percolated UIPT distribution. Recall that $\tau^*=\Delta_{M}(\sigma)\in \DFS_{M^*}$ is well-defined almost surely. 
	\begin{compactitem}
		\item[(i)] The tree $\tau^*$ is associated to a DFS of $M^*$ such that every time we are in Case (a) of Definition~\ref{def:DFS-algo-inf}, the next position of the chip $v=s(k+1)$ is chosen according to the rule (ii') of Definition~\ref{def:space-filling}. We call it the \emph{space-filling exploration} of $M^*$.
		\item[(ii)] Let $(t_i)_{i\in\ZZ}$ be the triangles of $M$ such that the associated vertices $t_i^*$ of $M^*$ are ordered by the time they are first visited in the space-filling exploration of $M^*$. Then $(\etavf^{-1}(t_i))_{i\in \Z}$ is increasing.
		\item[(iii)] For the space-filling exploration of $M^*$, we define the \emph{treatment} of edges as above the statement of Proposition~\ref{prop32}. Let $(e_i)_{i\in \Z}$ be the edges of $M$ such that the associated edges $e_i^*$ of $M^*$ are treated in this order \referee{ and $e_0$ is the root-edge. Then, for all $i\in \Z$, $\etav(i)=e_i$.} 
	\end{compactitem}
\end{prop}


Let $(M,\sigma)\in\Phi^\infty(w)$ for some $w\in\imK$. 
We will now explain that if we are given $M\in\imT$ and $\tau^*\in \DFS_{M^*}$, then it is easy to recover $w$. 
Indeed, one can first determine $\sigma$ as in Definition~\ref{def:Gamma}, and then, by Proposition~\ref{prop:dfs-ordering-inf}(iii), one can determine the order of treatment of the edges of $M^*$, and further determine the letter $w_i$ by examining the $i$th treated edge $e_i^*$ (with $e_0^*$ being the root of $M^*$) as follows. If $e_i^*$ is not in $\tau^*$, then $w_i=c$. If $e_i^*=\{u,v\}$ is in $\tau^*$ with $u$ the parent of $v$, then we consider the edges $e_i,e_\ell,e_r$ in clockwise order around the face of $M$ dual to $v$, and get $w_i=a$ if $e_r$ is treated before $e_\ell$, and $w_i=b$ otherwise.

\section{Discrete dictionary III: tree of clusters, envelope excursions, and pivotal points}\label{sec:dual-DFS}
In this section, we explain how to obtain information about the clusters of a percolated triangulation $(M,\sigma)$ in $\bmT_P$ or $\imT_P$ in terms of the associated walk $w=\bPhi^{-1}(M,\sigma)$. In particular we will describe 
the \emph{tree of clusters} (describing the nesting structures of the clusters), 
the \emph{envelope closing times} (corresponding to the completion of clusters), and 
the \emph{pivotal points}
in terms of the walk $w$.

\subsection{Tree of clusters, and its relation to the exploration tree}
\label{sec:tree-cluster}
We first define the \emph{tree of clusters} of a site-percolated near-triangulation $(M,\sigma)\in\bmT_P$. Recall that the white (resp.\ black) \emph{clusters} are the connected components of the subgraph of $M$ induced by the white (resp.\ black) vertices. Since $\sigma$ satisfies the root-interface condition, there is a single white cluster containing the white outer vertices of $M$, which we call the \emph{outer white cluster} of $(M,\sigma)$. The \emph{outer black cluster} is defined similarly, and the other clusters are called \emph{inner clusters} of $(M,\sigma)$. 

We say that a cluster $C$ is \emph{incident} to a percolation cycle or percolation path $\gamma$ if $\gamma$ goes through at least one edge incident to $C$. 
For instance, the outer clusters are the only clusters incident to the percolation path of $(M,\sigma)$. For an inner cluster $C$ of $(M,\sigma)$ there is a percolation cycle $\gamma(C)$ called \emph{outside-cycle of $C$} 
which is the percolation cycle incident to $C$ which separates $C$ from the root-face of $M$.

\begin{definition}\label{def:tree-clusters}
	The \emph{tree of clusters} of $(M,\sigma)\in\bmT_P$, denoted $\cltree(M,\sigma)$, is the graph whose vertices are the clusters, and whose edges joins vertices if they correspond to clusters incident to a common percolation cycle or percolation path.
\end{definition}

It is easy to see that the tree of clusters $T=\cltree(M,\sigma)$ is indeed a tree. The tree encodes the nesting structure of the percolation cycles as illustrated in Figure~\ref{fig:tree-of-clusters2}. Each vertex of $T$ corresponds to a cluster and we can talk about its \emph{color} (black or white). Note also that every edge of $T$ is bicolor and corresponds to a percolation cycle or percolation path.

\fig{width=\linewidth}{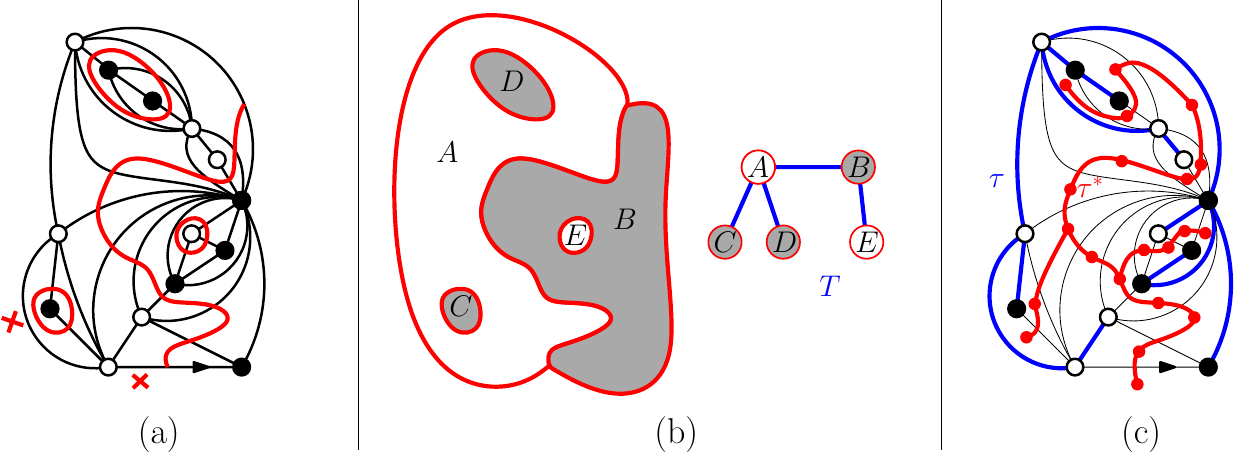}{(a) A percolated near-triangulation $(M,\sigma)\in\bmT_P$ and the percolation interfaces. (b) The clusters of $(M,\sigma)$ and the corresponding tree of clusters $T$. (c) The spanning tree $\tau^*=\dfs(M,\sigma)$ of $M^*$ and the spanning tree $\tau=\dfsdual(M,\sigma)$ of $M$. One can check that $T$ is obtained from $\tau$ by contracting every unicolor edge.}

We now explain the relation between the tree of clusters and the exploration tree.
\begin{definition}\label{def:dfsdual}
	Let $(M,\sigma)\in \bmT_P$ be a site-percolated near-triangulation, and let $\tau^*=\dfs(M,\sigma)\in \DFS_{M^*}$ be its exploration tree as defined in Section~\ref{sec:exploration-tree-from-walk}. We denote by $\dfsdual(M,\sigma)$ the spanning tree $\tau$ of $M$ which is dual to $\tau^*$ (that is, $\tau$ is made of the dual of the edges of $M^*$ not in $\tau^*$).
\end{definition}

\begin{prop}\label{prop:tree-clusters-dual-to-DFS}
	Let $(M,\sigma)\in \bmT_P$ be a site-percolated near-triangulation. The tree of clusters $T=\cltree(M,\sigma)$ is obtained from the spanning tree $\tau=\dfsdual(M,\sigma)$ of $M$ by contracting every unicolor edge (and forgetting the planar embedding). Moreover, for each vertex $v$ of $T$, the number of vertices in the cluster of $(M,\sigma)$ corresponding to $v$ is equal to the number of vertices of $\tau$ contracted to $v$. 
\end{prop}

Proposition~\ref{prop:tree-clusters-dual-to-DFS} is illustrated by Figure~\ref{fig:tree-of-clusters2}. 

\begin{proof} 
	Let $T'$ be the tree obtained from $\tau$ by contracting every unicolor edge.
	By Property (ii) of Theorem~\ref{thm:DFS-to-perco}, there is exactly one edge of $\tau$ across each percolation cycle or percolation path of $M$. Hence the path of $\tau$ between two vertices in a same cluster of $(M,\si)$ only goes through unicolor edges (because otherwise it would have to cross the same percolation cycle twice which is impossible).  Hence each cluster of $(M,\si)$ will contract to a single vertex of $T'$. Conversely if two vertices of $\tau$ belong to different clusters, they will not be contracted to the same vertex of $T'$.
	This creates a bijection between the vertices of $T'$ and the clusters of $(M,\si)$. Moreover, two clusters of $(M,\si)$ are adjacent if and only if the two corresponding vertices of $T'$ are adjacent. Hence, $T'=T$.
\end{proof}

\subsection{Tree of clusters as a function of the Kreweras walk}
\label{sec:tree-cluster-walk}
Motivated by Proposition~\ref{prop:tree-clusters-dual-to-DFS}, we would like to express the spanning tree $\dfsdual(M,\sigma)$ of $M$ as a function of the walk $w=\bPhi^{-1}(M,\sigma)$.

We first need to expand on the vocabulary about Kreweras walks established at the beginning of Section~\ref{sec:LR-decomposition-walk}.
Let $w=w_1w_2\ldots w_n\in \bmK$. We say that a matching $w_i,w_k$ is an \emph{ancestor} of a matching $w_{i'},w_{k'}$ if $w_i$ and $w_i'$ are either both $a$-steps or both $b$-steps, and $w_{i'},w_{k'}$ are enclosed by the matching $w_i,w_k$. We say that $w_i,w_k$ is the \emph{parent-matching} of $w_{i'},w_{k'}$ if $w_i,w_k$ is an ancestor of $w_i',w_k'$ and there is no ancestor of $w_{i'},w_{k'}$ enclosed by the matching $w_i,w_k$.

Let $w=w_1w_2\ldots w_n\in \bmK$ and let $\widetilde{w}=wab$ be the walk obtained from $w$ by appending an unmatched $a$-step and an unmatched $b$-step at the end of $w$ which we call the \emph{top steps}.
Let $A_w$ be the set of unmatched $a$-steps of $\widetilde{w}$, let $B_w$ be the set of unmatched $b$-steps of $\widetilde{w}$, and let $C_w$ be the set of $c$-steps of $\widetilde{w}$. Let $\widetilde{V}_w=A_w\cup B_w\cup C_w$, and let $V_w$ be the set obtained from $\widetilde{V}_w$ by removing the two top steps. We now define two functions $\matchingpar$ and $\treepar$ from $V_w$ to $\widetilde{V}_w$.
\begin{compactitem}
	\item For an unmatched $a$-step (resp.\ $b$-step) $w_i$, we set $\matchingpar(w_i)=w_j$, where $w_j$ is the next unmatched $a$-step (resp.\ $b$-step) in $\widetilde{w}$.
	\item For an unmatched $c$-step of type $a$ (resp.\ $b$) $w_i$, we set $\matchingpar(w_i)=w_j$, where $w_j$ is the next unmatched $c$-step of type $a$ (resp.\ $b$) in $\widetilde{w}$ if there is such a $c$-step, and $w_j$ is the first unmatched $a$-step (resp.\ $b$-step) in $\widetilde{w}$ otherwise.
	\item For a matched $c$-step $w_k\in C_w$ we consider the corresponding far-matching $w_i,w_k$. If $w_i,w_k$ has a parent-matching $w_{i'},w_{k'}$, then we set $\matchingpar(w_k)=w_{k'}$. Otherwise, if $c$ is of type $a$ (resp.\ $b$), then we set $\matchingpar(w_k)=w_j$, where $w_j$ is the first unmatched $a$-step ($b$-step) or unmatched $c$-step of type $a$ (resp.\ $b$) following $w_k$ in $\widetilde{w}$.
\end{compactitem}
An element of $\widetilde{V}_w$ is said to be \emph{white} if it is an $a$-step or a $c$-step of type $a$, and \emph{black} otherwise. For $v\in V_w$, we set $\treepar(v)=\matchingpar^{m}(v)$, where $m$ is the minimum positive integer such that $\matchingpar^{m}(v)$ has the same color as $\matchingpar^{m-1}(v)$ (note that if $\matchingpar^{k}(v)$ is a top step, then $\matchingpar^{k-1}(v)$ and $\matchingpar^{k}(v)$ have the same color, so $m\leq k$ and $\treepar(v)$ is well-defined).

\begin{definition}\label{def:dfsdualw}
	We denote by $\dfsdual(w)$ the tree with vertex set $\widetilde{V}_w$ and edge set $\{(\widetilde{a},\widetilde{b})\}\cup\{(v,\treepar(v))~|~v\in V_w\}$, where $\widetilde{a},\widetilde{b}$ are the top steps of $\widetilde{w}$. 
\end{definition}

The construction of $\dfsdual(w)$ is illustrated in Figure~\ref{fig:dfs-dual-theorem2}(b). We now define a correspondence $\lambdav$ between the vertices of $\dfsdual(w)$ and the vertices of $(M,\sigma)=\bPhi(w)$. For $x\in C_w$, we let $\lambdav(x)$ be the in-vertex of $M$ associated to the $c$-step $x$ by the mapping $\etavf$. 
For $x\in A_w$ (resp.\ $B_w$), we consider the edge $e$ of $M$ associated to $x$ by the mapping $\etae$. The edge $e$ is crossed by the percolation path of $(M,\sigma)$, and we let $\lambdav(x)$ be the white (resp.\ black) endpoint of $e$ (which is an outer-vertex). Lastly, for the top $a$-step (resp.\ $b$-step) $x$, we let $\lambdav(x)$ be the top-left (resp.\ top-right) vertex of $M$. It is clear from Properties (i-ii) of Theorem~\ref{thm:bij-extend}, that $\lambdav$ is a one-to-one correspondence between the vertices of $\dfsdual(w)$ and the vertices of $M$.

\fig{width=\linewidth}{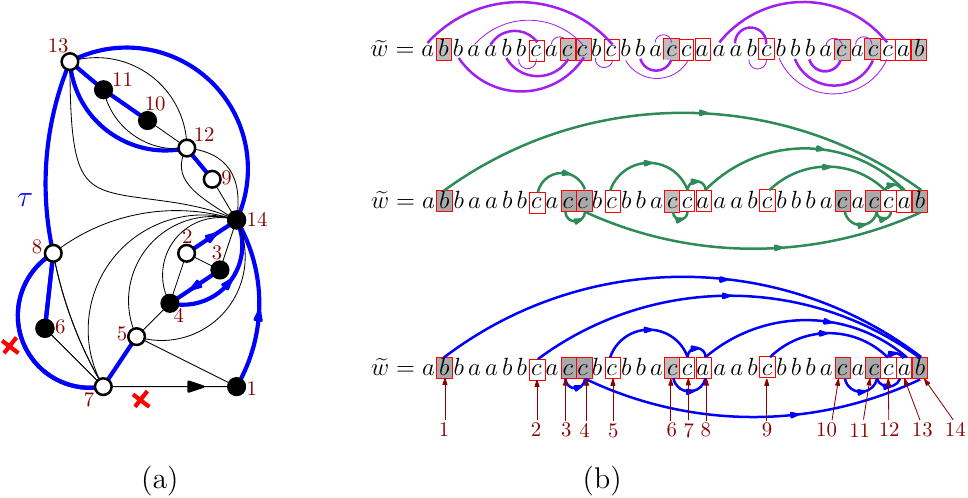}{Illustration of Theorem~\ref{thm:dfsdual} for $w=abbaabbcaccbcbbaccaaabcbbbacacc\in\bmK$. (a) The percolated near-triangulation $(M,\sigma)=\bPhi(w)$ and the spanning-tree $\tau=\dfsdual(M,\sigma)$. (b) The walk $\widetilde{w}=wab$ and the construction of $\dfsdual(w)$. The $c$-steps of type $a$ (resp.\ $b$) and the unmatched $a$-steps (resp.\ $b$-steps) are indicated by a \referee{white (resp.\ black)} box. The top row indicates the matching pairs $a,c$ (above $\widetilde{w}$) and the matching pairs $b,c$ (below $\widetilde{w}$). The far-matchings are drawn in bold lines while the close-matchings are drawn in thin lines.
	The middle row shows the arrows from each boxed step $w_k$ toward $\matchingpar(w_k)$. The bottom row shows the arrows from each boxed step $w_k$ toward $\treepar(w_k)$. One can check that the tree $\dfsdual(w)$ made of these arrows is isomorphic to the spanning-tree $\tau$ (the isomorphism $\lambdav$ is indicated by the labelling of the steps in $\widetilde{V}_w$ from 1 to 14 and the corresponding labelling of the vertices of $M$).}

\begin{thm}\label{thm:dfsdual}
	Let $w\in \bmK$ and let $(M,\sigma)=\bPhi(w)$. The tree $\dfsdual(w)$ (given by Definition~\ref{def:dfsdualw}) and the spanning tree $\dfsdual(M,\sigma)$ of $M$ (given by Definition~\ref{def:dfsdual}) are isomorphic trees. In fact, the mapping $\lambdav$ gives the isomorphism between these trees, and the color of vertices coincide through this isomorphism.
\end{thm}

Theorem~\ref{thm:dfsdual} is illustrated in Figure~\ref{fig:dfs-dual-theorem2}. It is proved in Section~\ref{sec:appendix-dual-dfs}.

\begin{remark}
	It is not hard to define a canonical embedding of $\dfsdual(w)$ so that $\dfsdual(w)$ and $\dfsdual(M,\sigma)$ coincide as plane trees but this is omitted as it is not relevant to our study. 
\end{remark}

\subsection{Tree of clusters for infinite maps}\label{sec:cluster-tree-inf}
In this section we briefly explain how to adapt Sections~\ref{sec:tree-cluster} and~\ref{sec:tree-cluster-walk} to the infinite volume setting. 
Let $(M,\sigma)\in\imT_P$ be sampled from the distribution of the percolated UIPT. We define the tree of clusters $\cltree(M,\sigma)$ as in definition~\ref{def:tree-clusters}. The spanning tree $\dfs(M,\si)$ of $M^*$ is well-defined almost surely, so we can define the spanning tree $\dfsdual(M,\sigma)$ of $M$ as in Definitions~\ref{def:dfsdual}.

We recall from \cite{angel-peeling} that almost surely all the clusters of $(M,\si)$ are finite. This implies that $\cltree(M,\sigma)$ is almost surely locally finite and one-ended. Observe also that the proof of Proposition~\ref{prop:tree-clusters-dual-to-DFS} carries to the infinite setting so we obtain:
\begin{prop}\label{prop:tree-clusters-dual-to-DFS-infty}
	The relation between $\cltree(M,\sigma)$ and $\dfsdual(M,\sigma)$ stated in Proposition~\ref{prop:tree-clusters-dual-to-DFS} holds in the infinite setting.
\end{prop}

We will now define $\dfsdual(w)$ similarly as in Definition~\ref{def:dfsdualw}. Let $w\in\imK$ be sampled with the uniform distribution. 
Let $V_w$ be the set of $c$-steps of $w$.
For all $w_k\in V_w$ we may define $\matchingpar(w_k)$ exactly as in Section~\ref{sec:tree-cluster-walk}: since $w$ is in $\imK$, the far-matching $w_i,w_k$ has a parent-matching $w_{i'},w_{k'}$ and we set $\matchingpar(w_k)=w_{k'}$. Since $w$ is sampled with the uniform distribution, it is easy to see that almost surely, for every $v\in V_w$ there exists a positive integer $m$ such that $p^{m-1}(v)$ and $p^m(v)$ have the same color. In this case, we define $\treepar(v)=p^m(v)$ for the minimum such $m$. Finally, we define $\dfsdual(w)$ as the tree with vertex set $V_w$ and edge set $\{(v,\treepar(v))~|~v\in V_w\}$. 

A vertex $x\in V_w$ is white (resp. black) if it is a $c$-step of type $a$ (resp. $b$), and we define $\lambda_{\op{v}}(x)$ to be the vertex of $(M,\si)=\Phi^\infty(w)$ associated to the $c$-step $x$ by the mapping $\etavf$. We now state the infinite volume analogues of Theorem~\ref{thm:dfsdual-inf}. The proof is given in Section~\ref{app:inf}.
\begin{thm}\label{thm:dfsdual-inf}
	Let $w\in \imK$ be sampled from the uniform distribution and let $(M,\sigma)=\Phi^\infty(w)$. 	The tree $\dfsdual(w)$ and the spanning tree $\dfsdual(M,\sigma)$ of $M$ are almost surely well-defined and are isomorphic trees. In fact, the mapping $\lambdav$ gives the isomorphism between these trees, and the color of vertices coincide through this isomorphism.
\end{thm}

\subsection{Envelope excursion of a percolation cycle}\label{subsec:env}
We now define subwalks of a Kreweras walk which correspond to the construction of clusters.
Let $(M,\sigma)$ be in $\bmT_P$ or in $\imT_P$, and let $\gamma$ be a percolation cycle. If $(M,\sigma)\in\bmT_P$, then the \emph{inside-region} of $\gamma$ is the region enclosed by $\gamma$ not containing the root-face of $M$. If $(M,\sigma)\in\imT_P$, then the \emph{inside-region} of $\gamma$ is the bounded region enclosed by $\gamma$. We say that a vertex or edge $x$ of $M$ is \emph{inside} $\gamma$ if $x$ is (entirely) in the inside-region of $\gamma$. 

As stated in Theorem~\ref{thm:DFS-to-perco} (Property (ii)) and Theorem~\ref{thm:dfs-inf}, every percolation cycle $\gamma$ of $(M,\sigma)$ has a unique edge $e^*$ of $M^*$ which is not in the exploration tree $\dfs(M,\sigma)$. We call $e^*$ the \emph{envelope edge} of $\gamma$. 
Recall from Theorem~\ref{thm:DFS} (ii), that the mapping $\etavf$ gives a one-to-one correspondence between the $c$-steps of $w$ and the in-edges of $M^*$ not in $\dfs(M,\sigma)$. We call \emph{envelope step} of $\gamma$ the $c$-step of $w$ corresponding to its envelope edge. Note that in the finite case  $(M,\sigma)\in\bmT_P$, the envelope edge is dual to an inner edge of $M$, hence the envelope step is always a \emph{matched $c$-step}. We call \emph{envelope excursion} a subwalk $w_iw_{i+1}\ldots w_k$ of $w$ such that $w_i$ is far-matched to $w_k$ and $w_k$ is an envelope step. The time $k$ is called the \emph{envelope closing time} of $\ga$. 
The following result is immediate from Theorems~\ref{thm:dfsdual} and~\ref{thm:dfsdual-inf}. 
\begin{cor} Let $w$ be in $\bmK$ or $\imK$.
	Let $w'=w_i\dots w_k$ be a subwalk of $w$, with $w_k=c$. The subwalk $w'$ is an envelope excursion if and only if $w_i$ is the far-match of $w_k$ and $\treepar(w_k)$ is of different color than $w_k$. 
	\label{cor9}
\end{cor}
\begin{definition}\label{def:lt-cluster}
	To a percolation cycle $\ga$ of $(M,\sigma)$, we associate a discrete looptree $\frk L(\ga)$ which, roughly speaking, corresponds to the boundary of the cluster enclosed by $\ga$. Precisely, let $C$ be the percolation cluster for which $\ga$ is the outside-cycle.
	Let $V$ be the set of vertices of $M$ which are in $C$ and incident to an edge crossed by $\ga$, and let $E$ be the set of edges of $M$ incident to a triangle crossed by $\ga$ and joining two vertices of $C$. We denote by $\frk L(\ga)$ the (unrooted) discrete looptree obtained from the submap $M_\ga$ of $M$ with vertex-set $V$ and edge-set $E$ by replacing every bridge of $M_\ga$ by a double edge. We define the root-edge $e_0$ of $\frk L(\ga)$ as follows: denoting by $e^*\in E(M^*)$ the envelope edge of $\ga$, by $e\in E(M)$ the edge dual to $e^*$, and by $v$ the endpoint of $e$ inside $\ga$, we define $e_0$ as the first edge of $\frk L(\ga)$ following $e$ in counterclockwise direction around $v$.
\end{definition}

\begin{lemma} \label{lem:envelope-matching-is-loop-building}
	Let $w\in\bmK$ and $(M,\sigma)=\bPhi(w)$, or let $w\in\imK$ and $(M,\sigma)=\Phi^\infty(w)$.
	Let $\gamma$ be a percolation cycle of $(M,\sigma)$, let $w'=w_iw_{i+1}\ldots w_k$ be the corresponding envelope excursion, and let $\wt w=w_iw_{i+1}\ldots w_{k-1}$. 
	\begin{compactitem}
		\item[(i)] The envelope step $w_k$ corresponds via $\etae$ to an edge of $M$ crossing $\gamma$, namely the dual $e$ of the envelope edge of $\gamma$. The envelope step $w_k$ also corresponds via $\etavf$ to a vertex of $M$ inside $\gamma$, namely the endpoint of $e$ inside $\gamma$. The step $w_i$ corresponds via $\etavf$ to a triangle crossed by $\gamma$.
		\item[(ii)] All the edges of $M$ which are either inside $\gamma$ or crossing $\gamma$ correspond via $\etae$ to steps in the envelope excursion $w'$. 
		All the vertices of $M$ inside $\gamma$, and all the triangles inside $\gamma$ or crossed by $\gamma$ correspond via $\etavf$ to steps in the envelope excursion $w'$. 
		\item[(iii)] The triangles crossed by $\gamma$ (or, dually, the vertices of $M^*$ on $\ga$) correspond via $\etavf$ to the spine steps of $\wt w$, that is, the $a$-steps and $b$-steps of $w'$ not enclosed in any close-matching inside $\wt w$. 
		\item[(iv)] 
		If $w_k$ is a $c$-step of type $a$ (resp. $b$), then we consider the discrete looptree $\frk L=\frk L_\ell(\wt w)$ (resp. $\frk L=\frk L_r(\wt w)$) given by Definition~\ref{def:LR-w}. 
		The root-edge of $\frk L$ is on a bubble $B$ of degree 2. Upon deleting the two edges of $B$ from $\frk L$, we get an isolated vertex and an unrooted map $\frk L'$, which we canonically root at the edge of $\frk L'$ following $B$ in counterclockwise direction around $\frk L$ (see Figure~\ref{fig:looptree-cycle}). Then $\frk L'$ is equal to   $\frk L(\ga)$.
	\end{compactitem}
\end{lemma}
\referee{We mention that some of the steps of the enveloppe excursion $w'$ corresponding to a percolation cycle $\gamma$ encode (via $\etavf$ and $\etae$) some vertices, triangles and edges which are outside of $\gamma$. See for instance Figure \ref{fig:looptree-cycle}.}

\fig{width=\linewidth}{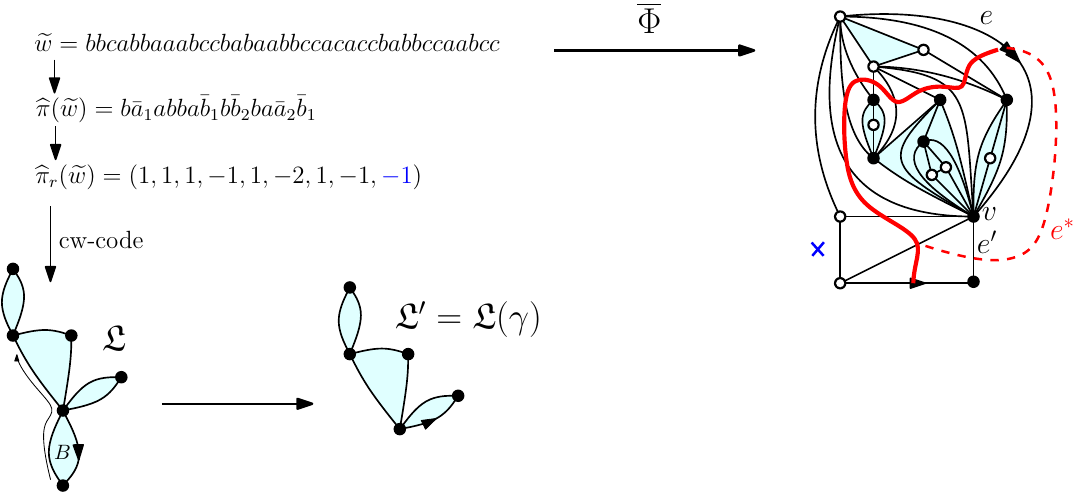}{Illustration of  Lemma~\ref{lem:envelope-matching-is-loop-building}(iv). 
	Here the envelope excursion of $\gamma$ is $w'=\wt{w} c$, where $\wt w=bbc  abbaaabcc b abaabbccacacc b a bbcc aabcc$. The construction of the looptree $\frk L$ is indicated on the left. The map $\bPhi(\wt w)$ is indicated on the right. Note that the edges $e$ and $e'$ would be identified in $(M,\si)$, and their dual $e^*$ is the envelope edge of $\ga$.
}

Lemma~\ref{lem:envelope-matching-is-loop-building} is proved in Section~\ref{sec:appendix-dual-dfs}, together with the following claim.

\begin{claim}\label{claim:disjoint-intervals}
	If $w'$, $w''$ are envelope excursions of $w$, then either $w',w''$ are disjoint or one is included in the other. 
\end{claim}


\subsection{Pivotal points of the percolation}\label{subsec:pivot}
Let $(M,\sigma)$ be in $\bmTP$ or in $\imTP$. Let $\Gam$ be the set of percolation cycles of $(M,\sigma)$.
For $v\in V(M)$, let $\Gam_v$ be the set of percolation cycles of the percolated map obtained from $(M,\sigma)$ by flipping the color of $v$. Let $\Gam\Delta \Gam_v$ denote the symmetric difference of $\Gam$ and $\Gam_v$. 
We say that $v$ is a \emph{pivotal point} of $(M,\sigma)$ if $\Gam\Delta \Gam_v$ contains at least three cycles. Equivalently, a pivotal point is a vertex of $M$, such that flipping the color of $v$ results in some splitting or merging of some percolation cycles.
We are mainly interested in \emph{macroscopic} pivotal points, which are pivotal points such that at least three cycles in $\Gam\Delta \Gam_v$ are macroscopic in the scaling limit (equivalently, they enclose a linear number of edges). 
\begin{definition}\label{def:pivot} 
	We call \emph{area} of a percolation cycle $\ga$ of $(M,\si)$ the number of vertices of $M$ in the inside-region of $\ga$. 
	Given $\alpha>0$, we say that a vertex $v\in V(M)$ is an $\alpha$-\emph{pivotal point} if there are at least three cycles in $\Gam\Delta \Gam_v$ with area at least $\alpha$.
\end{definition}
\begin{figure}[h!]
	\centering
	\includegraphics[scale=0.8]{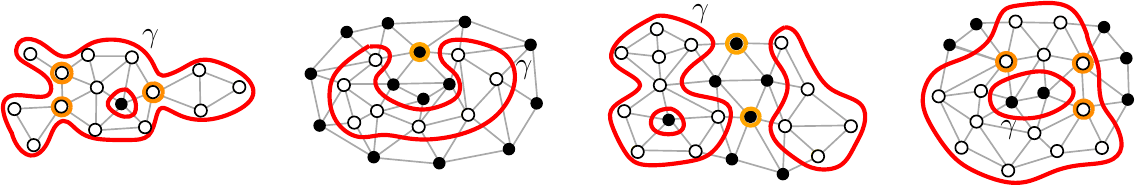}
	\caption{Pivotal points (marked in orange) of the four types 1, 2, 3, and 4 (from left to right).}\label{fig:piv}
\end{figure}

We say that a pivotal point is an $\alpha$-pivotal point of type 1 (resp.\ type 3) if there are two non-nested cycles in $\Gam_v\setminus \Gam$ (resp.\ $\Gam\setminus \Gam_v$) each with area at least $\alpha$. We say that a pivotal point is an $\alpha$-pivotal point of type 2 (resp.\ type 4) if there are two nested cycles in $\Gam_v\setminus \Gam$ (resp.\ $\Gam\setminus \Gam_v$) such that the areas $A,B$ of the inner and outer cycles satisfy $A\geq \al$ and $B-A\geq \al$. See Figure~\ref{fig:piv}. Observe that a pivotal point may be of multiple types. However, \nina{as explained in Remark \ref{rmk9},} by combining our scaling limit result in Section~\ref{sec:conv} with the result of Camia and Newman \cite[Theorem 2]{camia-newman-full} (which we recall in Lemma~\ref{lem:pivotal-continuous}) it follows that, for any fixed $c>0$ and $\eps>0$, with probability converging to 1 as $n\rta\infty$, each $\eps n$-pivotal point $v$ for which $|\etavf^{-1}(v)|<c n$ is an $\eps n$-pivotal point of a unique type. Furthermore, by the same results, with probability converging to 1 as $n\rta\infty$, for any $v$ for which $|\etavf^{-1}(v)|<c n$, the symmetric difference $\Gam\Delta \Gam_v$ contains at most 3 cycles of area at least $\eps n$.

Given a percolation cycle $\gamma$ of area at least $\alpha$ we say that an $\alpha$-pivotal point $v$ is \emph{associated with $\gamma$} if the condition for $\alpha$-pivotal of type 1 or 2 is satisfied and $\ga\in\Gam\setminus \Gam_v$, or if the condition for $\alpha$-pivotal of type 3 or 4 is satisfied for $\ga,\ga'\in\Gam\setminus \Gam_v$, $\ga,\ga'$ have area at least $\alpha$, and the envelope excursion of $\ga$ starts after the envelope excursion of $\ga'$.

We now describe a way of detecting the pivotal points associated with a cycle $\ga$ in terms of its envelope excursion.

Let $w\in\{a,b,c \}^\ZZ$. As in Section~\ref{sec:LR-decomposition-infinite}, we write $w=w^-w^+$, and we define $T$ and $\wh\pi(w^-)$ as in Definition~\ref{def:pi-w-inf}. We denote by $(\wh L_k,\wh R_k)_{k\in\ZZ^{\leq 0}}$ the walk on $\Z^2$ with steps given by $\wh\pi(w^-)$ and such that $(\wh L_0,\wh R_0)=(0,0)$. Let $\wh\eta_{\op{e}}:\ZZ^{\leq 0}\to E(M)$ be defined by $\wh\eta_{\op{e}}(k)=\eta_{\op{e}}(T(k))$. Note that the image of $\wh\eta_{\op{e}}$ is the set of edges crossed by the percolation path of $(M^-,\sigma^-)$ (not including the root-edge). Define $\wh\eta_{\op{v}}:\ZZ^{\leq 0}\to V(M)$ by letting $\wh\eta_{\op{v}}(k)$ be the unique vertex on the face
$\eta_{\op{vf}}(T(k))\in F(M)$ which is not an endpoint of $\wh\eta_{\op{e}}(k)$.
We point out that $\wh\eta_{\op{v}}$ is \emph{not} the composition of $\eta_{\op{vf}}$ and $T$ (in fact $\wh\eta_{\op{v}}$ is not injective, and its image is the set of vertices incident to the triangles on the percolation path of $(M^-,\si^-)$). If $w\in\bmK$ and $w^-$ is a given prefix of $w$, then we define $\wh\pi(w^-)$, $\wh\eta_{\op{e}}$, and $\wh\eta_{\op{v}}$ similarly. 

Recall the notation $\wt L,\wt R, \frk L_\ell(w^-), \frk L_r(w^-)$ of Section~\ref{sec:LR-decomposition-infinite}: the walk $\wt L$ (resp. $\wt R$) is obtained from $\wh L$ (resp. $\wh R$) by removing the 0-steps, and encodes a forested line $\frk L_\ell(w^-)$ (resp. $\frk L_r(w^-)$). The following result will be proved in Section~\ref{sec:appendix-dual-dfs}.
\begin{lemma}\label{prop:piv}
	In the setting above, if $i\in\ZZ^{\leq -1}$ is such that	
	\eqb
	\wh L_i<\wh L_{i-1},\textrm { and } 
	\wh L_i\geq \min\{ \wh L_j~|~ i<j\leq 0, \wh L_j\neq \wh L_{j-1} \},
	\label{eq:piv}
	\eqe
	then $\wh\eta_{\op{v}}(i)$ is a vertex on at least two bubbles of the forested line $\frk L_\ell(w^-)$ encoded by $\wh L$ (see Section~\ref{sec:LR-decomposition-infinite} for the definition of $\frk L_\ell(w^-)$). By symmetry, the same statement holds with $\wh R$ instead of $\wh L$.
	
	Conversely, let $\ga$ be a percolation cycle of $(M,\si)$, and suppose that $w$ is recentered in such a way that $w_0$ is the envelope step of $\ga$. Let $s>0$ be the length of the cycle $\ga$. 
	If $v$ is an $\al$-pivotal point of type 1 or 2 associated with $\ga$, then the following holds:
	\begin{compactitem}
		\item[(a)] there exists $i\in\{-s+1,\ldots,0 \}$ such that $v=\wh\eta_{\op{v}}(i)$ and either $\wh L_i<\wh L_{i-1}$ and the looptree encoded by the \emph{excursion of $\wh L$ ending at $i$} (that is, the walk $\wh L_m,\wh L_{m+1},\ldots,\wh L_i$, where $m=\min\{n~|~\forall j\in[n,i],\wh L_j\geq \wh L_i\}$) encloses at least $\alpha$ vertices of $M$, or the same holds with $\wh R$ instead of $\wh L$.
	\end{compactitem}
	For any percolation cycle $\gamma'$ of $(M,\si)$ there is a set $S_{\ga,\ga'}$ of at most 3 vertices such that if $v\notin S_{\ga,\ga'}$ is an $\al$-pivotal point of type 3 or 4 associated with $\ga$ such that flipping the color of $v$ merges the cycles $\ga$ and $\ga'$, then the following holds:
	\begin{compactitem}
		\item[(b)] there exists $i\in\{-s+1,\ldots,0 \}$ such that $v=\wh\eta_{\op{v}}(i)$ and either $\wh L_i<\wh L_{i-1}$ and  $i$ is a \emph{global running infimum} for $\wh L$ relative to time $-s$ (that is, $\wh L_i<\wh L_j$ for all $j\in \{-s,\ldots,i-1\}$), or the same holds with $\wh R$ instead of $\wh L$. 
	\end{compactitem}
\end{lemma}

\section{Proofs of the bijective correspondences}\label{sec:appendix}
\subsection{Proofs for Section~\ref{sec:exploration-tree-from-map}: link between percolation configurations and DFS trees.}\label{sec:appendix-DFS-Perco}
In this subsection, we prove Theorem~\ref{thm:DFS-to-perco} and Claim~\ref{claim:color-doesnt-matters}. 
We start with a basic claim about depth-first search processes.


\begin{claim}\label{claim:descendants}
	Let $G$ be a graph and let $v_0$ be a vertex. Let $T$ be a tree obtained by a DFS of $G$ starting at $v_0$. Consider an arbitrary step of the DFS of $G$. Let $u$ be the position of the chip at that step, and let $U$ be the set made of $u$ and all yet unvisited vertices of $G$. If a vertex $v\in U$ is reachable from $u$ by a path of $G$ using only vertices in $U$, then $v$ will be a descendant of $u$ in $T$. 
\end{claim}

\begin{proof}
	Consider a path $u=v_1,\ldots,v_k=v$ with all the vertices $v_i$ in $U$. By Claim~\ref{claim:characterization-DFStree}, for all $i\in[k-1]$, the vertices $v_{i}$ and $v_{i+1}$ are $T$-comparable. Since the vertices $v_1,\ldots,v_k$ are in $U$, they cannot be ancestors of $u$ in $T$. Hence, by induction on $i\in[k]$, the vertex $v_i$ is a descendant of $u$ in $T$. 
\end{proof}

Next, we prove Claim~\ref{claim:color-doesnt-matters}.
From now on, we let $M$, $M^*$ and $v_0$ be as in Definition~\ref{def:Perc-and-DFS-sets}. 

\begin{proof}[Proof of Claim~\ref{claim:color-doesnt-matters}]
	For a face $f$ of $M^*$, we denote by $v_f$ the first vertex incident to $f$ visited during the DFS.
	Consider a DFS $X$ of $M^*$ satisfying the hypotheses of Claim~\ref{claim:color-doesnt-matters}, and let $\tau^*$ be the associated spanning tree. 
	We want to show that $\tau^*=\Delta_M(\gao)$. For this it suffices to show that changing the choices of made during the DFS $X$ when the chip position $u$ is not equal to $v_f$ (where $f$ is the forward edge at $u$) does not affect the final DFS tree $\tau^*$.
	
	Consider a step in the DFS $X$ where there are several edges between the chip position $u$ and some unvisited vertices. Let $f$ be the forward face of $u$ and let us suppose that $u\neq v_f$. By Claim~\ref{claim:descendants}, all the vertices incident to $f$ are descendant of $v_f$ in $\tau^*$, so $u$ is a descendant of $v_f$ and all the vertices on the path $P$ from $v_f$ to $u$ in $\tau^*$ have already been visited.
	
	Let $e_2$ and $e_3$ be the forward edges at $u$ and let $u_2$ and $u_3$ be their endpoints. By hypothesis 
	there are several edges toward unvisited vertices, hence $u_2,u_3$ are yet unvisited vertices. We let $U_2$ and $U_3$ be respectively the set of unvisited vertices reachable from $u_2$ and from $u_3$ through unvisited vertices.
	We now observe that there is no path of unvisited vertices between $u_2$ and $u_3$, because such a path would have to cross $P$ (indeed $P$ has both endpoints on $f$ and $u_2,u_3$ are on different sides of $P$ at $u$). Therefore $U_2\cap U_3=\emptyset$ and the DFS $X$ will visit independently the two sets $U_2,U_3$: the chip will visit one set entirely, then backtrack to $u$, then visit the other set entirely. Thus, the choice of the DFS between the forward edges $e_2$ and $e_3$ at this step will not affect the final DFS tree $\tau^*$ (it will only change the order of visit of the vertices). 
\end{proof}

We now state another easy claim which will be used in the proof of Theorem~\ref{thm:DFS-to-perco}.

\begin{claim} \label{claim:color-matters}
	Let $\gao\in \Perc_M$ and let $f$ be a face of $M^*$ dual to an inner vertex of $M$. 
	We consider the DFS of $M^*$ corresponding to $\Delta_M(\gao)$ (see Definition~\ref{def:Delta}). Let $u$ be the first vertex of $M^*$ incident to $f$ encountered during this DFS (note that $f$ is the forward face of $u$). Then, the tree $\tau^*=\Delta_M(\gao)$ will contain exactly one of the forward edges at $u$. Namely, if $f$ is black (resp. white) $\tau^*$ will contain the left forward edge $e_2$ (resp. right forward edge $e_3$), but not the other forward edge.
\end{claim}

\begin{proof} 
	It is obvious from the definition of $\Delta_M$ that, if $f$ is black (resp. white), then $\tau^*$ will contain $e_2$ (resp. $e_3$). So we only need to show that $\tau^*$ does not contain both $e_2$ and $e_3$. Suppose by contradiction that both $e_2$ and $e_3$ are in $\tau^*$, so that $u$ has two children $u_1,u_2$ in $\tau^*$. 
	Since all the vertices of $M^*$ incident to $f$ are unvisited at the time the DFS arrives at $u$, Claim~\ref{claim:descendants} ensures that they will all be descendants of $u$ in $\tau^*$. So there must exist two adjacent vertices $v_1,v_2$ incident to $f$, with $v_1$ descendant of $u_1$ and $v_2$ descendant of $u_2$. In this case, $v_1,v_2$ are not $\tau^*$-comparable, which contradicts Claim~\ref{claim:characterization-DFStree}.
\end{proof}

\begin{proof}[Proof of Theorem~\ref{thm:DFS-to-perco}]
	We first prove that $\Delta_M$ is injective. Let $\gao_1,\gao_2$ be distinct inner colorings of $M$, and let $F$ be the set of faces of $M^*$ having different colors in $\gao_1$ and $\gao_2$. 
	Let $u$ be the first vertex of $M^*$ incident to a face in $F$ encountered during the DFS of $M^*$ corresponding to $\Delta_M(\gao_1)$.
	It is clear that $u$ is also the first vertex of $M^*$ incident to a face in $F$ encountered during the DFS corresponding to $\Delta_M(\gao_2)$. Note also that the parent-edge $e_1$ of $u$ is the same in $\Delta_M(\gao_1)$ and in $\Delta_M(\gao_2)$. Now, Claim~\ref{claim:color-matters} ensures that the trees $\Delta_M(\gao_1)$ and $\Delta_M(\gao_2)$ each contain a different forward edge of $u$. Hence $\Delta_M(\gao_1)\neq \Delta_M(\gao_2)$. Thus $\Delta_M$ is injective.
	
	Next we show that $\Delta_M\circ\Lambda_M=\textrm{Id}$.
	Let $\tau^*\in \DFS_{M^*}$, let $\gao=\Lambda_M(\tau^*)$, and let $\tau'=\Delta_M(\gao)$. 
	Suppose by contradiction that $\tau'\neq \tau^*$. Let $e$ be the first edge in $\tau'\setminus \tau^*$ added to $\tau'$ during the DFS corresponding to $\Delta_M(\gao)$. Let $u,v$ be the endpoints of $e$, with $u$ the parent of $v$ in $\tau'$. Observe that the path $P$ from $v_0$ to $u$ is the same in $\tau^*$ and in $\tau'$ (by the choice of $e$). Let $e_1$ be the parent-edge of $u$ in $\tau^*$ (or equivalently, in $\tau'$). Let $f$ be the forward edge at $u$ and let $e_2$ and $e_3$ be the left-forward and right-forward edges respectively. In other words, $e_1,e_2,e_3$ are the edges incident to $u$ in clockwise order, with $e_1$ in $P$ and $f$ is the face of $M^*$ between $e_2$ and $e_3$.
	By Claim~\ref{claim:characterization-DFStree} we know that the vertices $u,v$ are $\tau^*$-comparable. Moreover, $v$ cannot be an ancestor of $u$ in $\tau^*$ (because it is not on $P$), hence it is a descendant of $u$ in $\tau^*$. Let $Q$ be the path of $\tau^*$ from $u$ to $v$. Note that $Q\cup\{e\}$ form a cycle of $M^*$, and that $P$ and $f$ are on different sides of this cycle because the paths $P,Q\subset \tau^*$ cannot cross. In particular, $f$ is not incident to $v_0$ (hence $f$ is the dual of an inner face of $M$). In fact, $u$ is the first vertex incident to $f$ encountered during the DFS corresponding to $\Delta_M(\gao)$ (otherwise by Claim~\ref{claim:descendants}, $u$ would have to be a descendant in $\tau'$ of the first vertex incident to $f$). Thus by definition of $\Delta_M$, the face $f$ is black if $e=e_2$ and white if $e=e_3$. However, by definition of $\Lambda_M$, the face $f$ is white if $e=e_2$ and black if $e=e_3$. We reach a contradiction, hence $\tau'=\tau^*$ and $\Delta_M\circ\Lambda_M=\textrm{Id}$.
	
	Since $\Delta_M\circ\Lambda_M=\textrm{Id}$ and $\Delta_M$ is injective, we see that $\Delta_M$ is a bijection and that $\Lambda_M$ is the inverse mapping.
	We now need to prove the statements (i-iii) about percolation interfaces. Let $\sigma$, $\gao$, and $\tau^*$ be as in Theorem~\ref{thm:DFS-to-perco}. In the definition of $\tau^*=\Delta_M(\gao)$ (Definition~\ref{def:Delta}) we have used the convention that the faces of $M^*$ dual to outer vertices of $M $ are considered white. However, by Remark~\ref{rk:outer-coloring-does-not-matter}, we would have gotten the same tree $\tau^*$ by using the colors in $\sigma$ of the outer-vertices of $M$. Now consider the DFS of $M^*$ resulting from this convention. 
	It is clear from the definitions, that this DFS the chip will first visit all the vertices on the percolation path (from $v_0$ to its other end) before visiting any other vertex. 
	This proves (i). Similarly, for any percolation cycle $C$ of $(M,\sigma)$, we can consider the first time the DFS reaches a vertex $v$ of $M^*$ on the cycle $C$. It is clear from the definitions that in the next few steps of the DFS, the chip will follow the edges of the cycle $C$ starting at $v$ (without visiting any vertex not on $C$) until it reaches the second neighbor of $v$ on $C$. This proves (ii). Moreover, by definition, the direction in which the percolation interface $C$ is followed is such that the black faces are on the right and the white faces are on the left. This proves (iii).
\end{proof}


\subsection{Proofs for Section~\ref{subsec:finite}: the bijection $\Phi$ and its relation to the bijection from \cite{bernardi-dfs-bijection}.}\label{sec:appendix-bijection}
In this Section we prove Theorem~\ref{thm:bij-OB} using the results from \cite{bernardi-dfs-bijection}. Roughly speaking, our bijection $\Phi$ between $\mK$ and $\mT_P$ is obtained by composing the bijection $\Om$ obtained in \cite{bernardi-dfs-bijection} (and whose definition is recalled below) with the bijection $\Lambda_M$ described in Section~\ref{sec:exploration-tree-from-map} (Definition~\ref{def:Gamma}). 
The bijection $\Om$ was actually defined in the ``dual setting'', that is, in terms of near-cubic maps instead of near-triangulations, and we refer the reader to Section~\ref{sec:exploration-tree-from-map} for definitions about duality.

We first recall how the bijection $\Om$ was defined in \cite{bernardi-dfs-bijection}; see Figure~\ref{fig:bijection-Kreweras} for an example.
The bijection $\Om$ is between the set $\mK$ of Kreweras walks and a set $\mC_T$ of near-cubic maps with a marked spanning tree. 
More precisely, $\mC_T$ is the set of triples $(M^*,e^*,\tau^*)$, where $M^*$ is the dual of a near-triangulation $M\in\mT$, $e^*$ is an edge of $M^*$ incident to the root-vertex $v_0$, and $\tau^*$ is in $\DFS_{M^*}^{e^*}$ (see Definition~\ref{def:Perc-and-DFS-sets-e}). 
We call \emph{head-edge} the edge $e^*$ and we call \emph{head-vertex} the endpoint of $e^*$ distinct from $v_0$ (if $e^*$ is a self-loop, the head-vertex is $v_0$).

The definition of $\Om$ follows a similar scheme as the Definition~\ref{def:Phi} of $\Phi$, in the sense that the image $\Om(w)$ of a walk $w=w_1\cdots w_n\in\mK$ is defined as the result of applying successively some elementary construction steps $\Om_{w_1}$, $\Om_{w_2}$, \ldots, $\Om_{w_n}$. In our figures we adopt a convenient drawing convention (used in \cite{bernardi-dfs-bijection}) for the elements $(M^*,e^*,\tau^*)$ of $\mC_T$. 
Namely, we do not draw the map $M^*$ itself, but rather the map $\overline{M}^*$ obtained from $M^*$ by ``ungluing from $v_0$ all the non-root edges incident to $v_0$'': these edges become edges incident to special non-root vertices of degree 1 called \emph{buds} (so $M^*$ would be obtained from the drawn map $\overline{M}^*$ by gluing all the buds to $v_0$). See Figure~\ref{fig:bijection-Kreweras} for an example.

\fig{width=\linewidth}{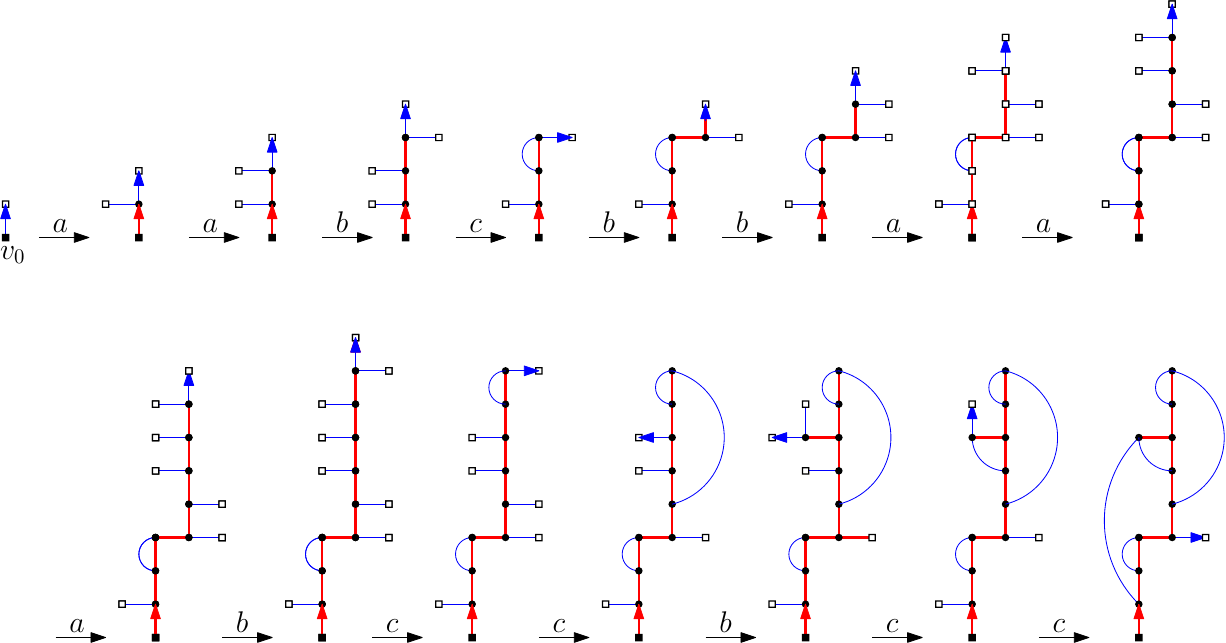}{The bijection $\Om$ defined in \cite{bernardi-dfs-bijection}, applied to the walk $w=aabcbbaaabccbcc$. Here the elements $(M^*,e^*,\tau^*)\in\mC_T$ are drawn according to the following conventions: (1) Instead of drawing $M^*$, we draw the the map $\overline{M}^*$ obtained from $M^*$ by ``ungluing from $v_0$ all the non-root edges incident to $v_0$'' these edges become incident to special vertices of degree 1 called \emph{buds} which we draw as white squares,
	(2) the DFS tree $\tau^*$ is indicated by bold red lines, (3) the root-edge $e_0^*$ and the head-edge $e^*$ are indicated by arrows: the arrow pointing toward a bud is $e^*$.}

We now give the precise definition of $\Om$ from \cite{bernardi-dfs-bijection}.
Let $M_0^*$ be the rooted map with one vertex and one self-loop, let $e_0^*$ be its root-edge, and let $\tau_0^*$ be the unique spanning tree of $M_0^*$. For $w=w_1\cdots w_n\in\mK$, the image $\Om(w)$ is defined as the triple $(M^*,e_0^*,\tau^*)=\Om_{w_n}\circ\cdots \circ\Om_{w_2}\circ \Om_{w_1}(M_0^*,e_0^*,\tau_0^*)$, where the mappings $\Om_a,\Om_b,\Om_c$ are defined in Definition~\ref{def:Om} and represented in Figure~\ref{fig:bij-rules-cubic}.

\fig{width=\linewidth}{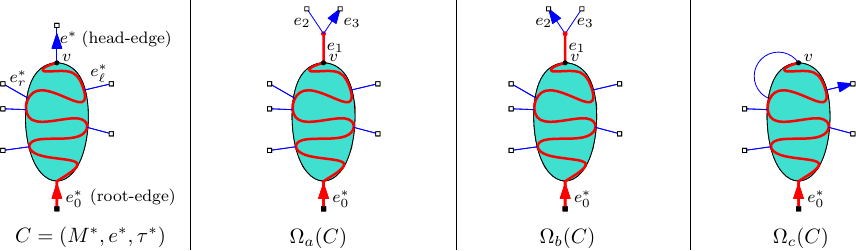}{The elementary steps $\Om_a,\Om_b,\Om_c$ applied to an element $C=(M^*,e^*,\tau^*)\in\mC_T$. The element $C$ is represented on the left, with $v$ being the non-root vertex incident to the head-edge $e^*$. The tree $\tau^*$ is not represented, except for the path $P$ of $\tau^*$ from the root-vertex $v_0$ to the head-vertex $v$ which is represented in bold red line. 
	As in Figure~\ref{fig:bijection-Kreweras}, the non-root edges incident to the root-vertex in $M^*$ are not represented as incident to $v_0$, but instead as edges incident to ``buds''. Since $\tau^*$ is in $\DFS_{M^*}^{e^*}$, any bud is adjacent to a vertex on the path $P$.}

\begin{definition}\label{def:Om}
	Let $C=(M^*,e^*,\tau^*)\in\mC_T$, let $v_0$ be the root-vertex of $M^*$ and let $v$ be the head-vertex.
	\begin{compactitem}
		\item The image $\Om_a(C)$ (resp.\ $\Om_b(C)$) is obtained from $C$ by replacing $e^*$ by a new vertex $u$ incident to three new edges $e_1,e_2,e_3$ in clockwise order around $u$, with $e_1$ joining $u$ to $v$, and $e_2,e_3$ joining $u$ to $v_0$ (hence the edges $e_2,e_3$ are drawn as incident to buds in our figures). The edge $e_1$ is added to the tree $\tau^*$, and the edge $e_3$ (resp. $e_2$) becomes the new head-edge.
		\item In order to define $\Om_c(C)$ we consider the edges $e_\ell^*,e_r^*$ preceding and following, respectively, the head-edge $e^*$ in counterclockwise order around $v_0$. The image $\Om_c(C)$ is only defined if the edges $e_\ell^*,e_r^*$ are both distinct from the root-edge of $M^*$. 
		In that case, we consider the non-root endpoints $v_\ell,v_r$ of $e_\ell^*,e_r^*$.  Since $\tau^*$ is in $\DFS_{M^*}^{e^*}$ the vertices $v_\ell,v_r$ are both ancestors of the head vertex $v$, so one is an ancestor of the other. If $v_\ell$ is an ancestor of $v_r$, then $\Om_c(C)$ is obtained by deleting $e^*$ and $e_\ell^*$ and replacing them by an edge between $v$ and $v_\ell$, while $e_r^*$ becomes the new head-edge (see Figure~\ref{fig:bij-rules-cubic}). 
		If $v_r$ is an ancestor of $v_\ell$, then $\Om_c(C)$ is obtained by deleting $e^*$ and $e_r^*$ and replacing them by an edge between $v$ and $v_r$, while $e_\ell^*$ becomes the new head-edge. 
	\end{compactitem}
\end{definition}

The following result is proved in \cite{bernardi-dfs-bijection}.

\begin{thm}[\cite{bernardi-dfs-bijection}]\label{thm:bij-original}
	The mapping $\Om$ is a bijection between $\mK$ and $\mC_T$.
\end{thm}

\fig{width=\linewidth}{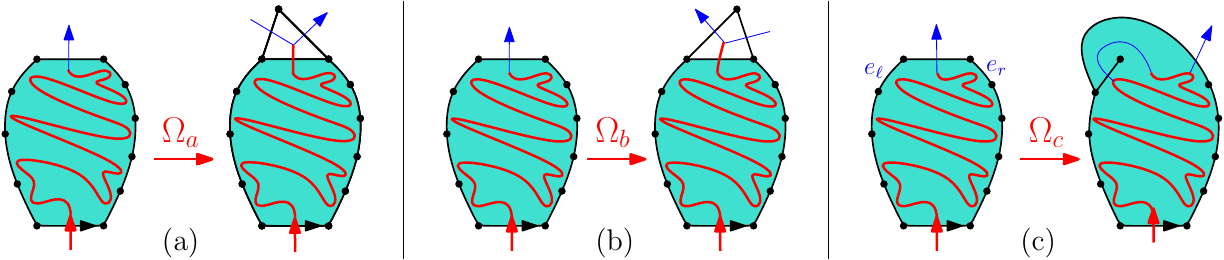}{The mapping $\Om_a,\Om_b,\Om_c$ presented in the dual setting (that is, in terms of the near-triangulation). Here, instead of representing $(M^*,e^*,\tau^*)\in\mC_T$ we represent the near-triangulation $M$, the edge $e^*$ (blue arrow) and the path of $\tau^*$ from the root-vertex $v_0$ to the head-vertex (this is the percolation path in the associated percolation configuration $\sigma=\Lambda_M^{e^*}(\tau^*)$).}

In order to make the relation between the mappings $\Om$ and $\Phi$ clearer, we represent the mapping $\Om_a,\Om_b,\Om_c$ in the dual setting in Figure~\ref{fig:bij-rules-triang-tree}. It is easy to see, by comparing Figures~\ref{fig:bij-rules-triang} and~\ref{fig:bij-rules-triang-tree} (and using Theorem \ref{thm:DFS-to-perco} to identify the percolation path of $(M,\sigma)$ with the path of the DFS tree $\tau^*$ from the root-vertex to the head-vertex), that the composition of the bijection $\Om$ with the mapping $\Lambda_M^{e^*}$ of Definition~\ref{def:Perc-and-DFS-sets-e}, is equal to the mapping $\Phi$. More precisely, for 
a walk $w\in\mK$, if $\Om(w)=(M^*,e^*,\tau^*)$, then $\Phi(w)=(M,\sigma)$, where $\sigma=\Lambda_M^{e^*}(\tau^*)$. By Theorem~\ref{thm:bij-original} and Corollary~\ref{cor:DFS-site-perco}, this shows that $\Phi$ is a bijection between $\mK$ and $\mT_P$. This concludes the proof of Theorem~\ref{thm:bij-OB}, as the other statements in this theorem are clear from the definition of $\Phi$.

\subsection{Proofs for Section~\ref{subsec:bij-chordal}: extending $\Phi$ to $\bmK$}\label{sec:appendix-bijection-extended}
Here we prove Theorem~\ref{thm:bij-extend} starting from Theorem~\ref{thm:bij-OB}.

\begin{proof}[Proof of Theorem~\ref{thm:bij-extend}]
	We first give an alternative description of the mapping $\bPhi$ in terms of $\Phi$.

	Let $w=w_1w_2\ldots w_n$ in $\bmK\setminus \mK$, and let $(M,\sigma)=\bPhi(w)$. We now describe an alternative construction of $(M,\sigma)$, which is represented in Figure~\ref{fig:example-bPhi-proof}. Let $i$ and $j$ be the number of unmatched $c$-steps in $w$ of type $a$ and $b$ respectively. Let $a^ib^j$ be the word made of $i$ consecutive $a$-steps, followed by $j$ consecutive $b$-steps, and let $w':=a^ib^jw$. 
	Clearly $w':=a^ib^jw$ is in $\mK$, so we can apply $\Phi$ to $w'$. Let $(M_0',\sigma_0')=\Phi(a^ib^j)$, and let $(M',\sigma')=\Phi(w')=\phi_{w_n}\circ\cdots\circ\phi_{w_1}(M_0',\sigma_0')$. 
	Let $(M'',\sigma'')$ be obtained from $(M',\sigma')$ by deleting the root-edge and all the edges corresponding to inner edges of $M_0'$, and setting the new root-edge of $M'$ to be the top-edge of $M_0'$ (oriented from the top-left vertex of $M_0'$ to the top-right vertex of $M_0'$). This operation is illustrated in Figure~\ref{fig:example-bPhi-proof} (right-column). The outer edges of $M''$ which correspond to outer edges of $M_0'$ are marked as inactive, and the others are marked as active.

	\fig{width=.7\linewidth}{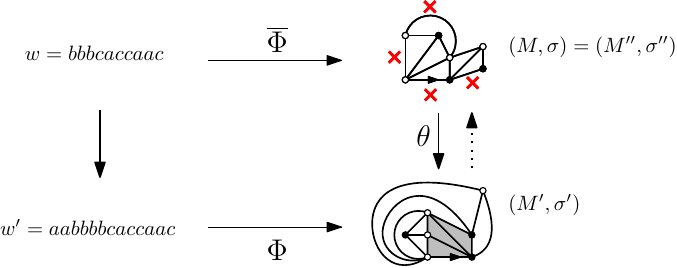}{The mapping $\bPhi$ in terms of the bijection $\Phi$. Here $w=bbbcaccaac$ has $i=2$ (resp. $j=1$) $c$-steps without matching $a$-step (resp. $b$-step). The pair $(M,\sigma)=\bPhi(w)$ was computed in Figure~\ref{fig:example-bPhi-proof}. The first $i+j=3$ triangles on the percolation path of $(M',\sigma')=\Phi(w')$ are indicated in gray.}
	
	We claim that $(M'',\sigma'')=(M,\sigma)$. In order to show this, let us define a mapping $\theta$ on $\bmTP$ (see Figure~\ref{fig:example-bPhi-proof}). Let $(\bar M,\bar \sigma)\in \bmTP$ with $i$ inactive left edges and $j$ inactive right edges. We define $\theta(\bar M,\bar\sigma)$ as the percolated near-triangulation obtained by gluing $(\bar M, \bar\sigma)$ and $(M_0',\sigma_0')=\Phi(a^ib^j)$ as follows: the root-edge of $(\bar M,\bar \sigma)$ is glued to the top-edge of $(M_0',\sigma_0')$, and then the inactive left and right edges of $(\bar M,\bar \sigma)$ are glued to the left and right edges of $(M_0',\sigma_0')$. 
	It is easy to see that $\theta(\bar M,\bar \sigma)$ is in $\mT_P$, and that the inner triangles of $(M_0',\sigma_0')$, are the first $i+j$ triangles on the percolation path of $\theta(\bar M,\bar \sigma)$.

	In order to show $(M'',\sigma'')=(M,\sigma)$, it suffices to show that $\theta(M,\sigma)=\theta(M'',\sigma'')$ (since $\theta$ can be inverted by removing the first $i+j$ triangles on the percolation path). Moreover, it is clear from the definition that $\theta(M'',\sigma'')=(M',\sigma')=\Phi(w')$. Hence it remains to show that $\theta(M,\sigma)=(M',\sigma')$.
	Now, it is easy to see from the definitions that the processes for constructing $(M,\sigma)=\bphi_{w_n}\circ\cdots\circ\bphi_{w_1}(M_0,\sigma_0)$ and $(M',\sigma')=\phi_{w_n}\circ\cdots\circ\phi_{w_1}(M_0',\sigma_0')$ are almost identical. The only difference is when treating the $c$-steps of $w$ without matching $a$-step or $b$-step. For a $c$-step without matching $a$-step (resp. $b$-step), the step $\bphi_c$ in the construction of $(M,\sigma)$ makes the top-edge an inactive left edge (right-edge), while the step $\phi_c$ in the construction of $(M',\sigma')$ glues the top-edge to a left edge (resp. right edge) of $(M_0',\sigma_0')$. It is easy to see that all of these extra ``edge gluings'' occurring in the construction of $(M',\sigma')$ can be delayed to the end of the construction, without affecting the final result. Thus $(M',\sigma')=\theta(M,\sigma)$. Hence $(M'',\sigma'')=(M,\sigma)$.

	We now argue that $\bPhi$ is a bijection. We first show injectivity. Suppose $w$ and $\tilde{w}$ are walks in $\bmK$ such that $\bPhi(w)=\bPhi(\tilde{w})$. The numbers $i$ and $j$ of unmatched $c$-steps of type $a$ and $b$ are the same in $w$ and $\tilde{w}$ since they are equal to the number of inactive left and right edges of $\bPhi(w)=\bPhi(\tilde{w})$. Moreover, the process described above for constructing $(M'',\sigma'')=\bPhi(w)$ from $(M',\sigma')=\Phi(a^ib^jw)$ is clearly injective (with inverse $\theta$). Hence $\Phi(a^ib^jw)=\Phi(a^ib^j\tilde{w})$. Since $\Phi$ is injective by Theorem~\ref{thm:bij-OB}, this implies $w=\tilde{w}$.
	
	We now show surjectivity. Let $(M,\sigma)\in\bmTP$, and let $i$ and $j$ be the number of inactive left and right edges of $(M,\sigma)$ respectively. We want to find a preimage of $(M,\sigma)$ by $\bPhi$. Let $(M',\sigma')=\theta(M,\sigma)$.
	Since $(M',\sigma')\in\mT_P$ there exists $w'\in \mK$ such that $\Phi(w')=(M',\sigma')$. 
	We also claim that $w'$ is of the form $a^ib^jw$ for some walk $w$ in $\bmK$. This property can be deduced from the fact that the percolation path of $(M',\sigma')$ start by $i+j$ triangles having $i+j+2$ distinct vertices, with the first $i$ triangles having 2 white vertices and the $j$ subsequent vertices having 2 black vertices. The ``distinct'' condition above is the key: it guarantees that the beginning of the percolation path of $(M',\sigma')$ are the first triangles that have been created during the construction of $\Phi(w')$, and remains part of the percolation path at every successive step of the construction. 
	Furthermore, it is easy to see that $w$ has $i$ unmatched $c$-steps of type $a$ and $j$ unmatched $c$-steps of type $b$ (indeed the $i+j$ first steps of $w'=a^ib^jw$ are matched since the $i+j$ first triangles on the percolation path of $(M',\sigma')$ are not incident to the left or right edges). Hence it is clear (from the alternative description of $\bPhi$ given above) that $\bPhi(w)=(M,\sigma)$, which proves the surjectivity of $\bPhi$.

	Lastly, the stated correspondences between the steps of $w$ and the vertices, faces and edges of $\bPhi(w)$ are clear from the definitions.
\end{proof}

\subsection{Proofs for Section~\ref{subsec:future-past}: alternative description of $\bPhi$ and future/past decomposition} \label{sec:appendix-future-past}
We now prove Proposition~\ref{prop:alternative-description} and Proposition~\ref{prop:future-past}.

\begin{proof}[Proof of Proposition~\ref{prop:alternative-description}]
	Let $w=w_1\ldots w_n\in\bmK$. We claim that the site-percolated maps $T_w$ and $\bPhi(w)$ are obtained by gluing some triangles in the exact same way: the only difference is that the gluings corresponding to the $c$-steps are made ``one at a time'' in the definition of $\bPhi(w)$, while they are all ``delayed until the end'' in the definition of $T_w$. Let us justify this statement briefly. 
	Let $k\in [n]$ and let $(M_k,\sigma_k)=\bPhi(w_1w_2\ldots w_k)$. If $w_{k+1}=a$ (resp. $b$), then it is clear that applying $\bphi_{w_{k+1}}$ to $(M_k,\sigma_k)$ is the same as gluing the root-edge of the brick $T_{a}$ (resp. $T_{b}$) to the top-edge of $(M_k,\sigma_k)$.
	
	Consider now the case where $w_{k+1}=c$ is a matched $c$-step of type $a$. Let $t_\ell$ (resp. $t_r$) be the last triangle on the percolation path of $(M_k,\sigma_k)$ incident to a left (resp. right) edge. Note first that the triangles $t_\ell$ are $t_r$ are the triangles created when treating the $a$-step and $b$-step matched to the current $c$-step. Hence, the left and right edges incident to $t_r$ and $t_l$ are active and the last triangle $t$ on the percolation path is equal to $t_r$ (this situation is represented in Figure~\ref{fig:bij-rules-triang}). Moreover, it is easy to check that the gluing and recoloring performed by applying $\bphi_c$ to $(M_k,\sigma_k)$ in this case is the same as gluing the brick $T_\ta$ to the top-edge $e$ of $(M_k,\sigma_k)$ (that is, gluing the root-edge of $T_{\ta}$ to $e$, and gluing the left and right sides of the brick $T_\ta$ to the sides of $t_\ell$ and $t_r$ respectively, all that while keeping the colors of the brick $T_\ta$). Symmetrically, if $w_{k+1}$ is a matched $c$-step of type $b$, then applying $\bphi_c$ to $(M_k,\sigma_k)$ is the same as gluing the brick $T_\tb$ to the top-edge of $(M_k,\sigma_k)$.
	
	Next consider the case where $w_{k+1}$ is an unmatched $c$-step of type $a$. In this case, $(M_k,\sigma_k)$ has no active left edge. Moreover the last triangle $t_r$ on the percolation path of $(M_k,\sigma_k)$ incident to a right edge was created when treating the $b$-step matched to the current $c$-step. Hence it is easy to check that the gluing and recoloring performed by applying $\bphi_c$ to $(M_k,\sigma_k)$ in this case is the same as gluing the brick $T_\tb$ to the top-edge $e$ of $(M_k,\sigma_k)$ (that is, gluing the root-edge of $T_{\ta}$ to $e$, and gluing the right side of the brick $T_\ta$ to the side of $t_r$, all that while keeping the colors of the brick $T_\ta$). The case where $w_{k+1}$ is an unmatched $c$-step of type $b$ is symmetric.
\end{proof}

\begin{proof}[Proof of Proposition~\ref{prop:future-past}] By Proposition~\ref{prop:alternative-description}, we can think of $\Phi(w)$ as obtained by gluing some bricks $T_a,T_b,T_\ta,T_\tb$. 
	Now, we know that the past near-triangulation $(P,\al)$ is obtained by gluing the bricks corresponding to steps in $u$ together, and it is easy to see that the future near-triangulation $(Q,\be)$ is obtained by gluing the bricks corresponding to steps in $v$ together. Lastly, the gluing of $(P,\al)$ with $(Q,\be)$ is easily seen to correspond to performing the remaining gluings between the bricks corresponding to $u$ and the bricks corresponding to $v$.
\end{proof}

\subsection{Proofs for Section~\ref{sec:spine}: the percolation path as a function of the Kreweras walk.}\label{sec:appendix-spine}
In this subsection we prove Theorem~\ref{thm:LR}. 
Let $w\in\bmK$. We consider the description of the percolated near-triangulation $(M,\sigma)=\bPhi(w)$ given by Proposition~\ref{prop:alternative-description}. Let $\tw$, $\tT_w$, and $T_w=\bPhi(w)$ be as in Section~\ref{subsec:future-past}. By definition, $T_w$ is obtained from $\tT_w$ by gluing the matching pairs of opening and closing sides of $\tT_w$. 
We now think of performing these gluings in two stages: first we glue the matching sides of $\tT_w$ corresponding to steps inside a common cone-excursion, and then we glue the other matching sides of $\tT_w$. We demote by $T_w^\circ$ the site-percolated triangulation obtained after the first stage. We will now describe $T_w^\circ$ and its relation to $\LR(M,\sigma)$.

Let $w_j$ be a spine step of $w$ which is close-matched, and let $w'$ be the cone excursion starting at $w_j$. Let $T^\circ_{w'}$ be the sub-triangulation of $T_w^\circ$ corresponding to the cone-excursion $w'$. 
By definition of $T^\circ_w$, all the matching sides of $T^\circ_{w'}$ are glued, so $T_{w'}^\circ=\bPhi(w')$. 
We now describe the sub-triangulation $T_{w'}^\circ=\bPhi(w')$. For concreteness, let us assume that $w_j$ is an $a$-step (the case of a $b$-step being symmetric).
Let $k$ be the height of the cone excursion $w'$. We claim that $T_{w'}^\circ$ has $k+2$ outer vertices of which exactly one is white. Moreover all the unicolor outer edges are inactive (equivalently, they correspond to closing edges of $T_w^\circ$). Lastly, the outer white vertex is incident to a unique inner triangle which corresponds to the brick $T_{w_j}$, and we call it \emph{spine-triangle} of $T_{w'}^\circ$. This situation is represented in Figure~\ref{fig:cone-excursion-image}.
These claims are direct consequences of Remark~\ref{rmk8}. Indeed, we can write $w'=aw''c$ where $w''$ is a walk in $\smK$ starting at $(0,0)$ and ending at $(1-k,0)$. By Remark~\ref{rmk8}, the associated percolated triangulation $T_{w''}=\bPhi(w'')$ has $k+1$ outer-vertices of which one is white. Hence the claims follow upon observing that $T_{w'}^\circ$ is obtained from $T_{w''}$ by gluing a brick $T_a$ (to the root-edge) and a brick $T_{\tb}$ (to the top-edge). We call \emph{shell} of $T_{w'}^\circ$ the map obtained from $T_{w'}^\circ$ by erasing the inner vertices and the inner edges, except for the edges of the spine triangle, 
and if there were no inner vertices replacing the unicolor edge  of the spine triangle by a double edge. In other words, the shell of $T_{w'}^\circ$ is made of the spine triangle attached along its unicolor edge to a simple unicolor cycle of length $k+1$; see Figure~\ref{fig:cone-excursion-image}.

\fig{width=.6\linewidth}{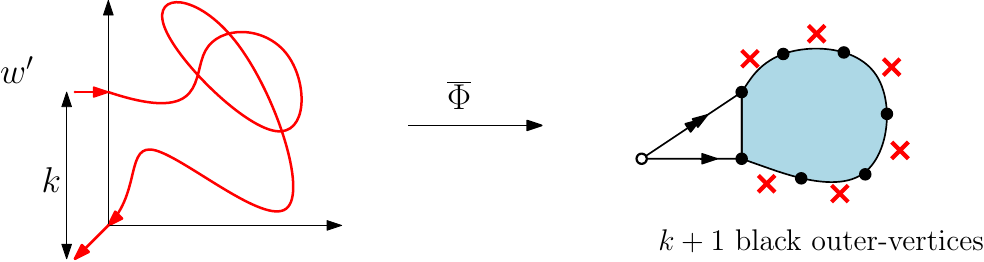}{A cone excursion $w'$ and its image $T_{w'}^\circ=\bPhi(w')$.} 

Now we can think of $T_w^\circ$ as obtained by gluing several bricks and sub-triangulations corresponding to the spine steps of $w$: one brick $T_a$ for each unmatched or far-matched spine $a$-step, one brick $T_b$ for each unmatched or far-matched spine $b$-step of $w$, and one sub-triangulation $T_{w'}^\circ=\bPhi(w')$ for each cone excursion $w'$ starting with a (close-matched) spine step of $w$. Hence, these sub-triangulations of $T_w^\circ$ are in one-to-one correspondences with the steps of $\hpi(w)$: the steps $a$ (resp. $b$) in $\hpi(w)$ correspond to the bricks $T_a$ (resp. $T_b$), while the letters $\ba_k$ (resp. $\bb_k$) correspond to sub-triangulations made of a spine triangle attached to a percolated near-triangulation with $k+1$ white (resp. black) outer vertices. This is represented in Figure~\ref{fig:proof-spine2}.

\fig{width=\linewidth}{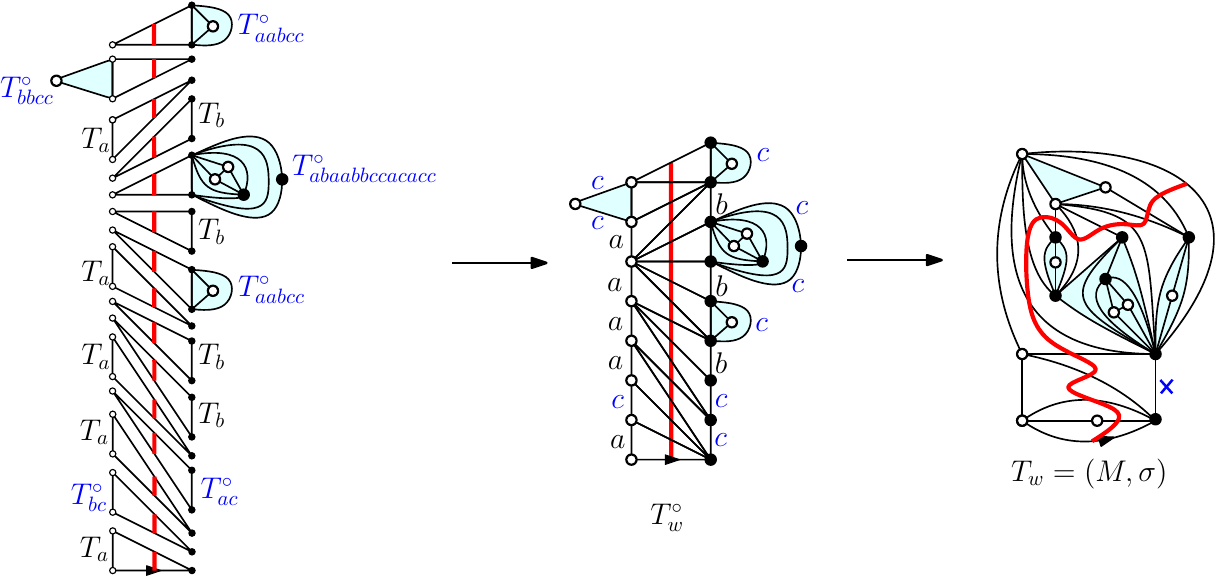}{The percolated near-triangulations $T_w^\circ$ and $T_w=\bPhi(w)$. Here, $w=a bc a a ca b b a aabcc b abaabbccacacc b a bbcc aabcc$ so that $\hpi(w)=a\ba_1a\bb_1abba\bb_1b\bb_2ba\ba_2\bb_1$ as in Figure~\ref{fig:LR-theorem2}.
}

From the above description, we see that the set of triangles on the percolation path of $T_w^\circ$ is made of the bricks $T_a$ and $T_b$ corresponding to the $a$-steps and $b$-steps of $\hpi(w)$, and of the spine triangles of the near-triangulations corresponding to the steps $\ba_k$ and $\bb_k$ of $\hpi(w)$. More precisely, the sub-triangulation of $T_w^\circ$ made of the triangles on the percolation path is $\Phi(\hpi_s(w))$ (see Definition~\ref{def:hpi} of $\hpi_s(w)$). Observe also that all the vertices on the left side (resp. right side) of $T_w^\circ$ are white (resp. black). Hence the vertices will not change color when gluing the matching sides of $T_w^\circ$ in order to obtain $T_w=\bPhi(w)$. Hence the triangles on the percolation path of $T_w^\circ$ and $T_w$ are the same. 
Also, for any cone excursion $w'$ corresponding to a letter $\ba_k$ or $\bb_k$ of $\hpi(w)$, the inner vertices of $T_{w'}^\circ$ will not be incident to the triangles on percolation path of $T_w$. So the map $\LR(T_w)$ is obtained by 
\begin{compactenum}
	\item replacing each sub-triangulation $T_{w'}^\circ=\bPhi(w')$ of $T_w^\circ$ corresponding to letters $\ba_k$ or $\bb_k$ of $\hpi(w)$ by its shell (which is a spine triangle attached to a cycle of length $k+1$),
	\item performing the gluing of the matching sides of $T_w^\circ$: the side edges of the shells are closing in $T_w^\circ$, while the side edges of the bricks $T_a,T_b$ are opening in $T_w^\circ$ (and the opening and closing sides are matched like a parenthesis system).
\end{compactenum}

It remains to prove that the counterclockwise code of the white cluster of $\LR(T_w)=\LR(\bPhi(w))$ is $\hpi_\ell(w)$  and the clockwise code of the black cluster of $\LR(T_w)=\LR(\bPhi(w))$ is $\hpi_r(w)$.
In order to do this proof, we can assume that $w$ is in $\mK$. Indeed, consider $w\in\bmK$ having $i'$ and $j'$ unmatched $c$-steps of type $a$ and $b$ respectively. Then it is clear that from the relation between $\bPhi(w)$ and $\Phi(a^{i'}b^{j'}w)$ established in the proof of Theorem~\ref{thm:bij-extend}, that the white (resp. black) cluster of $\LR(\bPhi(w))$ is equal to the white (resp. black) cluster of $\LR(\Phi(a^{i'}b^{j'}w))$.
So let us now assume that $w\in \mK$.
Note first that the letters $\ba_k$ (resp. $\bb_k$) of $\hpi(w)$ correspond to inner faces of degree $k+1$ of the white (resp. black) cluster of $\LR(T_w)$. Furthermore, we claim that the counterclockwise code of the white cluster $C_\ell$ of $\LR(T_w)$ is $\hpi_\ell(w)$ (see Definition~\ref{def:hpi} of $\hpi_\ell(w)$). Indeed, as we have established, the walk $\hpi(w)$ describes the succession of triangles along the percolation path of $T_w$. Hence the steps of $\hpi(w)$ which are in $\{a\}\cup\{\ba_k,k\geq1\}$ (hence correspond to triangles with a unicolor white edge) corresponds to the edges of $C_\ell$. And, moreover, the steps $\ba_k$ in $\hpi(w)$ correspond to the last edge along some face of degree $k+1$ of $C_\ell$. Similarly, the clockwise code of the black cluster $C_r$ of $\LR(T_w)$ is $\hpi_r(w)$ (see Definition~\ref{def:hpi} of $\hpi_r(w)$). This completes the proof of Theorem~\ref{thm:LR}.

\subsection{Proofs for Section~\ref{sec:exploration-tree-from-walk}: the exploration tree as a function of the Kreweras walk.}\label{sec:appendix-DFS-from-walk}
In this subsection, we prove Theorem~\ref{thm:DFS} and Proposition~\ref{prop32}. The proofs of these results are actually easier to state in terms of the bijection $\Om$ defined in Section~\ref{sec:appendix-bijection}. But we first need to extend the mapping $\Om$ (which is defined on $\mK$) into a mapping $\bOm$ defined on $\bmK$.

The bijection $\bOm$ is illustrated in Figure~\ref{fig:omega-bar}. Roughly speaking, it is obtained by translating the bijection $\bPhi$ in terms of cubic maps with a marked DFS tree. We now give a direct definition. 
Let $\bmC_T$ be the set of triples $(M^*,e^*,\tau^*)$, where 
\begin{compactitem}
	\item $M^*$ is the dual of a near-triangulation $M\in\mT$, 
	\item the edges of $M^*$ incident to the root-vertex $v_0$ are marked as either \emph{active} or \emph{inactive}. The root-edge is inactive and the active edges are consecutive around $v_0$ 
	\item $e^*$ is an active edge of $M^*$ (incident to $v_0$). The edge $e^*$ is called the \emph{head-edge} and its non-root endpoint is called the \emph{head-vertex} (if $e^*$ is a self-loop, the head-vertex is $v_0$).
	\item $\tau^*$ is a DFS tree in $\DFS_{M^*}$ such that all the non-root vertices of $M^*$ incident to active edges are ancestors of the head-vertex.
\end{compactitem}

We now define the mapping $\bOm$ from $\bmK$ to $\bmC_T$. An instance of $\bOm$ is represented in Figure~\ref{fig:omega-bar}(b).
Let $M_0^*$ be the rooted map with one vertex and one edge $e_0^*$ (which is self-loop which is both the root-edge and the head-edge), and let $\tau_0^*$ be the unique spanning tree of $M_0^*$. For $w=w_1\cdots w_n\in\bmK$, the image $\bOm(w)$ is defined as the triple $(M^*,e_0^*,\tau^*)=\bOm_{w_n}\circ\cdots \circ \bOm_{w_1}(M_0^*,e_0^*,\tau_0^*)$, where the mappings $\bOm_a,\bOm_b,\bOm_c$ are defined as follows. The mapping $\bOm_a$ and $\bOm_b$ are equal to $\Om_a$ and $\Om_b$ (see Definition~\ref{def:Om} and Figure~\ref{fig:bij-rules-cubic}) with the newly created edges incident to $v_0$ marked as active.
We now consider $C\in\bmC_T$ and want to define $\bOm_c(C)$. Let $e_\ell^*,e_r^*$ be respectively the edges preceding and following the head-edge $e^*$ in counterclockwise order around $v_0$. The image $\bOm_c(C)$ is only defined if at least one of the edges $e_\ell^*,e_r^*$ is active.
If the edges $e_\ell^*,e_r^*$ are both active, we define $\bOm_c(C)=\Om_c(C)$ (see Definition~\ref{def:Om} and Figure~\ref{fig:bij-rules-cubic}). 
If $e_\ell^*$ (resp. $e_r^*$) is inactive, then $\bOm_c(C)$ is obtained by deleting $e^*$ and creating an edge between the head-vertex $v$ and $v_0$, 
while $e_r^*$ (resp. $e_l^*$) becomes the new head-edge (see Figure~\ref{fig:omega-bar}(a)). 

It is immediate from Theorem~\ref{thm:bij-extend} and Theorem~\ref{thm:DFS-to-perco} that $\bOm$ is a bijection between $\bmK$ and $\bmC_T$. 

\fig{width=\linewidth}{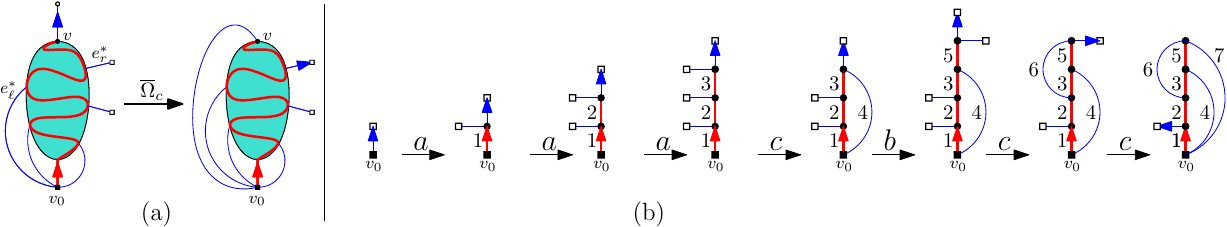}{The bijection $\bOm$. In this figure, a special convention was used to represent the edges of $M^*$ incident to the root-vertex $v_0$. Indeed only the inactive edges of $M^*$ are represented as incident to $v_0$; the active edges are instead represented as edges incident to special vertices of degree 1 called \emph{buds}. So the map $M^*$ is the map which would be obtained by identifying all the buds to $v_0$. (a) The mapping $\bOm_c$ when the edge $e_\ell^*$ is inactive (this represents the treatment of an unmatched $c$-step of type $a$). (b) The bijection $\bOm$ applied to the word $w=aaacbcc\in\bmK$. The edges are labeled by the order in which they are \emph{created} (that is to say, the order in which they become in-edges). 
}

\begin{proof}[Proof of Theorem~\ref{thm:DFS}]
	Let $w\in\bmK$, let $(M,\sigma)=\bPhi(w)$ and let $\tau^*=\dfs(M,\sigma)$. Let $w_{t_1},\ldots,w_{t_n}$ be all the $a$-steps and $b$-steps of $w$, and let $w^{(1)},\ldots,w^{(n)}$ be the prefixes of $w$ ending with the steps $w_{t_1},\ldots,w_{t_n}$ respectively. 

	Recall from property (ii) of Theorem~\ref{thm:bij-extend} that the $a$-steps and $b$-steps of $w$ correspond via $\etavf$ (and duality) to the non-root vertices of $M^*$. Let $v_1,\ldots,v_n$ be the vertices of $M^*$ corresponding to these steps. Note that the root-edge of $M^*$, which is in $\tau^*$, joins the root-vertex $v_0$ to $v_1$.
	By property (i) of Theorem~\ref{thm:LR}, for all $k\in[n]$ the length $h_k$ of $\pi(w^{(k)})$ is equal to the number of non-root vertices on the percolation path of $\bPhi(w^{(k)})$, which is a path $P_k$ of $M^*$ going from $v_1$ to $v_k$. Moreover, by property (i) of Theorem~\ref{thm:DFS-to-perco} the percolation path $P_k$ is contained in $\tau^*$ (since it is contained in $\dfs(\bPhi(w^{(i)}))$). Hence $h_k$ is the height of the vertex $v_k$ in $\tau^*$. 
	
	Thus $h_1,\ldots,h_n$ represent the respective height of the non-root vertices $v_1,\ldots,v_n$ of $\tau^*$. Recall from Definition~\ref{def:DFS-w} that the height-code of $\dfs(w)$ is $(0,h_1,\ldots,h_n)$. So in order to show that the rooted trees underlying $\tau^*$ and $\dfs(w)$ are equal, it suffices to show that the order $v_0,v_1,\cdots, v_n$ of the vertices of $M^*$ corresponds to a \emph{pre-order} of the vertices of the tree $\tau^*$, that is, an order for which any vertex $v$ precedes all its descendants in $\tau^*$, and all the descendants of $v$ precedes all the non-descendants of $v$ appearing after $v$. This property is easy to check from the definition of $\bOm$.
	
	This completes the proof of Theorem~\ref{thm:DFS} because the additional properties (i-ii) are direct consequences of the definitions. 
\end{proof}

\begin{proof}[Proof of Proposition~\ref{prop32}]
	Let $w\in\bmK$, let $(M,\si)=\bPhi(w)$, and let $(M^*,e^*,\tau^*)=\bOm(w)$. We denote by $v_0$ the root-vertex of $M^*$. 
	Recall that the functions $\etae$ and $\etavf$ indicate the order in which the in-edges and inner triangles of $M$ are \emph{created} during the bijection $\bPhi$. Equivalently, $\etae$ and $\etavf$ indicate the order in which the in-edges and non-root vertices of $M^*$ are \emph{created} during the related bijection $\bOm$. Here our convention is that we consider the in-edges of $M^*$ to be \emph{created} during $\bOm$ when they become in-edges (while the active edges of $M^*$, including the head-edge, are considered to be created after all the other edges); see Figure~\ref{fig:omega-bar}(b).
	
	It is clear from the definition of $\bOm$ that the vertices of $M^*$ are created according to a DFS of $M^*$ with associated DFS tree $\tau^*=\dfs(M,\si)$. We call this the \emph{creation-DFS} of $M^*$. Proving Claim~(i) of Proposition~\ref{prop32} amounts to showing that the creation-DFS and space-filling exploration of $M^*$ are equal. In other words, we need to prove that the creation-DFS follows rules (ii'-1) and (ii'-2) of Definition~\ref{def:space-filling}.
	
	Since we already know that the creation-DFS and space-filling exploration of $M^*$ give the same DFS tree $\tau^*$, we can unambiguously talk about the \emph{forward face}, and \emph{forward edges} of any non-root vertex of $M^*$. Let $u$ be a non-root vertex of $M^*$ and let $f$ be its forward face. 
	We denote by $v_f$ the first vertex of $M^*$ incident to $f$ visited during the creation-DFS. Note that $v_f$ is the ancestor in $\tau^*$ of all the other vertices of $M^*$ incident to $f$ (hence $v_f$ is also the first vertex incident to $f$ visited during the space-filling exploration). 
	
	Suppose first that $u=v_f$. By Claim~\ref{claim:color-matters}, only one forward edge of $u$ is in $\tau^*$: its forward left edge if $f$ is black, and its forward right edge otherwise. This shows that the creation-DFS follows rule (ii'-1) of Definition~\ref{def:space-filling}.
	
	It remains to prove that the creation-DFS follows the rule (ii'-2) of Definition~\ref{def:space-filling}. We first need to reformulate this rule.
	We consider the set $\mP$ all the paths of $\tau^*$ from a vertex incident to $f$ to another vertex incident to $f$ (and not incident to $f$ in between). Let $P\in\mP$. We say that a vertex or edge $x$ of $M^*\setminus P$ is \emph{enclosed by} $P$ if $P$ together with the interior of $f$ separates $x$ from the head-edge. Observe that all the vertices of $M^*$ enclosed by $P$ have an ancestor on $P$ (since $P$ separate them from $v_0$). Moreover one of the endpoints of $P$ is an ancestor of all the other vertices on $P$ (otherwise, the common ancestor $x$ on $P$ would have 2 children $x_1,x_2$ on $P$ and there would be an edge enclosed by $P$ between some descendant of $x_1$ and some descendant of $x_2$, which is impossible by Claim~\ref{claim:characterization-DFStree}). We now denote by $v_P^-$ and $v_P^+$ the endpoints of $P$ with $v_P^-$ the ancestor of $v_P^+$. We denote by $e_P^-$ the forward edge at $v_P^-$ enclosed by $P$, and we denote by $e_P^+$ the forward edge at $v_P^+$ not enclosed by $P$. The situation is illustrated in Figure~\ref{fig:order-DFS}(a). By the preceding discussion, $v_P^-$ is the ancestor of all the vertices on $P$ and enclosed by $P$.  Moreover, all the vertices enclosed by $P$ and incident to $f$ are descendants of $v_P^+$ (otherwise Claim~\ref{claim:characterization-DFStree} would be violated for some edge incident to $f$ and enclosed by $P$). Hence $e_P^-\notin\tau^*$. We will now prove the following claim.
	
	\begin{claim}\label{claimX}
		For any path $P\in\mP$, during the application of the bijection $\bOm$ the vertices and edges enclosed by $P$ are all created before the edge $e_P^+$.
	\end{claim}
	
	\fig{width=\linewidth}{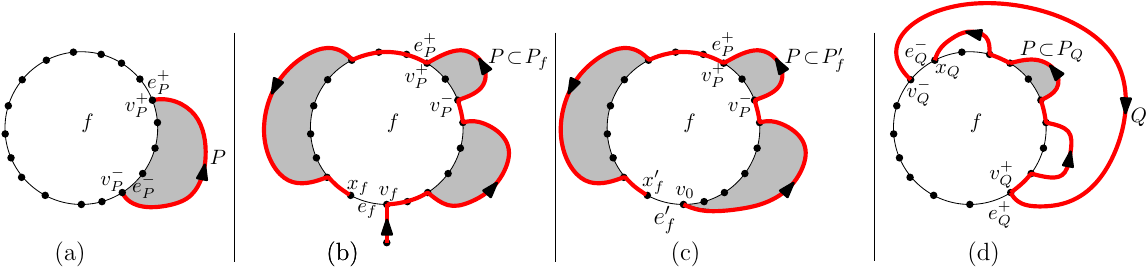}{Notation for the proof of Claim~\ref{claimX}. The paths $P\in\mP$ are indicated in bold and oriented from parents to children.   (a) A path $P\in\mP$, its endpoints $v_P^-,v_P^+$, and the enclosed region (dashed). (b) The case where $P$ is not enclosed by any path in $\mP$, and $f$ is a white face not incident to $v_0$. (c) The case where $P$ is not enclosed by any path in $\mP$, and $f$ is a white face incident to $v_0$. (d) The case where $P$ is enclosed by another path $Q\in\mP$.}
	
	Note that Claim~\ref{claimX} implies that the creation-DFS follows rule (ii'-2) of Definition~\ref{def:space-filling}. 
	Indeed, Claim~\ref{claimX} shows that rule (ii'-2) applies when the chip position is $u=v_P^+$. We now consider a path $P\in \mP$, and prove that Claim~\ref{claimX} holds for $P$ by induction on the number of paths in $\mP$ enclosing $P$.
	
	Let us first suppose that $P$ is not enclosed by any path in $\mP$, and that $f$ is not incident to $v_0$ (equivalently, $v_f\neq v_0$). The situation is represented in Figure~\ref{fig:order-DFS}(b). Let $e_f$ be the edge incident to $v_f$ and to $f$ and which is not in $\tau^*$. Let $x_f$ be the endpoint of $e_f$ distinct from $v_f$ (note that $x_f$ is a descendant of $v_f$ in $\tau^*$). Consider the path $P_f$ of $\tau^*$ between $v_f$ and $x_f$. We say that a vertex or edge of $M^*\setminus P^*$ is \emph{enclosed} by $P_f$ if it is separated from the head-edge by the cycle $P_f\cup\{e_f\}$.
	Since $P$ is not enclosed by any path in $\mP$, the path $P$ must be part of the path $P_f$, and all the vertices and edges enclosed by $P$ are enclosed by $P_f$. Moreover, it is easy to see from the definition of $\bOm$ that all the vertices and edges enclosed by $P_f$ are created before $e_f$ during the application of the bijection $\bOm$ (since the cycle $P_f\cup\{e_f\}$ separates them from the head-edge). Now, we still need to prove that the vertices and edges enclosed by $P$ are created before the edge $e_P^+$. Note that $e_P^+$ is either equal to $e_f$ (in which case we are done) or on the path of $P_f$ between $v_P^+$ and $x_f$. Suppose now that $e_P^+\neq e_f$. Let $V_P$ (resp. $V_P'$) be the strict descendants of $v_P^+$ enclosed by $P$ (resp. not enclosed by $P$). We know that either all the vertices in $V_P'$ are created before all the vertices in $V_P$ or the converse is true (because the creation-DFS is a DFS). Now, since the vertices in $V_P$ are created before $e_f$ which is incident to $x_f\in V_P'$, we conclude that all the vertices in $V_P$ are created before all the vertices in $V_P'$. Hence, all the vertices in $V_P$ are created before $e_P^+$. As mentioned earlier, the set $V_P$ includes all the vertices incident to $f$ and enclosed by $P$. Hence the edges incident to $f$ and enclosed by $P$ are all created before $e_P^+$. This implies that all the vertices and edges enclosed by $P$ will be created before the edge $e_P^+$ (because $P$ together with the edges incident to $f$ separate the head-edge from all the vertices and edges enclosed by $P$). 
	
	Next we consider the case where $P$ is not enclosed by any path in $\mP$, and $f$ is incident to $v_0$. This case is treated almost exactly as the previous one, except we need to adjust our definitions slightly; see Figure~\ref{fig:order-DFS}(c). Let $e_2,e_3$ be edges preceding and following $f$ in clockwise order around $v_0$. If $f$ is white (resp. black) we set $e_f'=e_2$ (resp. $e_3$). Note that $e_f'$ is not the root-edge of $M^*$ hence is not in $\tau^*$. We denote by $x_f'$ the non-root endpoint of $e_f'$, and we denote by $P_f'$ the path in $\tau^*$ from $v_f$ to $x_f$. With this notation (illustrated in Figure~\ref{fig:order-DFS}(c)), this case is treated exactly as the previous one, upon replacing $x_f,e_f,P_f$ by $x_f',e_f',P_f'$; see Figure~\ref{fig:order-DFS}(c).
	
	Next suppose that $P$ is enclosed by another path in $\mP$. Consider the path $Q\in\mP$ enclosing $P$ and not enclosing any other path $Q'\in\mP$ enclosing $P$. The situation is represented in Figure~\ref{fig:order-DFS}(d). Let $x_Q$ be the endpoint of $e_Q^-$ distinct from $v_Q^-$. Recall that $e_Q^-\notin\tau^*$ and that $x_Q$ is a descendant of $v_Q^+$. We consider the path $P_Q$ of $\tau^*$ from $v_Q^+$ to $x_Q$. It is clear that the path $P$ must be part of the path $P_Q$. By the induction hypothesis, we also know that the edges enclosed by $Q$ are all created before $e_Q^+$. We contend that $e_Q^-$ is created just before $e_Q^+$. Indeed, by definition of $\bOm$, the edge created just before $e_Q^+$ must be incident to $f$ and join a descendant of $v_Q^+$ to an ancestor of $v_Q^+$; and $e_Q^-$ is the only such edge. From this point, we follow the same line of argument as above, with $v_Q^-,x_Q,e_Q^-,P_Q$ playing the role of $v_f,x_f,e_f,P_f$. This concludes the proof of Claim~\ref{claimX} and of Claim (i) of Proposition~\ref{prop32}. 
	
	It remains to prove Claim (ii) of Proposition~\ref{prop32}, which is a statement about the order in which the in-edges of $M^*$ are created during the bijection $\bOm$. We need to prove that this order of creation is equal to the order of treatment of the in-edges defined for the space-filling exploration of $M^*$. 
	We will prove that the creation-order and treatment-order coincide, and that moreover for any in-edge $e$ of $M^*$, the head-vertex $u$ at the time of the creation of $u$ is equal to the chip position at the time of the treatment of $e$. First observe that, from Claim (i) of Proposition~\ref{prop32}, we already know that this property holds when restricting these orders to the set of edges in the tree $\tau^*$. We now extend it to the set of all in-edges.
	
	Let $e$ be an in-edge of $M$ not in $\tau^*$. Let $\{u,v\}$ be the endpoints of $e$, with $u$ the ancestor of $v$ in $\tau^*$. It is clear from the definition of $\bOm$ that when $e$ is created, the head-vertex is at $v$. Now consider a time during $\bOm$ at which the head-vertex is $v$, and $e$ has not yet been created. 
	It is clear that if all the other edges incident to $v$ have already been created (that is to say, are in-edges), then the next step during the bijection $\bOm$ will be to create the edge $e$. Next we suppose that another edge $e'$ of $M^*$ incident to $v$ has not yet been created. Then the next step in the bijection $\bOm$ will be either to create $e$ or to create $e'$. 
	We consider the set $E'$ of edges of $M^*$ distinct from $e$ which have not yet been created. We want to show that $e$ is the next edge to be created if and only if there is a path in $E'$ between $v$ and an ancestor of $u$. Let us first suppose that $e$ is the next edge to be created. This means that after $e$ is created, $e'$ becomes the head-edge. If the head-vertex stays a descendant of $v$ until the end of the bijection $\bOm$, then there is a path between $v$ and $v_0$ in $E'$, namely the path of $\tau^*$ from $v$ to the head-vertex of $M^*$ followed by the head-edge of $M^*$. If the head vertex cease to be a descendant of $v$ during the bijection $\bOm$, then the edge $e''\in E'$ created just before this event joins a descendant $v'$ of $v$ to an ancestor of $u$. Hence the path made of the path of $\tau^*$ from $v$ to $v'$ followed by $e''$ is a path in $E'$ joining $v$ to an ancestor of $u$.
	Thus if $e$ is the next edge to be created, then there is a path in $E'$ between $v$ and an ancestor of $u$.
	Suppose conversely that there is a path $P$ in $E'$ between $v$ and an ancestor of $u$. 
	We want to prove that $e$ is the next edge to be created. Suppose by contradiction that this is not the case. Then all the descendants of $v$ and the incident edges will be created before $e$ is created. In particular, since the path $P\subseteq E'$ exists, there is an edge $e''\in E'$ between  a descendant of $v$ and an ancestor of $u$ which will be created before $e$ is created. 
	Thus, it is easy to see that just before $e$ is created, there is no active edge incident to the vertices on the path $Q$ of $\tau^*$ between $u$ and $v$ (roughly speaking, the edge $e''$ prevents the existence of such active edges on one side of $Q$, while the edge $e$ prevents the existence of such active edges on the other side of $Q$). This contradicts the rule of $\bOm_c$ when creating $e$. We reach a contradiction, hence if there is a path in $E'$ between $v$ and an ancestor of $u$, then $e$ is the next edge to be created.
	
	We have proved that the order of creation of the in-edges of $M^*$ follows the same rule as the order of treatment of the edges during the space-filling exploration. This proves Claim (ii) of Proposition~\ref{prop32}.
\end{proof}

\subsection{Proofs for Section~\ref{sec:dual-DFS}: dual of the DFS tree and pivotal points in terms of the walk}\label{sec:appendix-dual-dfs}
In this subsection we prove Theorem~\ref{thm:dfsdual} as well as Lemmas~\ref{lem:envelope-matching-is-loop-building} and~\ref{prop:piv}, and Claim~\ref{claim:disjoint-intervals}.


\begin{proof}[Proof of Theorem~\ref{thm:dfsdual}]
	We prove the result by induction on the length of the walk $w\in\bmK$. The result is clearly true when $w$ is the empty word. We now consider a non-empty walk $w=w'w_n$, where $w_n$ is the last step. Let $(M,\sigma)=\bPhi(w)$, and let $(M',\sigma')=\bPhi(w')$. By the induction hypothesis, we know that the trees $\tau'=\dfsdual(M',\sigma')$ and $\wttau'=\dfsdual(w')$ are isomorphic, with $\lambdav$ giving the isomorphism between vertices. We want to show that the trees $\tau=\dfsdual(M,\sigma)$ and $\wttau=\dfsdual(w)$ are isomorphic. 
	We consider the three cases $w_n=a$, $b$ or $c$, and check that the changes from $\tau'$ to $\tau$ are isomorphic to the changes from $\wttau'$ to $\wttau$. The changes from $\tau'$ to $\tau$ are easy to understand and are represented in Figure~\ref{fig:bij-rules-dfsdual2}. We now analyze the changes from $\wttau'$ to $\wttau$. Below we denote by $\wta$ and $\wtb$ the top-steps of $\wtw:=wab$, and by $\wta'$ and $\wtb'$ the top-steps of $\wtw':=w'ab$.
	
	\fig{width=\linewidth}{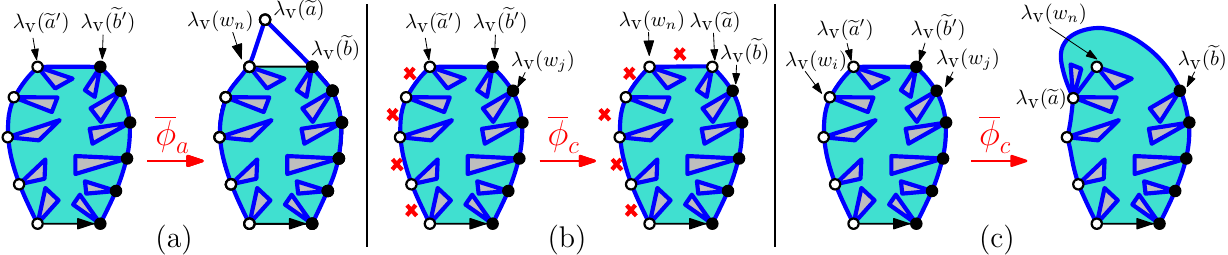}{The changes induced on the tree $\dfsdual(M,\sigma)$ by the mappings $\bphi_a$ and $\bphi_c$. (a) The changes from $\tau'$ to $\tau$ when applying $\bphi_a$. (b) The changes from $\tau'$ to $\tau$ when applying $\bphi_c$ in the case of an unmatched $c$-step of type $a$. (c) The changes from $\tau'$ to $\tau$ when applying $\bphi_c$ in the case of a matched $c$-step of type $a$.}

	We first consider the case $w_n=a$. In this case $w_n$ is an unmatched $a$-step. Hence the vertices of $\wttau'$ having $\wta'$ as their parent, have $w_n$ as their parent in $\wttau$, while $w_n$ has $\wta$ as its parent. Nothing else changes from $\wttau'$ to $\wttau$. It is clear that the changes are the same from $\tau'$ to $\tau$ (see Figure~\ref{fig:bij-rules-dfsdual2}(a)). So $\tau$ and $\wttau$ are isomorphic, with $\lambdav$ giving the isomorphism. The case where $w_n$ is a $b$-step is symmetric.
	
	Next, we consider the case where $w_n$ is an unmatched $c$-step of type $a$. Let $w_j$ be the $b$-step matched to $w_n$. The changes from $\wttau'$ to $\wttau$ are as follows.
	\begin{compactitem} 
		\item The vertex $w_j$ of $\wttau'$ is deleted (it corresponds to an unmatched $b$-step of $\wtw'$ but a matched $b$-step of $w$). 
		The vertices of $\wttau'$ which had $w_j$ as their parent in $\wttau'$ now have $\wtb$ as their parent in $\wttau$. So, in some sense, the vertex $w_j$ of $\wttau'$ is replaced by $\wtb$.
		\item The vertex $w_n$ is added. The vertices of $\wttau'$ which had $\wta'$ as their parent in $\wttau'$ now have $w_n$ as their parent in $\wttau$. So, in some sense, the vertex $\wta'$ of $\wttau'$ is replaced by $w_n$.
		\item The vertex $w_n$ has $\wta$ as its parent in  $\wttau$. Moreover, the vertices of $\wttau'$  which had $\wtb'$ as their parent in $\wttau'$ (these are matched $c$-steps of type $b$ following $w_j$) now have $\wta$ as their parent in $\wttau$ (because if $w_k$ had $\wtb'$ as its parent in $\wttau'$, then in $\wt w$ one has $p(w_k)=w_n$ and $\vec{p}(w_k)=\wta$). 
		So, in some sense, the vertex $\wtb'$ of $\wttau'$ is replaced by $\wta$. 
	\end{compactitem}
	It is easy to see that these changes from $\wttau'$ to $\wttau$ are isomorphic to the changes from $\tau'$ to $\tau$; see Figure~\ref{fig:bij-rules-dfsdual2}(b). The case where $w_n$ is an unmatched $c$-step of type $b$ is symmetric.
	
	Lastly we consider the case where $w_n$ is a matched $c$-step. Let $w_i$ and $w_j$ be the steps matching $w_k$, with $i<j$. Let us suppose for concreteness that $w_i=a$  (the case $w_i=b$ being of course symmetric). The changes from $\wttau'$ to $\wttau$ are as follows.
	\begin{compactitem} 
		\item The vertices $w_i$ and $w_j$ of $\wttau'$ are deleted.
		The vertices of $\wttau'$ which had $w_i$ (resp. $w_j$) as their parent in $\wttau'$, have $\wta$ (resp. $\wtb$) as their parent in $\wttau$. So, in some sense, the vertices $w_i$ and $w_j$ of $\wttau'$ are replaced by $\wta$ and $\wtb$ respectively.
		\item The vertex $w_n$ is added and its parent is  $\wta$. Moreover, the vertices of $\wttau'$ distinct from $w_i$ which had $\wta'$ as their parent in $\wttau'$ (these are matched $c$-steps of type $a$ whose far-match appears after $w_i$), have $w_n$ as their parent in $\wttau$.  So, in some sense, the vertex $\wta'$ is replaced by $w_n$.
		\item The vertices of $\wttau'$ distinct from $w_j$ which had $\wtb'$ as their parent in $\wttau'$ (these are matched $c$-steps of type $b$ whose far-match appears after $w_j$), have $\wta$ as their parent in $\wttau$. 
		In summary, the children of $\wta$ are $w_n$ together with the vertices of $\wttau'$ which were the children of either $w_i$ or $\wtb'$. So in some sense, the vertex $\wta$ is obtained by merging $w_i$ and $\wtb'$. 
	\end{compactitem}
	It is easy to see that these changes from $\wttau'$ to $\wttau$ are isomorphic to the changes from $\tau'$ to $\tau$; see Figure~\ref{fig:bij-rules-dfsdual2}(c). So $\tau$ and $\wttau$ are isomorphic, with $\lambdav$ giving the isomorphism. 
\end{proof}

\begin{proof}[Proof of Lemma~\ref{lem:envelope-matching-is-loop-building}]
	The listed properties are easy consequences of Theorem~\ref{thm:LR} (for $(M,\sigma)\in\bmT_P$) and Theorem~\ref{thm:spinelt-inf} (for $(M,\sigma)\in\imT_P$). Let us suppose for concreteness that $w\in\bmK$, and that the envelope step $w_k$ is a $c$-step of type $a$ (that is, $w_i=a$); the other cases are treated in the exact same manner. 
	
	Let $w''=w_1\ldots w_{k-1}$
	be the prefix of $w$ ending just before the envelope step $w_{k}$, and let $(M'',\sigma'')=\bPhi(w'')$. This percolated near-triangulation is represented in Figure \ref{fig:proof-enveloppe-interval}.
	Let $v_\ell$ be the top-left vertex of $(M'',\sigma'')$, let $e_{\textrm{top}}$ be the top-edge and let $e_\ell$ be the left edge incident to $v_\ell$. Let also $t_{\textrm{top}}$ and $t_\ell$ be the inner triangles of $M''$ incident to $e_{\textrm{top}}$ and $e_\ell$. 
	Since the steps $w_i$ and $w_k$ form a far-matching, it is clear that $\etavf(i)=t_\ell$, and that $\bphi_{w_k}(M'',\sigma'')$ is obtained from $(M'',\sigma'')$ by gluing the edges $e_{\textrm{top}}$ and $e_\ell$ together. Hence, $\etae(k)$ is the edge $e$ obtained by gluing $e_{\textrm{top}}$ and $e_\ell$ together, and $\etavf(k)=v_\ell$. 
	By definition, the envelope edge $e^*$ (which is the dual of $e$) is the only edge of $\ga$ which is not in $\tau^*$. Hence the path $P=\ga\setminus e^*$ is the path $\tau^*$ between the endpoints $t_\ell^*$ and $t_{\textrm{top}}^*$ of $e^*$. Thus Property (i) of Lemma~\ref{lem:envelope-matching-is-loop-building} holds. Moreover, $P$ is the part of the percolation path of $(M'',\si'')$ from the triangle $\etavf(i)=t_\ell$ to the top-edge, hence by Theorem~\ref{thm:LR} the triangles crossed by $P$ are in one-to-one correspondence via $\etavf$ to the spine steps of $w''$ appearing after $w_i$. It is clear that these steps are exactly the spine steps of $w_iw_{i+1}\cdots w_{k-1}$. Hence Property (iii) holds. 
	
	Next, observe that all the vertices, edges and faces of $M$ inside $\ga$ or crossed by $\ga$ are in-edges of $\bPhi(w_1w_2\cdots w_k)$ (since all the edges of $\ga$ are in-edges, and the top-edge is outside $\ga$). Hence these vertices, edges and faces correspond via $\etavf$ or $\etae$ to steps in $w_1w_2\ldots w_{k}$. Observe now that $t_\ell^*$ is the ancestor in $\tau^*$ of all the vertices of $M^*$ on $\ga$ or inside $\ga$. Hence, none of the vertices, edges and faces inside $\ga$ or crossed by $\ga$ are in-edges of $\bPhi(w_1w_2\cdots w_{i-1})$. Thus Property (ii) holds. 
	
	It remains to prove Property (iv). It is clear from the definitions that $\frk L(\ga)$ corresponds to a part of the white cluster $W$ of $\LR(M'',\si'')$. Precisely, it corresponds to the part of $W$ ``rooted'' at the vertex $v_\ell$. By Theorem~\ref{thm:LR}, $\LR(M'',\si'')=\LR(w'')$, hence the white cluster $W$ is encoded by $\frk L_\ell(w'')$. Moreover, it is clear from the definition of $\frk L_\ell$ that the looptree of $W$ rooted at $v_\ell$ correspond to the steps of $\hpi_\ell(w'')$ appearing after $w_i$ (because it is the last unmatched $a$-step of $w''$). Property (iv) then follows without difficulty.
\end{proof}

\fig{width=.5\linewidth}{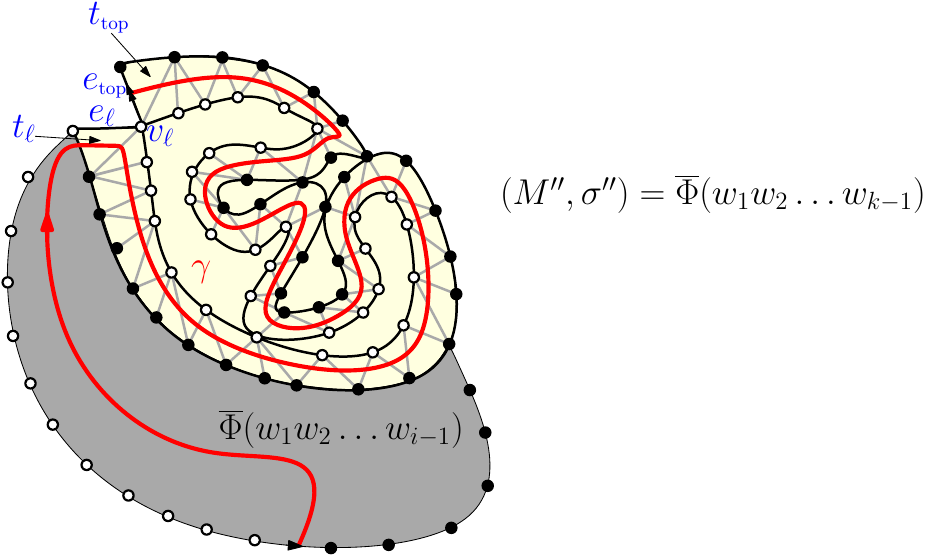}{Notation for the proof of Lemma~\ref{lem:envelope-matching-is-loop-building}. The figure shows the map $(M'',\sigma'')=\bPhi(w'')$, where $w''=w_1\ldots w_{k-1}$. The part of $(M'',\si'')$ corresponding to $\bPhi(w_1\ldots w_{i-1})$ is represented in gray.
	The map $\bphi_{w_k}(M'',\sigma'')$ would be obtained from $(M'',\sigma'')$ by gluing the edges $e_{\textrm{top}}$ and $e_\ell$ together.}

\begin{proof}[Proof of Claim~\ref{claim:disjoint-intervals}]
	Suppose by contradiction that there exist overlapping envelope excursions $w'=w_i\ldots w_k$ and $w''=w_{i'}\ldots w_{k'}$ with $i<i'<k<k'$. Since the matchings $w_i,w_k$ and $w_{i'},w_{k'}$ are crossing, the steps $w_i,w_{i'}$ are different (one is an $a$-step and the other is a $b$-step). Let $w_s,w_t$ (resp $w_{s'},w_{t'}$) be the parent-matching of $w_i,w_k$ (resp. $w_{i'},w_{k'}$). Since $w',w''$ are envelope excursions, both $w_s,w_t$ and $w_{s'},w_{t'}$ are close-matchings, hence they do not cross any other matchings. Thus $s,s'<i$, $k'<t,t'$, and either the matching $w_s,w_t$ encloses $w_{s'},w_{t'}$ or the matching  $w_{s'},w_{t'}$ encloses $w_s,w_t$. But if $w_s,w_t$ encloses $w_{s'},w_{t'}$, then $w_s,w_t$ cannot be the parent-matching of $w_{i},w_k$ (because the matching made of $w_{t'}$ and its far-match encloses $w_{i},w_k$).
	Similarly if $w_{s'},w_{t'}$ encloses $w_{s},w_{t}$, then $w_{s'},w_{t'}$ cannot be the parent-matching of $w_{i'},w_{k'}$. We reach a contradiction. 
\end{proof}

\begin{proof}[Proof of Lemma~\ref{prop:piv}]
	We first prove the first part of Lemma~\ref{prop:piv} using Theorem~\ref{thm:spinelt-inf}. As in Section~\ref{sec:LR-decomposition-infinite}, we denote by $\wt L$ and $\wt R$ the walks obtained from $\wh L$ and $\wh R$ by removing all the 0-steps. We also denote $\frk L_\ell:=\frk L_\ell(w^-)$ and $\frk L_r:=\frk L_r(w^-)$ the forested lines they encode. The condition \eqref{eq:piv} for a time $i$ thus correspond to a time $\wt i<0$ such that 
	\begin{equation}\label{eq:pivtilde}
		\wt L_{\wt i}<\wt L_{\wt i-1} \textrm{ and } \wt L_{\wt i}\geq \min\{L_j~|~\wt i<j\leq 0\}.
	\end{equation}
	Recall also that the negative steps of $\wt L$ correspond to the completion of a bubble of a discrete looptree in the forested line $\frk L_\ell$. In fact, Condition \eqref{eq:pivtilde} corresponds to the completion of a bubble which is either not incident to the semi-infinite path of $\frk L_\ell$ (if $L_{\wt i}>\min\{L_j~|~\wt i<j\leq 0\}$), or incident to the semi-infinite path at a point where at least two bubbles are attached (if $L_{\wt i}=\min\{L_j~|~\wt i<j\leq 0\}$). In both cases, the vertex $\wh\eta_{\op{v}}(i)$ is the vertex of $M$ where the bubble is attached, hence it is incident to at least two bubbles of $\frk L_\ell$.
	
	Next we prove Property (a) about pivotal points of type~1 or~2. The situation is illustrated in Figure~\ref{fig:proof-pivotal}(a). Let us assume for concreteness that the percolation cycle $\ga$ is the outside-cycle of a \emph{white} cluster $C$. If $v$ is a pivotal point of type~1 associated with $\ga$, then $v$ is a cut-point of the white looptree at level~0, denoted by $\frk L_0$, of the forested line $\frk L_\ell$. Recall that $\frk L_0$ is encoded (via $\ccwcode$) by the non-zero steps of $\wh L_{|[-s,0]}$. Hence, the last time $i\in[-s+1,0]$ corresponding to the completion of a bubble of $\frk L_\ell$ attached to $v$ satisfies Property (a) (with $\wh L_i< \wh L_{i-1}$). 
	Similarly, if $v$ is a pivotal point of type~2 associated with $\ga$, then $v$ is a cut-point of the black forested line, surrounded by the part of the percolation path corresponding to $\ga$. Since this part of the black forested line is encoded by $\wh R_{|[-s,0]}$, we can find $i\in[-s+1,0]$ satisfying Property (a) (with $\wh R_i< \wh R_{i-1}$). 
	
	
	\fig{width=.8\linewidth}{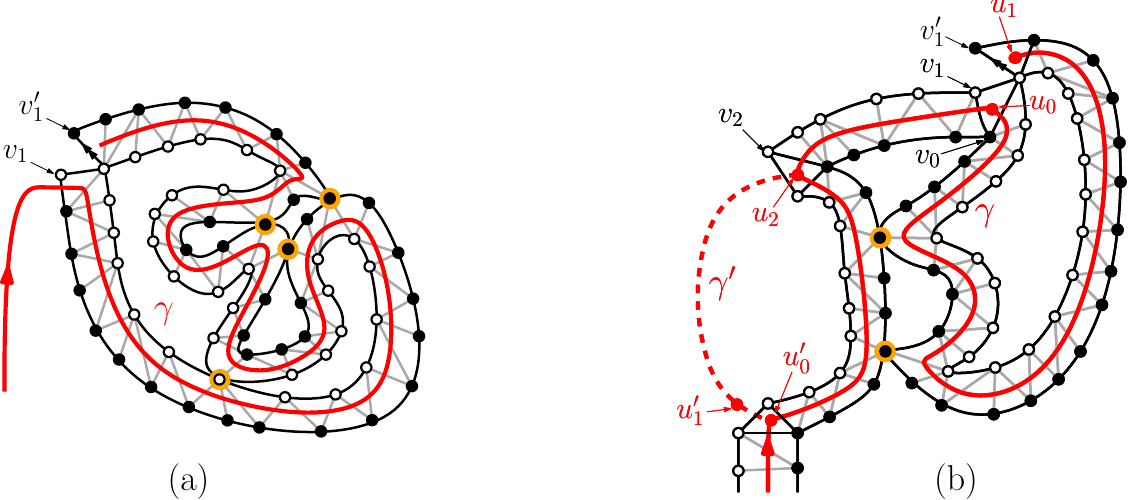}{(a) Pivotal points (indicated in orange) of type 1 (white) and type 2 (black) associated with $\ga$. We represented the portion of the percolation path of $(M^-,\si^-)$ going through $\ga$. (b) Pivotal points of type 3 in $V_{\ga,\ga'}$ (the situation for pivotal points of type 4 would be identical except for the dotted part of $\ga'$ which would surround $\ga$ clockwise). We represented the portion of the percolation path $P^-$ of $(M^-,\si^-)$ going through $\ga'$ and $\ga$. Note that the figure illustrates the situation in $(M^-,\si^-)$, as opposed to $(M,\si)$: in $(M,\si)$ the vertices $v_1$ and $v_1'$ would be identified, and some of the vertices would be of a different color.}
	
	It remains to prove Property (b) about pivotal points of type 3 and 4. The situation is illustrated in Figure~\ref{fig:proof-pivotal}(b). We again assume that the percolation cycle $\ga$ is the outside-cycle of a \emph{white} cluster $C$.
	Let $\ga'$ be a percolation cycle which is the outside-cycle of a cluster $C'$, such that the envelope excursion of $\ga'$ starts before the envelope excursion of $\ga$. Let $V_{\ga,\ga'}$ be the set of vertices of $(M,\si)$ such that flipping the color of $v$ connects the two clusters $C,C'$; these are pivotal points of type 3 or 4 (depending on the color of $C'$) associated with $\ga$. We want to show that all but 3 vertices in $V_{\ga,\ga'}$ satisfy Property (b). A crucial point is to show that these vertices are incident to some edges $e,e'$ of $M^-$ (as opposed to $M$) crossed by $\ga$ and $\ga'$ respectively.

	We introduce some notation which is illustrated in Figure~\ref{fig:proof-pivotal}(b). Let $\tau^*=\dfs(M,\si)$ and let $P=\ga\cap\tau^*$. We know that $P$ is a path in $\tau^*$ containing every edge of $\ga$ except one. Let $u_0,u_1\in V(M^*)$ be the endpoints of $P$ with $u_0$ the ancestor of $u_1$. We define $P',u_0',u_1'$ similarly for $\ga'$. Let $e_0$ be the parent edge of $u_0$ in $\tau^*$, and let $v_0,v_1\in V(M)$ be the endpoints of the edge $e_0^*$, with $v_0$ on the right of the edge $e_0$ oriented toward $u_0$.
	Let $u_2$ be the first common ancestor of $u_1$ and $u_1'$ in $\tau^*$. Let $f_2\in F(M^*)$ be the forward face of $u_2$, and let $v_2=f_2^*\in V(M)$.
	Let $P^-$ be the percolation path of $(M^-,\si^-)$. Observe that the vertices on $P^-$ are the ancestors of $u_1$ and that $P'\cap P^-$ goes from $u_0$ to $u_2$.

	Let $v$ in $V_{\ga,\ga'}\setminus \{v_0,v_1,v_2\}$. We will show that $v$ satisfy Property (b). 
	Since $v\in V_{\ga,\ga'}$, it is incident to an edge $e$ of $M$ crossing $\ga$ and to an edge $e'$ of $M$ crossing $\ga'$. Dually, there are edges $e_*\in\ga,e_*'\in \ga'$ incident to the face $v^*$ of $M^*$. Moreover, $v\neq v_1$ implies that $e_*$ is not the envelope edge of $\ga$ so $e_*\in \ga\cap P^-$. We now show that $e_*'\in \ga'\cap P^-$. Indeed, let us assume by contradiction that $e_*'\notin P^-$. In this case, the endpoints $x',y'$ of $e_*'$ are descendants of $u_2$ which are ancestors of $u_1'$, while the endpoints $x,y$ of $e_*$ are descendants of $u_2$ which are ancestors of $u_1$. The vertices $x,y,x',y'$ are all incident to the face $v^*$. This is represented in Figure~\ref{fig:proof-piv}(a).
	
	\fig{width=.75\linewidth}{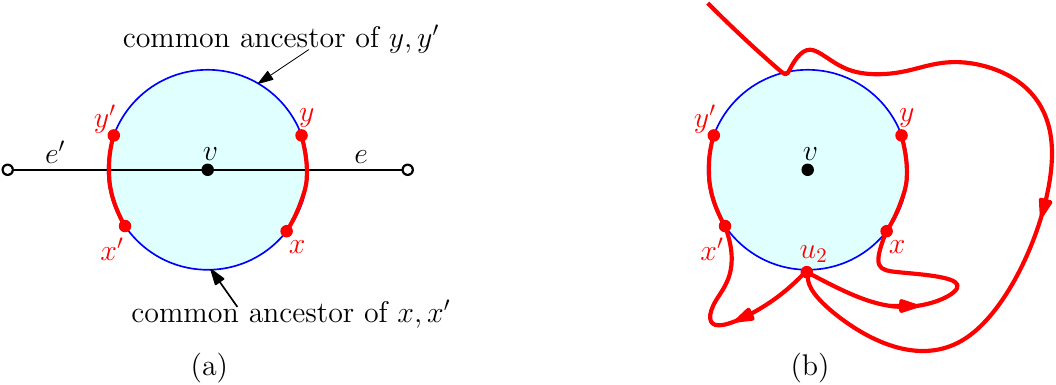}{Proof of Property (b) by contradiction. (a) The face $v^*$ of $M^*$ and the vertices $x,y,x',y'$. The part of the boundary of $v^*$ between $x,x'$ (resp. $y,y'$) must contain a common ancestor of $x,x'$ (resp. $y,y'$), hence an ancestor of $u_2$. (b) Branches of $\tau^*$ toward $x,x',y,y'$.}
	
	Since $\tau^*$ is a DFS tree, the vertices $x,y$ must be separated from the vertices $x',y'$ by some ancestors of $u_2$ around the face $v^*$. Since we need two separators around $v^*$, it is easy to see that planarity imposes that one of them is $u_2$, and that $v^*$ is the forward face at $u_2$. This is represented in Figure~\ref{fig:proof-piv}(b). This means $v=v_2$, which is a contradiction.
	Thus, $v$ is incident to an edge $e$ of $M^-$ crossing $\ga\cap P^-$ and an edge $e'$ of $M^-$ crossing $\ga'\cap P^-$. 
	Thus, $v$ is on the black forested line $\frk L_r$, and it disconnects $v_0$ from the semi-infinite line of $\frk L_r$. 
	Hence $v$ is a cut point of $\frk L_r$ and corresponds to a negative step of $\wh R$.
	More precisely, let us consider the first triangle $t$ incident to $v$ and crossed by $\ga\cap P^-$, and the last triangle $t'$ incident to $v$ and crossed by $\ga'\cap P^-$. The triangles $t,t'$ corresponds to some steps $i,i'$ of the walk $(\wh L_k,\wh R_k)$, which satisfy $\wh \eta_{\op{v}}(i)=v$, $i'\leq -s< i\leq 0$ (because the last $s$ steps correspond to the triangles on $\ga$), 
	and $R_{i}=R_{i'-1}>\min\{R_j~|~i'\leq j<i\}$ (because the black edge incident to $t'$ starts a bubble of $\frk L_r$ which ends with the black edge incident to $t$). Hence Property (b) holds (with $R_i<R_{i-1}$).
\end{proof}

\subsection{Description of the crossing events in terms of Kreweras walks}\label{sec:crossing-discrete}

\nina{In this section we will introduce crossing events $E_{\op{b}}(v),E_{\op{w}}(v)$ for a percolated near-triangulation $(M,\si)$ and describe these crossing events in terms of the associated Kreweras walk. The continuum analog of these descriptions will be given in Section~\ref{sec:crossing} (see~\eqref{eq:crossing2} and~\eqref{eq:crossing4}). The convergence from discrete to continuum is stated in Section~\ref{sec:conv} and proved in Section~\ref{app:conv}.}

Let $(M,\si)$ be a percolated triangulation with a simple boundary \nina{of length $h+2$ such that the only white outer vertex is the root-vertex}. Let $a_1, a_2,a_3,a_4$ be distinct outer edges of $M$ appearing in this order in clockwise direction around the root-face of $M$, with $a_1$ being the root-edge. For $i,j\in[4]$, we denote by $V(a_i,a_{j})$ 
the set of outer vertices of $M$ situated between $a_i$ and $a_j$ in clockwise order around the root-face (including one endpoint of $a_i$ and one endpoint of $a_j$).

\fig{width=.9\linewidth}{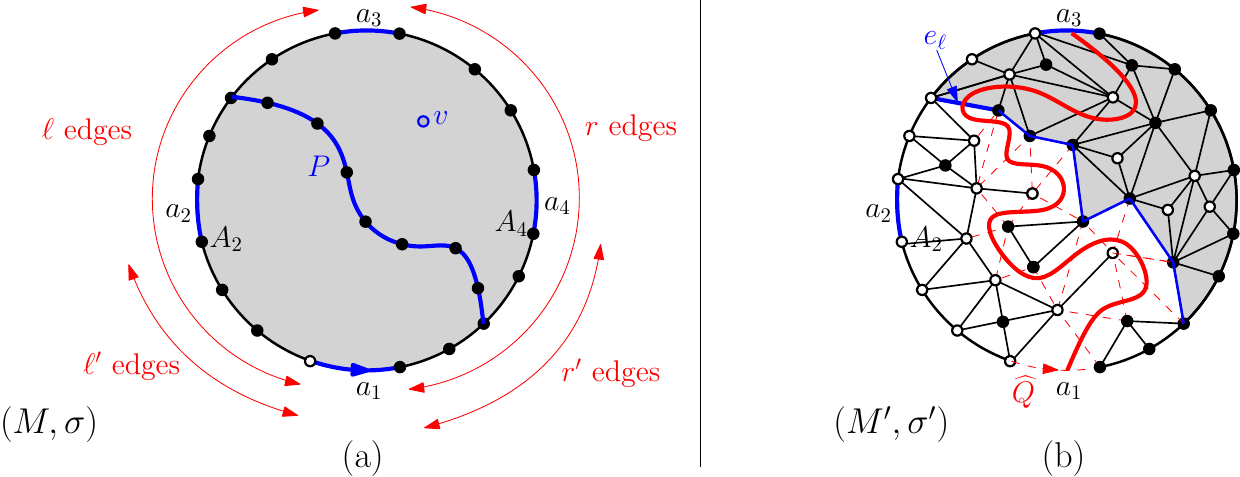}{(a) The event $E_\ob(v)$ for a percolated near-triangulation $(M,\si)$ and a vertex $v$, with a path $P$ satisfying Conditions (a-c). (b) Illustration of characterization of $E_\ob(v)$ given by Fact~\ref{fact:a}. The edge $e_\ell$ and the map $M_\ob$ (the dashed edges crossed by $\hQ$ are not in $M_\ob$). The 1-component of $M_\ob$ containing $a_3$ is represented in gray, while the path $P$ discussed in the proof of Fact~\ref{fact:a} is indicated in blue.}

For a vertex $v\in V(M)$, we denote by $E_\ob(v)$ the event (depending on $(M,\si)$ and $a_2,a_3$) that there exists a simple path $P$ on $M$ such that 
\begin{compactitem}
	\item[(a)] $P$ has one endpoint in $V(a_2,a_3)$ and one endpoint in $V(a_3,a_1)$, 
	\item[(b)] all the non-endpoint vertices of $P$ are inner black vertices of $(M,\si)$,
	\item[(c)] either $v\in P$ or $v$ is on the same side of $P$ as the edge $a_3$.
\end{compactitem}
The event $E_\ob(v)$ is illustrated in Figure~\ref{fig:crossing-discrete}(a). We define the (symmetric) event $E_\ow(v)$ in the same way except that Conditions (a-b) are replaced by 
\begin{compactitem}
	\item[(a')] $P$ has one endpoint in $V(a_3,a_4)$ and one endpoint in $V(a_1,a_3)$,
	\item[(b')] all the non-endpoint vertices of $P$ are inner white vertices of $(M,\si)$.
\end{compactitem}


\nina{We will now describe the crossing events $E_{\op{b}}(v),E_{\op{w}}(v)$ in terms of the Kreweras walk.} By Remark~\ref{rmk8}, \nina{percolated maps $(M,\si)$ of the form we consider here are} 
in bijection, via $\bPhi$, with Kreweras walks from $(0,0)$ to $(0,-h)$ staying in the quadrant $\{(x,y)~|~x\geq 0,y\geq -h\}$ (in this bijection, the $h$ right edges of $(M,\si)$ are considered inactive).
Let $w=w_1\ldots w_n=\Phi^{-1}(M,\si)\in\smK$.

Let $a_1,a_2,a_3,a_4$ be distinct outer edges of $M$ appearing in this order in clockwise direction around the root-face of $M$, with $a_1$ being the root-edge.
Recall the notation $V(a_i,a_j)$, and the definition of the crossing events $E_\ob(v)$ and $E_\ow(v)$ (depending on $(M,\si)$ and $a_2,a_3,a_4$). The event $E_\ob(v)$ is illustrated in Figure~\ref{fig:crossing-discrete}(a). 
We denote by $\ell'$ and $\ell$ the number of vertices in $V(a_1,a_2)$ and $V(a_1,a_3)$ respectively, \emph{minus 1}. Similarly, we denote by $r'$ and $r$ the number of vertices in $V(a_4,a_1)$ and $V(a_3,a_1)$ respectively, \emph{minus 1}. 

We first characterize the events $E_\ob(v)$ and $E_\ow(v)$ in terms of percolation paths. Let $(M,\si')$ be the percolated near-triangulation obtained from $(M,\si)$ by recoloring all the vertices in $V(a_1,a_3)$ in white.
We denote by $\hQ$ the percolation path of $(M,\si')$. We consider the path $\hQ$ as oriented from $a_1$ to $a_3$.
Let $e_\ell$ (resp. $e_r$) be the first edge crossed by $\hQ$ which is incident to a vertex in $V(a_2,a_3)$ (resp. $V(a_3,a_4)$). 
Let $M_\ob$ (resp. $M_\ow$) be the map obtained from $M$ by deleting all the edges crossed by $\hQ$ strictly before reaching $e_\ell$ (resp. $e_r$).

Recall that the \emph{1-components} of a connected graph are the maximal sets of vertices which are not separated by \emph{cut-points} (vertices whose deletion produces an non-connected graph), and the edges between these vertices.

\begin{fact}\label{fact:a}
	The event $E_\ob(v)$ (resp. $E_\ow(v)$) occurs if and only if $v$ and $a_3$ are in the same 1-component of $M_\ob$ (resp. $M_\ow$).
\end{fact}

Let us briefly sketch the justification of Fact~\ref{fact:a}. If $v$ and $a_3$ are in the same 1-component of $M_\ob$, then there is a path $P$ satisfying the Conditions (a-c) of $E_\ob(v)$ which is made of $e_\ell$ and some unicolor black inner edges incident to the triangles crossed by the portion of the percolation path $\hQ$ from $a_1$ to $e_\ell$. This is illustrated in Figure~\ref{fig:crossing-discrete}(b). Conversely, if $v$ and $a_3$ are not in the same 1-component of $M_\ob$, then the white vertices incident to the edges crossed by the portion of the percolation path $\hQ$ from $a_1$ to $e_\ell$ prevent the existence of a path path $P$ satisfying Conditions (a-c).

Let $A_2$ (resp. $A_4$) be the endpoint of $a_2$ (resp. $a_4$) in $V(a_1,a_2)$ (resp. $V(a_4,a_1)$). It is easy to see from Fact~\ref{fact:a} that $E_\ob(A_4)$ (resp. $E_\ow(A_2)$) occurs if and only if $e_\ell$ appears before (resp. after) $e_r$ along $\hQ$. Hence $E_\ob(A_4)$ occurs if and only if $E_\ow(A_2)$ does not occur.

Next, we characterize the event $E_\ob(v)$ in terms of the past/future decomposition of $(M,\si)$ at the edge $a_3$, which is represented in Figure~\ref{fig:crossing-discrete2}(a). 
Let $\etae$ be the bijection between $[n]$ and the non-top edges of $(M,\si)$ given by Definition~\ref{def:eta-extend}. Let $t_3=\etae^{-1}(a_3)$, let $w^{-}=w_1\ldots w_{t_3-1}$, and let $w^+=w_{t_3}\ldots w_n$.  Note that $a_3$ is the top-edge of the past map $(M^-,\si^-)=\bPhi(w^-)$ and becomes inactive at time $t_3$. Hence all the right edges of $(M^-,\si^-)$ are inactive, while all the left edges are active; see Figure~\ref{fig:crossing-discrete2}(a).
We consider the decomposition $w^+=cv^1c v^2 c\cdots cv^k$ with $v^i$ in $\smK$, separated by $c$-steps without matching steps in $w^+$ (note that $w^+$ starts with a separating $c$-step). We denote by $(M^+,\si^+)$ the \emph{future near-triangulation} obtained by gluing together the maps $(Q_i,\be_i)=\bPhi(v^i)$.
Since all right edges of $(M^-,\si^-)$ are inactive, the separating $c$-step of $w^+$ are all of type $b$. Hence for all $i\in[k]$ the top-left vertex of $(Q_{i-1},\be_{i-1})$ is identified with the black endpoint of the root-edge of $(Q_{i},\be_{i})$, with the notational convention $(Q_{0},\be_{0})=(M^-,\si^-)$. 
For $i\in[k+1]$, let $v_i^+$ and $e_i^+$ be the top-left vertex and top-edge respectively of $(Q_{i-1},\be_{i-1})$.
Note that the vertices $v^+_1,\ldots,v^+_k$ are all white in $(M,\si')$. This shows that the left vertices of the map $(M^-,\si^-)$ are white in $(M,\si')$, hence that the percolation path of $(M^-,\si^-)$ is the same as the percolation path $\hQ$ of $(M,\si')$; see Figure~\ref{fig:crossing-discrete2}(a). In particular, the edge $e_\ell$ considered in Fact~\ref{fact:a} is an inner edge of $(M^-,\si^-)$ which is incident to a left vertex $v^-$ of $(M^-,\si^-)$ which will be identified with one of the vertices in $\{v^+_1,\ldots, v^+_k\}$. More precisely, $v^-$ will be identified with $v_{i_*}^+$, where $a_2$ is a right edge of $(Q_{i_*},\be_{i_*})$; see Figure~\ref{fig:crossing-discrete2}(b).
We can now state another characterization of $E_\ob(v)$.

\begin{fact} \label{fact:b}
	The event $E_\ob(v)$ occurs if and only if either $v=v^+_{i_*}$ or $\etavf^{-1}(v)\in [t_\ell^-, t_\ell^+]$, where $t_\ell^-=\etae^{-1}(e_\ell)$ and $t_\ell^+=\etae^{-1}(e_{i_*}^+)$. 
\end{fact}

\fig{width=.8\linewidth}{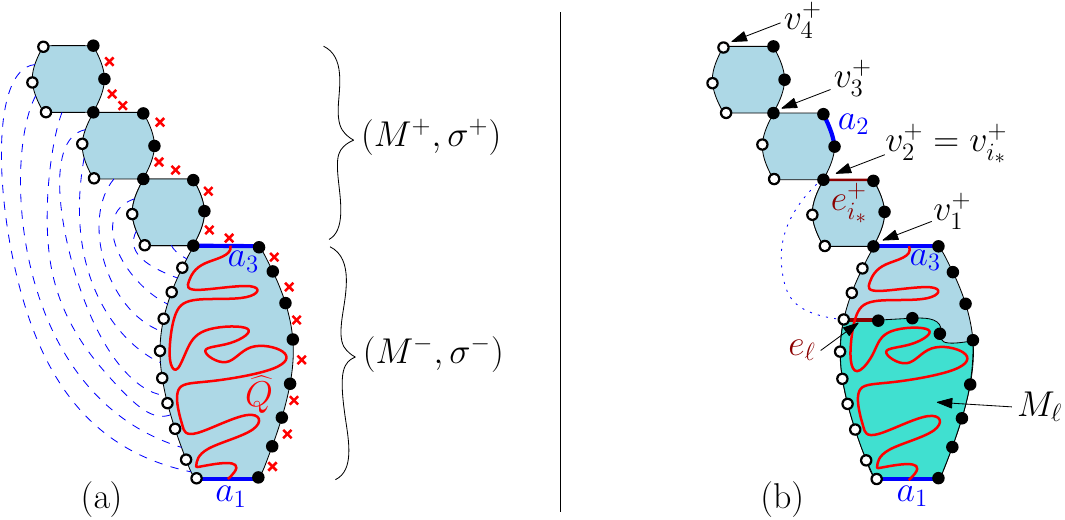}{(a) Past/future decomposition of $(M,\si)$ at time $t_3-1$. Here $k=3$. (b) Illustration of Fact~\ref{fact:b}.}

Let us justify Fact~\ref{fact:b} using Fact~\ref{fact:a}. Let $\eps$ be the outer edge following $e_{i_*}^+$ in counterclockwise order around $M$.
Note that $\etavf^{-1}(v^+_{i^*} ) = \etae^{-1}(\eps)>t_3$.
For a vertex $v\neq v^+_{i_*}$ such that $\etavf^{-1}(v)>t_3$ (equivalently, such that $v\in M^+$), it is easy to see that $v$ and $a_3$ are in the same connected component of $M_\ob$ if and only if $v$ is in $(Q_{i},\be_{i})$ for $i< i_*$. Indeed, the other vertices are separated by the cutpoint $v^+_{i_*}$ because the path $\hQ$ crosses all the triangles of $M^-$ incident to left edges of $M^-$. Thus for a vertex $v\neq v^+_{i_*}$ such that $\etavf^{-1}(v)>t_3$, the event $E_\ob(v)$ occurs if and only if $\etavf^{-1}(v)\in (t_3, t_\ell^+]$.
For a vertex $v$ such that $\etavf^{-1}(v)\leq t_3$ (equivalently, such that $v$ is an in-vertex of $M^-$, or the top-right vertex of $M^-$), it is easy to see that $E_\ob(v)$ occurs if and only if $v$ is \emph{not} an in-vertex of $(M_\ell,\si_\ell)=\bPhi(w_1\ldots w_{t_\ell^--1})$. Indeed, the percolation path of $(M_\ell,\si_\ell)$ is the part of $\hQ$ occurring before $e_\ell$, and all the right active vertices of $(M_\ell,\si_\ell)$ will remain distinct black vertices of $(M,\si)$. Hence the right active edges of $(M_\ell,\si_\ell)$ together with $e_\ell$ will form a path $P$ satisfying the Conditions (a-c) of the event $E_\ob(v)$; see Figure~\ref{fig:crossing-discrete2}(b). Thus, for a vertex $v$ such that $\etavf^{-1}(v)\leq t_3$, the event $E_\ob(v)$ occurs if and only if $\etavf^{-1}(v)\in [t_\ell^-, t_3]$.

We will now give a characterization of $E_{\ob}(v)$ in terms of $w$. For this we need to express $t_\ell^-$ and $t_\ell^+$ in terms of $w$.
For $i\in \{0,\ldots,n\}$, let $(M^{(i)},\si^{(i)}):=\bPhi(w_1\ldots w_i)$. Let $(L_i,R_i)_{i\in\{0,\ldots,n\}}$ represent the coordinates of the lattice walk $w$, \emph{with the convention that the walk starts at $(L_0,R_0)=(0,h)$ and ends at (0,0)} (note that this convention is different than the one used in Section~\ref{sec:bijection}, but aligns better with the mating-of-trees convention of Section~\ref{sec:crossing}). 
Note that $L_i$ is the number of (active) left edges of $M^{(i)}$ while $R_i=h+r_i-r_i'$, where $r_i$ and $r_i'$ are the number of active and inactive right edges of $M^{(i)}$ respectively. It is easy to check the following fact.

\begin{fact}\label{fact:c}
	Let $k\in[h]$ and let $b_k$ be the $k^{\textrm{th}}$ right edge of $(M,\si)$ in clockwise order around the root-face, starting from the top-edge.
	Then, $\etae^{-1}(b_k)=\inf\{i~|~R_i<k\}.$ 
\end{fact}

Note that $a_3=b_\ell$, hence by Fact~\ref{fact:c},
\begin{equation}\label{eq:t3-discrete}
	t_3=\inf\{t~|~R_t< \ell\}.
\end{equation}
Let $\textrm{cut}(t_3)$ be the set of indices in $[t_3,n]$ corresponding to the separating $c$-steps of $w^+=cv^1c v^2 c\cdots cv^k$. Note that $t\in \textrm{cut}(t_3)$ if and only if $\etae(t) \in\{e^+_1,\ldots,e^+_k\}$. Hence it is easy to see that $t_\ell^+=t_2$, where
\begin{equation}\label{eq:t2-discrete}
	t_2=\sup\{t\in \textrm{cut}(t_3)~|~R_t\geq \ell'\}.
\end{equation}
Moreover, observing that $t_\ell^--1$ is the first time the vertex $v^+_{i_*}$ is created (as the top-left vertex of $(M_\ell,\si_\ell)$), we get $t_\ell^--1=t_2'$, where
\begin{equation}\label{eq:t2p-discrete}
	t_2'=\inf\{t~|~L_{t'}\geq L_{t_2},~\forall t'\in[t,t_2]\}.
\end{equation}
Summarizing, we get the following characterization of $E_\ob(v)$ in terms of $w$.

\begin{fact}\label{fact:d}
	Let $v$ be a non-top vertex of $(M,\si)$. 
	The event $E_\ob(v)$ occurs if and only if either $v=v^+_{i_*}$ or $\etavf^{-1}(v)\in [t_2'+1, t_2]$, where $t_2$ and $t_2'$ are given by~\eqref{eq:t2-discrete} and~\eqref{eq:t2p-discrete} respectively. 
\end{fact}

We can also get a similar characterization for $E_\ob(A_4)$.
\begin{fact}\label{fact:e}
	The event $E_\ob(A_4)$ occurs if and only if $\ds \inf_{t\in [0,t_2']}R_t\geq h-r'.$
\end{fact}
Indeed, by Fact~\ref{fact:c}, $\inf_{t\in [0,t_2']}R_t< h-r'$ if and only if $\etae^{-1}(a_4)\leq t_2'$. The latter condition occurs if and only if $a_4$ is an inactive right edge of $(M_\ell,\si_\ell)$, hence if and only if $A_4$ and $a_3$ are \emph{not} in the same 1-component of $M_\ob$.\\


It remains to characterize $E_\ow(v)$ in terms of $w$. We first translate the characterization given by Fact~\ref{fact:a} using the past/future decomposition of $(M,\si)$ at time $t_3-1$ (see Figure~\ref{fig:crossing-discrete2}(a)). Let $t_r^-=\etae^{-1}(e_r)$, and let $(M_r,\si_r)=\bPhi(w_1\ldots w_{t_r^--1})$. 
Let $j_*$ be the smallest value $j\in[k+1]$ such that $v^+_j$ is identified to a left vertex of $(M_r,\si_r)$. We now state our second characterization of $E_\ow(v)$.

\begin{fact} \label{fact:f}
	The event $E_\ow(v)$ occurs if and only if either $v=v^+_{j_*}$ or $\etavf^{-1}(v)\in [t_r^-, t_r^+]$, where $t_r^-=\etae^{-1}(e_r)$ and $t_r^+=\etae^{-1}(e_{j_*}^+)$ if $j_*\in[k]$, and $t_r^+=n$ if $j_*=k+1$. 
\end{fact}
The proof of Fact~\ref{fact:f} is similar to the proof of Fact~\ref{fact:b} and is left to the reader. We now express $t_r^\pm$ in terms of $w$. 
First, it is easy to check that the number of left edges of $(M^-,\si^-)$ which are also left edges of $(M_r,\si_r)$ is $\ds \inf_{s\in[t_r^-,t_3-1]}(L_s)$. Hence it is easy to check that $t_r^+=t_4'$, where
\begin{equation}\label{eq:t4p-discrete}
	t_4'=\inf\{t\in \textrm{cut}(t_3)~|~L_t\leq \inf_{s\in[t_r^-,t_3-1]}(L_s)\}.
\end{equation}
We now characterize $t_r^-$. By Fact~\ref{fact:c}, $a_4=b_{h-r'}=\etae(t)$ for $t=\inf\{i~|~R_i<h-r'\}$. 
Let us denote by $\textrm{anfr}(t_3)$ the set of indices $t\in [t_3]$ such that either $t=t_3$ or $w_t$ is a spine step of $w^-$. Recall from Theorem~\ref{thm:LR} that $e\in E(M)$ is crossed by the percolation path $\hQ$ if and only if $\etae^{-1}(e)\in \textrm{anfr}(t_3)$. 
We now claim that 
$t_r^-=t_4$, where
\begin{equation}\label{eq:t4-discrete}
	t_4=\inf\{t\in \textrm{anfr}(t_3)~|~R_t<h-r'\}.
\end{equation}
Indeed, with this definition, $w_{t_4-1}$ is the $c$-step ending the maximal cone excursion $w'$ of $w^-$, such that $a_4$ is a right edge of $\bPhi(w')$. To summarize, we obtain the following characterization of $E_\ow(v)$.

\begin{fact}\label{fact:g}
	Let $v$ be a non-top vertex of $M$. 
	The event $E_\ow(v)$ occurs if and only if either $v=v^+_{j_*}$ or $\etavf^{-1}(v)\in [t_4, t_4']$, where $t_4$ and $t_4'$ are given by~\eqref{eq:t4-discrete} and~\eqref{eq:t4p-discrete} respectively. 
\end{fact}


\subsection{Proofs for the infinite volume results}\label{app:inf}

In this subsection we prove Theorems~\ref{thm:UIPT},~\ref{thm:dfs-inf},~\ref{thm:dfs-inf2}, and~\ref{thm:dfsdual-inf} as well as Lemma~\ref{prop:normal} and Proposition~\ref{prop:dfs-ordering-inf}.

\begin{proof}[Proof of Lemma~\ref{prop:normal}] First we argue that $w^+$ a.s.\ has infinitely many cut-times. It is clear that the set $S\subseteq \NN$ of cut-times of $w^+$ is a \emph{renewal process}. That is to say, there exist i.i.d.\ random variables $\{X_i\}_{i\in\NN}$ supported on $\ZZ^{>0}=\{1,2,\dots \}$ such that $S=\{\sum_{i=1}^nX_i\,:\,n\in\NN\}$. We need to show that $|S|=\infty$ a.s., or equivalently that each $X_i$ is finite a.s.
	If the probability $q$ that $X_1$ is finite was less than 1, then $|S|$ would be a geometric variable (of parameter $1-q$) and would have a finite expectation. Hence, it suffices to show $\E[|S|]=\infty$. Now, observe that for all $k\geq 0$ the event $k \in S$ has the same probability as the event that a uniformly random walk with steps $\{(-1,0),(0,-1),(1,1)\}$ starting at $(0,0)$, starts with a step $(1,1)$ and stays in the quadrant $\{(i,j)~|~i,j>0\}$ for at least $k-1$ additional steps. By \cite[Theorem 1]{dw-cones}, the probability that a random walk stays inside a cone for a long period of time is asymptotically the same as for the limiting Brownian motion. By using the exponent for cone excursions of the Brownian motion \cite[equation (4.3)]{shimura-cone}, we get that $\P[k\in S]\sim Ck^{-3/4}$ for some constant $C$ (here the cone angle is $2\pi/3$). This implies $\E[|S|]=\infty$ as desired.

	The proof that $w^+$ a.s.\ has infinitely many split-times is done in the exact same way. Again the set of split-times is a renewal process, and we estimate the probability that a given time is a split-time by using the Brownian motion exponent for cone excursions. 
\end{proof}

We now embark on the proof of the main results. We first need to define a notion of \emph{local convergence} of several functionals of the percolated triangulations.
For $n\in\N_+$, let $(M_n,\sigma_n)$ be a pair made of a rooted loopless triangulation $M_n$ of the sphere with $n+2$ vertices (hence $3n$ edges), and a percolation configuration $\sigma_n$ of $M_n$ such that the root-edge goes from a white vertex to a black vertex.
By Corollary~\ref{cor:chordal-case}, $(M_n,\sigma_n)$ is encoded bijectively (through $\bPhi$) by a walk $\widetilde{w}^n\in\smK$ of length $3n-1$, such that every step is fully-matched, except one $c$-step which is only matched to a single $a$-step or $b$-step. Let $\frk e_n\in E(M_n)$ be a non-top edge of $M_n$ chosen uniformly at random among the non-top edges (equivalently, $\frk e_n$ is in the image of $\etae:[3n-1]\to E(M)$).

We now consider some additional functionals of $(M_n,\sigma_n,\frk e_n)$. We denote by $w^n$ the walk obtained by recentering $\widetilde{w}^n$ in such a way that the  $0$th step corresponds to the edge $\frk e_n$, and we adjust the domain of $\etae$ correspondingly, that is to say, $w^n=w^n_{-k}w^n_{-k+1}\dots w^n_{0}\dots w^n_{3n-k-2}$ for some $k\in\{0,\dots,3n-2 \}$ chosen such that $\etae^{-1}(\frk e_n)=0$. 
Let $\tau_n^*=\dfs(M_n,\sigma_n)$ be the DFS tree of the dual map $M^*_n$ of $M_n$ (as defined in Section~\ref{sec:exploration-tree-from-walk}). Let $h^n$ be the height code of $\tau_n^*$ recentered (in time and height) in such a way that $h^n=(h^n_{-k'},\dots,h^n_{2n-k'})$ with $h^n_0=0$, where $k'$ is the number of $a$ and $b$ steps in $w^n_{-k}w^n_{-k+1}\dots w^n_{-1}$. We obtain a tuple $(w^n,M_n,\sigma_n,\frk e_n,\tau_n^*,h^n)$. 
We say that this tuple converges \emph{locally} if the walks $w^n$ and $h^n$ converge on compact sets of indices, and if for any $R\in\N$, we have convergence of $M_n$, $\sigma_n$, and $\tau_n^*$ restricted to the ball of radius $R$ of $M_n$ centered at $\frk e_n$.

Observe that the triple $(M_n,\sigma_n,\frk e_n)$ is chosen uniformly at random if and only if $w^n$ is chosen uniformly at random from the walks of length $3n-1$ in $\smK$ with a single unmatched $c$-step, recentered at a uniformly random time. 

\begin{lemma}\label{prop17} 
	Consider the random tuple $(w^n,M_n,\sigma_n,\frk e_n,\tau_n^*,h^n)$ defined as above, with $(M_n,\sigma_n,\frk e_n)$ chosen uniformly at random. Then the following hold.
	\begin{compactitem}
		\item[(i)] The tuple 
		$(w^n,M_n,\sigma_n,\frk e_n,\tau_n^*,h^n)$ converges locally in law to a limiting tuple $(w,M,\sigma,\frk e,\tau^*,h)$.
		\item[(ii)] In the limiting tuple, $w\in\{a,b,c \}^\ZZ$ has steps chosen uniformly and independently at random, $M$ has the law of a loopless UIPT rooted at the undirected edge $\frk e$, $\sigma$ is an instance of critical site-percolation on $M$ (that is, the color of the vertices are uniformly random), and $\tau^*\subset E(M^*)$.
		\item[(iii)] In the limiting tuple, the walk $w$ determines the tuple $(M,\sigma,\frk e,\tau^*,h)$. Furthermore, $(M,\sigma)=\Phi^\infty(w)$.
		\item[(iv)] In the limiting tuple, $\tau^*$ is a one-ended spanning tree of the dual map $M^*$ which is in $\DFS_{M^*}$. Moreover $\Lambda_{M}(\tau^*)=\sigma$. Lastly, $\tau^*$ satisfies the properties (ii-iii) of Theorem~\ref{thm:dfs-inf}, $\tau^*$ satisfies the properties of Proposition~\ref{prop:dfs-ordering-inf}, and $\sigma=\Lambda_{M^*}(\tau^*)$.
		\item[(v)] The rooted plane trees $\tau^*$ and $\dfs(w)$ have the same underlying rooted tree (although their planar embeddings may differ).
		\item[(vi)] The assertions of Theorem~\ref{thm:dfsdual-inf} about $\dfsdual(M,\sigma)$ and $\dfsdual(w)$ are true.
	\end{compactitem}
\end{lemma}

\begin{proof}
	\referee{We start by proving the assertions (i) and (iii).}
	The walk $w^n$ converges locally in law to a walk $w$ with uniform and independent increments; this can be proved by proceeding as in \cite[Section 4.2]{shef-burger}. 
	Define $(M',\sigma')=\Phi^\infty(w)$. First we will argue that, given any $R\in\N_+$, there a.s.\ exists some (random) $K=K(R,w)\in\N_+$ such that $w_{-K}\dots w_K$ determines $(M',\sigma')$ restricted to the ball of radius $R$ around the root of $M'$. 
	Since $w$ is uniformly random in $\mK^\infty$, it is easy to see that $M'$ is almost surely locally finite. Hence the ball of radius $R+1$ of $M'$ is finite. Hence there exists $K'=K'(R,w)\in\N_+$ such that the preimage by $\etae$ of all the edges in this ball is in $[-K',K']$.
	Hence there exist $\ell_\pm\in\N_+$ such that $w_{-K'}\dots w_{K'}$ is a subword of $w^-_{\ell_-}\dots w^-_1w^+_1c\dots cw^+_{\ell_{+}}$. The subword $w^-_{\ell_-}\dots w^-_1w^+_1c\dots cw^+_{\ell_{+}}$ determines $(M',\sigma')$ restricted to the radius $R$ ball, which implies that some appropriate $K$ can be found a.s., since we can choose $K$ such that $w^-_{\ell_-}\dots w^-_1w^+_1c\dots cw^+_{\ell_{+}}$ is a subword of $w_{-K}\dots w_{K}$. By Proposition~\ref{prop:future-past} and the infinite volume construction in Section~\ref{sec:bij-inf}, the pair $(M_n,\sigma_n)$ is determined from the walk $w^n$ in a local way in exactly the same way as $\Phi^\infty(w)$ is obtained from $w$. These observations imply that $(w^n,M_n,\sigma_n,\frk e_n)$ converges locally to $(w,M,\sigma,\frk e)$, where $(M,\sigma)=\Phi^\infty(w)$ and $\frk e$ is the (undirected) root-edge of $(M,\sigma)$.
	
	In the finite volume case and for $e\in E(M_n)$ with dual $e^*\in E(M_n^*)$, we have $e^*\in\tau^*_n$ if and only if $w_{\etae^{-1}(e)}\in\{a,b \}$. This implies the convergence of $\tau^*_n$ to a limit $\tau^*$ if we view $\tau^*_n$ and $\tau^*$ as subsets of $E(M^*_n)$. Moreover $w$ determines $\tau^*$. Let $h$ be the spine length sequence of $w$ (as defined in Section~\ref{sec:DFS-infinite}). We claim that $h^n$ converges locally to $h$. Indeed, each value $h_i^n$ of $h$ depends on a finite portion of $w$. Moreover, it is easy to see using Proposition~\ref{thm:DFS}, that each value $h_i^n$ of $h^n$ is obtained from $w^n$ in exactly the same way as $h_i$ is obtained from $w$. So $h^n$ converges to $h$. Combining this with the results of the previous paragraph, we obtain assertions (i) and (iii) of the lemma.

	
	\referee{Next we prove assertion (ii).}
	Recall from Section~\ref{sec:bij-inf} that the UIPT is defined as the local limit of a uniformly random rooted triangulation of size $n$ (local limit around its root-edge). Now, taking a uniformly random rooted triangulation of size $n$ and re-rooting it at a uniformly random non-top edge again gives a uniformly random triangulation. Moreover, the total variation distance between $\frk e_n$ and a uniformly random edge goes to 0, so the total variation distance between $M_n$ and a uniformly random rooted loopless triangulation goes to 0. Hence, it is immediate from assertion (i) of the theorem that the map $M$ has the law of the loopless UIPT. We verified above that $w$ has steps chosen uniformly and independently at random, and the remaining claims of assertion (ii) are immediate since $(w,M,\sigma,\frk e,\tau^*,h)$ is the local limit of $(w^n,M_n,\sigma_n,\frk e_n,\tau_n^*,h^n)$. In particular, the claim that $\sigma$ is uniform follows since for any $R\in\N$ the ball of radius $R$ in $M_n$ centered at $\frk e_n$ converges in law as $n\rta\infty$, so the probability that the root-edge of $M_n$ (which has non-random colors) is contained in this ball converges to 0.
	\referee{This proves  assertion (ii) of the lemma.}
	
	\referee{We now embark on the proof of assertion (iv).}
	So far $\tau^*$ has been defined as a subset of $E(M^*)$. Now we will argue that $\tau^*$ is a spanning tree of $M^*$. Since $\tau^*$ is the limit of a tree, we know that it does not contain any cycles. To conclude that $\tau^*$ is a tree, we need to show that it is a.s.\ connected. 
	We denote by $(M^-,\sigma^-)$ the past site-percolated near-triangulation encoded by the past word $w^-$. 
	We denote by $M_-^*$ the map obtained from the dual of $M^-$ by deleting the vertex of infinite degree (which corresponds to the infinite face of $M^-$).
	Let $\tau_-^*=\{e\in\tau^*\,:\,\etae^{-1}(e)<0 \}$ be the restriction of $\tau^*$ to edges $M_-^*$.
	It is easy to see from the definition of $M^-$ that $\tau_-^*$ is a spanning tree of $M_-^*$. Indeed, consider the pairs $(P_i,\alpha_i)=\Phi(w^-_i)$ corresponding to the decomposition $w^-=\dots w^-_2 w^-_1$ with $w^-_i$ in $\mK$. It is easy to see that the restriction $\tau_{-i}^*$ of $\tau_-^*$ to $(P_i,\alpha_i)$ is a spanning tree of $P^*_i$ and the dual of the root-edge of $P_i$ is always in $\tau^*_-$, so the trees $\tau_{-i}^*$ are all connected.
	\referee{Thus $\tau^*$ is a spanning tree of $M^*$.}
	
	\referee{We now argue that $\tau^*$ is one-ended and satisfies the other claimed properties.}	
	It is clear that the subtree $\tau_-^*$ is one ended, and we consider the infinite path $P$ of $\tau^*$ from $\infty$ to the vertex $t^*$ of $M^*_-$, where $t$ is the triangle of $M^-$ incident to the top-edge.
	Now, recall that the future site-percolated triangulation $(M^+,\sigma^+)$ is defined in terms of the decomposition $w^+=w_1^+c w_2^+c\dots$ with $w_i^+$ in $\smK$. The restriction $\tau_{i}^*$ of $\tau^*$ to $(Q_i,\beta_i)=\bPhi(w_i^+)$ is a spanning tree of $Q^*_i$ and the dual of the root-edge of $Q_i$ is always in $\tau^*$, attached to $P$. This can be seen by observing that an edge $e\in E(M^*)$ is in $\tau^*$ if and only if $w_{\etae^{-1}(e^*)}\in\{a,b \}$, so by using that $(Q_i,\beta_i)=\bPhi(w_i^+)$ for $w_i^+\in\smK$ we see that the only edge on the boundary of $Q_i$ which is dual to an edge in $\tau^*$ is the root-edge. So the set $\tau^*$ is obtained from $\tau_-^*$ by gluing the finite subtrees $\tau_i^*$ to the path $P$. Hence it is a one-ended spanning tree. Moreover, since the restrictions $\tau_{-i}^*$, $\tau_{i}^*$ of $\tau^*$ are DFS trees of $P_i^*$ and $Q_i^*$, it is not hard to see that all the edges of $M^*$ join $\tau^*$-comparable vertices. In other words, $\tau^*$ is in $\DFS_{M^*}$. It is clear that $\Lambda_{M}(\tau^*)=\sigma$, since  $\Lambda_{M}(\tau_n^*)=\sigma_n$ and  $(M,\sigma,\tau^*)$ is the local limit of $(M_n,\sigma_n,\tau_n^*)$. 
	Furthermore the percolation path $P$ of $(M^-,\sigma^-)$ is contained in $\tau^*$. Hence it is clear that $\tau^*$ satisfies Properties (ii-iii) of Theorem~\ref{thm:dfs-inf} (since it is satisfied by $\tau^*_n$).

	Next, we show that $\tau^*$ satisfies the properties of Proposition~\ref{prop:dfs-ordering-inf}. 
	We first need to show that there exists a DFS of $M^*$ which follows rule (ii') of Definition~\ref{def:space-filling}, and that the associated tree is $\tau^*$. 
	By Proposition~\ref{prop32}, applying the DFS rule $(ii')$ of Proposition~\ref{prop:dfs-ordering-inf} to $(P^*_i,\al_i)$ (resp. $Q_i^*,\be_i)$) gives the DFS tree $\tau_{-i}^*$ of $P_i^*$ (resp. $\tau^*_i$ of $Q_i^*$). By applying Proposition~\ref{prop:dfs-ordering-inf} to subwords of the form $w_{k,l}=w^-_{-k}\dots w^-_1w^+_1 c\dots c w^+_{\ell}$ in $\mK$ we see that applying the DFS rule $(ii')$ of Proposition~\ref{prop:dfs-ordering-inf} to the submap of $M^*$ corresponding to $(P_i,\al_i)$ (resp. $(Q_i,\be_i)$) still gives $\tau_{-i}^*$ (resp. $\tau^*_i$). Now it is easy to see that we can concatenate all the DFS corresponding to these submaps into a DFS of $M^*$ which follows rule (ii') of Definition~\ref{def:space-filling}. Moreover, by applying Proposition~\ref{prop32}, to subwords of the form $w_{k,l}$ in $\mK$ we see that $\tau^*$ satisfies the properties of Proposition~\ref{prop:dfs-ordering-inf}.
	\referee{This concludes the proof of assertion (iv).}

	
	We obtain assertion (v) by using that the analogous result holds in the finite setting (Theorem~\ref{thm:DFS}),
	and that both $h$ and $\tau^*$ are determined from $w$ in a local way. Similarly, assertion (vi) follows easily from the fact that  $(w,M,\sigma,\frk e,\tau^*,h)$ is the local limit of $(w^n,M_n,\sigma_n,\frk e_n,\tau_n^*,h^n)$. 
\end{proof}

Observe that Properties (ii-iii) of Lemma~\ref{prop17} proves Theorem~\ref{thm:UIPT} except for its last sentence (uniqueness of $w\in\imK$ such that $\Phi^\infty(w)=(M,\si)$). Moreover, Property (iv) of Lemma~\ref{prop17} shows that for $(M,\si)\in\imTP$ chosen according to the percolated UIPT distribution, there exists almost surely $\tau^*\in\DFS_{M^*}$ such that $\Lambda_{M}(\tau^*)=\si$ (but it does not imply the uniqueness of $\tau^*$). Lemma~\ref{prop17} also proves Theorems~\ref{thm:dfs-inf},~\ref{thm:dfs-inf2}, and~\ref{thm:dfsdual-inf} as well as Proposition~\ref{prop:dfs-ordering-inf}, except for Theorem~\ref{thm:dfs-inf}(i) (uniqueness of $\tau^*$). 

\xin{We prove the last sentence of Theorem~\ref{thm:UIPT} and Theorem~\ref{thm:dfs-inf}(i) in  the appendix.
	But we emphasize that these two statements  are not needed for the rest of the paper.}

\section{The mating-of-trees correspondence}\label{sec:dictionary}

In this section, after recalling basic concepts of $\SLE_6$ and $\sqrt{8/3}$-LQG in Sections~\ref{sec:sle} and~\ref{sec:lqg}, respectively, 
we will review the mating-of-trees theorem (Theorem~\ref{thm:mot}) in the infinite volume case in Section~\ref{sec:mot}. Then in Sections~\ref{sec:mot}-\ref{sec:cont-piv}, we present the continuum analogue of the future/past decomposition and spine-looptree decomposition and use them to construct branching $\SLE_6$, $\CLE_6$, and $\LQG$ pivotal measure in the mating-of-trees framework, as the continuum analogs to the DFS tree, percolation cycles, and the counting measure on percolation pivotal points, respectively. In Section~\ref{sec:disk}, we present the finite volume variants of these constructions. Finally, in Section~\ref{sec:crossing}, we elaborate on the mating of trees aspect of percolation crossing events.
The reader is advised to look at Table \ref{table-dictionary}
while reading this section, to see the discrete-continuum correspondences. In some sense, this section can be viewed as a mini-survey for the mating-of-trees theory for $\gamma=\sqrt{8/3}$ and $\kappa=6$ based on \cite{wedges} and related papers (see references in this section). 
However, Sections~\ref{sec:cont-cle},~\ref{sec:cont-piv} and~\ref{sec:crossing}, which are the most crucial parts to this paper, are \emph{new} to the best of our knowledge.\footnote{Section~\ref{sec:cont-cle} has some overlap with an unpublished manuscript of Gwynne and Miller (in particular, the concept of envelope interval), but neither of the two perspectives we take on $\CLE_6$ loops here is explicitly explored in that manuscript.}

\subsection{Schramm-Loewner evolutions}
\label{sec:sle}

\xin{In this subsection we recall several variants  of  SLE$_6$, including chordal and radial $\SLE_6$, branching $\SLE_6$, and space-filling $\SLE_6$ defined on the disk, as well as some of their whole-plane analogs. To clarify the notation, we will use $\etab$ to represent the space-filling $\SLE_6$ since this is the most relevant  variant of $\SLE_6$ for this paper. We will use  $\wh\etab$ to represent a chordal or  radial $\SLE_6$. We will use $\wh\etab^{z}$ to denote the branch towards the point $z$ of a branching $\SLE_6$.
}

The ``classical'' SLE$_6$ is a family of random continuous curves modulo parametrization indexed by a triple $(D,x,y)$, where $D\subset \C$ is a simply connected domain whose boundary $\bdy D$ is a continuous curve, and $x\neq y$ are two points such that $x\in \bdy D, y\in \ol D$. For each triple $(D,x,y)$, $\SLE_6$ indexed by $(D,x,y)$ is a random curve in $D$ starting at $x$ and ending at $y$. $\SLE_6$ is characterized by the following three properties:
\begin{itemize}
	\item {\it Conformal invariance:} Given $(D,x,y)$ and $(D',x',y')$ and a conformal map $\varphi:D\to D'$ such that $\varphi(x)=x'$ and $\varphi(y)=y'$, suppose $\wh\etab$ has the law of an $\SLE_6$ indexed by $(D,x,y)$, then $\varphi\circ \xin{\wh\etab}$ has the law of an $\SLE_6$ indexed by $(D',x',y')$.
	
	\item {\it Domain Markov property:} Let $(D,x,y)=(\bbH,0,\infty)$ or $(D,x,y)=(\D,1,0)$, where $\bbH$ is the upper half plane and $\D$ is the unit disk centered at 0. 
	Suppose $\wh\etab$ is sampled from $\SLE_6$ indexed by $(D,x,y)$ and is under the so-called capacity parametrization\footnote{Capacity parametrization is only to make sure that $\wh \etab $ is progressively adapted, namely, for any time $t>0$ in the domain of definition of $\wh \etab$, the number $t$ can be determined purely by the segment $\wh \etab|_{[0,t]}$ modulo parametrization. Any progressively adapted parametrization will work here.
	} (see for example \cite{lawler-book}). Then for any $t>0$, conditioning on $\wh\etab([0,t])$, the law of $\wh\etab([t,\infty)$ is the $\SLE_6$ indexed by $(D_t,\wh\etab(t),y )$, where $D_t$ is the connected component of $D\setminus\xin{\wh\etab}([0,t])$ containing $y$.
	\item {\it Target invariance:} Given $(D,x,y)$ and $(D,x,y')$ with $y\neq y'$, suppose $\wh\etab$ and $\wh\etab'$ are sampled from the $\SLE_6$ indexed by $(D,x,y)$ and $(D,x,y')$, respectively.
	Let $\tau$ be the first time $\wh\etab$ separates $y$ and $y'$; more precisely, $$\tau=\inf\{t\ge 0: \textrm{$y$ and $y'$ are in different connected components of $D\setminus \wh\etab([0,t])$}\}.$$ Similarly, let $\tau'$ be the first time $\xin{\wh\etab'}$ separates $y$ and $y'$. Then the curves $\wh\etab|_{[0,\tau]}$ and $\wh\etab'|_{[0,\tau']}$ modulo parametrization have the same law. 
\end{itemize}
It is proved by Schramm \cite{schramm0} that if a family of curves satisfies the first two properties and that $W_t$ denotes the \emph{driving function for the Loewner chain} encoding the curve indexed by $(\bbH, 0,\infty)$ with the capacity parametrization, then there exists a $\kappa>0$ such that $W_t$ has the law of a linear Brownian motion with variance $\kappa$.
Moreover, the driving function of the random curve indexed by $(\D,1,0)$ as a radial Loewner chain has the law of $e^{iW_t}$. These two curves are called the \emph{chordal $\SLE_{\kappa}$ on $(\bbH,0,\infty)$} and the \emph{radial $\SLE_\kappa$ on $(\D,1,0)$}, respectively. 
Finally, if the family further satisfies the target invariance property, then $\kappa=6$. By the Riemann mapping theorem and conformal invariance, this uniquely specifies the law of $\SLE_6$ indexed by any $(D,x,y)$. When $y\in \bdy D$ (respectively, $y\in D$), we call it the \emph{chordal $\SLE_{6}$ on $(D,x,y)$} (respectively, the \emph{radial $\SLE_{6}$ on $(D,x,y)$}).

Although the Loewner chain perspective on $\SLE_6$  is important and fruitful, we will not give more details on this perspective but only refer to \cite{lawler-book}, since one can understand much of the mating-of-trees perspective on $\SLE_6$, including everything in this paper, without going into Loewner evolutions. 
Here we recall two topological properties of $\SLE_6$ which \emph{are} important to this paper. First, both chordal and radial SLE$_6$ are curves whose trace has zero Lebesgue measure, but which creates ``bubbles'' (bounded simply connected domains) by hitting its past and the domain boundary infinitely often \cite{schramm-sle}. Second, recall the domain $D_t$ in the definition of the \emph{domain Markov property} above. Both of the two arcs on $\bdy D_t$ from $\wh\etab(t)$ to $y$ are simple curves almost surely (equivalently,  all the cut points of $D_t$ separate $\wh\etab(t)$ and $y$). In fact, by the so-called SLE duality, before hitting $\bdy D$, the two curves evolve as variants of $\SLE_{8/3}$ \cite{dubedat-duality}.

Let us now define the \emph{branching $\SLE_6$} which was first introduced in \cite{shef-cle}.
Let  $D$ be a domain, let  $x\in \bdy D$ and let $y_1,\dotsm, y_k$ be distinct points in $\ol D\setminus\{x\}$. 
In light of the target invariance and domain Markov property of $\SLE_6$, the
$\SLE_6$ on $(D,x,y_i)$ ($1\le i\le k$) can be coupled together in such a way that the $k$ curves agree before the set $\{y_1,\dots,y_k \}$ is separated into two complementary connected components of the curve. After the separation time, the remaining parts of the $k$ curves evolve independently in each component in the same fashion,
until all the $k$ points are in $k$ different connected components. Then we finish the $k$ curves by running independent $\SLE_6$ curves in these $k$ components targeted at $y_1,\dotsm, y_k$, respectively. Fixing a countable dense subset $Y$ of $D$, we can extend the collection of branching coupled $\SLE_6$ curves on $(D,x,z)$ to all $z\in Y$. Write the SLE$_6$ on $(D,x,z)$ in this coupling as $\wh \etab^z$ for all $z\in Y$. For any $z \in \ol \D$, by taking a limit, we can almost surely obtain a
curve from $x$ to $z$, which we denote by $\wh \etab^{z}$. By locality of SLE$_6$, this curve has the law of an SLE$_6$ on $(D,x,z)$. The collection of curves $\{\wh \etab^z\}_{z\in \ol D\setminus\{x\} }$ is the \emph{branching $\SLE_6$} on $D$ rooted at $x$. Given any fixed $z$, the almost surely well-defined curve $\wh\etab^z$ has the marginal law of $\SLE_6$ on $(D,x,z)$.
We call $\etab^z$ the \emph{branch} of $\{\wh \etab^z\}_{z\in \ol D\setminus\{x\} }$ targeted at $z$.  \xin{As stated in \cite[Section~4.2]{shef-cle}, the law of branching SLE$_6$ does not depend on the choice of the dense set $Y$.} 
In Section~\ref{sec:disk}, we will give a version of branching SLE$_6$ on $\D$ which is almost surely defined for all points $z\in\ol D\setminus \{x\}$ simultaneously. 

More recently, a \emph{space-filling} version of SLE$_6$ has been constructed. The construction of such an object as a continuous curve is a major achievement of the theory of imaginary geometry \cite{ig1,ig2,ig3,ig4} and an essential input to the mating-of-trees theory \cite{wedges} (see Theorem~\ref{thm:mot}). We start from the chordal variant. Intuitively, a \emph{chordal space-filling $\SLE_6$ on $(D,x,y)$} $\etab$ is obtained from running a chordal $\SLE_6$ on $(D,x,y)$ such that whenever a bubble is disconnected from the target $y$,
the curve immediately fills the bubble by a space-filling $\SLE_6$ curve before heading to $y$. To make it precise, let us start from a sample of branching $\SLE_6$ $\{\wh \etab_z\}_{z\in \ol D\setminus \{x\}}$ rooted at $x$. Given two rationals $q\neq q'$ in $D$, we can run the curves $\wh\etab^{q}$ and $\wh\etab^{q'}$ until the split, which is exactly the time when they separate $q$ and $q'$ into two complementary connected components, one of which contains $y$. We write $q\prec q'$ if $q'$ and $y$ are in the same component. Then $\prec$ defines a total ordering on $\Q^2\cap D$. It is proved in \cite{ig4} that with probability 1, there exists a continuous curve $\etab$ such that $\etab$ visits $\Q^2\cap D$ in the order is given by $\prec$.  This provides a rigorous definition of $\etab$.  
\xin{In fact, as explained in \cite[Section 4.3]{ig4}, the ordering $\prec $ is equivalent to a construction using the so-called flow lines in the imaginary geometry. Then using tools from imaginary geometry, the existence of a continuous curve $\etab$ is proved in this context.}

In this paper, we will also consider a variant of space-filling SLE$_6$ where the initial and terminal point of the curve is the same. More precisely, let $D$ be a domain as above, and let $x\in \bdy D$ and $y_n\in \bdy D$ be such that $y_n$ approaches $x$ from the left side as $n\to\infty$. 
\xin{Then the   ordering $\prec_n$ defining the chordal space-filling $\SLE_6$ on $(D,x,y_n)$ weakly converges to a random ordering $\prec$ on $\Q^2\cap D$. Again this ordering $\prec$  can be constructed using imaginary geometry, from which one can see that $\prec$ defines a random  continuous space-filling curve on $D$ starting and ending at $x$ (see  \cite[Section 4.3]{ig4}).} We call this curve the \emph{counterclockwise space-filling $\SLE_6$ on $(D,x,x)$}. 
If $y_n$ approaches $x$ from the right side, then the curve obtained from the same procedure is called the \emph{clockwise space-filling $\SLE_6$ on $(D,x,x)$}.

Given a sample $\etab$ of the chordal space-filling $\SLE_6$ on $(D,x,y)$, 
we can obtain a chordal $\SLE_6$ on $(D,x,y)$ by skipping the times $t$ 
when $y$ and $\etab(t)$ are not on the boundary of the same connected component of $D\setminus \etab([0,t])$.
Conversely, we may construct an instance of a chordal space-filling SLE$_6$ $\etab$ on $(D,x,y)$ by filling in the bubbles created by a chordal $\SLE_6$ $\wh\etab$ on $(D,x,y)$ by independent space-filling SLE$_6$ curves starting and ending at the point where the bubble is enclosed.
The orientation of the SLE$_6$ curves in each bubble is opposite to the orientation of the boundary of the bubble as defined by the order of visit by $\etab$. \xin{To see this, one can check that the ordering $\prec_{x,y}$ defining $\etab$ restricted to rational points in each bubble of $\wh\etab$
	gives the total ordering that defines the SLE$_6$ curves in that bubble with the desired  orientation.}

As in the discrete, we will also consider the infinite volume version of these $\SLE_6$ curves. Let $\wh\etab^R$ be the radial $\SLE_6$ indexed by $(R\D,R,0)$. 
The \emph{whole-plane SLE$_6$}\footnote{
	More precisely, 
	this curve should be called \emph{whole-plane SLE$_6$ from $\infty$ to $0$}. By using M\"obius transform over $\C\cup \{\infty\}$, we may define \emph{whole-plane SLE$_6$ from $x$ to $y$ for any 
		$x\neq y\in \C\cup \{\infty\}$. }
}
is defined by taking the local limit around $0$ of $\wh \etab^R$ as $R\to\infty$. 
More precisely, for any fixed $R_0>0$, the law of $\wh\etab\cap R_0\D$ is the weak limit of the law of $\wh\etab^R\cap R_0\D$. 
Similarly, \emph{whole-plane branching $\SLE_6$} is defined as the local limit around $0$ of the branching $\SLE_6$ on $R \D$ rooted at $R$, and the \emph{whole-plane space-filling $\SLE_6$} is defined as the local limit around $0$ of the chordal space-filling $\SLE_6$ on $R\D$ from $R$ to $-R$. (See Figure~\ref{fig:pastwedge} for an illustration of branching and space-filling $\SLE_6$ near 0.) \xin{The existence of these  infinite volume limits  can easily be justified via imaginary geometry. See \cite[Lemma 2.3]{hs-Euclidean}
	and its proof.}  

Given an instance of the whole-plane branching $\SLE_6$, one can define a total ordering $\prec $ on $\Q^2$ similarly as above (in the case $D\subsetneq \C$) with $\infty$ playing the role of the target point $y$. This ordering almost surely defines a curve which has the law of the whole-plane space-filling $\SLE_6$. Therefore, both in the disk case and in the whole-plane case, the space-filling $\SLE_6$ is a function of the branching $\SLE_6$. As will be explained in Remark~\ref{rmk:branch}, in this coupling of space-filling $\SLE_6$ and branching $\SLE_6$, the former also determines the latter almost surely via an explicit function. Therefore, they encode the same amount of information, and they can be thought of as two different ways of representing the full scaling limit of planar critical percolation. A third way called $\CLE_6$ will be discussed in Sections~\ref{sec:cont-cle}--\ref{sec:disk}. 
A fourth way called quad crossings will be alluded in Section~\ref{sec:crossing}.

\subsection{Liouville quantum gravity}
\label{sec:lqg}

Liouville quantum gravity is a theory of random surfaces which is our continuum analogue of random planar maps. We first recall that the Gaussian free field (GFF) \cite{shef-gff} is a random planar Gaussian distribution which is defined as follows. Let $D\subsetneq\C$ be a simply connected domain, let $C_0^\infty(D)$ be the set of $C^\infty$ functions compactly supported in $D$, and let $H_0(D)$ be the Hilbert space closure of $C_0^\infty(D)$ equipped with the \emph{Dirichlet inner product}:
\eqbn
(f,g)_\nabla = \frac{1}{2\pi} \int_D \nabla f(z)\cdot \nabla g(z)\,dz,\qquad f,g\in C_0^\infty(D).
\eqen
Let $(\phi_n)_{n\in\N}$ be an orthonormal basis for $H_0(D)$.\footnote{Precisely, $(\phi_n)_{n\in\N}$ constitute a basis in the sense of Hilbert spaces, that is, every element of $H_0(D)$ can be written uniquely as an \emph{infinite} linear combination $\sum_{n} c_n\phi_n$ with $\sum_{n} c_n^2<\infty$.} The \emph{zero-boundary GFF} $\gff$ on $D$ can be expressed as a random linear combination of these basis elements
\eqbn
\gff = \sum_{n=1}^{\infty} \alpha_n\phi_n
,\qquad \textrm{where $\alpha_n$ are independent standard normal variables}.
\eqen
This sum does not converge in $H_0(D)$ but does converge almost surely in $H^{-1}(\D)$, the Sobolev space with index $-1$ (\cite[Section~4.2]{dubedat-coupling}).
In particular, $\gff$ is a random distribution. 
For $f\in H_0(D)$, and $c_n=(\phi_n,f)_\nabla$ for $n\in\N$,  the random series $\sum_n  c_n \al_n$ converges in $L^2$ to a Gaussian variable with mean zero, which we denote by $(\gff,f)_\nabla$.
For all $f,g\in H_0(D)$ the covariance between $(\gff,f)_\nabla$ and  $(\gff,g)_\nabla$ is $(f,g)_\nabla$.  By integration by parts, for each $f\in C_0^\infty(D)$ 
we have $(\gff,\Delta f)=-2\pi (\gff,f)_\nabla$ almost surely.
The \emph{free-boundary GFF on $D$} may be defined similarly by replacing $(\phi_n)_{n\in\N}$ by an orthonormal basis for the Hilbert space closure ${H}(D)$ of the set of functions $f$ in $C^\infty(D)$ with finite \emph{Dirichlet energy} $(f,f)_\nabla$.

Let $\gff$ be a GFF in some domain $D\subset\C$. For $\gamma\in(0,2)$, $\gamma$-Liouville quantum gravity ($\gamma$-LQG) may be heuristically defined as the Riemannian manifold with metric tensor given by the Euclidean metric tensor times $e^{\gamma \gff}$. This heuristic definition of $\gamma$-LQG does not make literal sense since $\gff$ is a distribution and not a function. However, one can make sense of the \emph{area measure} (\xin{see e.g.\ \cite[Proposition 1.1]{shef-kpz}})
\begin{equation}\label{eq:area}
	\mub=\mub_{\gff}=e^{\gamma\gff} dxdy
\end{equation}
by regularizing the field. Consider the area measure $\mub_\eps=e^{\gamma \gff_\eps}dxdy$, where $dxdy$ is the Lebesgue area measure and $\gff_\eps$ is a ``regularized version'' of $\gff$ (for example $\gff_\eps=(\gff,f_{\eps,z})+0.5\gamma\ln(\eps)$, where $f_{\eps,z}$ is a positive smooth function supported on $B_\eps(z)$ such that $\int_{B_\eps(z)}f_{\eps,z}(w)\,dw=1$).
One can show that for many definitions of $\gff_\eps$ this area measure converges as $\eps\rta 0$, and that the limiting measure is independent of the exact definition of $\gff_\eps$; \xin{see e.g.\ \cite{berestycki-gmt-elementary}}. For a free-boundary GFF $\gff$ in some domain $D\subsetneq\C$, the field $\gff$ also induces a \emph{length measure} along $\partial D$. 
If $\partial D$ contains a line segment $\ell$, then the LQG boundary length measure along $\ell$ is defined heuristically as the measure \xin{(see e.g.\ \cite[Section~6.3]{shef-kpz})}
\begin{equation}\label{eq:length}
	\nub=\nub_{\gff}=e^{\ga \gff/2}d\ell,
\end{equation} 
where $d\ell$ is the Euclidean length (although, as for the LQG area, a regularization procedure is needed to make this definition rigorous). It is clear that the above measures can also be defined for certain other fields, for example fields which can be written as the sum of a GFF and a continuous function (possibly with logarithmic singularities). 
See \cite{shef-kpz,rhodes-vargas-review,berestycki-gmt-elementary,shef-zipper,wedges,aru-survey} for further details.

Consider the collection of pairs $(D,\gff)$, where $D\subset\C$ is simply connected and $\gff$ is a distribution defined on $D$. We say that two pairs $(D,\gff)$ and $(\wt D,\wt \gff)$ are \emph{$\ga$-equivalent} if there is a conformal map $\phi:\wt D\to D$ such that 
\begin{equation}\label{eq:Q}
	\wt \gff=\gff\circ\phi+Q\log|\phi'| \qquad \textrm{where} \; Q=2/\gamma+\gamma/2. 
\end{equation}
A $\ga$\emph{-LQG surface} is a \emph{$\ga$-equivalence class} of pairs $(D,\gff)$ for $\gff$ a Gaussian free field or a related kind of distribution. 
The LQG area measure is preserved when applying the coordinate change formula, in the sense that almost surely, for any open set $A\subset\wt D$ it holds that $\mub_{\gff}(\phi(A))=\mub_{\wt \gff}(A)$. The LQG boundary measure is preserved in a similar way. 
\xin{To see why $Q=2/\gamma+\gamma/2$ makes the LQG  area and length preserved, we first note that if $h$ were a smooth function then $Q$ should be replaced by $2/\gamma$ in order to preserve the area and length. The extra $\gamma/2$ in $Q$ comes from the fact that $h$ is rough and we have to perform a regularized limit. See e.g.\ \cite[Section 2]{shef-kpz}.}

In this paper \xin{we focus on the scaling limit of uniform triangulations, which corresponds to LQG with $\ga=\sqrt{8/3}$}.
we will primarily be considering three $\ga$-LQG surfaces for $\ga=\sqrt{8/3}$: the \emph{$\sqrt{8/3}$-LQG cone}, the \emph{$\sqrt{8/3}$-LQG sphere}, and the \emph{$\sqrt{8/3}$-LQG disk}, which we now define. 

Consider a zero-boundary GFF $\gff$ on $\D$. In short, the $\sqrt{8/3}$-LQG cone is the local limit around a point sampled from the LQG area measure of the $\sqrt{8/3}$-LQG surface $(\D,\gff+C)$ for $C\to +\infty$. 
However, \emph{convergence as LQG surfaces} is a subtle notion of convergence,
as $\mub_{\gff+C}$ clearly does not have limit as a random measure. 
To obtain a meaningful limit, we follow \cite[Appendix A]{wedges}.
Conditioning on $\gff$, let $z_0\in \D$ be sampled from the probability measure $\mub_{\gff}/\mub_{\gff}(\D)$. For $C>0$, let $R$ be chosen such that $\mub_{\gff^C} (\D)=1$, where $\gff^C(z)=\gff(R^{-1}(z+z_0)))+C-Q\log R$. \xin{(Note that $\gff^C$ is defined as $\wt \gff$ in~\eqref{eq:Q} with $\phi(z)=R^{-1}(z+z_0)$.)}
Then $\gff^C$ converges in law as a random distribution to a random distribution $\gff$ on $\C$, which can be written as the so-called whole-plane GFF plus a \xin{random} radially symmetric continuous function. 
Moreover, $\mub_{\gff^C}$ converges in law to the random measure $\mub_{\gff}$  on $\C$ associated with $\gff$ as in~\eqref{eq:area} such that $\mub_{\gff}(\D)=1$.
The morale behind this procedure is that one specifies a representative in the equivalence class (that is, the \emph{embedded} surface) consistently before taking the limit.
In the example of the quantum cone studied here we choose the representative which gives unit mass to $\D$, and in the limit we obtain an embedded surface which also gives unit mass to $\D$. We could have chosen another representative, whose limit would be a different representative of the same $\sqrt{8/3}$-LQG surface. In this paper, when we refer to a \emph{$\sqrt{8/3}$-LQG cone}, we always refer to the particular representative $(\C,\gff)$ which gives unit mass to $\D$. \xin{This is merely for concreteness.
	Results in this paper  either does not depend on  the embedding, or once it is stated clearly for one embedding, it can be transferred to a statement  under another embedding straightforwardly according to~\eqref{eq:Q}.}

The \emph{unit boundary length $\sqrt{8/3}$-LQG disk} can be constructed similarly by a limit \cite[Proposition 5.10]{wedges}. Consider a smooth bounded domain $D$ with a linear segment $\mathrm{L}$ of $\partial D$. Let $\gff$ be a GFF on $D$ with free-boundary conditions on $\mathrm{L}$ and zero-boundary conditions on $\partial D\setminus \mathrm{L}$.
\xin{Then \eqref{eq:length} defines the LQG boundary measure  $\nub$ on $\mathrm{L}$.}
Fix $C,\eps>0$ and condition $\gff$ on the event that  $\{\sqrt{C}\leq \nub_{\gff}(\mathrm{L})\leq \sqrt{C}(1+\eps) \}$. Let $\wh \gff=\gff-\sqrt{8/3}\log C$, which gives 
$1\leq \nub_{\xin{\wh \gff}}(\mathrm{L})\leq 1+\eps$. Then $(D,\wh \gff)$ converges as a $\sqrt{8/3}$-LQG surface to the unit boundary length $\sqrt{8/3}$-LQG disk when we first send $C\rta\infty$ and then send $\eps\rta 0$. 
In this limit, the segment \xin{$\bdy D\setminus \mathrm{L}$} will collapse to a boundary marked point on the surface. For concreteness we will specify a representative $(\D,\gffd)$ of the unit boundary length $\sqrt{8/3}$-LQG disk in Section~\ref{sec:disk}, where $1$ is the marked point. Given  $(\D,\gffd)$ and a random variable $C$ independent of $\gffd$, the LQG surface $(\D, \gffd+C)$ is called the \emph{$\sqrt{8/3}$-LQG disk with boundary length $L:=e^{\frac12\sqrt{8/3} C}$}. For a fixed $A>0$, conditioning on $\mub_{\gffd+C}(\D)=A$, the resulting surface is called 
the \emph{$\sqrt{8/3}$-LQG disk with boundary length $L$ and area $A$}. 
\xin{The legitimacy of the degenerate conditioning on $\mub_{\gffd+C}(\D)=A$ is explained below \cite[Defition 4.21]{wedges} via the general theory of regular conditional probability and a scaling argument.} 

The \emph{unit area $\sqrt{8/3}$-LQG sphere} can be constructed by a limiting procedure in a similar manner as the disk \xin{\cite[Appendix~A.2]{wedges}}. We will provide such a limiting construction in Section~\ref{sec:disk}, where we will also specify a particular representative $(\C\cup \{\infty\},\gffs)$ for concreteness.

\subsection{Mating of trees}
\label{sec:mot}
In this subsection we present the \emph{mating-of-trees} construction by Duplantier, Miller, and Sheffield \cite{wedges}, which can be thought as the continuum analogue of (the different versions of) our bijection~$\Phi$. This construction relates LQG surfaces to planar Brownian motion. 

We will use the following terminology for Brownian motion. 
For $\beta>0$, a \emph{two-sided linear Brownian motion with variance $\beta$} is a stochastic process $(\Xb_t)_{t\in \R}$ indexed by $\R$ such that $(\Xb_t)_{t\ge 0}$ and $(\Xb_{-t})_{t\ge 0}$ are two independent standard linear Brownian motions multiplied by $\sqrt\beta$. For $\alpha\in (-1,1)$ and $\beta>0$ a stochastic process $\wt\Zb=(\wt\Lb_t,\wt\Rb_t)_{t\in\R}$ is called a \emph{planar two-sided Brownian motion with correlation $\alpha$ and variance $\beta$} if
$\wt \Lb+\wt \Rb$ and $\wt \Lb- \wt \Rb$ are two independent two-sided linear Brownian motions with variance $2(1+\alpha)\beta$ and $2(1-\alpha)\beta$, respectively. 

Let us now recall the mating-of-trees construction for LQG cones. Consider a pair $(\gff,\etab)$, where $\gff$ is the field on $\C$ associated with a $\sqrt{8/3}$-LQG cone and $\etab$ is a whole-plane space-filling SLE$_6$ which is independent with $\gff$ as a curve modulo parametrization. 
We parametrize $\etab$ such that $\etab(0)=0$ and $\mub(\etab([s,t]))=t-s$ for any $s<t$, where $\mub=\mub_h$ is the area measure associated with $\gff$.
At any time $t$ the boundary of $\etab((-\infty,t])$ can uniquely be written as the union of two semi-infinite non-crossing curves $\xi^{\op{L}}_t,\xi^{\op{R}}_t$ ending at $\etab(t)$, such that $\xi^{\op{L}}_t,\xi^{\op{R}}_t$ evolve continuously in $t$. The curves $\xi^{\op{L}}_t,\xi^{\op{R}}_t$ are called the \emph{left frontier} and the \emph{right frontier}, respectively, of $\etab((-\infty,t])$, with the convention that the set $\etab((-\infty,t])$ is on the right of the left frontier when this frontier is oriented toward $\etab(t)$; see Figure \ref{fig:pastwedge}.
By SLE duality, for any $t\in\R$ each of the two frontiers of $\etab((-\infty,t])$ is an SLE$_{8/3}$-like curve.\footnote{More precisely, these frontiers have the law of a whole-plane SLE$_{8/3}(-2/3)$ \cite{ig4}.} 
By \cite{shef-zipper,wedges} we can define an LQG length measure along such curves by using the definition of the LQG boundary measure. More precisely, we consider a conformal map from the complement of the SLE$_{8/3}$-type curve such that the curve is mapped to a straight line, and we define the LQG length measure by considering the pullback of the LQG boundary length measure as defined in Section \ref{sec:lqg}. Let $\Zb=(\Lb_t,\Rb_t)_{t\in\R}$ be defined such that $\Lb_t$ (respectively, $\Rb_t$) is the length of the left (respectively, right) frontier of $\etab((-\infty,t])$, relative to the length of the frontier at time 0. We call $\Zb$ the \emph{boundary length process} of $(\gff,\etab)$.
\begin{theorem}\xin{$\mathrm{(\cite[Theorem~8.1]{wedges})}$}
	Let $(\gff,\etab)$ and $\Zb$ be as above.
	The stochastic process $\Zb$ is well-defined, and there exists $\beta>0$ such that $\Zb$ has the law of a planar two-sided Brownian motion with correlation $1/2$ and variance $\beta$.
	Furthermore, the pair $(\gff,\etab)$ modulo rotation about the origin is measurable with respect to the $\sigma$-algebra generated by $\Zb$.
	\label{thm:mot}
\end{theorem}
\begin{remark}\label{rmk:variance}
	The variance of the planar two-sided Brownian motion $\Zb$ is not evaluated explicitly in \cite{wedges}. Throughout the paper we use $\beta$ to denote this unknown constant. 
	The third named author is working on a project with Gwynne and Remy on evaluating the constant explicitly by combining the constructive field theory \cite{dkrv-lqg-sphere} and the mating-of-trees approaches to LQG.
\end{remark}

\begin{figure}
	\centering
	\includegraphics[scale=1]{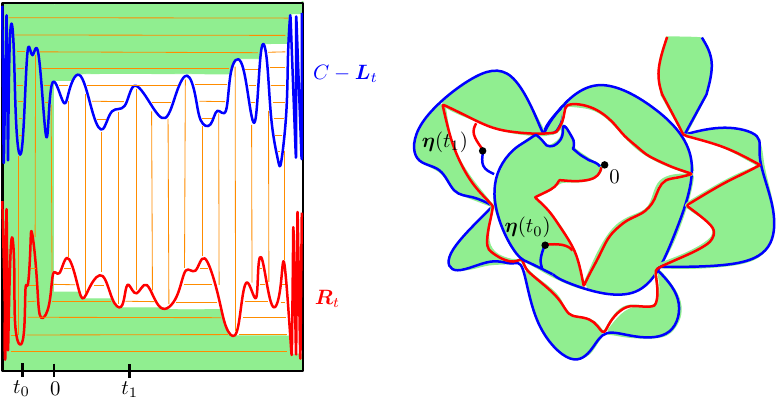}
	\caption{
		Left: A segment of the Brownian motion $\Zb$ in Theorem~\ref{thm:mot}. The $\Lb$-coordinate is drawn upside down by considering the graph of $C-\Lb$ for a large constant $C$. 
		Right: The whole plane space-filling $\SLE_6$ drawn at times \xin{$t_0<0<t_1$} (recall that $\etab(0)=0$). 
		\xin{Here the red (resp., blue) curves represent the right (resp., left) frontier of $\etab$ at times $t_0$, $0$, and $t_1$.}
	}
	\label{fig:mating}
\end{figure}
Theorem \ref{thm:mot} gives a construction of LQG and SLE often referred to as the \emph{mating-of-trees}. This name comes from the following interpretation of the construction. Observe that the left frontier of $\etab(-\infty, t]$ at different $t$'s are rays from $\etab(t)$ to $\infty$ that merge upon intersection. This produces a ``spanning tree'' of $\C$ rooted at $\infty$. Another ``dual spanning tree'' can be obtained by consider the right frontier. Then $\etab$ can be thought of as the Peano curve snaking in between the two spanning trees. 
The LQG boundary length and area measures endow the two trees with a metric-measure structure. 
In this point of view, Theorem~\ref{thm:mot} says that the two trees are both continuum random trees (CRT) in the sense of Aldous \cite{aldous-crt1,aldous-crt2,aldous-crt3}. 
Moreover, their contour function $\Lb$ and $\Rb$ form a two-sided planar Brownian motion with correlation $1/2$. (See Figure~\ref{fig:mating}.)
However, this pair of CRT's should not be confused  with the continuum analogue of the DFS tree that we will introduce in Section~\ref{sec:dictionary-branch}. 

\xin{\begin{remark}[Recover $(\gff,\etab)$ from $\Zb$]
		We have explicitly described how $\Zb$ is determined by $(\gff,\etab)$. The fact that $(\gff,\etab)$ is determined by $\Zb$ modulo rotation was first proved in~\cite[Section~8]{wedges} without giving an explicit procedure. Such a procedure is recently provided through~\cite{gms-tutte}. The idea is to consider each segment $\etab[i/n, (i-1)/n]$  of $\etab$ as a cell on $\C$, 
		whose adjacency relation gives a planar map called the mated-CRT map. It is easy to see that the mated-CRT map is determined by $\Zb$ explicitly. The main result of~\cite{gms-tutte} asserts that under the so-called Tutte embedding, which only depends on the map itself, the (normalized) vertex counting measure of the mated-CRT map along with the ordered cells converges to $(\mu_{\gff},\etab)$ in probability. Since the LQG area measure  $\mu_{\gff}$ determines $\gff$~\cite{bss-lqg-gff}, we recover $(\gff,\etab) $ from $\Zb$.
\end{remark}}

There are also ``finite volume versions'' of Theorem \ref{thm:mot} corresponding to the LQG disk and the LQG sphere \cite[\xin{Appendix A}]{wedges}  \cite{sphere-constructions}, which we will present in Section~\ref{sec:disk}. 
We also remark that there is an analogue of Theorem \ref{thm:mot} (and of its finite volume variants) for space-filling SLE$_\kappa$ and the $\ga$-LQG cone for arbitrary $\kappa>4$ and $\ga=4/\sqrt{\kappa}$-LQG.
The correlation of the planar Brownian motion in the general case is $-\cos((4\pi)/\kappa)$ \cite{wedges,kappa8-cov}.

We will now present the continuum analogue of the future/past decomposition presented in Section~\ref{subsec:future-past}. 
Let $(\gff,\etab)$ and $\Zb$ be as in Theorem~\ref{thm:mot} .
We call the closed set $\etab(-\infty,0] \subset\C$ the \emph{past wedge} of $(\gff,\etab)$ (relative of 0). The open set $\C\setminus \etab((-\infty,0]))$ is called the \emph{future wedge} (relative to 0). See Figure~\ref{fig:pastwedge} for an illustration. 
It is clear that the past wedge and the closure of the future wedge have the same law. However, they are presented in a asymmetric manner since, in analogy to the discrete, it is more instructive to think of the past wedge as a connected set with a curve in it (see Section~\ref{sec:dictionary-fl}) and the future wedge as a chain of simply connected domains. To illustrate this point of view for the future wedge,
let us define the set 
$$\cut(0)=\{t>0: \Lb_s>\Lb_t\;\textrm{and}\; \Rb_s>\Rb_t\;\textrm{for all}\; s\in (0,t) \}.$$ 
The connected components in the future wedge are exactly the collection the interiors of $\etab([s,t])$ for $(s,t)$ a connected component of $(0,\infty)\setminus \cut(0)$. The set $\cut(0)$ has the law of the range of a stable subordinator, hence it has a local time \cite{bertoin-sub}. Restricting $(\gff,\etab)$ to the components in the future wedge and parametrizing them by the local time of $\cut(0)$, the future wedge can be viewed as a Poisson point process of $\sqrt{8/3}$-LQG surfaces decorated with space-filling curves from one boundary point to another. Conditioning on the area and length of these surfaces, they are independent $\sqrt{8/3}$-LQG disks \cite{wedges} decorated with chordal space-filling $\SLE_6$. We refer to \cite{gwynne-miller-char,gwynne-miller-sle6} for a detailed account of the mating-of-trees theory for these curve-decorated $\LQG$ surfaces. We need some theory from these works in Section~\ref{sec:cont-cle}, which we will review then.	
\begin{figure}
	\centering
	\includegraphics[scale=1]{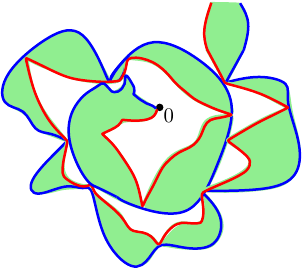}\quad
	\includegraphics[scale=1]{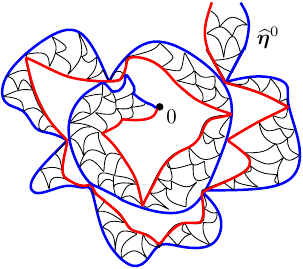}\quad
	\includegraphics[scale=1]{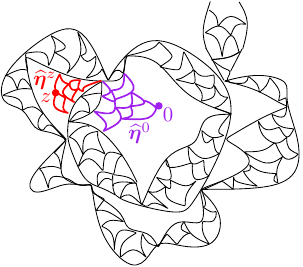}
	\caption{
		Left: The past wedge is shown in green, and the blue (respectively, red) curve from 0 to $\infty$ is the left (respectively, right) frontier of $\etab((-\infty,0])$. The blue and red curves starting from $\etab(t_0)$ and $\etab(t_1)$, respectively, show the left and right frontier of $\etab((-\infty,t_0 ])$ and $\etab((-\infty,t_1 ])$ (with $t_0<0<t_1$).
		Middle: The branch $\wh\etab^0$ of the branching SLE$_6$ $\taub^*$ from $\infty$ to 0 has the law of a whole-plane SLE$_6$. The frontier of this curve (in blue and red) defines the boundary of the past wedge.
		Right: The whole-plane SLE$_6$ $\wh\etab^0$ (respectively, $\wh\etab^z$) is the concatenation of the black curve and the purple (respectively, red) curve. The two branches $\wh\etab^0$ and $\wh\etab^z$ of the branching SLE$_6$ $\taub^*$ agree until the target points 0 and $z$ lie in different complementary components of the curve.
	}
	\label{fig:pastwedge}
\end{figure}

\subsection{Branching SLE$_6$}
\label{sec:dictionary-branch}
We now present the continuum analogue of the exploration tree defined in Section~\ref{sec:DFS}.

Let us first introduce some notations related to cone excursions. For $t_1 < t_2$, a one-dimensional path $(\wt\Lb_t)_{t\in[t_1,t_2]}$ is called an \emph{excursion} if $\wt\Lb_{t_1} = \wt\Lb_{t_2} < \wt\Lb_t$ for all $t\in (t_1,t_2)$. A
two-dimensional path $(\wt\Lb_t,\wt\Rb_t)_{t\in[t_1,t_2]}$ is called a \emph{cone excursion} if $\wt\Lb_t > \wt\Lb_{t_2}$ and $\wt\Rb_t > \wt\Rb_{t_2}$ for all
$t\in (t_1, t_2)$, and if either $\wt\Lb_{t_1} = \wt\Lb_{t_2}$ or $\wt\Rb_{t_1} = \wt\Rb_{t_2}$. In particular, note that one of the coordinates $\wt\Lb$
and $\wt\Rb$ define an excursion, while the other coordinate attains a running infimum at time $t_2$.

Let $(\gff,\etab)$ and $\Zb$ be as in Theorem \ref{thm:mot}. A \emph{cone interval} is an interval $I=[t_1,t_2]$ such that $\Zb|_I$ is a cone excursion.  
By definition of cone intervals and elementary properties of $\Zb$, it is easy to see  that with probability 1, for any pairs of cone intervals $I=[s,t]$ and $J=[s',t']$, we have 
\begin{equation}\label{eq:cone-int}
	I \cap J =\emptyset\quad\textrm{or}\quad I \subset J \quad\textrm{or}\quad J \subset I. 
\end{equation}

If $\Lb_{t_1}>\Lb_{t_2}$ (respectively, $\Rb_{t_1}>\Rb_{t_2}$), then $I$ is called a \emph{left} (respectively, \emph{right}) cone interval and $\Zb|_I$ is called a \emph{left} (respectively, \emph{right}) cone excursion. 
If $s>t$ is such that for all $t'\in(t,s)$, $\Lb_{t'}>\Lb_{s}$ and $\Rb_{t'}>\Rb_{s}$, then we say that $s$ is an \emph{ancestor} of $t$. For $u\in \R$, we say that a  time $t<u$ is \emph{ancestor-free} relative to time $u$ if it has no ancestors in $[t,u]$. 
The set of ancestor-free times relative to time $u$ is denoted by $\ans(u)$. Equivalently,
$$\ans(u)=\{t<u~|~\textrm{there exists no cone interval } [t_1,t_2]\subset(-\infty,u) \textrm{ such that }t\in(t_1,t_2)\}.$$
The set $\etab(\ans(0))$ is the trace of a continuous curve inside the past wedge relative to 0, denoted by $\wh\etab^0$ (see Figure~\ref{fig:pastwedge}). The curve $\wh\etab^0$ has the law of a whole-plane $\SLE_6$ from $\infty$ to 0. In each connected component of the interior of the past wedge, $\wh \etab^0$ evolves as a chordal $\SLE_6$ conditioned on  \xin{touching every point on} the domain boundary (see \cite{ig4} for a precise definition of this variant of $\SLE_6$). 

Given any point $z\in \C$, let $t^z=\sup\{t\in \R:\etab(t)=z\}$. Then $\etab(\ans(t_z))$ is the trace of a continuous curve in the past wedge relative to time $t_z$. This curve is from $\infty$ to $z$, and we denote it by $\wh\etab^z$. See Figure~\ref{fig:pastwedge}.
\begin{definition}\label{def:branch}
	The collection of curves $\taub^*:=\{\wh \etab^z\}_{z\in \C}$ is called the \emph{branching $\SLE_6$} associated with $(\gff,\etab)$. The curve $\etab^z$ is called the \emph{branch} of $\taub^*$ targeted at $z$.
\end{definition}
\begin{remark}
	\label{rmk:branch}
	As explained in \cite[\xin{Section~10.2}]{wedges}, for a fixed point $z\in \C$, one can get $\wh\etab^z$ almost surely by running $\etab$ from $-\infty$ to $z$ and skipping all the times $t$ when $\etab(t)$ and $z$ are not at the boundary of the same connected component of $\C\setminus \etab((-\infty,t])$. 
	In particular, $\taub^*$ modulo parametrization can be determined by $\etab$ modulo parametrization without any reference to $\gff$.
	Moreover, combined with \cite{ig4}, the law of $\taub^*$ is the whole-plane branching $\SLE_6$ as defined in Section~\ref{sec:sle}. 
	To be precise, $\taub^*$ is a \emph{version} of the whole-plane branching $\SLE_6$ as a stochastic process indexed by $\C$ in the sense that 
	they share the same finite marginal distribution. 
	Finally, the coupling of $\etab$ and $\taub^*$ coincides with the one described in Section~\ref{sec:sle}. In particular, in this coupling, modulo parametrization $\etab$ and $\taub^*$ determine each other. 
\end{remark}
Now we assume that $u$ is a fixed deterministic time. As explained in \cite[\xin{Section~10.2}]{wedges}, the set $-\ans(0):=\{t:-t\in\ans(0) \}$ has the law of a range of a stable subordinator. In particular, 
it is possible to define the \emph{local time} $\ellb^u=(\ellb^u_t)_{t\le 0}$ of $\ans(u)$, such that $\ellb^u(u)=0$ and $\ellb^u$ is a continuous and non-decreasing process which is constant on intervals disjoint from $\ans(u)$. This local time is only uniquely defined up to a multiplicative constant, which we will fix in Section \ref{sec:walk} (see Remark~\ref{rmk:constant}). 
Note that by our convention that $\ellb^u(u)=0$, $\ellb^u_t<0$ for all $t<u$. By definition, $d\ellb^u$ induces a measure supported on $\ans(u)$. Let  
\begin{equation}\label{eq:Tu}
	\Tb^u_t=\inf \{s\le u: \ell^u_s>t\}\qquad \textrm{for }t\le 0.
\end{equation} 
We call $(\Tb^u_t)_{t\le 0}$ the \emph{inverse} of $\ellb^u$. 
Now for the point $z=\etab(u)$, it is almost surely the case that $t_z=u$.
We parametrize $\wh\etab^z$ by requiring $\wh\etab^z(t)=\etab(\Tb^u_t)$.
This parametrization is called \emph{quantum natural parametrization} of $\wh\etab^u$ \cite[\xin{Section~6.5}]{wedges}. 

For each fixed $u$, we have defined the quantum natural parametrization of the branch of $\taub^*$ targeted at $\etab(u)$. This definition works well if $u$ is replaced by another backward stopping time of $\Zb$. However, we will not define the quantum natural parametrization for all branches of $\taub^*$ simultaneously. For an arbitrary random time the set $\ans(u)$ is challenging to understand. For example, for any fixed $z\in \C\setminus\{0\}$, the time $t_z$ is not a backward stopping time and $\Zb|_{[-\infty,t_z]}$ does not evolve as a Brownian motion. So the above definition via local time of the range of stable subordinators fails. 
However, in Section~\ref{sec:cont-cle}, we do manage to extend the definition of local time of $\ans(t)$ for a particular kind of random times called the envelope stopping time, and use it to define a parametrization of $\CLE_6$.
There is yet another point of view on the quantum natural parametrization, which is as a Gaussian multiplicative chaos over the Minkowski content of $\SLE$ curves 
(see \cite{benoist-lqg-chaos,hs-Euclidean}). 
This perspective may allow the definition of the quantum natural parametrization for all branches of $\taub^*$. But we will not pursue this direction as it is not needed for our work.

\subsection{Looptree and forested lines}
\label{sec:dictionary-fl}
We now present the continuum analogue of the discrete looptrees and the spine-looptrees decomposition defined in of Section \ref{sec:spine}. \xin{This subsection is based on \cite[Section~10.2]{wedges}.}

Let $(\gff,\etab)$ and $\Zb$ be as in Theorem \ref{thm:mot}. Fix a time $u\in \R$. Here and in the rest of the paper, whenever the reference time $u$ is clear from the context, we will drop the $u$-dependence of objects viewed from $u$. For example, we will write $\ellb, \Tb,\wh\etab$ respectively for the local time $\ellb^u$, its inverse $\Tb^u$  and the curve $\wh \etab^{z}$ for $z=\etab(u)$.
With this notation, define $\wh\Zb=(\wh\Lb_t,\wh\Rb_t)_{t\leq 0}$ by $\wh\Zb_t=\Zb_{\Tb_t}$. By \cite[\xin{Section~10.2}]{wedges}, if $u$ is a deterministic time, then (the c\'adl\'ag modification of) the processes $(\wh\Lb_{-t})_{t\ge 0}$ and $(\wh\Rb_{-t})_{t\ge 0}$ are independent $3/2$-stable L\'evy processes with only positive  jumps. We call $\wh\Zb$ \emph{the L\'evy process relative to time $u$}, although for a random $u$ the associated $\wh\Zb$ is not necessarily distributed as a L\'evy process. In this section we assume that $u$ is a deterministic time and describe the spine-looptrees decomposition of the past wedge relative to time $u$. In Section~\ref{sec:cont-cle}, the time $u$ will be a particular random variable.

We start by explaining the concept of \emph{looptree} as introduced in  \cite{curien-kortchemski-looptree-def}. 
A looptree is a metric space which is a continuum analogue of the discrete looptrees described in Section \ref{sec:LR-decomposition-map}. 
Heuristically, a looptree can be thought as the metric space obtained by arranging circles of various lengths, called \emph{bubbles}, in a tree-like fashion. However, there are additional difficulties in the continuum, and one cannot really define looptrees by gluing some bubbles together. As we recall now, looptrees can be associated to L\'evy excursions in a way which is analogous to the clockwise-code of discrete looptree given by Lemma~\ref{lem:cwcode}.

We call \emph{backward infimum time for $\wh\Lb$} a time $t\leq 0$, such that $\wh \Lb_{t}<\wh \Lb_{s}$ for all $s\in(t,0]$.
To a backward infimum time $t$ for $\wh \Lb$, we associate a time 
\begin{equation}\label{eq:st}
	\frk s_t:=\inf\{ s\leq t\,:\,\wh\Lb_{s'}\ge \wh\Lb_{t}\,\,\forall s'\in[s,t] \}. 
\end{equation}
Let $\ul{\wh\Lb}_t=\inf_{s\in [t,0]} \wh\Lb_s$ for all $t\leq 0$. The set $\{(\frk s_t,t) : t\; \textrm{is a backward infimum time for}\; \wh\Lb \;\textrm{and}\; \frk s_t\neq t\}$ is almost surely equals the set of connected components of $(-\infty,0)\setminus\{s:\ul{\wh\Lb}_s={\wh\Lb}_s \}$. Let $t$ and $\frk s_t$ be as above, and let $\Xb :=(\wh\Lb-\ul{\wh\Lb})|_{[\frk s_t,t]}$ (we drop the dependence of $\Xb$ on $t$ for simplicity). By the fluctuation theory of L\'evy process \cite{bertoin-book}, $\Xb$ has the law of a $\frac32$-stable L\'evy excursion with a random duration $[\frk s_t,t]$ (see \cite{curien-kortchemski-looptree-def} for more details on stable L\'evy excursions).

Given the excursion $\Xb$, one consider  the equivalence relation on $[\frk s_t,t]$ defined by $t_1\sim_{\Xb} t_2$ if and only if 
\begin{equation}\label{eq:d=0}
	\Xb_{t_1}=\Xb_{t_2}=\inf_{s\in [t_1,t_2]} \Xb_{s}.
\end{equation}
It was shown in \cite{curien-kortchemski-looptree-def}, that one can almost surely associate to the excursion $\Xb$ a metric $d_{\Xb}$ on the quotient space $([\frk s_t, t]/\sim_{\Xb}, d_{\Xb})$ which makes it a compact metric space. This metric space is called the  (metric) \emph{looptree} associated with $\Xb$ and is denoted by $\frk L_{\Xb}$. Let $\pi_{\Xb}$ be the quotient map from $[\frk s_t,t]$ to $\frk L_{\Xb}=[\frk s_t,t]/\sim_{\Xb}$. We call $\pi_{\Xb}(\frk s_t)=\pi_{\Xb}(t)$ the \emph{root} of $\frk L_{\Xb}$.

\begin{figure}
	\centering
	\includegraphics[scale=1]{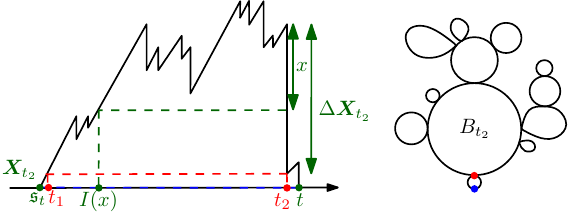}
	\caption{\xin{Left: The figure illustrates the excursion $\Xb$ and various times associated to it. 
			Right: The figure illustrates the looptree $\frk L_{\Xb}$ corresponding to $\Xb$.  The blue point is the root of $\frk L_{\Xb}$. The red point corresponds to $\pi_{\Xb}(t_1)=\pi_{\Xb}(t_2)$, which is the base point of $B_{t_2}$.}}
	\label{fig:looptree-I}
\end{figure}

We will not give the explicit formula of $d_{\Xb}$ (see \cite[Equation (2.5)]{curien-kortchemski-looptree-def}) since it is quite notationally involved and is not relevant to our later discussion.\footnote{In fact, the paper \cite{curien-kortchemski-looptree-def} consider stable excursions $\Xb^{\mathrm{exc}}$ with only positive jumps. To put our excursion $\Xb$ in their framework and define $d_{\Xb}$, one takes $\Xb^{\mathrm{exc}}$ to be the c\'adl\'ag modification of the time reversal of $\Xb$ and then plugs $\Xb^{\mathrm{exc}}$ into \cite{curien-kortchemski-looptree-def}. The equivalence relation $\sim_{\Xb}$ is extracted from \cite[Equation (2.5)]{curien-kortchemski-looptree-def} under this identification.} But $\frk L_{\Xb}$ as a compact connected topological space (which we would call the \emph{topological looptree}) is explicitly prescribed by~\eqref{eq:d=0} (in other words, the metric space $(\frk L_{\Xb},d_{\Xb})$ is homeomorphic to  $[\frk s_t, t]/\sim_{\Xb}$). 
We will give a variational characterization of $d_{\Xb}$ now.
For all $s\in [\frk s_t,t]$, let $\Delta \Xb_s=\Xb_{s^-}-\Xb_s$. Suppose  $t_2$ is a jumping time for $\Xb$ so that $\Delta\Xb_{t_2}>0$. Let $t_1$ be the almost surely unique time $t_1$ such that $\eqref{eq:d=0}$ holds. For each $x\in (0,\Delta \Xb_{t_2}]$, 
let $I(x)=\sup\{s<t_2: X_{s} -X_{t_2}=\Delta \Xb_{t_2}-x\}$. Then $I(\Delta \Xb_{t_2})=t_1$. By convention we assume $I(0)=t_2$. Then by the definition of $\sim_{\Xb}$, the mapping $\pi_{\Xb} \circ I:[0,\Delta \Xb_{t_2}] \to \frk L_{\Xb}$ defines a cycle embedded in $\frk L_{\Xb}$, which we call the \emph{bubble} of $\frk L_{\Xb}$ associated with the jump $\Delta\Xb_{t_2}$ and denote by $B_{t_2}$. \xin{See Figure~\ref{fig:looptree-I}.}
We call $\pi_{\Xb}(t_1)=\pi_{\Xb}(t_2)$ the \emph{base point} of $B_{t_2}$. Moreover, the metric of $\frk L_t$ restricted to $B_{t_2}$ is simply given by 
\begin{equation}\label{eq:dx}
	d_{\Xb}(\pi_{\Xb}\circ I(x), \pi_{\Xb} \circ I(y)) =\min \{ |x-y|, \Delta\Xb_{t_2} -|x-y|\},
\end{equation} which means that $d_{\Xb}$ restricted to $B_{t_2}$ is a metric cycle of length $\Delta \Xb_{t_2}$. Moreover, $\pi_{\Xb}\circ I$ specifies this isometry after quotienting the endpoints of $[0,\Delta \Xb_{t_2}]$.
This associates every jump of $\Xb$ with a bubble in $\frk L_{\Xb}$. One way to characterize $d_{\Xb}$ on $\frk L_{\Xb}$ is to say that it is the \emph{smallest} metric on the topological looptree $\frk L_{\Xb}$ such that $(\frk L_{\Xb},d_{\Xb})$ is a geodesic metric space (i.e., every pair of points can be joined by an arc
isometric to a compact interval of the real line) and~\eqref{eq:dx} holds.
Here by \emph{smallest} we mean given any metric $d$ with the said properties on the $\frk L_{\Xb}$ and two points $p_1,p_2$ on $\frk L_{\Xb}$, we have $d(p_1,p_2)\ge d_{\Xb}(p_1,p_2)$ 

It is elementary to see that, almost surely for $\Xb$, it does not exist three times $\frk s_t\le t_1<t_2<t_3\le \frk t$ such that $d_{\Xb}(t_1,t_2)=d_{\Xb}(t_2,t_3)=0$.
Equivalently, almost surely the cardinality of $\pi_{\Xb}^{-1}(p)$ is $\le 2$ for all $p\in \frk L_{\Xb}$.
We call $p$ a \emph{double point} of $\frk L_{\Xb}$ if $\frk L_{\Xb}\setminus\{p\}$ has at least two connected components.
We denote the set of doubles points of $\frk L_{\Xb}$ by $\dbl_{\Xb}$.
By definition, $\pi_{\Xb}(t_2)$ is  a base point of a bubble in $\frk L_{\Xb}$ if and only $\Delta X_{t_2}>0$. 
This countable set of base points of bubbles only occupy a tiny portion of $\dbl_{\Xb}$.
In fact, a non-atomic measure will be defined on $ \dbl_{\Xb}$ in Section~\ref{sec:cont-piv}.

Given $\frk s_t\le t_1<t_2\le t$ such that $\Xb$ have jumps at $t_1,t_2$, let $B_1,B_2$ be the corresponding bubbles on $\frk L_{\Xb}.$ 
We say that $B_2$ is an \emph{ancestor} of $B_1$ in the looptree $\frk L_{\Xb}$ if and only if $\Xb_{t_2}<\inf_{s\in[t_1,t_2^-]}\Xb_s$. 
This defines a tree structure on bubbles. Intuitively, the looptree is obtained by ``gluing'' together its bubbles according to this tree structure. However in the looptree $\frk L_{\Xb}$ the gluing is not so straightforward as one can check that almost surely any two bubbles are disjoint. The gluing is achieved by filling the gaps between bubbles with ``dust" (which are points not on any bubbles) so that we obtain a compact connected space.

Although our discussion above is for the particular L\'evy excursion $\Xb$ on $[\frk s_t,t]$ and its looptree $\frk L_{\Xb}$, the discussion directly extends to any process $\Xb'$ whose law is a $\frac32$-stable L\'evy excursion with only negative jumps. Hence we can associate a looptree $\frk L_{\Xb'}$ to any such excursion. In particular, given a jump time $t_2$ and $\Xb_{t_1}=\Xb_{t_2}$ as in~\eqref{eq:d=0}, running $\Xb$ backward from $t_2$ to $t_1$, we obtain a countable collection of L\'evy excursions away from the running infimum of $\Xb$ relative to time $t_2^-$, each of which gives a subspace of $\frk L_{\Xb}$ which itself is a looptree with root on the bubble corresponding to $\Delta\Xb_{t_2}$.

\begin{figure}
	\centering
	\includegraphics[scale=1]{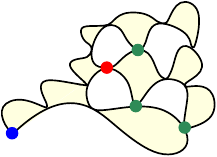}
	\caption{\xin{The figure illustrates a looptree $\frk L_{\Xb}$ embedded in $\C$.  The blue point is the root. The red point corresponds to a point $\wh \etab(t_1)=\wh\etab(t_2)$ for $t_1\sim_{\Xb}t_2$.  For each green marked point, although it is a single point on $\C$, it corresponds to two points on  $\frk L_{\Xb}$ touching each other under the embedding $\phi_{\Xb}$.  The yellow bubbles are independent $\sqrt{8/3}$-LQG disks filling $\frk L_{\Xb}$ to make it an LQG looptree.}}
	\label{fig:looptree-z}
\end{figure}

So far $\frk L_{\Xb}$ is only an compact metric space determined by the excursion $\Xb$.
We now describe  $\frk L_{\Xb}$ as an embedded topological space in the $\sqrt{8/3}$-LQG cone. By Theorem~\ref{thm:mot} the boundary length process $Z$ is well-defined as a continuous function. \xin{Recall that  $\Xb=(\wh\Lb-\ul{\wh\Lb})|_{[\frk s_t,t]}$.}
By the definition of $Z$ as the boundary length process,
$t_1\sim_{\Xb}t_2$ implies that  $\wh\etab(t_1)=\wh\etab(t_2)$. Therefore there exists a unique continuous map $\phi_{\Xb}$ from $\frk L_{\Xb}$ to $\C$ such that $\wh\etab(s)=\phi_{\Xb}\circ \pi_{\Xb}(s)$ for all $s\in [\frk s_t,t]$. In other words, the curve $\wh\etab([\frk s_t,t])$ is the image of $\frk L_{\Xb}$ under the embedding $\phi_{\Xb}$. Below we describe the geometric properties of $\frk L_{\Xb}$ under this embedding.

Consider a jump $\Delta \Xb_s$ of $\Xb$ where $s\in (\frk s_t,t)$. By definition of $\wh \Zb$, the jump $\Delta \Xb_s$ corresponds to a cone excursion of $\Zb$. Precisely,  denoting $\Tb_{s^-}:=\lim_{r\to s,r<s}\Tb_r$, the walk $\Zb|_{[\Tb_{s^-},\Tb_s]}$ is 
a left cone excursion. In particular, the set $\etab([\Tb_{s^-},\Tb_s])$ is the closure of a simply connected domain whose boundary is a Jordan curve. 
Reasoning as above, the image by $\phi_X$ of the bubble $B_s$ of $\frk L_{\Xb}$ corresponding to $\Delta\Xb_s$ is the boundary of the set $\etab([\Tb_{s-},\Tb_s])$ and the base point is mapped to $\etab(\Tb_{s-})=\etab(\Tb_s)$. Moreover, $\phi_{\Xb}$ restricted to $B_s$ is an isometry of $\partial \etab([\Tb_{s-},\Tb_s])$ parametrized by the $\sqrt{8/3}$-LQG boundary measure. In particular, the total length of this bubble equals the size $\Delta \Xb_s$ of the jump. Since $\frk L_{\Xb}$ is embedded in $\C$, we call the interior of $\etab([\Tb_{s^-}, \Tb_s])$ the \emph{interior} of $B_s$. 
The field $\gff$ for the $\sqrt{8/3}$-LQG cone restricted to the interior of $B_s$
defines a $\sqrt{8/3}$-LQG surface with total area $\Tb_s-\Tb_{s^-}$. Conditioning on boundary lengths and areas,
the collection of $\sqrt{8/3}$-LQG surfaces inside each bubble of $\frk L_{\Xb}$
are independent $\sqrt{8/3}$-LQG disks. We call the looptree $\frk L_{\Xb}$ along with the $\sqrt{8/3}$-disks inside each bubble the \emph{LQG looptree} associated with $\Xb$ and still denote it by $\frk L_{\Xb}$. 

By our convention of left and right frontiers, $\wh \etab|_{[\frk s_t,t]}$
can be viewed as a curve tracing $\frk L_{\Xb}$ in counterclockwise direction. For any $s\in [\frk s_t,t]$, by the definition of $\wh \Lb$, the quantum length of the counterclockwise arc from $\wh \etab(s)$ to the root of $\frk L_{\Xb}$ on the outer boundary of $\wh\etab([\frk s_t, s])$ equals $\Xb_s$. We call $\Xb$ the \emph{counterclockwise code} of $\frk L_{\Xb}$. 
For $\wh\Rb$, we can similarly define the backward infimum times, the countable collection of L\'evy excursion, and their associated looptrees. The only difference is that now the curve $\wh\etab$ traces these looptrees in clockwise direction and hence the L\'evy excursion should be considered the \emph{clockwise code} of the looptree.

\bigskip

\begin{figure}
	\centering
	\includegraphics[scale=1]{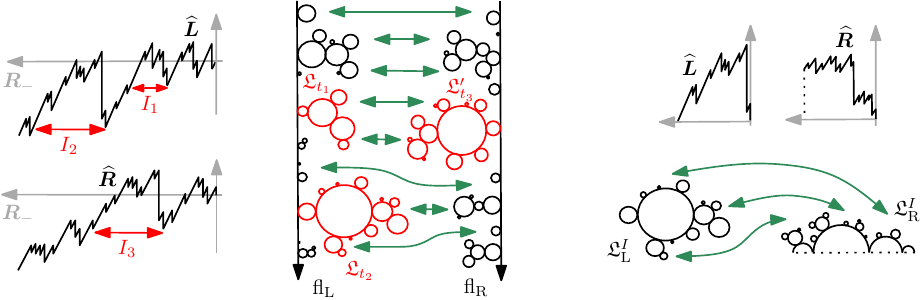}
	\caption{Left: The figure illustrates \xin{$\fll$ and $\flr$} in the spine-looptrees decomposition of the past wedge relative to time~$u$. This is to be compared with the discrete encoding in Figures~\ref{fig:LR-theorem2} and~\ref{fig:LR-inf}. Right: The figure illustrates $(\wh\Lb,\wh\Rb)$ for a $\CLE_6$ loop restricted to an envelope interval and the associated looptrees.
	}\label{fig:past-wedge}
\end{figure}

We now turn our attention to the concept of forested line. Recall that $\wh\Lb$ is the first coordinate of the L\'evy process relative to time $u$.
We define an equivalence relation $\sim_{\Lb}$ on $(-\infty,0]$ in the same spirit as $\sim_{\Xb}$ above. 
For each $t_1 \le t_2\le 0$, $t_1\sim_{\Lb} t_2$ if and only if $\wh \Lb_{t_1}=\wh \Lb_{t_2}=\inf_{s\in [t_1,t_2]} \wh \Lb_{s}$.
Restricted to any interval $[\frk s_t, t]$ such that $t$ is a backward infimum time of $\wh\Lb$ and $\frk s_t\neq t$, the equivalence relation $\sim_{\Lb}$ is the same as $\sim_{\Xb}$.
Let $\fll$ be the quotient space $(-\infty,0]/\sim_{\Lb}$ and let $\pi_{\Lb}$ be the quotient map. Then $\pi_{\Lb}([\frk s_t,t])$ is the looptree $\frk L_{\Xb}$ as a subspace of $\fll$.
Given $x\in (-\infty, 0]$, let $I_{\Lb}(x):=\sup\{s: \wh\Lb_s=x \}$. Then $\pi \circ I_{\Lb}$ defines a continuous ray in $\fll$ emanating from $\pi_{\Lb}(0)$ which consists of the roots of all the looptrees $\frk L_{\Xb}$ as well as points that are not in any of these looptrees. Let $d_{\Lb}$ be the \emph{smallest} metric (in the same sense as for $d_{\Xb}$ above) on $\fll$ such that
\begin{compactitem}
	\item[1.] the embedding from each $\frk L_{\Xb}$ to $\fll$ is isometric, and
	\item[2.] the embedding $\pi_{\Lb}\circ I_{\Lb}$ from $(-\infty,0]$ to $\fll$ is isometric.
\end{compactitem}
The topology of $\fll$ easily implies that such $d_{\Lb}$ exists. Replacing $\wh\Lb$ with $\wh\Rb$, we can define the equivalence relation $\sim_{\Rb}$, the topological space $\flr=(-\infty,0]/\sim_{\Rb}$, the quotient map $\pi_{\Rb}$, and the metric $d_{\Rb}$. \xin{We  reiterate that  $\sim_{\Lb}$ and  $\sim_{\Rb}$, hence  $\fll$ and  $\flr$, depend on $u$.} See Figure~\ref{fig:past-wedge} for an illustration of $\fll$ and $\flr$.

\begin{definition}
	The metric space $\fll$ under $d_{\Lb}$ is called the \emph{left forested line} relative to time $u$. The metric space $\flr$ under $d_{\Rb}$ is called the \emph{right forested line} relative to time $u$. 
\end{definition}

As in the looptree case, both $\fll$ and $\flr$ are naturally embedded in the $\sqrt{8/3}$-LQG cone through~$\wh\etab$.
Note that for any $t_1<t_2\le 0$, if $t_1\sim_{\Lb} t_2$ or $t_1\sim_{\Rb} t_2$, then $\etab (t_1) =\etab(t_2)$. Therefore there exist unique embeddings $\phi_{\Lb}:\fll\to \C$ and $\phi_{\Rb}: \flr\to \C$ such that $\wh \etab= \phi_{\Lb}\circ \pi_{\Lb}=\phi_{\Rb}\circ \pi_{\Rb}$. The curve $\wh \etab$ can be viewed as the image of $\fll$ and $\flr$ under $\phi_{\Lb}$ and $\phi_{\Rb}$, respectively. The embedding $\phi_{\Lb}$ restricted to each individual looptree $\frk L_{\Xb}$ coincide with the embedding $\phi$ defined above, which allows us to define the LQG surface structure of $\fll$ by thinking of each $\frk L_{\Xb}$ as an LQG looptree. The same holds for $\flr$. 
The image of the ray $\pi_{\Lb}\circ I_{\Lb}$ (respectively, $\pi_{\Rb}\circ I_{\Rb}$) under $\phi_{\Lb}$ (respectively, $\pi_{\Rb}$) is the left (respectively, right) frontier of the past wedge relative to $u$, with the pushforward of the metric on the ray agreeing with the $\sqrt{8/3}$-LQG boundary length on the frontier. 
The curve $\wh \etab$ goes from $\infty$ to $0$ and traces bubbles in $\fll$ (respectively, $\flr$) in counterclockwise (respectively, clockwise) direction. Therefore, we may call $\wh \Lb$ (respectively, $\wh \Rb$) the \emph{counterclockwise} (respectively, \emph{clockwise}) \emph{code} of $\fll$ (respectively, $\flr$).

We now discuss the relation between the double points of $\fll$, $\flr$, and $\wh \etab$. Recall that $\wh \etab$ is the image of $\fll$ under $\phi_{\Lb}$, but that the embedding $\phi_{\Lb}$ is not injective. 
To simplify the notation, for each $s\in(-\infty,0]$ we denote the point $\pi_{\Lb}(s)$ (respectively, $\pi_{\Rb}(s)$) on $\fll$ (respectively, $\flr$) by $\fll(s)$ (respectively, $\flr(s)$). Since SLE$_6$ has no triple points almost surely \cite[Remark 5.3]{miller-wu-dim}, there are no points  on $\fll$ (respectively, $\flr$) whose pre-image under $\pi_{\Lb}$ (respectively, $\pi_{\Rb}$) has cardinality strictly larger than 2.
Let $\dbl_{\Lb}$ (respectively, $\dbl_{\Lb}$) be the set of points on $\fll$ (respectively, $\flr$) whose pre-image under $\pi_{\Lb}$ (respectively, $\pi_{\Rb}$) has cardinality exactly~2.  
Then $p\in \dbl_{\Lb}$ if and only if $p$ is a double point of a looptree in $\fll$ or the unique point of the looptree on the infinite line. 
Using the fact that $\Lb$ and $\Rb$ are independent $\frac32$-stable process, we have that for $t_1<t_2\le 0$, $\fll(t_1) = \fll(t_2)$ implies $\flr(t_1) \neq \flr(t_2)$ almost surely, while $\flr(t_1) = \flr(t_2)$ implies $\fll(t_1) \neq \fll(t_2)$ almost surely.
Furthermore, 
$\fll(t_1) \neq \fll(t_2)$ and $\wh\etab(t_1)=\wh \etab(t_2)$ if and only if $\flr(t_1)=\flr(t_2)$, and the same statement holds if we swap left and right. 
Let $\dbl_{\wh \etab}$ be the set of double points of $\wh \etab$, that is, points whose pre-image under $\wh\etab^{-1}$ is has cardinality 2.
Then 
\begin{equation}\label{eq:double}
	\dbl_{\wh \etab}=\phi_{\Lb} (\dbl_{\Lb} ) \cup \phi_{\Rb} (\dbl_{\Rb} )\quad
	\textmd{and}\quad
	\phi_{\Lb} (\dbl_{\Lb} ) \cap \phi_{\Rb} (\dbl_{\Rb})=\emptyset.
\end{equation}

We conclude this section by the promised spine-looptrees decomposition of the past wedge. Cutting the past wedge $\etab( (-\infty,u])$ along the ``spine'' $\wh\etab$ decomposes the wedge into a left part and a right part. The left (respectively, right) part is the forested line $\fll$ (respectively, $\flr$) where each bubble is attached with a $\sqrt{8/3}$-LQG quantum disk according to the boundary measure such that each looptree in $\fll$ becomes a LQG looptree. We call the left (respectively, right) part \emph{the left (respectively, right) LQG forested line relative to $u$}. Conversely, the past wedge can be obtained by gluing together these two LQG forested lines by identifying  the points $\fll(s)$ and $\flr(s)$ to $\wh\etab(s)$ in such a way that  $\wh\etab$ becomes the interface between $\fll$ and $\flr$.

\subsection{Envelope intervals and CLE$_6$}
\label{sec:cont-cle}	
The \emph{conformal loop ensemble} CLE$_6$ $\Gab$ is a random countable collection of \emph{loops} \cite{shef-cle}. The \emph{loops} are closed SLE$_6$-like curves which cannot cross themselves and each other. However, they may touch each other and be nested, and their union form a dense subset of the considered domain. They describe the scaling limit of the collection of percolation cycles for critical percolation on the triangular lattice \cite{camia-newman-full}. 
Let $(\gff,\etab)$ and $\Zb$ be as in Theorem \ref{thm:mot}.  In this section we describe the $\CLE_6$ 
associated with $(\gff,\etab)$ and $\Zb$ in analogy to Section~\ref{subsec:env}.

Given $\Zb$, an \emph{envelope interval} is a cone interval (see Section 6.4) $[t_1,t_2]$ with the particular property that for some $\eps>0$ there are no cone intervals $J$ such that $[t_1+\eps, t_2-\eps] \subset J \subsetneq [t_1,t_2]$. It can be shown via elementary Brownian motion argument that with probability 1, any interval $[t_1,t_2]$ with the property that $\Lb_{t_1}=\Lb_{t_2} =\min_{t\in [t_1,t_2]} \Lb_t$ and $\Rb_{t}\ge \Rb_{t_2}$ 
for all $t\in [t_1,t_2]$ is a cone interval. The same holds if $\Lb$ is replaced by $\Rb$. By~\eqref{eq:cone-int} and the preceding,
for any fixed finite interval $I$, the intersection of all cone intervals enclosing $I$ is almost surely a cone interval. 
Therefore, there exists a  smallest cone interval, which we denote by $\env(I)$, 
such that almost surely $I\subsetneq \env(I)$ and hence $\env(I)$ is an envelope interval.
Conversely, any envelope interval is of the form $\env(I)$ for some interval $I=[s,t]$ with $s,t\in \Q$. In particular, there are countably many  envelope intervals. 

In this section we explain how to associate a loop $\gab$ with an envelope interval $[t_1,t_2]$. We will provide two perspectives, both of which are useful in our proofs in Section~\ref{app:conv}.
The first perspective is to view $\gab$ as a concatenation of two $\SLE_6$ curves, one in past wedge and one in the future wedge relative to a fixed time $u$.
This perspective allows us to provide a definition of $\CLE_6$ via Theorem~\ref{thm:mot}. The second perspective is to view $\gab$ as an embedded looptree, which will be useful for the study of pivotal measures in Section~\ref{sec:cont-piv}.\\
\begin{figure}
	\centering
	\includegraphics[scale=1]{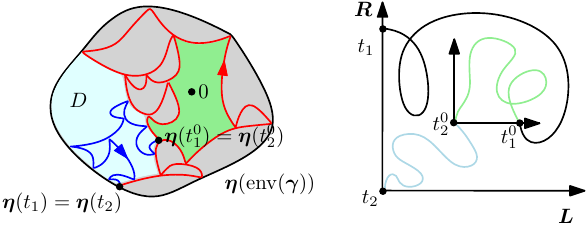}
	\caption{
		Left: construction of a $\CLE_6$ loop $\gab$ as the concatenation of a past segment $\gab_1$ (red curve) and a future segment $\gab_2$ (blue curve) relative to 0. The light blue region $D$ is the component of the future wedge containing the future segment. The green region is the connected component of $\etab(\env(\gab))\setminus \gab$ containing 0.
		Right: The cone excursions corresponding to $\env(\gab)$ and $[t_1^0,t_2^0]$. Since in this example these two intervals are of opposite types, we have $0=\etab(0)\in \reg(\gab)$.}\label{fig:cle-2parts} 
\end{figure}

We start with the first perspective, which is illustrated in Figure \ref{fig:cle-2parts}. 
Let $[t_1,t_2]$ be an envelope interval. Without loss of generality, we assume that $0\in (t_1,t_2)$, and we consider the past wedge and the future wedge relative to time 0.
Recall the time set $\cut(0)$ defined in Section~\ref{sec:mot}. Let
\begin{equation}\label{eq:t2}
	t_2^0=\sup\{ t< t_2: t\in \cut(0) \}.
\end{equation}
Then by the definition of envelope intervals, we have $t_2\in \cut(0)$ and $t_2^0<t_2$. Hence, the interior of $\etab([t_2^0,t_2])$, which we denote by $D$,
must be a connected component of the future wedge relative to time $0$.  It is shown in \cite{gwynne-miller-char} that conditioning on the total boundary lengths of the two arcs of $\partial D$ from $\etab(t_2^0)$ to $\etab(t_2)$, $(D,\gff,\etab|_{[t_2^0,t_2]})$ is a $\sqrt{8/3}$-LQG disk with specified boundary lengths, decorated with a chordal space-filling $\SLE_6$ on $D$ from $\etab(t_2^0)$ to $\etab(t_2)$. By definition of $\Zb$, the subwalk $\Zb|_{[t_2^0,t_2]}$ describes the LQG length evolution of the left and right arcs of the boundary of $D\setminus \etab([t_2^0,t])$ between the points $\etab(t)$ and $\etab(t_2)$ as $t$ runs from $t_2^0$ to $t_2$.

Consider the set $\ans(t_2)\cap [t_2^0,t_2]$, which are the ancestor-free time relative to $t_2$ in $[t_2^0,t_2]$. As explained in  \cite[Section 7.4]{gwynne-miller-char}, which is based on \cite{wedges}, although $t_2$ is a random time, it is possible to define a local time $\ellb$ on $\ans(t_2)\cap [t_2^0,t_2]$ in the same manner as in Section~\ref{sec:dictionary-branch} for a fixed time. Namely, there exists a process $\ellb$ on $[t_2^0,t_2]$ 
such that $\ellb(0)=0$ and $\ellb$ is a continuous and non-decreasing process which is constant on intervals disjoint from $\ans(t_2)$. Again $\ellb$ is uniquely determined up to a multiplicative constant and we fix it to be consistent with the convention when $t_2$ is a fixed time. More precisely, for any rational $u\in [t^0_2,t_2]$, the Stieltjes measures $d\ellb^u$ and $d\ellb$ agree on $\ans(u)\cap \ans(t_2)\cap [t_2^0,t_2]$. Let $\Tb_t=\inf \{s\le t_2: \ellb_s>t\}$ be the inverse of $\ellb$ as in Section~\ref{sec:dictionary-branch}, and for each $\ellb_{t_2^0} \le t\le 0$, 
let $\wh \etab(t)=\etab(\Tb_t)$. Then, $\gab_2:=\wh \etab|_{[\ellb_{t_2^0} , 0]}$ is the classical chordal $\SLE_6$ corresponding to the space-filling chordal $\SLE_6$ $\etab|_{[t_2^0,t_2]}$. The $\SLE_6$ curve $\gab_2$ is the segment of the curve $\gab$ lying in the future wedge relative to $0$.

To complete the definition of $\gab$ we need to define its segment lying in the past wedge.
Let
\begin{equation}\label{eq:t1}
	t_1^0=\sup \{ t\le 0: \Lb_t=\Lb_{t_2^0}\;\textrm{or}\; \Rb _{t}=\Rb_{t_2^0} \}.
\end{equation}
Since $t_1^0$ is a backward stopping time for $\Zb$, the two processes $\Zb|_{(-\infty,t_1^0]}$ and $\Zb|_{(-\infty,0]}$ have the same law. Note that $\ans (0)\cap (-\infty, t_1^0]= \ans(t_2)\cap (-\infty,t_2^0)$. 
We can extend the definition of the local time $\ellb$ defined above for $[t_2^0,t_2]$ to the entire $(-\infty,t_2]$ in such a way that $\ellb$ is a continuous and non-decreasing process which is constant on intervals disjoint from $\ans(t_2)$ (in particular $\ell_{t_1^0}=\ell_{t_2^0}$). The multiplicative constant for $\ellb$ is fixed by requiring that $d\ellb$ agrees with $d\ellb^0$ on $\ans (0)\cap (-\infty, t_1^0]$. We also extend $\Tb_t=\inf \{s\le t_2: \ellb_s>t\}$ and $\wh\etab(t)=\etab(\Tb_t)$ to all $t\in (-\infty,t_2]$. The curve $\wh\etab$ is a called the \emph{quantum natural parametrization} of the branch of $\taub^*$ targeted at $\etab(t_2)$. 

We call the curve $\gab=\wh \etab|_{[\ellb_{t_1}, 0 ]}$   the \emph{$\CLE_6$ loop associated with the envelope interval $[t_1,t_2]$}. We write $[t_1,t_2]$ as $\env(\gab)$ and call $t_2$ the \emph{envelope closing time} of $\gab$, which is the last time when $\etab$ visits
$\gab$. We say that the \emph{orientation} of $\gab$ is clockwise (respectively, counterclockwise) if $\gab$ traces the boundary of the unbounded component of $\C\setminus \gab$ in clockwise (respectively, counterclockwise) direction. 
The \emph{region enclosed by $\gab$}, denoted by $\reg(\gab)$ is defined to be the union of all bounded connected components of $\C\setminus \gab$ whose boundary is traced by $\gab$ in the orientation of $\gamma$.
As promised, the segment $\gab_1:=\wh \etab|_{[\ellb_{t_1}, \ellb_{t_1^0} ]}$ is the segment of $\gab$ lying in  the past wedge relative to 0, while  $\gab_2:=\wh \etab|_{[\ellb_{t_2^0},0]}$ is the segment lying in  the future wedge.
If we choose another rational $u\in (t_1,t_2)$, the definition of $\ellb$ and $\gab$ almost surely stay the same, although the decomposition of $\gab$ into future and past segments could change. 

\begin{figure}
	\centering
	\includegraphics[scale=1]{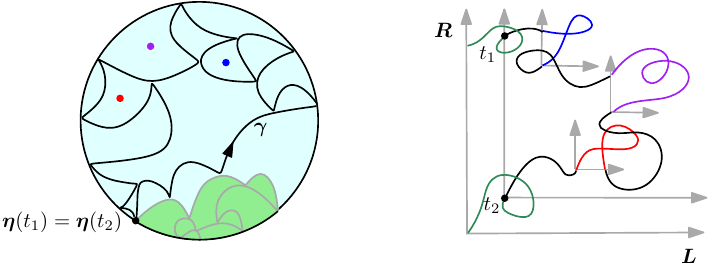}
	\caption{
		Left: A $\CLE_6$ loop $\gab$. The light blue region is filled by $\etab(\env(\gab))$. 
		The curve $\gab$ is oriented in counterclockwise direction and $\env(\gab)$ is a right cone interval.
		The three colored points are three different reference points one can use to decompose $\gab$ into a future and a past segment.
		The red point is inside $\reg(\gab)$ while the other two are not. 
		Right: The cone excursion $[t_1,t_2]=\env(\gab)$ is the middle segment in non-green colors. Each of the three smaller cone excursions with colors corresponds to the connected component of $\etab(\env(\gab))\setminus\gab$ containing the marked point in the left figure with the same color. }
	\label{fig:env1}
\end{figure}

Note that $[t_1^0,t_2^0]$ is a cone excursion inside $\env(\gab)$ containing 0. Moreover, $[t_1^0,t_2^0]$ is \emph{maximal} inside $\env(\gab)$ 
in the sense that if $J$ is a cone interval with $[t_1^0,t_2^0]\subsetneqq J\subset \env(\gab)$, we must have $J=\env(\gab)$.
In general, for a deterministic time $u\in [t_1,t_2]$, the maximal cone interval $[t_1^u,t_2^u]$ inside $[t_1,t_2]$ containing $u$ can be almost surely found in the same manner as $[t_1^0,t_2^0]$ (see~\eqref{eq:t2} and~\eqref{eq:t1}). The interior of $\etab([t_1^u,t_2^u])$ is the connected component of $\etab(\env(\gab))\setminus \gab$ containing $\etab(u)$. 
If $\env(\gab)$ is a left (respectively, right) cone interval, then $\gab$ is clockwise (respectively, counterclockwise) oriented.
If $[t_1^u,t_2^u]$ is a left (respectively, right) cone interval, then $\gab$ visits $\bdy \etab([t_1^u,t_2^u])$ in counterclockwise (respectively, clockwise)  direction.
In particular, for any rational $u\in \env(\gab)$, $\etab(u)\subset \reg (\gab )$ if and only if $[t_1^u,t_2^u] $ and $\env(\gab)$ are of different types.
See Figure \ref{fig:env1} for an illustration of $\gab$ and $\env(\gab)$.

\begin{definition}\label{def:CLE}
	Associating each envelope interval of $\Zb$ with a loop as above gives a countable collection of loops on $\C$ which we call the $\CLE_6$ associated with $(\gff,\etab)$ and denote by $\Gab$. Moreover, each curve has a parametrization which we still call the \emph{quantum natural parametrization}.
\end{definition}

The discrete analogue of this construction of $\CLE_6$ is explained in the proof of Lemma~\ref{prop21}. In Section \ref{sec:disk}, it will be shown that the CLE$_6$ of Definition \ref{def:CLE} coincides (modulo parametrization) with the CLE$_6$ classically defined as in \cite{shef-cle}.

Let $\gab'\neq \gab$ be two $\CLE_6$ loops. 
By~\eqref{eq:cone-int}, we have almost surely
\begin{equation}\label{eq:env}
	\env(\gab) \cap \env(\gab') =\emptyset\quad\textrm{or}\quad\env(\gab) \subset \env(\gab') \quad\textrm{or}\quad\env(\gab') \subset \env(\gab). 
\end{equation}
We claim that
\begin{equation}\label{eq:intersection}
	\env(\gab) \subset \env(\gab') \quad \textrm{implies}\quad \gab\cap \gab' \subset \bdy\etab(\env(\gab)).
\end{equation}
In fact, as explained above Definition~\ref{def:CLE}, any connected component of $\etab(\env(\gab')) \setminus \gab'$ is  the interior of the image by $\etab$ of a maximal cone intervals. 
In light of~\eqref{eq:cone-int},  $\env(\gab) \subset \env(\gab') $ implies 
that  $\etab(\env(\gab))$  is contained in the closure of a connected component of $\etab(\env(\gab')) \setminus \gab'$. Therefore, any touching of $\gab$ and $\gab'$ must be on $\bdy\etab(\env(\gab))$.
\begin{figure}
	\centering
	\includegraphics[scale=1]{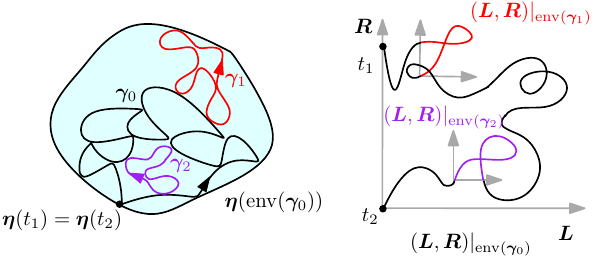}
	\caption{
		Left: Three $\CLE_6$ loops $\gab_0,\gab_1,\gab_2$, where $\gab_1,\gab_0$ are unnested, while $\gab_2,\gab_0$ are nested. The blue region is $\etab(\env(\gab_0))$.
		Right:  Cone excursions corresponding to $\env(\gab_i)$ for $i=0,1,2$. Both $\env(\gab_1)$ and $\env(\gab_2)$ are subintervals of $\env(\gab_0)$. 
	}
	\label{fig:2loops}
\end{figure}

Given $\gab,\gab'\in \Gab$ 
such that $\gab\neq \gab'$, we say that $\gab$ and $\gab'$ are \emph{nested} if either $\gab\subset \reg(\gab')$ or $\gab'\subset\reg(\gab)$. Otherwise, we say that  $\gab$ and $\gab'$ are \emph{unnested}. Recall the relation  explained above between the orientation of $\CLE_6$ loops and the type of the cone excursions. 
Suppose $\gab\cap \gab'\neq \emptyset$. If $\etab(\env(\gab)) \subset \reg(\gab')$, then the orientations of $\gab$ and $\gab'$ are opposite. 
If $\etab(\env(\gab)) \subset \etab(\env(\gab'))\setminus\reg(\gab')$, then the orientations of $\gab$ and $\gab'$ are the same. 
Therefore when $\gab\cap \gab'\neq\emptyset$, 
they have the same orientation if and only if they are unnested. See Figure~\ref{fig:2loops}.

\bigskip

Now we build on the previous description of $\CLE_6$ to establish our second perspective on $\CLE_6$, which is based on looptrees.
Let $\gab$ be a $\CLE_6$ loop. Without loss of generality, we assume $0\in \env(\gab)$. Recall $\ellb, \Tb,\wh\etab$ defined above as the local time, its inverse, and the branch in $\taub^*$ relative to $t_2$.
As in Section~\ref{sec:dictionary-fl}, for all $t\le 0$, let $\wh\Zb_t=(\wh\Lb_t,\wh\Rb_t):=\Zb_{\Tb_t}$ be the \emph{L\'evy process} relative to time $t_2$

From the discussion above, if $t_2$ is the envelope closing time of a loop $\gab$ with $0\in \env(\gab)$, then $t_2$ is the right endpoint of a connected component of $[0,\infty) \setminus \cut(0) $. In fact, the converse is also true. Let $t_2$ be the right endpoint of a connected component of $[0,\infty) \setminus \cut(0)$, let $t^0_2$ be the left endpoint of the same component, let $t_1^0$ be defined by~\eqref{eq:t1}, and let $t_1=\inf\{ t\le t_1^0: \Lb_t\ge \Lb_{t_1^0} \textrm{ and }\Rb_{t}\ge \Rb_{t_1^0} \}$. Then $[t_1,t_2]$ is an envelope interval containing 0.

We now specify an enumeration of $\CLE_6$ loops whose envelope interval contains 0. In fact, the particular enumeration is not important to us, because we will be discussing properties which holds almost surely for all $\CLE_6$ loops simultaneously.
However, the following one has the advantage that the law of $\wh \Zb_t$ relative to each $t_2$ is explicit, which allows us to extend the constructions in Section~\ref{sec:dictionary-fl} to $\CLE_6$ loops. Fix $\eps>0, k\in \N$. Let $(t^0_2,t_2)$ be the $k$th component of $[0,\infty)\setminus \cut(0)$ (reading forward) such that $\Lb_{t^0_2} +\Rb_{t^0_2}-\Lb_{t_2}-\Rb_{t_2}\ge \eps$. Let $t_1$ be such that $[t_1,t_2]$ is an envelope interval. By varying $\eps$ and $k$, we exhaust all envelope intervals containing $0$. Doing the same for all $q\in \Q$ in place of $0$ provides an exhaustion of the CLE$_6$ $\Gab$.

We first recall the construction of $\gab$ from the interval  $[t_1,t_2]$ obtained for fixed $\eps,k$.  
Let $t_1^0$ be as in~\eqref{eq:t1}, and let $s_0=\ellb_{t_1^0}=\ellb_{t_2^0}$. 
Let $s=\inf\{t\le 0: \wh \Zb_{t'} \in \Zb_{t_2}+[0,\infty)^2,\;\forall \; t'\in [t,0] \}$. 
Then $s=\ellb_{t_1}$. We call $-s>0$ the \emph{quantum natural length} of $\gab$. The curve $\wh \etab|_{[s,0]}$ is $\gab$ equipped with its quantum natural parametrization starting and ending at $\etab(t_1)=\etab(t_2)$.

Next, we describe the law of $\wh\Zb$ (still for fixed $\eps,k$). 
Our choice of $t_2$ implies that $t_1^0$ is a backward stopping time for $\Zb$ and hence the walk $\Zb|_{(-\infty,t_1^0]} -\Zb_{t_1^0}$ recentered at $t_1^0$ has the same law as $\Zb|_{(-\infty,0]}$. 
Thus $\wh \Zb|_{(-\infty, s_0)} - \Zb_{t_1^0}$ has the same law as $\wh\Zb^0|_{(-\infty,0)}$ after time shift, where $\wh\Zb^0$ is the L\'evy process relative to 0, and is independent of $\Zb|_{[t_1^0,\infty)}$. At time $s_0$, $\wh \Zb$ has a jump from $\wh \Zb_{s_0^-} = \Zb_{t_1^0}$ to $\wh \Zb_{s_0} = \Zb_{t_2^0}$. 

It remains to describe the law of $\wh\Zb|_{[s_0,0]}$, which is a process from $\Zb_{t_2^0}$ to $\Zb_{t_2}$. 
Let $(\ell,r)=\Zb_{t^0_2} -\Zb_{t_2}$. Conditioned on $\ell$ and $r$, the pair $(\etab([t_2^0,t_2]),\gff)$ has the law of a $\sqrt{8/3}$-LQG disk decorated with an independent chordal space-filling $\SLE_6$, such that the boundary length of the two boundary arcs is $\ell$ and $r$, respectively. The law of $\wh\Zb|_{[s_0,0]}$ given $\ell,r$ is explicitly given in \cite[Theorem 1.2]{gwynne-miller-sle6}. The law can be described as the process $\wh \Zb^0+\Zb_{t_2^0}$ until it hits the boundary of the quadrant $\Zb_{t_2}+(0,\infty)^2$, conditioned on the exit location being the corner $\Zb_{t_2}$. Again, this is a conditioning on a zero probability event, but \cite{gwynne-miller-sle6} make sense of it by giving its
Radon-Nikodym derivative with respect to the unconditioned process when bounded away from the exit time from the quadrant.
This concludes our description on $\wh \Zb$.

By inspecting the explicit description of $\wh\Zb$ and 
using basic properties of stable process, one can check that conditioned on $\Zb_s$, the law of $\wh \Zb|_{[s,0]}$ is mutually absolutely continuous with respect to a more standard L\'evy type process which we describe now. For concreteness we assume now that $[t_1,t_2]$ is a right cone interval (the discussion extends to the other case by left/right symmetry).  
Let $\frk t'$ be a negative random variable with unbounded support and let $r$ be a positive  constant.
Let $\Yb$ be a $\frac32$-stable L\'evy excursion with only negative jumps on $[\frk t',0]$ and 
$\frk t=\sup \{t \le 0: \Yb_{t'} < r \;\forall \; \frk t' \le t'<t \}$. Conditioning on $\frk t$, let $\Xb$ be a $\frac32$-stable L\'evy excursion with only negative jumps on $[\frk t,0]$ independent of $\Yb$. 
Upon conditioning on $\{\Rb_s=r\}$, the law of $\wh \Zb|_{[s,0]}$ and $(\Xb,\Yb)|_{[\frk t,0]}$ are mutually absolutely continuous.

Let $\frk L_{\Xb},\frk L_{\Yb}$ be the pair of looptrees encoded by $\Xb_{[\frk t,0]}$ and $\Yb|_{[\frk t',0]}$. 
Then $\wh\Lb|_{[s,0]}$ defines a looptree $\frk L_{\op L}$ which is mutually absolutely continuous with respect to $\frk L_{\Xb}$, and which we call the \emph{left looptree} of $\gab$. 
All the discussion on the looptree $\frk L_{\Xb}$ in Section~\ref{sec:dictionary-fl} extends to $\frk L_{\op L}$. In particular, the region $\reg(\gab)$ enclosed by $\gab$ is the LQG looptree corresponding to $\frk L_{\op L}$, which is obtained by gluing independent $\sqrt{8/3}$-LQG disks to each bubble of $\frk L_{\op L}$ and embedding it in the $\sqrt{8/3}$-LQG cone via $\etab$.
Moreover, $\wh\Rb|_{[s,0]}$ defines a compact metric space $\frk L_{\op R}$ which is mutually absolutely continuous with respect to $\frk L_{\Yb}$ restricted to $[\frk t,0]$ (more precisely, its image under the quotient map). We abuse notation and call $\frk L_{\op R}$ the \emph{right looptree} of $\gab$, although $\frk L_{\op R}$ is not really a looptree but rather part of one. See Figure~\ref{fig:past-wedge} for an illustration of $\frk L_{\op L}$ and $\frk L_{\op R}$.
As in the spine-looptrees decomposition of the past wedge in Section~\ref{sec:dictionary-fl}, the $\CLE_6$ loop $\gab$ decomposes $\etab([t_1,t_2])$ into $\frk L_{\op L}$ and $\frk L_{\op R}$. 
In the remaining part of the paper, similarly as in Definition \ref{def:lt-cluster}, we use the symbol $\frk L(\gab)$ to denote  the looptree of $\gab$ whose LQG looptree structure is given by  $\gff$ restricted to $\reg(\gab)$. In other words, for a right (respectively, left) cone interval, we denote  $\frk L(\gab)=\frk L_{\op L}$ (respectively,  $\frk L(\gab)=\frk L_{\op R}$).

\subsection{Pivotal points and pivotal measure}
\label{sec:cont-piv}

Given the $\CLE_6$ $\Gab$ of Definition~\ref{def:CLE},
a point $p\in \C$ is called a \emph{double point} of a $\CLE_6$ loop $\gab\in \Gab$ if $\gab$ visits $p$ at least twice. 
The set of double points of $\gab$ is denoted by $\dbl_{\gab}$. A point $p\in \C$ is called a \emph{pivotal point} of $\Gab$ if it is a double point of a loop or the intersection of two loops.
The set of pivotal points of $\Gab$ is denoted by $\piv$.
Let us recall a result of Camia and Newman  about pivotal points.
\begin{lemma}[\cite{camia-newman-full}, Theorem 2]\label{lem:pivotal-continuous}
	Almost surely, the sets $\dbl_{\gab}$ for $\gab\in \Gab$, and the sets $\gab\cap \gab'$ for $\gab\neq \gab'\in \Gab$ are all disjoint. \nina{Furthermore, almost surely no point in $\dbl_{\gab}$ is a triple point of $\gab$ for $\gab\in\Gab$.}
\end{lemma}
\begin{remark}\label{rmk:CN}
	Lemma~\ref{lem:pivotal-continuous} is proved for pivotal points of the $\CLE_6$ as defined in the classical \cite{shef-cle}. We can apply it to the collection $\Gab$ of Definition~\ref{def:CLE}, since as mentioned above, we will show in Section~\ref{sec:disk} that these two notions of $\CLE_6$ agree modulo parametrization.
	Alternatively, we expect that Lemma~\ref{lem:pivotal-continuous} could be proved from Theorem~\ref{thm:mot} purely based on properties of $\frac32$-stable process.
\end{remark}

In this section we will first define natural measures on double points of looptrees and then use it to define the \emph{LQG pivotal measure} on $\piv$ associated with $(\gff,\etab)$. 

Let $\Xb$ be a $\frac32$-stable L\'evy excursion on $[\ell,0]$ with only negative jumps. At this point of general discussion, $\ell$ could be an arbitrary negative random variable but later we will take it to be  the quantum natural length of a $\CLE_6$ loop $\gab$ or one of the values $\frk t'$ or $\frk t$ defined at the end of Section~\ref{sec:cont-cle} (for the related processes $\Xb,\Yb$). Let $\frk L_{\Xb}$ be the looptree associated with $\Xb$. Recall the equivalence relation $\sim_{\Xb}$ on $[\ell ,0]$ defined by~\eqref{eq:d=0}. Unlike in Section~\ref{sec:cont-cle}, we now write the quotient map from $[\ell,0]$ to $\frk L_{\Xb}$ as $\pi_{\Xb}$ to indicate the dependence on $\Xb$.

For any $s\in [\ell,0]$, let
\begin{equation}\label{eq:foward}
	\Ab_{\Xb}(s)=\{ t\in(s,0]\,:\,\inf_{t'\in[s,t]} \Xb_{t'}= \Xb_t\}
\end{equation} 
be the set of forward running infima of $\Xb$ relative to time $s$. By the fluctuation theory of L\'evy process with only negative jumps, for any fixed $s$, the law of $\Ab_{\Xb}(s)$ is the range of a stable subordinator of Hausdorff dimension $\frac13$. Therefore, 
one can almost surely define a local time for $\Ab_{\Xb}(s)$.
Let $\pb_{\Xb}(s)$ be the Stieltjes measure  of this local time,\footnote{The local time $\pb_{\Xb}(s)$ is only defined up to multiplication by a constant. We will set this constant in~\eqref{eq83}. Recall the related discussion for the set of ancestor-free times in Section~\ref{sec:dictionary-branch}.} which is supported on $\Ab_{\Xb}(s)$.

Let $\dbl_{\Xb}$ be the set of double points of $\frk L_{\Xb}$, that is, points on $\frk L_{\Xb}$ with more than one pre-image under $\pi_{\Xb}$.
Then 
$$\dbl_{\Xb}(s):=\pi_{\Xb} (\Ab_{\Xb} (s))\subset \dbl_{\Xb} ~\textrm{ and }~ \dbl_{\Xb}=\cup_{s\in [\ell,0]} \dbl_{\Xb}(s).$$
For each fixed $s$, let $\nub_{\Xb}(s)$ be the push-forward of the measure $\pb_{\Xb}(s)$ onto $\frk L_{\Xb}$ by $\pi_{\Xb}$.
Note that $\nub_{\Xb}(s)$ represents the set of double points of $\frk L_{\Xb}$ separating $\pi_X(s)$ from the root. 
It is clear that for any fixed $s\neq s'$, 
the measures $\nub_{\Xb}(s)$ and $\nub_{\Xb}(s')$ almost surely agree on $\dbl_{\Xb}(s)\cap \dbl_{\Xb}(s')$. Therefore, it is tempting to define a measure $\nub$ on $\dbl_{\Xb}$ such that 
when restricting to any fixed $s$, it agrees with $\nub_{\Xb}(s)$. However, there are two caveats here. 
The first and obvious one is that $\nub_{\Xb}(s)$ is only almost surely well defined for any fixed $s$ while $[\ell, 0]$ is uncountable. The second one is more subtle and serious. It turns out that $\nub_{\Xb}$ is extremely big, in the sense that almost surely,
\begin{equation}\label{eq:big}
	\nub_{\Xb}(\pi_{\Xb}([t_1,t_2]) \cap \dbl_{\Xb})=\infty
\end{equation}
for any $\ell\le t_1<t_2\le 0$.
Therefore we consider exhaustions of $\dbl_{\Xb}$ when studying $\nub_{\Xb}$. In the following we introduce two such exhaustions for different purposes.

The first exhaustion is depending on $\Xb$. Let $\{Q_n\}$ be a sequence of increasing subsets of $[\ell,0]\cap\Q$ with finite cardinality such that $\cup_n Q_n=[\ell,0]\cap\Q$. We claim that 
\begin{equation}\label{eq:cover}
	\dbl_{\Xb}=\cup_{n=1}^\infty\cup_{q\in Q_n}\dbl_{\Xb}(q).
\end{equation}
Although being intuitive, it is not obvious that~\eqref{eq:cover} holds. However, it will become an immediate consequence of Lemma~\ref{lem:cover} stated below.
Since $\nub_{\Xb}$ is almost well-defined on each $\cup_{q\in Q_n}\dbl_{\Xb}(q)$ by putting $\{\nub_{\Xb}(q)\}_{q\in Q_n}$ together, we have a measure $\nub_{\Xb}$ defined on the entire double point set $\dbl$.

The second exhaustion is extrinsic and requires the LQG structure of $\frk L_{\Xb}$. More precisely, conditioning on $\frk L_{\Xb}$, we glue independent $\sqrt{8/3}$-LQG quantum disks to all bubbles of $\frk L_{\Xb}$ so that the boundary measure of the disk agrees with the metric on the bubble, making $\frk L_{\Xb}$ an LQG looptree. Given a double point $p\in \frk L_{\Xb}$, 
the set $\frk L_{\Xb}\setminus\{p\}$ has two connected components.
For $\eps>0$, we say that $p$ is 
\emph{$\eps$-significant} if for the closure of each connected component of $\frk L_{\Xb}\setminus\{p\}$, the sum of quantum area of the $\sqrt{8/3}$-LQG disks associated with all bubbles  is at least $\eps$. Let $\dbl_{\Xb,\eps}$ be the set of double points of $\frk L_{\Xb}$ which are at least $\eps$-significant. It is clear that $\dbl_{\Xb}=\cup_{\eps} \dbl_{\Xb,\eps}$ almost surely. 
The next lemma, whose proof is postponed to Section~\ref{sec:piv-def}, justifies~\eqref{eq:cover}.
\begin{lemma}	\label{lem:cover}
	Let $\{Q_n\}_{n\in \N}$ be the set of rationals as above. For any fixed $\eps>0$, there almost surely exists an $N\in\N$ such that
	\eqb
	\dbl_{\Xb,\eps}\subset \bigcup_{q\in Q_N}\dbl_{\Xb}(q).
	\label{eq:cover1} 
	\eqe 
\end{lemma}
We see that $\nub_{\Xb}(\dbl_{\Xb,\eps})<\infty$ almost surely for any $\eps>0$ by using this lemma and that $\nub_{\Xb}(\dbl_{\Xb}(q))<\infty$ for any rational $q$.

Both exhaustions of $\dbl_{\Xb}$ are useful. 
From the first exhaustion, we rigorously defined the measure $\nub_{\Xb}$ on the entire $\dbl_{\Xb}$ via local time on the running infima of L\'evy excursions. This is convenient for the study of $\nub$ in Section~\ref{subsub:pivot} and~\ref{subsec:flip} via the well established theory of L\'evy processes. 
The second exhaustion has the advantage of not relying on $\Xb$, which will be convenient for the study of $\nub_{\Xb}$ via GFF and SLE once the LQG looptree is embedded to $\sqrt{8/3}$-LQG cone (see Remark~\ref{rmk:pivot}). 
In this paper we mainly use the LQG exhaustion of $\dbl_{\Xb}$ as an objective rather than a tool. But the GFF/SLE perspective of $\nub$ will be a crucial input 
in \cite{hs-quenched}.
\bigskip

Now we turn our attention to the particular looptrees associated with $\CLE_6$ loops. We retain the notions and setup introduced in Section~\ref{sec:cont-cle} for the looptree perspective on $\CLE_6$. 
Let $\gab$ be a CLE$_6$ loop chosen from the exhaustion of $\Gab$ defined there. Let $s<0$ be such that $-s$ is the quantum natural length of $\gab$. Let $u$ be the envelope closing time of $\gab$, and let $(\wh \Lb,\wh \Rb)$ be the L\'evy process relative to the envelope closing time of $\gab$. As in Section~\ref{sec:cont-cle}, we assume for concreteness that $\env(\gab)$ is a right cone interval (the other case being symmetric). Let $\frk L_{\op L}$ and $\frk L_{\op R}$ be the left and right looptree of $\gab$. 

As in Section~\ref{sec:cont-cle}, we use $\Xb|_{[\frk t,0]}$ and $\Yb|_{[\frk t',0]}$ to denote the two independent $\frac32$-stable excursions such that $(\wh \Lb,\wh\Rb)|_{[s,0]}$ and $({\Xb},\Yb)|_{[\frk t,0]}$ are mutually absolutely continuous. Therefore we can reweigh the law of $({\Xb},\Yb)|_{[\frk t,0]}$ and then couple this pair with $(\wh \Lb,\wh\Rb)|_{[s,0]}$ such that 
$s=\frk t$ and $(\wh \Lb,\wh\Rb)|_{[s,0]}= ({\Xb},\Yb)|_{[\frk t,0]}$. Under this coupling, $\frk L_{\op L}=\frk L_{\Xb}$ and $\frk L_{\op R}$ is a subspace of $\frk L_{\Yb}$.
Therefore the measure $\nub_{\Xb}$ on $\dbl_{\Xb}$ turns into a measure $\nub_{\Lb}$ on $\dbl_{\Lb}$. Furthermore, the measure $\nub_{\Yb}$ restricted to $\dbl_{\Yb}\cap \frk L_{\op R}$ defines a measure $\nub_{\Yb}$.
To describe $\dbl_{\Yb}\cap \frk L_{\op R}$ purely in terms of $\wh \Rb$, let $\Ab_{\Rb}(s)$ be defined as in~\eqref{eq:foward} with $\wh \Rb$ in place of $\Xb$.
For $p\in \frk L_{\op R}$, we call $p$ a \emph{double point} of $\frk L_{\op R}$ if $\pi_{\Yb}^{-1}(p)$ 
has more than one point in $[s,0]$. Let $\dbl_{\Rb}$ be the set of double points of $\frk L_{\op R}$. Then $\dbl_{\Yb}\cap \frk L_{\op R}=\dbl_{\Rb} \cup \Ab_{\Rb}(s)$ and $\dbl_{\Rb} \cap \Ab_{\Rb}(s)\neq \emptyset$ almost surely (see the right side of Figure \ref{fig:past-wedge}, where the set $\Ab_{\Rb}(s)$ is represented as lying on the dotted line). Although we are using $\Xb,\Yb$, it is clear from the construction that the sets $\dbl_{\Lb}, \dbl_{\Rb},\Ab_{\Rb}(s)$ and the measures $\nub_{\Lb}$ and $\nub_{\Rb}$ are explicitly determined by $\wh \Lb,\wh\Rb$ without external reference. 

Let $\pi_{\Lb}$ and $\pi_{\Rb}$ be the quotient maps of from $[s,0]$ to $\frk L_{\op L}$ and $\frk L_{\op R}$, respectively. 
As in Section~\ref{sec:dictionary-fl}, 
there exist unique embeddings $\phi_{\Lb}$ of $\frk L_{\op L}$
and $\phi_{\Rb}$ of $\frk L_{\op R}$
such that $\phi_{\Lb} \circ \pi_{\Lb}(t) = \gab(t)$ and $\phi_{\Rb} \circ \pi_{\Rb}(t) = \gab(t)$ 
for all $t\in [s,0]$. Then~\eqref{eq:double} holds with $\gab$ in place of $\wh\etab$. Moreover, $\phi_{\Rb} (\Ab_{\Rb}(s))$ equals the intersection of $\gab$ and the boundary of $ \etab(\env(\gab))$.
The pushforward of $\nub_{\Lb}$ under $\phi_{\Lb}$ and $\nub_{\Rb}$ under $\phi_{\Rb}$ define a measure $\ol\nub_{\gab}$ supported on $\dbl_{\gab}\cup (\gab\cap \bdy\etab(\env(\gab)))$. In particular, by restriction, $\ol\nub_{\gab}$ induces a measure $\nub_{\gab}$ on $\dbl_{\gab}$.

Recall~\eqref{eq:intersection}. For any $\CLE_6$ loops $\gab'\neq \gab$, the intersection $\gab\cap \gab'$ belongs to either $\gab\cap \bdy\etab(\env(\gab))$ or 
$\gab'\cap \bdy\etab(\env(\gab'))$, depending on whether $\env(\gab)\subset \env(\gab')$ or vice versa. Therefore by restriction,
$\ol\nub_{\gab}$ induces a measure $\nub_{\gab,\gab'}$ on $\gab'\cap \gab$ for any $\env(\gab)\subset \env(\gab')$. By convention, we set $\nub_{\gab,\gab'}=0$ if $\env(\gab)\cap \env(\gab')=\emptyset$ or $\env(\gab')\subset \env(\gab)$ (recall~\eqref{eq:env}).
\begin{definition}\label{def:pivotal}
	Summing $\nub_{\gab}$ over $\gab\in \Gab$ and $\nub_{\gab,\gab'}$ over all $\gab\neq \gab'$, we obtain a measure $\nub$ supported on $\piv$. We call $\nub$ the \emph{LQG pivotal measure }of $\Gab$ associated with $(\gff,\etab)$.
\end{definition}
Since $\ol\nub_{\gab}$ restricted to $ \gab\cap \bdy\etab(\env(\gab))$ is almost surely finite for any $\gab\in \Gab$, the measure
$\nub$ restricted to the intersection of two distinct loops is finite.
However,~\eqref{eq:big} implies that the $\nub$-mass of any open set is infinite. 
Again, to say anything meaningful about $\nub$, one needs to take an approximate exhaustion of $\piv$. 
We now define one which is closely related to the second exhaustion for double points of looptrees. 

Given $z\in \piv$, a new loop configuration $\Gab_z$ can be obtained as follows. By Lemma~\ref{lem:pivotal-continuous}, almost surely $z$ is either in $\dbl_{\gab}$ for some $\gab\in\Gab$, or in  $\gab\cap \gab'$ for some $\gab,\gab'\in\Gab$ but not both.
If $z\in \dbl_{\gab}$ for some $\gab\in\Gab$, 
we split $\gab$ into two loops intersecting at $z$, so that the orientation of the two new loops are consistent with $\gab$.
If $z\in \gab\cap \gab'$ for some $\gab,\gab'\in\Gab$ such that $\gab\neq \gab'$, we obtain $\Gab_z$ by merging $\gab$ and $\gab'$ into one loop. Recall that $\gab$ and $\gab'$ have the same orientation if and only if they are unnested.
The new loop after merging $\gab$ and $\gab'$ can be singly orientated
in a way which is consistent \referee{with} the orientation of both $\gab$ and $\gab'$. 
We call \emph{flipping the color} of $z$, the operation of changing $\Gab$ to $\Gab_z$.
Let $\Gab\Delta\Gab_z$ denote the symmetric difference between $\Gab$ and $\Gab_z$. Almost surely, $\Gab\Delta\Gab_z$ always consists of 
exactly three loops, each of which encloses a region with a positive $\mub_{\gff}$-area.
Let $\sig(\Gab,z)$ be the minimum over the three areas. We call $\sig(\Gam,z)$ the \emph{significance} of $z$. Intuitively, the higher the significance of $z$, the more dramatic change it will cause when flipping the color of $z$. 	Let $z\in \piv$ and $\eps>0$. We say that $z$ is \emph{$\eps$-significant} if $\sig(\Gab,z)\geq\eps$. The set of $\eps$-significant points is denoted by $\piv_\eps$. It is clear from the definition that $\piv =\cup_{\eps>0} \piv_{\eps}$.

\begin{definition}\label{def:piv-eps}
	By restriction, $\nub$ induces a Borel measure on $\C$ supported on $\piv_{\eps}$, which is denoted by $\nub_\eps$. We call $\nub_\eps$ the \emph{$\eps$-pivotal measure} of $\Gab$ associated with $(\gff, \etab)$.
\end{definition}
We now argue that the measure $\nub_{\eps}$ is locally finite. 
Given a loop $\gab\in \Gab$, let $\frk L_{\op L}$ and $\frk L_{\op R}$ be its left and right looptrees, respectively.
Let $\dbl_{\Lb,\eps}$ (respectively, $\dbl_{\Rb, \eps}$) be the set of $\eps$-significant double points of $\frk L_{\op L}$ (respectively, $\frk L_{\op R}$) defined in the second exhaustion of double points of looptrees. 
Then
\begin{equation}\label{eq:piv-eps}
	\piv_\eps \cap \dbl_{\gab} = \phi_{\Lb} (\dbl_{\Lb,\eps})\cup \phi_{\Rb} (\dbl_{\Rb,\eps}).
\end{equation}
In particular, $\nub(\piv_{\eps}\cap \dbl_{\gab})<\infty$. Recall also that $\nub(\gab\cap \gab')<\infty$ for any $\gab\neq \gab'$. Now, given a bounded domain $D$ and $\eps>0$, there are only finitely many loops in $\Gab$ that have nonempty intersection with $D$ and LQG area larger than $\eps>0$ \cite{camia-newman-full}. Therefore, $\nub_{\eps}(D)<\infty$. So $\nub_{\eps}$ is locally finite.

\begin{remark}\label{rmk:pivot}
	In \cite{gps-pivotal}, the authors introduced another exhaustion of $\piv$ based on the so-called four-arm events. Under this exhaustion, they constructed a random local finite Borel measure $\lambda$. The construction is via establishing the scaling limit of the counting measure over pivotal points of critical site percolation on the regular triangular lattice. 
	In \cite{natural}, it is proved that the $3/4$-dimensional Minkowski content of $\piv_{\eps}$ exists and defines $\lambda$. 	In \cite{hs-quenched}, the second and third author of this paper prove that $\nub $ and $ e^{\frac14\sqrt{8/3} \gff}\,d\lambda$ agree up to a multiplicative constant, where the constant $\frac14$ comes from the KPZ relation applied to the set of pivotal points. 
	Our scaling limit result for pivotal measure is a random triangulation version of the result in \cite{gps-pivotal} for the triangular lattice.
\end{remark}
We conclude the section by introducing the notion of \emph{type}s for $\CLE_6$ pivotal points in analogy to Section~\ref{subsec:pivot}. Given a loop $\gab\in \Gab$, recall that $\gab$ decomposes $\etab(\env(\gab))$
into two looptrees, where $\frk L(\gab)$ is the one corresponding to $\reg(\gab)$. 
We say that a double point of $\gab$ is a pivotal point of \emph{type 1} if it comes from a double point of $\frk L(\gab)$ and of \emph{type 2} otherwise. If a point $p$ is a point of intersection of two distinct unnested (respectively, nested) $\CLE_6$ loops $\gab,\gab'\in \Gab$, we call $p$ a pivotal point of \emph{type 3} (respectively, \emph{type 4}). It is straightforward to see that this definition of types agrees with the one in Section~\ref{subsec:pivot} based on color flipping (see Figure~\ref{fig:piv}).

\subsection{Constructions for the disk and sphere}
\label{sec:disk}

In Section~\ref{sec:mot}-\ref{sec:cont-piv}, we explained the mating-of-trees theory for the $\sqrt{8/3}$-LQG cone. Starting from $(\gff,\etab)$ and $\Zb$, we defined the associated branching $\SLE_6$, $\CLE_6$, and LQG pivotal measure. In this section, we first present the disk variant of Theorem~\ref{thm:mot}. In this case, $\CLE_6$ is a well-known and well-studied subject in the literature (see for example \cite{camia-newman-full,shef-cle}). In this section we will also explain that the disk variant of Definition~\ref{def:CLE} agrees with the more classical construction in \cite{shef-cle}. Finally, we will briefly describe the sphere variant of Theorem~\ref{thm:mot}.

Given the foundation we laid for the $\sqrt{8/3}$-LQG cone case in Sections~\ref{sec:mot}-\ref{sec:cont-piv}, we will not start from scratch for the disk case. Instead, we embed the $\sqrt{8/3}$-LQG disk into the $\sqrt{8/3}$-LQG cone and trivially extend everything in Sections~\ref{sec:mot}-\ref{sec:cont-piv} to the disk case. Let $(\gff,\etab)$ and $\Zb$ be as in Theorem~\ref{thm:mot}. 
Let $\tau$ be the first time $t\in\R$ such that the following three conditions are satisfied:
	(i) the quantum boundary length of the component $D$ of $\C\setminus\etab((-\infty,t])$ containing 0 is smaller than 1, 
	(ii) there is a cone interval $[\frk s,\frk t]$ such that $\etab([\frk s,\frk t])=\ol D$, and 
	(iii) if $[\frk s',\frk t']$ is the smallest cone interval containing $[\frk s,\frk t]$ then $[\frk s,\frk t]$ and $[\frk s',\frk t']$ are in opposite directions, i.e., $[\frk s',\frk t']$ is a left (resp.\ right) cone interval if $[\frk s,\frk t]$ is a right (resp.\ left) cone interval.
For concreteness, we condition on the event that $[\frk s,\frk t]$ is a right cone interval so that $\Lb_{\frk s}=\Lb_{\frk t}$. Let $H=\Rb_{\frk s}-\Rb_{\frk t}$ and $A=\frk t-\frk s$. 
According to \cite[Proposition~5.1]{sphere-constructions} and its proof, conditioning on $H$, the law of $(D,\gff|_D)$ as an LQG surface is absolutely continuous with respect to the $\sqrt{8/3}$-LQG disk with boundary length $H$, with Radon-Nikodym derivative proportional to $A$ (see the more precise statement below~\eqref{eq:reweighting-A}).
This reweighing by $A$ is an instance of the inspection paradox for renewal processes.
Moreover, conditioning on $(D,\gff|_D)$, 
\begin{equation}\label{eq:uniform}
	\textrm{$\frk t$ (or equivalently, $-\frk s$) is uniform in $(0,A)$},
\end{equation}
and the curve $\etab|_{[\frk s,\frk t]}$ modulo parametrization is the counterclockwise space-filling $\SLE_6$ in $D$ starting and ending at $\etab(\frk s)=\etab(\frk t)$ independent of $\frk t$.

In light of the above discussion, from now on we work under the reweighed probability measure
\begin{equation}\label{eq:reweighting-A}
	d\wt\P:=\frac{c}{A}\1_{\Rb_{\frk s}>\Rb_{\frk t}}d\P,
\end{equation}
where $\P$ is the probability measure for $(\gff,\etab)$ and $c$ is a renormalizing constant.
Let $\phib$ be the conformal map from $D$ to $\D$ such that $\phib(0)=0$ and $\phib(\etab(\frk s))=1$. Let $\gffd=\gff\circ \phib^{-1}+Q\log|(\phib^{-1})'|-\frac{2}{\gamma}\log H$, where $\gamma=\sqrt{8/3}$ and $Q$ is as in~\eqref{eq:Q}. 
Then, under the new measure $\wt\P$, the LQG-surface $(\D,\gffd)$ is a representative (in the sense of the $\gamma$-equivalence of LQG-surfaces) of the unit boundary length $\sqrt{8/3}$-LQG disk as defined in Section \ref{sec:lqg}.
Moreover, under $\wt\P$, the field $\gffd$ satisfies the following two properties.
\begin{itemize}
	\item \emph{Bulk re-rooting invariance}: 
	Let $(z,\gffd)$ be coupled such that conditioning on $\gffd$, the law of $z$ is a point on $\D$ sampled from $\mub_{\gffd}$. 
	Let $\varphi_z$ be the M\"obius transform from $\D$ to $\D$ that maps $0$ to $z$ and preserves $1$. 	
	Then $\gffd\circ \varphi_z+Q\log|\varphi_z'| \overset{d}{=} \gffd$. 
	\item \emph{Boundary re-rooting invariance}: 
	Let $(z,\gffd)$ be coupled such that conditioning on $\gffd$, the law of $z$ be a point on $\bdy \D$ sampled from $\nub_{\gffd}$. 
	Let $\varphi_z$ be the M\"obius transform from $\D$ to $\D$ that maps $1$ to $z$ and preserves $0$. 	
	Then $\gffd\circ \varphi_z+Q\log|\varphi_z'| \overset{d}{=} \gffd$. 
\end{itemize}
The bulk re-rooting invariance follows from~\eqref{eq:uniform} and the boundary one follows from \cite[\xin{Proposition~A.8}]{wedges}. 
These re-rooting invariances characterize $\gffd$ among the representatives of the unit boundary length $\sqrt{8/3}$-LQG disk. Indeed, if $(\D,\gff')$ is a representative of the unit length $\sqrt{8/3}$-LQG disk, and $\gff'$ satisfies the two re-rooting invariance properties listed above, then $\gff'$ and  $\gffd$ are equal in law as random fields.

Let $\etad$ be the image of $\etab|_{[\frk s,\frk t]}$ under $\phib$. Then, modulo parametrization, $\etad$ is a counterclockwise space-filling $\SLE_6$ on $\D$ starting and ending at $1$ which is independent of $\gffd$. Let $\frk a:=\mub_{\gffd}(\D)=A/H^2$.
We reparametrize $\etad$ such that $\etad(0)=1$ and $\mub_{\gffd} (\etad([0,t])) =t $ for all $0\le t\le \frk a$. Let $\Zd$ be the right cone excursion on $[0,\frk a]$ defined by 
\begin{equation}\label{eq:Zd}
	\Zd_t:=(\Ld,\Rd)=H^{-1}(\Zb_{tH^2+\frk s}-\Zb_{\frk s}) \qquad \forall \; t\in [0,\frk a]. 
\end{equation}
For $t\in (0,\frk a)$, we call the clockwise (reps. counterclockwise) arc from $\etad(t)$ to $1$ on the boundary of $\D\setminus \etad[0,t]$ the left (respectively, right) frontiers of $\etad([0,t])$. 
Then $\Ld_t$ and $\Rd_t$ equals the $\nub_{\gffd}$-length of the left and right, respectively, frontiers of $\etad([0,t])$. Therefore we call $\Zd$ the boundary length process of $(\gffd,\etad)$.
The law of $\Zd=(\Zd_t)_{t\in [0,\frk a]}$ can be described as the Brownian motion $\Zb$ starting at $(0,1)$ and conditioned to exit the first quadrant $[0,\infty)^2$ at the origin at time $\frk a$. See for example \cite[Section 3]{sphere-constructions} for a rigorous meaning of this zero-probability conditioning. We call the law of $\Zd$ \emph{the Brownian cone excursion with correlation-$\frac12$ starting from $(0,1)$ with variance $\beta$}, or \emph{Brownian cone excursion from $(0,1)$} for short. The duration $\frk a$ of $\Zd$ is a random variable with inverse Gamma distribution. 

The following theorem is the disk variant of Theorem~\ref{thm:mot}.
\begin{theorem}[\xin{Theorem~2.1} \cite{sphere-constructions}]\label{thm:mot-disk}
	Let $(\D,\gffd)$ be the representative of the unit boundary length $\sqrt{8/3}$-LQG disk satisfying the above mentioned bulk and boundary re-rooting invariance properties.
	Let $\etad$ be a counterclockwise space-filling $\SLE_6$ on $\D$ starting and ending at $1$ which is independent of $\gffd$ modulo parametrization.
	Let $\etad$ be parametrized by its quantum area and $\Zd$ be its boundary length process. Then $\Zd$ is a Brownian cone excursion from $(0,1)$.
	Moreover, the pair $(\gff,\etab)$ is measurable with respect to the $\sigma$-algebra generated by $\Zd$.
\end{theorem}

The measurability statement in Theorem~\ref{thm:mot-disk} is inherited from the corresponding statement in Theorem~\ref{thm:mot}. Here there is no need to consider $(\D,\gffd)$ modulo the rotation about the origin as the boundary point 1 has been marked. This is the disk version of mating-of-trees theorem in \cite{sphere-constructions}. Theorem~\ref{thm:mot-disk} still holds if counterclockwise is replaced by clockwise and $(0,1) $ is replaced $(1,0)$. This corresponds to reweighing the law $\P$ of $(\gff, \etab)$ by $cA^{-1}\1_{\Lb_{\frk s}>\Lb_{\frk t}}$.

\bigskip

Recall the branching $\SLE_6$ $\taub^*$, the $\CLE_6$ $\Gab$ and the pivotal measure $\nub$ associated with $(\gff,\etab)$. We can define the corresponding objects for $(\gffd,\etad)$. There are two ways to do it. The first is to start from Theorem~\ref{thm:mot-disk} and repeat everything. Note that $\taub^*,\Gab$ and $\nub$ are described as almost sure explicit functions of $\etab$ and $\Zb$. Our construction of $(\gffd,\etad,\Zd)$ is done by reweighing and restriction of $(\gff,\etab,\Zb)$. Therefore, these explicit functions can be applied to $(\etad,\Zd)$ with straightforward adaption. This construction does not require a reference to $(\gff,\etab,\Zb)$. We omit the details here. Another equivalent construction is to simply map the restriction of $\taub^*,\Gab$ and $\nub$ from the domain $D$ to the unit disk $\D$ through the conformal map $\phib$. Now we will elaborate on the detail of this construction.

We start with $\taub^*$. Let $w\in \ol D\setminus\{\etab(\frk s)\}$. Since $[\frk s,\frk t]$ is a cone interval, it is necessarily the case that $\etab(\frk s)$ is the entrance point of the branch of $\taub^*$ targeted at $w$.
Let $z=\phib(w)$, and let $\etadh^z$ be the image of the segment on $\wh \etab^w$ from $\etab(\frk s)$ to $w$. Modulo parametrization, $\taud:=\{\etadh^z\}_{z\in\D\setminus\{1\}}$
is called the \emph{branching $\SLE_6$ associated with $(\gffd, \etad)$}. As in Remark~\ref{rmk:branch}, for any fixed $z\in \D$, $\etadh^z$ is an $\SLE_6$ on $(\D,1,z)$. Moreover,
$\taud$ is a version of branching $\SLE_6$ on $\D$ rooted at $1$ defined in Section~\ref{sec:sle}, and $(\etad,\taud)$ are coupled as in Section~\ref{sec:sle}.
Given a fixed $u>0$, on the event that $u<\frk a$, one can still almost surely parametrize the branch of $\taud$ targeted at $z=\etad(u)$ by its quantum natural time. In fact, let $w=\phib^{-1}(z)$. Since $\frk s$ is a backward stopping time for $\Zb$, the quantum natural parametrization can be defined almost surely for $\wh\etab^{w}$ via the local time on ancestor-free times. Let $s_0<0$ be such that $\wh \etab^{w}(s_0)=\etab(\frk s)$. The quantum natural parametrization of $\etadh^z$ is given by 
\begin{equation}\label{eq:taud}
	\etadh^z(s)=\phib\circ \wh\etab^z(s_0+sH^{3/2})\qquad \forall s\ge 0\; \textrm{such that }s_0+sH^{3/2}\le 0,
\end{equation}
where the scaling exponent $\frac32=2\times \frac34$ comes from the fact that the time of the Brownian motion is scaled by $H^2$ as~\eqref{eq:Zd}, while the Hausdorff dimension of ancestor-free time is $\frac34$.

We now turn our attention to $\CLE_6$. For each $\gab\in \Gab$ such that $\gab\subset \ol D$, $\gab_{\op D}:=\phib\circ \gab$ is a loop on $\ol \D$. Suppose $\gab$ is under its quantum natural parametrization, then as in~\eqref{eq:taud}, we parametrize $\gab_{\op D}$ by $\gab_{\op D}(s)=\phib\circ \gab(s_0+sH^{3/2})$ where $-s_0$ is the quantum natural length of $\gab$ and $s_0\le s_0+sH^{3/2}\le 0$. We call this parametrization the \emph{quantum natural parametrization} of $\gab_{\op D}$. Let $\Gabd$ be the collection of such loops $\gab_{\op D}$. 
Then $\Gabd$ is called the \emph{$\CLE_6$ associated with $(\gffd, \etad)$}

Starting from a sample of $\taud$, Sheffield \cite{shef-cle} defined the so-called \emph{$\CLE_6$ on $\D$ rooted at $1$}, which we denote by $\Gabd'$. Combining \cite{shef-cle} and~\cite{ig4}, viewed as curves modulo parametrization, $\Gabd'$ can be expressed as an explicit function of $\taud$ and $\etad$; \xin{see \cite[Section~3.6]{peano-survey}}.
Now we explain that modulo parametrization $\Gabd=\Gabd'$ almost surely, thus addressing the issue left in Remark~\ref{rmk:CN}. Let us first consider how to use $\taud$ to describe the outermost loop in $\Gabd$ whose region contains 0. By definition, it amounts to find the loop $\gab\in \Gab$ with the biggest envelope interval such that $0\in \reg(\gab)\subset \ol D$. 

Recall the bijection between envelope intervals and $\CLE_6$ loops in Section~\ref{sec:cont-cle}. 
Let $\env(\gab)=[t_1,t_2]$ and $t_2^0$ and $t_1^0$ be defined as in~\eqref{eq:t2} and~\eqref{eq:t1}, respectively. Then $[t^0_1,t^0_2]\subset [t_1,t_2]\subset [\frk s, \frk t]$. 
Since we are on the event that $[\frk s, \frk t]$ is a right cone interval, $\gab$ must be counterclockwise thus $\env(\gab)$ is a right cone interval. Therefore, $[t^0_1,t^0_2]$ is a left cone interval. Let $I$ be the maximal (i.e., largest) cone interval inside $[\frk s,\frk t]$ containing 0 and recall that $\env(I)$ is the smallest cone interval containing $I$. Then there exists a loop $\gab'\subset \ol D$ such that $\env(I)=\env(\gab')$. Since $\env(\gab')\subset \env(\gab)$ by definition, if $\gab'\neq \gab$, we must have $\gab'\subset \reg(\gab)$ and thus $\env(\gab')\subset [t_0^1,t_0^2]$, which contradicts the maximality of~$I$. Therefore $\gab'=\gab$ and hence $[t_0^1,t_0^2]$ can be identified as the maximal left cone interval inside $[\frk s,\frk t]$ and $\env(\gab)$ is the smallest cone interval containing $[t^1_0,t^2_0]$.
Let $\ellb$ be the local time of $\ans(0)$. 
Then $\ellb_{t_1^0} =\ellb_{t_2^0}$ is the first time after $\ellb_{\frk s}$ when $\wh\etab^0$ finishes tracing a bubble containing 0 in counterclockwise direction, and $\etab|_{[t_1^0,t_2^0]}$ fills this bubble. 

Recall the decomposition of $\gab$ it into a past segment $\gab_1$ and a future segment $\gab_2$ inside the past and future wedge relative to time $0$. In the previous paragraph we explained how to find the endpoint $p=\etab(t_2^0)$ of $\gab_1$, or equivalently, the staring point of $\gab_2$, in a way which only depends on topological properties of $\wh \etab^0$ in $D$. To find the other endpoint $q$ of $\gab_1$ and $\gab_2$, we let $D'$ be the component of the future wedge relative to time $0$ whose boundary contains $p$. Then $q$ is the other intersection point of the left and right frontier of $\etab((-\infty,0])$.

From the discussion above, after applying the conformal map $\phib$, the outermost loop $\gab_{\op D}$ in $\Gabd$ containing 0 can be recovered from $(\taud,\etad)$ as follows. Run branch $\etadh^0$ of $\taud$ targeted at $0$ until the first time
when this branching finishes tracing a bubble containing 0 in counterclockwise direction. Let $p$ be the last point on this bubble visited by $\etad$.
Let $F$ be the closure of the set of points visited by $\etad$ after its last visit of $p$. Then $p$ must be on the boundary of a unique connected component $\Delta$ of the interior of $F$. Moreover, the exists a unique point $q\neq p$ on the boundary of $\Delta$ such that $q$ is an intersection of the two boundary arcs of $F$ between $p$ and $1$. Now we first trace $\etadh^0$ from $q$ to $p$ and denote this 
path by $\gab_{\op D,1}$. Restricting $\etad$ to $\Delta$, we obtain a space-filling chordal $\SLE_6$ curve on $(\Delta,p,q)$. By skipping bubbles, we obtained a chordal $\SLE_6$ curve on $(\Delta,p,q)$, which we denote by $\gab_{\op D,2}$. Then $\gab_{\op D}$ is the concatenation of $\gab_{\op D,1}$ and $\gab_{\op D,2}$. Combining \cite{shef-cle} and~\cite{ig4}, the concatenation of $\gab_{\op D,1}$ and $\gab_{\op D,2}$ also almost surely gives the outermost loop in $\Gabd'$ containing $0$; \xin{also see \cite[Remark 3.28]{peano-survey} for a related construction.}  

To summarize, we have proved that the outermost loop in $\Gabd$ and $\Gabd'$ containing 0 are almost surely identical as curves modulo parametrization. 
Replacing $0$ by any other rational, the same argument implies that 
the collection of outermost loops of $\Gabd$ and $\Gabd'$ are almost surely the same as collection of curves modulo parametrization.
Recall that in \cite[\xin{Section~4.3}]{shef-cle} $\Gabd'$ is constructed recursively by further considering outermost loops inside each complementary component of the union of outermost loops. On the $(\D,\gffd)$ side, each such component is the interior of the range of $\etad$ restricted to a certain cone interval, inside which we can repeat the argument above. This shows that $\Gabd=\Gabd'$ almost surely as collection of curves modulo parametrization.

\begin{remark}
	Since $\etad$ and $\taud$ determine each other, $\Gabd$ is also a function of $\taud$. In fact, the original paper \cite{shef-cle} defined $\CLE_6$ only in terms of branching $\SLE_6$ without an explicit reference to $\etad$. Our construction above is adapted from the original construction in \cite{shef-cle} since this makes the agreement of $\Gabd$ and $\Gabd'$ more transparent. It is also possible to describe the branching $\SLE_6$ as a function of $\CLE_6$, so that $\CLE_6$ also contains the same amount of information as branching and space-filling $\SLE_6$ thus can be viewed as a third representation of the scaling limit of critical planar percolation. We will not review these better known constructions in detail since they are not needed for the rest of the paper.
\end{remark}

We define the pivotal points for $\Gabd$ similarly as we did for $\Gab$ in the whole-plane setting. We denote the set of pivotal points of $\Gabd$ by $\pivd$. The pivotal measure $\nud$ is simply the pushforward of $\nub$ under $\phib$ appropriately rescaled as a function of $H$ (similarly as in~\eqref{eq:Zd} and~\eqref{eq:taud}). For each $\gab_1,\gab_2\in \Gabd$ such that $\gab_1\neq \gab_2$,
let $\wt\gab_1 $ and $\wt \gab_2$ be the pre-image of $\gab_1$ and $\gab_2$ under $\phib^{-1}$. Let 
\begin{equation}\label{eq:piv-scale}
	\nud_{\gab_1} =H^{\frac12}\nu_{\wt \gab_1} \qquad\textrm{and}\qquad \nud_{\gab_1,\gab_2} =H^{\frac12}\nu_{\wt \gab_1,\wt \gab_2}.
\end{equation}
The exponent $\frac12=\frac32\times \frac13$ comes from the fact that the quantum natural length scales like $H^{3/2}$ 
while the set $A_{\Xb}(s)$ defined in~\eqref{eq:foward} has dimension $\frac13$.
The measure \emph{LQG pivotal measure associated with $(\gffd,\etad)$}, which is denoted by $\nud$,
is the summation over all the measures $\nud_{\gab_1}$ and $ \nud_{\gab_1,\gab_2}$ as in Definition~\ref{def:pivotal}. 
The set of $\eps$-significant point $\piv_{\op D,\eps}$ is defined to be $\phib(\piv_{\eps'})$, where $\eps'=H^2\eps$ since the function $\sig(\cdot)$ is defined via the $\mub_{\gff}$-mass, which scales as $H^2$. Then $\{\piv_{\op D,\eps}\}_{\eps>0}$ provides an exhaustion for $\pivd$. Let $\nub_{\op D,\eps}$ be the Borel measure induced by restricting $\nud$ to $\piv_{\op D,\eps}$. Then $\nub_{\op D,\eps}$ is almost surely finite.\\

We conclude our discussion on the disk variant of mating of trees by stating the fixed-area version of Theorem~\ref{thm:mot-disk}. Fix a constant $\frk m>0$.
We call the law of $\Zd$ in Theorem~\ref{thm:mot-disk} conditioned on $\frk a =\frk m$ the \emph{Brownian cone excursion with duration $\frk m$ from $(0,1)$}. 
This conditioning is another easy Brownian motion exercise (details can be found in \cite[\xin{Section 3}]{sphere-constructions}). Then Theorem~\ref{thm:mot-disk} readily gives the following result.
\begin{cor}[\cite{sphere-constructions}]\label{cor:fix-area}
	Fix a constant $\frk m>0$. Let $(\D,\gffd)$
	be the $\sqrt{8/3}$-LQG disk with length $1$ and area $A$. Let $\etad$ and $\Zd$ be defined in the same way as in Theorem~\ref{thm:mot-disk}. Then 
	$\Zd$ is a Brownian cone excursion with duration $\frk m$ starting from $(0,1)$.
	Moreover, $(\gff,\etab)$ is measurable with respect to the $\sigma$-algebra generated by $\Zd$.
\end{cor}
Here we abuse notation and still use $\gffd,\etad,\Zd$ to denote the objects in Theorem~\ref{thm:mot-disk} under conditioning. The field $\gffd$ is still characterized (among the representatives of the  $\sqrt{8/3}$-LQG disk with length $1$ and area $A$) by the bulk and boundary re-rooting invariance.
All the mating of trees theory for the random area case in Theorem~\ref{thm:mot-disk} extends to the fixed area case with straightforward modifications.

In the rest of this subsection, we will describe 
the sphere version of Theorem~\ref{thm:mot}. Again, in order to take advantage of the foundation we have already laid, we will think of the sphere as the disk with boundary length 0. Suppose $(\gffd,\etad)$ and $\Zd$ are as in Theorem~\ref{thm:mot-disk}. For $C\in \R$ and $\eps>0$, we define the event $E_{C,\eps}=\{ 1 \leq \mub_{\gffd-C}(D)\leq (1+\eps)\}$. On the event $E_{C,\eps}$, let $r=\inf\{r'\ge 0: \mub_{\gffd-C} (r'\D) = \frac12\}$. Similarly as in the construction of $\sqrt{8/3}$-cone in Section~\ref{sec:lqg}, let $\gffd^C$ be the field on $r^{-1} \D$ defined by $\gffd^C(z)=\gffd(rz)+C+Q\log r$ so that $(r^{-1}\D,\gffd^C)$ is equivalent to $(\D,\gffd-C)$ with $\mub_{\gffd^C}(\D)=\frac12$. 
According to \cite{sphere-constructions}, by first sending $C\rta\infty$ and then $\eps\rta0$,
the law of the random measure $\mub_{\gffd}$ conditioned on $E_{C,\eps}$ converges weakly to the random measure supported on the entire $\C$, which can be viewed as the $\sqrt{8/3}$-LQG measure of a random field $\gffs$. Namely, $\mub_{\gffd}$ equals $e^{\sqrt{8/3}\gffs}dxdy$ in the sense of~\eqref{eq:area}. We call the $\sqrt{8/3}$-LQG surface represented by $(\C \cup\{\infty\}, \gffs)$ the 
\emph{unit area $\sqrt{8/3}$-LQG sphere}. This surface has the topology of a sphere with a marked point at $\infty $, since the disk boundary collapse into a single point $\infty$.
\begin{remark}\label{rmk:sphere}
	There are various constructions of the unit area $\sqrt{8/3}$-LQG sphere in the literature. A limiting construction similar to the one in Section~\ref{sec:lqg} for the disk case is provided in \cite{wedges}, as well as an explicit description of the field in polar coordinates. More constructions related to the mating-of-trees are considered in \cite{sphere-constructions}. Our limiting construction above is immediate from the mating-of-trees perspective (see Theorem~\ref{thm:sphere}). In \cite{dkrv-lqg-sphere}, the field for the unit area $\sqrt{8/3}$-LQG sphere is constructed as a conformal field theory (CFT) on the Riemann sphere. In \cite{ahs-sphere}, it is shown that the CFT construction is equivalent to the ones in \cite{wedges,sphere-constructions}. 
	Finally, there is a random metric space called the \emph{Brownian map}, which is the Gromov-Hausdorff limit of the random triangulations considered in our paper (seen as metric spaces), as well as many other similar random planar maps model \cite{legall-uniqueness,miermont-brownian-map}. Recently, Miller and Sheffield~\cite{lqg-tbm1,lqg-tbm2,lqg-tbm3} managed to conformally embed the Brownian map into the Riemann sphere to obtain an instance of the unit area $\sqrt{8/3}$-LQG sphere.
\end{remark}

Now we investigate $\etad$ and $\Zd$ as we perform the above conditional limit. Since $\etad$ modulo parametrization is independent of $\gffd$, the limit of $\etad$ does not depend on the conditioning. The limiting curve $\etas$ modulo parametrization is exactly the whole-plane space-filling $\SLE_6$. 
We parametrize $\etas$ such that $\etas(0)=\infty$ and $\mub_{\gffs}(\etab([0,t])) =t$ for all $t\in[0,1]$.
Let $\Zs:=(\Ls,\Rs)$ be the boundary length process defined as the process $\Zb$ in Theorem~\ref{thm:mot} but with $(\gff,\etab)$ replaced by $(\gffs, \etas)$. 
Let $$\Zb^{\op D, C}_t= e^{-\frac12\sqrt{8/3}C} \Zd_{te^{\sqrt{8/3}C}} \qquad \forall t \in [0,\mub_{\gffd-C} (\D)]. $$

By definition of $E_{C,\eps}$, the process $\Zs$ is the weak limit of $\Zb^{\op D,C}$ conditioned on $E_{C,\eps}$ as $C\to\infty$ and $\eps\to 0$. 
Note that the time/space rescaling makes $\Zb^{\op D,C}$ a Brownian cone excursion with variance $\beta$ starting from $(0,e^{-\frac12\sqrt{8/3}C})$ instead of $(0,1)$. The conditioning on $E_{C,\eps}$ means that the duration of this excursion is in $[1,1+\eps]$. Therefore it is elementary to see that the limiting process $\Zs$ exists (see for example \cite[\xin{Section~3}]{sphere-constructions}) and can be interpreted as $\Zb|_{[0,1]}$ conditioned on $\{\Zb_t\in (0,\infty)^2\; \forall t\in (0,1)\textrm{ and } \Zb_0=\Zb_1=(0,0) \}$. We call the law of $\Zs$ the \emph{Brownian excursion in the first quadrant of duration 1}. 

Now we are ready to state the sphere variant of Theorem~\ref{thm:mot}.
\begin{theorem}\label{thm:mot-sphere}
	Let $(\gffs,\etas)$ be defined as above and let $\Zs$ be its boundary length process. Then $\Zs$ has the law of the Brownian excursion in the first quadrant with duration 1. Moreover, $(\gffs,\etas)$ modulo rotations about the origin is measurable with respect to $\Zs$.
\end{theorem}
As in the disk case, there are two ways to carry out our constructions in Sections~\ref{thm:mot}--\ref{sec:cont-piv} in the sphere case. If we replace the field $\gff$ by $\gffs$ in these constructions, the construction of the branching $\SLE_6$ $\taus$, the $\CLE_6$ $\Gabs$, and the LQG pivotal measures $\nus,\nus_{\eps}$ associated with $(\gffs,\etas)$ can be defined in the same way with small modifications. On the other hand, they can also be constructed by taking the weak limit of the disk case as above. We leave the details to the reader.

The laws of $\Gabs$ (on $\mathbb{S}=\C\cup \{\infty\}$) and of $\Gab$ (on $\C$) are the same modulo parametrization. The only difference lies in the parametrization introduced by $\gffs$ and $\gff$. 
The same comment applies to $\taus$. As in the disk case, $\nus_{\eps}$ is a finite measure for all $\eps>0$.

\subsection{Percolation crossing events and chordal SLE$_6$}\label{sec:crossing}

Macroscopic crossing events is a family of classical observables for percolation. Suppose $D$ is a simply connected domain as in Section~\ref{sec:sle}, 
with the triangular lattice of mesh size $n^{-1}$ on top of it. Consider also the Bernoulli-$\frac12$ site-percolation $\sigma_n$ on the triangular lattice. Given four distinct points $A_1,A_2,A_3,A_4$ clockwise aligned on $\bdy D$, one considers the \emph{black crossing event}, which is the event that there is a black crossing between the clockwise arcs $A_4A_1$ and $A_2A_3$ on $\bdy D$. It is proved by Smirnov \cite{smirnov-cardy} that as the mesh size goes to zero, the probability of this crossing event converges to an explicit formula predicted by Cardy \cite{cardy-formula}, in which the limiting probability only depends on $(D,A_1,A_2,A_3,A_4)$ through its cross-ratio, hence is conformally invariant. 
\nina{A crucial observation (going back at least to Schramm \cite{schramm0}) is that crossing events are functions of the percolation interfaces, hence the conformal invariance of their probabilities is a shadow of the conformal invariance of $\SLE_6$ (and $\CLE_6$). In this section, we explain this observation and its realization in the context of mating of trees on a $\sqrt{8/3}$-LQG disk. }
As a byproduct, we give the mating-of-trees representation of chordal branches of the branching $\SLE_6$,
as well as the future/past decomposition relative to a boundary point.

\begin{figure}
	\centering
	\includegraphics{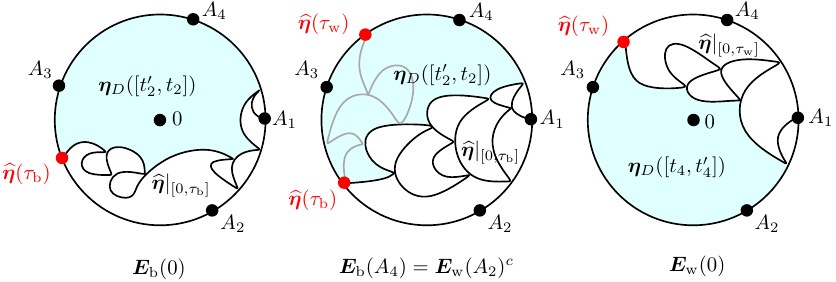}
	\caption{Illustration of the three events $\Eb_{\op b}(0)$, $\Eb_{\op b}(A_4)=\Eb_{\op w}(A_2)^c$, and $\Eb_{\op w}(0)$.
	}
	\label{fig:crossing}
\end{figure}
Suppose $\wh\etab$ is a chordal $\SLE_{6}$ on $(D,A_1,A_3)$, which is the scaling limit of the percolation interface from $A_1$ to $A_3$ if we assign black (respectively, white) boundary condition on the counterclockwise (respectively, clockwise) arc $A_1A_3$ on $\bdy D$. 
Let $\tau_{\op b}$ be the first time $\wh\etab$ hits the clockwise arc $A_2A_3$ and define 
\begin{equation}\label{eq:crossing}
	\Eb_{\op b}(v):=\{\textrm{$v$ and $A_3$ are in same connected component of $D\setminus \wh \etab([0,\tau_{\op b}])$}\} \qquad \forall v\in \ol D.
\end{equation}
This event is represented in Figure \ref{fig:crossing}. In the context of critical percolation on the triangular lattice, one can define the event $E^\triangle_{\op b}(v)$ that $\{A_3,v\}$ and the clockwise arc $A_1A_2$ are \emph{separated by a black crossing}, 
that is to say, 
there is a black path separating $D$ into two components, 
one containing $\{A_3,v\}$ and the other containing $A_1A_2$. 
It is elementary to check that the event $E^\triangle_{\op b}(v)$ (for both $v=A_4$
and $v\in D$) agrees with the discrete analogue of $\Eb_{\op b}(v)$ defined using the percolation interface from $A_1$ to $A_3$ in place of $\wh\etab$. 
\nina{Using this, it is not hard to see that in the setting of critical percolation on the triangular lattice, for each $v\in D\cup\{A_4 \}$, $(\1_{\Eb_{\op b}(v)},\wh\etab)$ is the joint scaling limit of $\1_{E^\triangle_{\op b}(v)}$ and the discrete percolation interface started from $A_1$.}

Let $\tau_{\op w}$ be the first time $\wh\etab$ hits the clockwise arc $A_3A_4$ and define
\begin{equation}\label{eq:crossing3}
	\Eb_{\op w}(v):=\{\textrm{$v$ and $A_3$ are in same connected component of $D\setminus \wh \etab([0,\tau_{\op w}])$}\} \qquad \forall v\in \ol D.
\end{equation}
Then by symmetry $\Eb_{\op w}(A_2)$ is the scaling limit of the event that there is a white crossing between the clockwise arcs $A_1A_2$ and $A_3A_4$ (for critical percolation on the triangular lattice). Furthermore, for $v\in D$, $\Eb_{\op w}(v)$ is the scaling limit of the event that there is a white path separating $\{A_3,v \}$ and the clockwise arc $A_4A_1$; see Figure \ref{fig:crossing}.

Now we consider the setup in Theorem~\ref{thm:mot-disk}, where $(\D,\gffd)$ is a particular representative of the unit length $\sqrt{8/3}$-disk with an independent counterclockwise space-filling $\SLE_6$ $\etab$ on $(\D,1)$, parametrized by $\mub_{\gffd}$-mass. Let $\Zd=(\Ld,\Rd)$ be the boundary length process of $(\gffd,\etad)$. Let $\ell,\ell',r,r'$ be four constants such that $0<\ell'<\ell, 0<r'<r$ and $\ell+r=1$. We specify $\{A_i\}_{1\le i\le 4}$ on $\bdy \D$ by setting $A_1=1$ and requiring that the $\nub_{\gffd}$-length 
of the clockwise and counterclockwise arc $A_1A_3$ is $\ell$ and $r$ respectively, while the $\nub_{\gffd}$-length of the clockwise arc $A_1A_2$ (respectively, $A_4A_1$) is $\ell'$ (respectively, $r'$). Define 
\[
t_{3}=\inf\{t\ge 0: \Rd_t<\ell\}.
\]
Then $t_3$ is the almost surely unique time that $\etad$ visits $A_3$. Let $\etadh$ be the branch of $\taud$ targeted at $A_3$. Then conditioning on $\gffd$, the curve $\etadh$ is a chordal $\SLE_6$ on $(\D, A_1, A_3)$. 
The space-filling curve $\etad|_{[0,t_3]}$ can be obtained by considering the curve $\etadh$, and filling a bubble (i.e., a complementary component) of $\etadh$ by a space-filling curve immediately after enclosure if and only if the bubble does not share a non-empty boundary arc with the clockwise arc $A_1A_3$.
The set $\etad|_{[0,t_3]}$ is the disk version of the past wedge relative to time $t_3$. The curve $\etadh$ can be considered as the spine inside this past wedge. 

The set $\D\setminus \etab([0,t_3])$ is the disk version of the future wedge relative to $t_3$, which is exactly the union of the connected components of $\D\setminus\etadh$ that share a non-empty boundary arc with the clockwise arc $A_1A_3$. These components are filled by $\etad|_{[t_3,\frk a]}$,
where $\frk a=\mub_{\gffd}(\D)$. Let $\cut(t_3)$ be the set of times at which both processes $\Ld$ and $\Rd$ achieve a running infimum relative to time $t_3$, that is,
$$\cut(t_3)=\{t\in[t_3, \frk a]: \Ld_{t'}> \Ld_t\;\textrm{and}\; \Rd_{t'}> \Rd_t, \; \forall \; t'\in [t_3,t) \}.$$
Then, similarly as in the case of the $\sqrt{8/3}$-LQG cone, the connected components of $\D\setminus \etab([0,t_3])$ are filled by $\etad$ restricted to the connected components of $[t_3,\frk a]\setminus \cut(t_3)$. Let 
\begin{equation}\label{eq:t3'}
	t_2=\sup\{t\in \cut(t_3): \Rd_t> \ell'\}
	\qquad \textrm{and}\qquad t_2'=\inf\{t< t_2: \Ld_{t'} > \Ld_{t_2}, \;\forall\; t'\in [t,t_2] \}.
\end{equation}
Then $t_2$ is the left endpoint of the connected component $I$ of $[t_3,\frk a]\setminus \cut(t_3)$ such that $A_2\in \etad(I)$.
Furthermore, $[t_2',t_2]$ is a cone interval containing $t_3$.
In fact, $\etad([t_2',t_2])$ is the closure of the connected component of $\D\setminus \wh \etab([0,\tau_{\op b}])$ with $A_3$ on its boundary, 
where $\tau_{\op b}$ is as in~\eqref{eq:crossing} with $D,\wh \etab$ replaced by $\D,\etadh$. 
Moreover, $t'_2$ (respectively, $t_2$) is the first (respectively, last) time $\etad$ visits $\etadh(\tau_{\op b})$.
Therefore we get the following mating-of-trees description for crossing events:
\begin{equation}\label{eq:crossing2}
	\Eb_{\op b}(0)=\{0\in [t_2',t_2]\}\qquad\text{and}\qquad \Eb_{\op b}(A_4)=\Big\{ \inf_{t\in [0,t_2']} \Rd_t > 1-r'\Big\}. 
\end{equation}

Next we give the mating-of-trees representation of $\Eb_{\op w}(A_2)$ and $\Eb_{\op w}(0)$. 
Since $\Eb_{\op w} (A_2)$ (respectively, $\Eb_{\op b} (A_4)$) is the event that $\wh\etab$ hits the arc $A_2A_3$ (respectively, $A_3A_4$) before the arc $A_3A_4$ (respectively, $A_2A_3$), the event $\Eb_{\op w} (A_2)$ is complementary to $\Eb_{\op b}(A_4)$, thus can be treated by~\eqref{eq:crossing2}. 
To treat $\Eb_{\op w}(0)$, we need to consider the past wedge $\etad([0,t_3])$.
As in the infinite volume case in Section~\ref{sec:dictionary-branch}, we say that a time $0\le t<t_3$ is \emph{ancestor-free} relative to time $t_3$ if it has no ancestors in $[t,t_3]$ (equivalently $t$ is not contained in the interior of a cone interval $I\subseteq [0,t_3)$),
and we let $\ans(t_3)$ be the ancestor-free times relative to $t_3$. Let
\begin{equation}\label{eq:t4}
	t_4= \inf\{t\in \ans(t_3):\Rd_t\le 1-r'\}.
\end{equation}
Then $\etab(t_4) =\wh\etab_D(\tau_{\op w})$, where $\tau_{\op w}$ is as in~\eqref{eq:crossing3} with $D,\wh \etab$ replaced by $\D,\etadh$.
Note that $t_4$ is the last time $\etad$ visits $\etadh(\tau_{\op w})$,
which is the position where $\etadh$ hits the clockwise arc $A_3A_4$.
Let 
\begin{equation}\label{eq:t4'}
	t'_4= \inf\{t\in \cut(t_3): \Ld_t <\inf_{s\in[t_4,t_3]}\Ld_{s}\}.
\end{equation} 
Then $\etad([t_4,t_4'])$ is the closure of the connected component of $\D\setminus \wh \etab([0,\tau_{\op w}])$ with $A_3$ on its boundary. 
Therefore we have 
\begin{equation}\label{eq:crossing4}
	\Eb_{\op w}(0)=\{0\in [t_4,t_4'] \}. 
\end{equation}
\referee{Note that the characterization~\eqref{eq:crossing2} and~\eqref{eq:crossing4} of the crossing events parallel the ones established for their discrete counterparts  in Section~\ref{sec:crossing-discrete} (see Facts~\ref{fact:d},~\ref{fact:e} and~\ref{fact:g}).}

Above we have considered two natural crossing events $\Eb_{\op b}(v)$ and $\Eb_{\op w}(v)$. It is equally natural to consider the events that there is a black or white crossing separating $\{A_i,0\}$ and $\{A_1,A_{5-i}\}$ for $i=2,3$. 
As pointed out in Remark~\ref{rmk:branch}, $\Gabd$ determines $\etad$.
It is proved in \cite[\xin{Section 4}]{ig4} that the law of the $\CLE_6$ $\Gabd$ associated with $(\gffd,\etad)$ modulo parametrization is invariant under all M\"obius transformations of $\D$, even if the root $1$ is not preserved. By the boundary re-rooting invariance of $(\D,\gffd)$, $\Gabd$ almost surely determines a clockwise and a counterclockwise space-filling $\SLE_6$ on $(\D, A_i)$ for $i=1,2,3,4$ by re-rooting. By Theorem~\ref{thm:mot-disk}, this induces eight Brownian cone excursions coupled together as boundary length processes, each of which has its own advantage in describing certain above mentioned crossing events. However, none of them will determine all the crossing events simultaneously in a simple manner, due to the fact that introducing the root breaks the rotational symmetry.

It is not surprising that $\CLE_6$ determines crossing events in a more symmetric manner \cite{camia-newman-full,gps-pivotal}. In fact, using $\CLE_6$, one can determine the scaling limit of crossing events even when $(D,A_1,A_2,A_3,A_4)$ is replaced by any $(D',A'_1,A'_2,A'_3,A'_4)$ where $D'\subset D$ and $A_i'\in \bdy D'$ for $i=1,2,3,4$.
The collection of all such crossing events provides a fourth description of the full scaling limit of critical planar percolation
which is called the \emph{quad crossings} and was first introduced in \cite{ss-planar-perc}. The mating-of-tree framework is suitable for encoding certain partial information about quad crossings, but is not
suitable for encoding all quad crossing information simultaneously.
This is a shortcoming of the mating-of-trees approach to the scaling limit of random percolated triangulations.
The second and third named authors' work \cite{hs-quenched}, based on this paper as well as \cite{ghs-metric-peano,ghss18, natural}, overcome this issue and eventually established a much stronger scaling limit result for percolated random triangulations. See Remark~\ref{rmk:CLE} and~\ref{rmk:LDP} for further comments.

\section{Convergence of percolated triangulations to $\CLE_6$ on $\sqrt{8/3}$-LQG}\label{sec:conv}
In this section we state our main convergence results. The proofs are given in Section~\ref{app:conv}. 

We first specify a notion of convergence for parametrized curves. For a metric space $(B,d_B)$, some intervals $I_1,I_2\subset\R$, and some functions $f_1:I_1\to B$ and $f_2:I_2\to B$, we define
\eqb
d_{\op{p}}(f_1,f_2) = \inf_{\psi} \sup_{t\in I_1} \Big(|\psi(t)-t|+d_{B}( f_1(t),f_2(\psi(t)) )\Big),
\label{eq35}
\eqe
where the infimum is taken over all increasing bijections $\psi:I_1\to I_2$. 
Let $d_{\BB S^2}$ denote the spherical metric on $\C$ obtained by considering stereographic projection of the Riemann sphere to $\C\cup \{\infty\}$. When we discuss convergence of curves below we will often use the metric $d_{\op{p}}$ with $(B,d_B)=(\C,d_{\BB S^2})$. We say that a sequence of curves $\xi_n:I_n\to \C$ \emph{converge as parametrized curves} in $(\C,d_{\BB S^2})$, if they converge for the metric $d_{\op{p}}$. If $B=\R^m$ for some $m\in\N_+$ we assume $d_B$ is equipped with the Euclidean metric, except for $B=\C$, where we assume $d_B=d_{\BB S^2}$.

\subsection{Scaling limit for infinite percolated triangulation (whole-plane setting)} \label{sec:cone}
Let us give a quick preview of this subsection. Let $(\gff, \etab,\Zb)$ be as in Theorem~\ref{thm:mot}. 
We will consider a sequence of uniformly random words $w^n\in\mK^\infty$, coupled with $(\gff, \etab,\Zb)$ in such a way that the associated properly rescaled lattice walks $Z^n$ converge almost surely to $\Zb$. 
By Theorem~\ref{thm:UIPT}, each walk $w^n$ corresponds to a percolated UIPT $(M_n,\sigma_n)$.
We then use the curve $\etab$ (and its discrete analogues $\etae$ and $\etavf$) to define an embedding $\phi_n$ of $(M_n,\sigma_n)$ in $\C$. Roughly speaking, the vertices of $(M_n,\si_n)$ are drawn along the curve $\etab$ in such a way that about $n$ vertices appear along each piece of $\etab$ parametrized by a unit interval. We then use $\phi_n$ to define the pushforward of the counting measure of the vertex set $V(M_n)$ onto $\C$,
an embedding of the exploration tree $\tau^*$,
an embedding of the percolation loops
and the pushforward of the counting measure of $\eps n$-pivotal points onto $\C$.
Our main result (Theorem~\ref{thm1}) is that under suitable normalization, these quantities converge jointly to their continuum counterparts defined in terms of $(\gff,\etab)$ in Section~\ref{sec:dictionary}.

We now give the precise definitions. Let $w\in \{a,b,c \}^\ZZ$ be a bi-infinite word, and let $Z=(L_k,R_k)_{k\in\Z}$ be the associated bi-infinite Kreweras walk on $\Z^2$ satisfying $Z_0=(0,0)$. 
Recall the constant $\beta>0$ mentioned in Remark~\ref{rmk:variance}.
For $n\in\N$, let $Z^n=(Z^n_t)_{t\in\R}$ be the rescaled version of $Z$ defined by 
\begin{equation}\label{eq:Zn}
	Z^n_t:=\frac 12 \sqrt{\beta/n} Z_{\lfloor 3nt\rfloor }
\end{equation}
\xin{Let $\Zb$ be as in Theorem~\ref{thm:mot}.} Then it is easy to see that if $w$ is uniformly random in $\{a,b,c \}^\ZZ$, then $Z^n$ converges in law to $\Zb$ in local uniform topology.

\xin{Throughout this section  we consider a coupling  of $(Z^n)_{n\in \N_+}$ and \emph{$(\gff,\etab,\Zb)$} such that $Z^n$ almost surely converges to $\Zb$ uniformly on compact sets, where $(\gff, \etab,\Zb)$ is as in Theorem~\ref{thm:mot}. The existence of such a coupling is ensured by Skorokhod's representation theorem. Under this coupling,  for $n\in\N_+$, we let $w^n$ be the bi-infinite word associated with $Z^n$. Note that in this coupling $Z^n$ does \emph{not} come from a common word.}

Let $(M_n,\sigma_n)=\Phi^\infty(w^n)$. Recall from Theorem~\ref{thm:UIPT}, that $(M_n,\sigma_n)$ is defined almost surely and has the law of the percolated loopless UIPT. We now define an embedding $\phi_n$ of $(M_n,\sigma_n)$ in $\C$.
Recall the bijections $\etae:\Z\to E(M_n)$ and $\etavf:\Z\to V(M_n)\cup F(M_n)$ (for simplicity we omit the dependence in $n$ in the notation of the functions $\etae$ and $\etavf$) in Definition~\ref{def:eta}. 
For each $n\in\N_+$, we define an embedding $\phi_n:V(M_n)\cup E(M_n)\to\C$ of $M_n$ as follows. For a vertex $v\in V(M_n)$ we define
\begin{equation}\label{def-phi}
	\phi_n(v) = 
	\etab\Big( \frac{1}{3n} \etavf^{-1}(v) \Big)\in \C.
\end{equation}
We also define the embedding $\phi_n$ on edges, by defining the image of $e=\{u,v\}\in E(M_n)$ to be
\begin{equation}\label{def-phie}
	\phi_n(e)=\frac{\phi_n(u)+\phi_n(v)}{2} \in\C.
\end{equation}

The following proposition, which is an immediate consequence of Lemma \ref{prop6}, shows that $\phi_n$ is a reasonable embedding.
\begin{prop}
	The following convergence holds in probability as $n\rta\infty$
	\eqbn
	\sup_{u,v\in V(M_n)\,:\,\{u,v \}\in E(M_n) } d_{\BB S^2}(\phi_n(u),\phi_n(v))
	\rta 0.
	\eqen
\end{prop}	
In fact, although this will not be proven in the present article, the preceding asymptotic result holds when the supremum is taken over pairs of vertices  $u,v$  of $M_n$ at graph distance $o(n^{1/4})$.

Next, we define the \emph{vertex measure} $\mu_n$ on $\C$ as the renormalized counting measure of (the embedding of) $V(M_n)$, where we assign mass $1/n$ to each vertex. In other words, for any Borel set $A\subset\C$,
\begin{equation}\label{eq:mu-def}
	\mu_n(A) = \frac{1}{n}\cdot\# \{ v\in V(M_n)\,:\, \phi_n(v)\in A \},
\end{equation}
where $\#$ represents the cardinality of a set.

The \emph{embedded space-filling percolation exploration} $\eta_n$ is a piecewise linear path in $\C$ which visits the edges of $M_n$ in the order they are treated in the space-filling exploration of $(M_n,\si_n)$ described in Section~\ref{sec:DFS-infinite}. 
In other words, 
$\eta_n:\R\to\C$ is the parametrized curve defined by
\begin{equation}\label{def-eta}
	\eta_n(t):=\phi_n(\etae(3tn)) \qquad \textrm{if } 3tn\in\ZZ,
\end{equation}
and by linear interpolation for other values of $t\in\R$. 
\old{Moreover, in any unit interval $I$, $\eta_e(I)$ contains $3n$ edges and hence $n+o_n(1)$ vertices. Therefore the $\mu_h$-mass of $\eta_e(I)$ equals 1.}

We now define notions related to the convergence of the DFS tree $\tau^*_n=\Delta_{M_n}(\sigma_n)$ (introduced in Section~\ref{sec:DFS-infinite}) 
toward the branching $\SLE_6$ $\taub^*=\{\wh\etab\}_{z\in \C}$ associated with $(\gff,\etab)$ (see Definition~\ref{def:branch}). 
Recall that although the branches of $\taub^*$ are simultaneously defined for all points on $\C$ as curves modulo parametrization, their quantum natural parametrization is only defined when
$$t_z:=\sup\{ t\in\R\,:\,\etab(t)=z \}$$ 
is a backward stopping time or an envelope closing time. Therefore the convergence we will consider is in the sense of finite marginals (see the definition below) and not in terms of contour function or Gromov-Hausdorff distance, which is commonly used in the context of continuum random trees.

In order to talk about finite marginal convergence, one needs to select finitely many points on the $\sqrt{8/3}$-LQG cone in the continuum and on the percolated UIPT in the discrete. There are various ways of doing so that work equally well. We will consider a way based on $\CLE_6$ loops and percolation cycles. 
We postpone the precise description to the statement of Theorem~\ref{thm1} and for now we introduce the following generic notation.
For a tuple $P\in\C^k$, let $\taub^*|_{P}$ denote the subtree consisting of the branches $\{\etab^z\}_{z\in P\cup\{0\}}$ of $\taub^*$. 

In the discrete,
for an edge $e\in E(M_n)$, we call \emph{path of $\tau^*_n$ toward $e$} the path of $\tau^*_n$ from $\infty$ to $u$, where $u\in V(M_n^*)$ is the endpoint of the dual edge $e^*=(u,v)$ which is the ancestor of the other endpoint $v$.
Recall from Definition~\ref{def:pi-w-inf} and Theorem \ref{thm:spinelt-inf} that the set of times $(T(k))_{k\in\ZZ^{\leq 0}}$ are defined so that the edges $\etae(T(k))$ are dual to the set of edges on the branch of $\tau^*_n$ toward $\etae(0)$ (equivalently the percolation path of the past map $(M_n^-,\si_n^-)$). 
Let $(T^n_s)_{s\leq 0}$ be given by
\begin{equation}\label{eq:Tn0}
	T_s^n=(3n)^{-1}T_{ \lfloor sn^{3/4} \rfloor }.
\end{equation} Define the \emph{embedded percolation exploration toward 0}, $\wh\eta^0=(\wh\eta^0_n(s))_{s\leq 0}$, by 
\begin{equation}\label{eq:eta-def}
	\wh\eta^0_n(s)
	:=\eta_n(T^n_s)
	=\phi_n( \etae(T_{sn^{3/4}})),\qquad\text{for}\,\,sn^{3/4}\in\ZZ^{\leq 0},
\end{equation}
and by linear interpolation for other $s$.
For any $z\in\C$, we define $\wh\eta^z_n$ in the exact same way, except that we consider the branch of $\tau^*_n$ toward the edge  $\etae(\lfloor n t_z \rfloor)$. We call $\wh\eta^z_n$ the \emph{embedded percolation exploration toward} $z\in\C$.

We now define an embedding of the full exploration tree $\tau^*_n$. 
For an edge $e^*$ of $M_n^*$ dual to an edge $e$ of $M$, we define the embedding $\phi_n(e^*)=\phi_n(e)$. We define the embedding of $\tau^*_n$ into $\C$, by drawing a straight line between $\phi_n(e^*_1)$ and $\phi_n(e_2^*)$ for each pair of edges $e_1^*,e_2^*$ of $\tau^*_n$ incident to a common vertex of $M_n^*$. Note that the embedding of $\tau^*_n$ is the union of the embedded percolation explorations $\wh\eta^z_n$ (but it could have some edge crossings). 

Fix $k\in \N_+$. For each $n\in \N_+$, consider a tuple $P_n=(p_1^n,\ldots,p_k^n)\in\C^{k}$ and the embedded ``subtree'' $\tau_n^*|_{P_n}$, which is the union of the branches $\wh\eta^{0}_n$ and $\wh\eta^{p_1^n}_n,\ldots,\wh\eta^{p_k^n}_n$. 
We say that the sequence of trees 
$\tau_n^*|_{P_n}$ \emph{converges} 
if there is a tuple $(p_1,\ldots,p_k)\in\C^{k}$ such that, for all $i,j\in[k+2]$ with the convention $p_{k+1}^n=p_{k+1}=0$ and $p_{k+2}^n=p_{k+2}=\infty$, the branch $\tau_{i,j,n}^*$ of $\tau_n^*|_{P_n}$ between $p_i^n$ and $p_j^n$ converges as a parametrized curve to a curve between $p_i$ and $p_j$ for the metric space $(\C,d_{\BB S^2})$, when $\tau_{i,j,n}^*$ is considered as a curve parametrized in such a way that it takes time $n^{-3/4}$ to trace each edge as in \eqref{eq:eta-def}.

Next, we discuss the convergence of the percolation cycles to the $\CLE_6$ $\Gab$ associated with $(\gff,\etab)$ (Definition~\ref{def:CLE}). 
In order to talk about convergence of individual cycles, it will be convenient to enumerate $\CLE_6$ loops in $\Gab$ in such a way that the same enumeration rule can also produce an enumeration of the percolation cycles of $(M_n,\sigma_n)$. There are various ways of doing so which would work equally well. Here we choose an enumeration based on the area enclosed by the cycles.

In the continuum, for each $\gab\in \Gab$ recall that $\reg(\gab)$ is the region enclosed by $\gab$. 
To simplify notation, we write $\mub_{\gff}(\reg(\gab))$ as $\area(\gab)$. 
Let $\gab'\in \Gab$ be the $\CLE_6$ loop with the smallest enclosed region 
such that $0\in \reg(\gab')$ and $\area(\gab')\geq 1$. For each $\gab\in \Gab$, let $\sma(\gab,\gab')$ be the $\CLE_6$ loop with the smallest enclosed region that encloses both 
$\gab$ and $\gab'$. Define the \emph{value} of $\gab$ by
\eqb
\op{val}(\gab):= \area (\sma(\gab,\gab'))+\area(\gab)^{-1}, 
\label{eq:value0}
\eqe
Note that for each $r>0$, the number of $\CLE_6$ loops with value smaller than $r$ is finite and a.s.\ no two cycles have the same value. The $\CLE_6$ loops can be enumerated so that their values are increasing. We denote this enumeration by $\Gab=\{\gab_j\}_{j\in \N_+}$.

In the discrete, for $\gam$ a percolation cycle of $(M_n,\si_n)$, we denote by  $\op{area}_n(\gam)$ the number of vertices enclosed by $\gam$, divided by $n$. 
Let $\gam'_n$ be the smallest percolation cycle such that area$_n(\gam'_n)\geq 1$, and such that both endpoints of the edge $\etae(0)$ are enclosed by $\gam'_n$. For a percolation cycle $\gam$ of $(M_n,\si_n)$, we let $\sma(\gam,\gam'_n)$ be the smallest common ancestor of $\gam$ and $\gam'_n$ in $\cltree(M_n,\sigma_n)$ (Definition~\ref{def:tree-clusters}) when we identify a percolation cluster with its outside-cycle (this is well-defined since $\cltree(M_n,\sigma_n)$ is one-ended almost surely, see the Section \ref{sec:cluster-tree-inf}). Define the \emph{value} of $\gam$ by
\eqb
\op{val}(\gam):=\op{area}_n(\sma(\gam,\gam'_n))+\op{area}_n(\gam)^{-1}.
\label{eq:value}
\eqe
As in the continuum, for each $r>0$, the number of percolation cycles with value smaller than $r$ is finite. Therefore
the cycles can be enumerated so that their values are non-decreasing, with draws resolved in an arbitrary way. We denote the enumerated percolation cycles of $(M_n,\si_n)$ by $\{\gam^n_j \}_{j\in \N_+}$.

We now define a parametrized embedding of the percolation loops.
Recall that each percolation cycle $\gam$ of $(M_n,\si_n)$ crosses a set of bicolor edges of $M_n$ separating two clusters.
Let $|\gam|$ denote the number of edges on the percolation cycle $\gam$.
For $j\in\N_+$, let $T^j(0)<\dots<T^j(|\gam_j^n|-1)$ be the set of times $t\in\Z$ such that $\etae(t)$ is an edge crossed by the percolation cycle $\gam_j^n$. Let $T^{n,j}=(T_s^{n,j})_{s\in[0,n^{-3/4}|\gam_j^n|)}$ be defined by the following rescaling, 
\begin{equation}\label{eq:T-cle}
	T_s^{n,j}=(3n)^{-1}T^j(\lfloor sn^{-3/4} \rfloor)
\end{equation}
We define the embedding $\ga^n_j:[0,n^{-3/4}(|\gam_j^n|-1)]\to\C$ of $\gam^n_j$ by setting 
\begin{equation}\label{eq:def-gab}
	\ga^n_j(s):=\eta_n( T_{s}^{n,j} )=\phi_n\big(\etae\big(T^j(sn^{-3/4})\big)\big)\qquad\text{for}\,\, sn^{3/4}\in\{0,\dots,|\gam_j|-1 \},
\end{equation}and by linear interpolation for other $u$. We are slightly abusing notation here by using the same symbol for the percolation cycle $\gam_j^n$ of the triangulation and its parametrized embedding into $\C$.

Lastly we turn our attention to the pivotal measure. Recall the LQG pivotal measure $\nub$ given by Definition~\ref{def:pivotal}  and its restriction $\nub_{\eps}$ to the set $\piv_\eps$ of $\eps$-pivotal points given by Definition~\ref{def:piv-eps}. Recall also the four types of pivotal points defined at the end of Section~\ref{sec:cont-piv}. 
We now define a partition of $\nub_{\eps}$ as in Definition~\ref{def:pivotal} based on Lemma~\ref{lem:pivotal-continuous} and the notion of types.
For $j\in \N_{+}$, let $\nub_j^{\eps,1}$ (respectively, $\nub_j^{\eps,2}$) denote $\nuloc_\eps$ restricted to the $\eps$-significant pivotal points of $\gab_j$ of type 1 (respectively, 2). Recall the measure $\nub_{\gab,\gab'}$ of Definition~\ref{def:pivotal}. For $i\neq j\in \N_{+}$, write $\nuloc_{\gab_j,\gab_i}$ as $\nub_{j,i}$ for simplicity. Recall that no cutoff $\eps$ is needed because the measure $\nub_{j,i}$ is finite. By our convention, $\nub_{j,i}=0$ if $\env(\gab_i)\subset \env(\gab_j)$. We do not introduce further notation to distinguish pivotal points of type 3 and type 4, since it can be determined by considering whether the two loops involved are nested or not. 
Note that by definition,
\begin{equation}\label{eq:gabeps-sum}
	\nub_\eps=\sum_{j\in\N_+}(\nub_j^{\eps,1}+\nub_j^{\eps,2})+\sum_{(j,i)\in S_\eps}\nub_{j,i},
\end{equation}
where $S_\eps$ is the set of pairs $(i,j)\in\N_+^2$ such that $\gab_i$, $\gab_j$ have area at least $\eps$ and either they are unested, or they are nested and the area between them is at least $\eps$. 

Now we define the corresponding discrete pivotal measures. 
Given $j\in \N_{+}$, let $\cP_j^{\eps,1}\subset V(M_n)$ (respectively, $\cP_j^{\eps,2}\subset V(M_n)$ ) be the set of $\eps n$-pivotal points $v\in V(M_n)$ of type 1 (respectively, type 2) associated with $\gam^n_j$, and let
\eqb
\nu^{\eps,k}_{j,n}:=n^{-1/4}\sum_{v\in\cP^{\eps,k}_j} \delta_{\phi_n(v)} \qquad \textrm{for }k\in\{1,2\}.
\label{eq:piv-dis}
\eqe
Let $i\neq j\in \N_{+}$. If the envelope closing time of $\gam^n_j$ is larger than that of $\gam^n_i$, then we define $\cP_{j,i } \subset V(M_n)$ as the set of pivotal points between $\gam^n_i$ and $\gam^n_j$. 
Otherwise, we define $\cP_{j,i } =\emptyset$. 
Finally, let 
\eqb
\nu_{j,i,n}:=n^{-1/4}\sum_{v\in\cP_{j,i}} \delta_{\phi_n(v)}.
\label{eq:piv-dis2}
\eqe

\begin{theorem}
	Consider the setting above, 
	where $(\gff,\etab,\Zb)$ is as in Theorem~\ref{thm:mot} and the sequence $(Z^n)$ of rescaled random walks is coupled with $(\gff,\etab,\Zb)$ in such a way that 
	$(Z^n)_{n\in\N}$ converges almost surely to $\Zb$. 
	Let $\Gab=\{\gab_i \}_{i\in \N_{+}}$, $\taub^*=(\wh\etab^z)_{z\in\C}$, and $\nub^{\ep,1}_{j}$, $\nub^{\ep,2}_{j},\nub_{i,j}$ be the CLE$_6$, the branching SLE$_6$, and the LQG pivotal measures associated with $(\gff,\etab)$.
	Then the following quantities converge jointly in probability as $n\rta\infty$. 
	\begin{compactitem}
		\item[(i)] Area measure: the vertex counting measure $\mu_n$ restricted to any ball converges in the weak topology to the $\sqrt{8/3}$-LQG area measure $\mub_{\gff}$.
		
		\item[(ii)] Space-filling percolation exploration: $\eta_n=(\eta_n(t))_{t\in\R}$ converges uniformly to the space-filling SLE$_6$ $\etab$.
		
		\item[(iii)] Percolation cycles: the embedded percolation cycles $\ga_1^n,\ga_2^n,\dots$ converge to the CLE$_6$ loops $\gab_1,\gab_2,\dots$ as parametrized curves in $(\C,d_{\BB S^2})$. For all $i,j\in\N_+$, $\1_{\ga_i^n\subset\reg(\ga_j^n)}$ converges to  $\1_{\gab_i\subset\reg(\gab_j)}$, and $\area_n(\ga^n_j)$ converges to $\area(\gab_j)$.
		
		\item[(iv)] Pivotal measures: for any fixed $\ep>0$ and $i,j\in\N_+$, the pivotal measures $\nu^{\ep,1}_{j,n}$, $\nu^{\ep,2}_{j,n}$, $\nu_{i,j,n}$, and $\nu_\eps^n$ converge in the weak topology to $\nub^{\ep,1}_{j}$, $\nub^{\ep,2}_{j}$, $\nub_{i,j}$, and $\nub_\eps$. \nina{Furthermore, $\1_{\nu^{\ep,1}_{j,n}(\C)=0}$ converges to $\1_{\nub^{\ep,1}_{j}(\C)=0}$ and the analogous statement holds for the other three measures.}
	\end{compactitem}
	Given any $m\in \N_+$, let $E^m_n\subset E(M_n)$ (respectively, ${\bf E}^m\subset\C$) denote the set of edges of $M_n$ (respectively, points of $\C$) which are in the inside-region of at least one of the percolation cycles $\gam^n_1,\gam^n_2,\dots,\gam^n_m$ (respectively, $\gab_1,\gab_2,\dots,\gab_m$). 
	For any $k\in\N_+$, 
	we may sample edges $e_1^n,\dots,e_{k}^n\in E(M_n)$ uniformly and independently at random from $E_n^m$ and extend the coupling
	in such a way that $P_n := (\phi_n(e_1^n), \dots, \phi_n(e_{k}^n))$ converges in probability to a tuple $P$ of $k$ points independently sampled from $\mub_{\gff}$ restricted to ${\bf E}^m$.
	\begin{compactitem}
		\item[(v)] exploration tree finite marginals: the subtree $\tau^*_n|_{{P_n}}$ converges to $\taub^*|_{P}$. Furthermore, for any fixed $t\in\R$ and $z:=\etab(t)$ the curve $\wh\eta^z_n$ converges uniformly to $\wh\etab^z$. 
	\end{compactitem}
	\label{thm1}
\end{theorem}

\begin{remark}
	Observe that (as explained right after Theorem~\ref{thm2}) Theorem~\ref{thm1} gives convergence of many interesting observables of the percolated UIPT which do \emph{not} depend on the particular embedding we chose. 
\end{remark}

\begin{remark} \label{rmk:KPZ}
	The scaling $n^{-1}$ for area in \eqref{eq:mu-def}, $n^{-3/4}$ for the length of percolation loops in \eqref{eq:def-gab}, $n^{-1/2}$
	for the magnitude of the random walk in \eqref{eq:Zn} (equivalently, length of the frontiers of the near-triangulation), and $n^{-1/4}$ for pivotal measure in \eqref{eq:piv-dis} and \eqref{eq:piv-dis2}, 
	correspond to the scaling of $H^2,H^{3/2},H,H^{1/2}$ for their continuum analogs in \eqref{eq:taud}, \eqref{eq:Zd}, \eqref{eq:piv-scale} (with $H^2$ corresponding to the area). 
	The ratios $4:3:2:1$ between the scaling exponents is an instance of the KPZ relation\footnote{ \nina{The KPZ relation is named after  Knizhnik, Polyakov, and Zamolodchikov \cite{kpz-scaling}. We refer to \cite{shef-kpz} and reference therein for further details.}} 
	between the scaling exponents of planar fractals in Euclidean geometry and LQG geometry. The planar fractals involved here are: $\SLE_6$ curves, $\SLE_{8/3}$ curves, and pivotal points of $\CLE_6$.
\end{remark}

\begin{remark}\label{rmk:lt}
	Beside the convergence results stated in the above theorem, it will be clear from the proof that the discrete looptrees $\frk L(\ga^n_1),\frk L(\ga^n_2),\dots$ (Definition~\ref{def:lt-cluster}) converge in the Gromov-Hausdorff topology\footnote{More precisely, we are considering here the discrete looptrees $\frk L(\ga^n_i)$ as compact metric spaces $(V_i^n,d_{\op{gr}}^n)$, where $V_i^n$ is the vertex set of $\frk L(\ga^n_i)$ and $d_{\op{gr}}^n$ is the graph distance rescaled by $n^{1/2}$.}
	to the corresponding continuum looptrees $\frk L(\gab_1),\frk L(\gab_2)\ldots$.
	This convergence result is immediate from \cite[Theorem 4.1]{curien-kortchemski-looptree-def} and results in Section~\ref{app:conv}, since the convergence of the walk encoding a looptree implies the convergence of the looptree itself in the Gromov-Hausdorff topology, provided the height of the tree of bubbles times the renormalization factor for the walk goes to zero. The constraint is satisfied in our case since the height of the tree of bubbles is of order $n^{1/4}=o(n^{1/2})$, as it is encoded by an excursion of a walk with increments that are in the domain of attraction of a $3/2$-stable random variable \cite{Duquesne03,k13}. 
\end{remark}

\nina{\begin{remark}
		\label{rmk9}
		In this remark we will explain that by combining Theorem \ref{thm1} and \cite[Theorem 2]{camia-newman-full} (which is recalled in Lemma~\ref{lem:pivotal-continuous}) we get that 
		for any fixed $c>0$ and $\eps>0$, with probability converging to 1 as $n\rta\infty$, each $\eps n$-pivotal point $v$ for which $|\etavf^{-1}(v)|<c n$ is an $\eps n$-pivotal point of a unique type. Furthermore, by the same results, with probability converging to 1 as $n\rta\infty$, for any $v$ for which $|\etavf^{-1}(v)|<c n$, the symmetric difference $\Gam\Delta \Gam_v$ contains at most 3 cycles of area at least $\eps n$.
		We will only justify the latter assertion since the former assertion can be proved in the same manner. We will argue by contradiction, and we assume that there is a $\delta>0$ such that for arbitrarily small $n$, with probability at least $\delta$ we can find a pivotal point $v_n\in V(M_n)$ for which $|\etavf^{-1}(v_n)|<c n$ and the symmetric difference $\Gam\Delta \Gam_{v_n}$ contains more than 3 cycles of area at least $\eps n$. Since the domain $\etab([-c,c])$ is bounded we can find a subsequence of $n\in\N$ such that with probability at least $\delta$, $v_n$ exists and $\phi_n(v_n)$ converges to some point $z\in\C$ along this subsequence. By Theorem \ref{thm1}(i,iii) it holds that the symmetric difference $\Gab\Delta \Gab_z$ contains more than 3 CLE$_6$ loops of area at least $\eps$. This contradicts Lemma~\ref{lem:pivotal-continuous} and therefore allows us to conclude the argument. We remark that the exact same argument works in the finite-volume setting of Section \ref{sec:finite}, except that it is not necessary with any constraints on $|\etavf^{-1}(v_n)|$ in this setting since the domain is bounded.
\end{remark}}

\subsection{Scaling limit for finite percolated triangulations (disk and sphere settings)} \label{sec:finite}

In this section, we state our main scaling limit result for finite volume maps. 

We start with the setup for the disk. In the continuum, fix a constant $\frk m>0$, and let $(\gffd,\etad,\Zd)$ be defined as in Corollary~\ref{cor:fix-area}:
$(\D,\gffd)$ is the representative of the $\sqrt8/3$-LQG disk with length 1 and area $\frk m$ specified by the bulk and boundary re-rooting invariance;
$\etad$ is the counterclockwise space-filling $\SLE_6$ on $(\D,1)$ parametrized in such a way that $\etad(0)=\etad(\frk m)=1$ and $\mub_{\gffd}(\etab([0,t]))=t$ for all $t\in [0,\frk m]$;
$\Zd$ is the boundary length process of $(\gffd,\etad)$, 
which is a Brownian cone excursion of duration $\frk m$ from $(0,1)$. Let $\ub$ be the almost surely unique time such that $\etad(\ub)=0$. Recall that by \eqref{eq:uniform}, conditioning on $(\gffd,\etad)$, the time $\ub$ is uniform on $(0,\frk m)$.

We now define the discrete counterparts of $\gffd,\etad,\Zb$ and $\ub$.
Let $\{h_n\}_{n\in N_{+}}$ and $\{m_n\}_{n\in N_{+}}$ be two sequences of non-negative integers such that 
$$
\lim_{n\to\infty} \frac{h_n}2\sqrt{\beta/n}=1 \qquad\textrm{and}\qquad \lim_{n\to\infty} m_n/n=\frk m,
$$
where $\beta$  is the constant defined in Remark~\ref{rmk:variance}, and we assume $3m_n+2h_n<3\frk m n$.
Let $M_n$ be a uniformly chosen triangulation with a simple boundary having $h_n+2$ outer vertices and $m_n$ inner vertices. Let $\frk e_n$ be the root-edge of $M_n$. 
Conditioning on $M_n$, let $\frk v_n$ be a uniformly sampled inner vertex and let $\sigma_n$ be the site-percolation on $M_n$ where each inner vertices is colored black or white independently with probability $\frac12$. Moreover every outer vertex of $M_n$ is colored black except the origin of $\frk e_n$ which is colored white, so that $(M_n,\si_n)\in \bmT_P$. 
Recall Corollary~\ref{cor:chordal-case} about the bijection $\bPhi$. Let $Z=\bPhi^{-1}(M_n,\sigma_n)\in \smK$. 
Moreover all the $a$-steps and $b$-steps in $Z$ are matched. All the $c$-steps have a matching $a$-step, and all but $h_n$ $c$-steps have a matching $b$-step. Recall that the walk $Z$ starts at $(0,0)$, ends at $(0,-h_n)$ and has $3m_n+2h_n$ steps. As in \eqref{eq:Zn}, let
\begin{equation}\label{eq:Zdn}
	Z^n_t:=\frac {1}{2} \sqrt{\beta/n} (Z_{\lfloor 3nt\rfloor }+(0,h_n)),
\end{equation}
so that $Z^n$ starts at $(0,\sqrt{\beta/n}\cdot h_n/2)$ and ends at $(0,0)$. Let $u_n$ be such that $\etavf(3nu_n)=\frk v_n$ (see Definition~\ref{def:eta}). 

By Lemma~\ref{lem:uniform}, there exists a coupling of $(M_n,\sigma_n,\frk v_n)_{n\in \N_+}$ and $(\gffd,\etad, \ub)$ such that $\lim_{n\to\infty}Z^n=\Zd$ in uniform topology and $\lim_{n\to \infty}u_n=\ub$ almost surely. In the disk case we will work under such a coupling throughout this section. 

Now let us define the embeddings $\phi_n:V(M_n)\cup E(M_n)\to\D$ by \eqref{def-phi} and \eqref{def-phie} as in the UIPT case.
In this case the image of $\phi_n$ is inside $\ol\D$. 
We still define the measure $\mu_n$ by \eqref{eq:mu-def}, where the support of $\mu_n$ becomes $\ol \D$.
Let the space-filling percolation exploration $(\eta_n(t))_{t\in[0,\frk m]}$ still be defined as in the UIPT case (see \eqref{def-eta}) 
with the range of $t$ being $[0,\frk m]$ instead of $\R$. Note that under the coupling above $\phi_n(\frk e_n)\rta 1$ and $\phi_n(\frk v_n)\rta 0$ as $n\to \infty$.

In the disk case, besides the area measure, it is also natural to consider the boundary measure. For $n\in \N_+$, let $\nu_n$ be the uniform probability measure on the outer vertices of $M_n$ embedded into $\D$ via $\phi$. Note that $\nu_n$ is a measure on $\ol{\D}$, not on $\partial\D$. 

Recall the branching $\SLE_6$ $\taud=\{\etadh^z \}_{z\in \D\setminus \{1\} }$ associated with $(\gffd,\etad)$ (see Section  \ref{sec:disk}) and the exploration tree $\tau^*_n=\dfs(M_n,\sigma_n)$ (see Section \ref{sec:exploration-tree-from-walk}).
For each $z\in\ol\D\setminus\{1\}$, let $\wh\eta^z_n$ be the embedded percolation exploration from $\frk e_n$ to the edge $\etae^{-1}(\lceil 3nt_z\rceil )$ with $t_z=\sup\{t\in\R\,:\,\etad(t)=z \}$.
For $k\in \N_+$ and a $k$-tuple $P\in (\D\setminus \{1\})^k$, let $\taud|_{P}$ be the subtree $\{\etadh^z\}_{z\in P\cup\{0\}}$ of $\taud$, and let $\tau^*_n|_{{P}}$ be the embedded subtree of $\tau^*_n$ defined as in the UIPT case. We adopt the same notion of finite marginal convergence of $\tau^*_n$.

Recall from Section \ref{sec:disk} the $\CLE_6$ $\Gabd$ and the pivotal measure $\nub_{\op D,\eps}$ associated with $(\gffd,\etad)$. 
We define an enumeration of the $\Gabd$ which we again denote by $\gab_1,\gab_2,\dots$. But since $\mub_{\gffd}(\D)<\infty$, we do not rely on $\val$ in \eqref{eq:value} but simply require that $\area(\gab_j)$ is decreasing.
Similarly, for the percolation cycles of $(M_n,\perc_n)$, we enumerate them as $\gam^n_1,\gam^n_2,\dots$ in such a way that $\op{area}_n(\gam^n_j)$ is non-increasing, with ties broken in an arbitrary manner. 
The embeddings $\ga_j^n$ of the percolation cycles and the times $(T_s^{n,j})_{s\in [0,n^{-3/4}|\gam^n_j|)}$ are defined as in \eqref{eq:def-gab} and \eqref{eq:T}.
We still define the type of a pivotal point of $\Gabd$ and the measure $\nub^{\eps,1}_{j}$, $\nub^{\eps,2}_{j}$, and $\nub_{i,j}$ supported on pivotal points of various types in the same way as in the $\sqrt{8/3}$-LQG cone case. 
In the discrete, let $\nu^{\eps,1}_{j,n}$, $\nu^{\eps,2}_{j,n}$, and $\nu_{i,j,n}$ be defined as in~\eqref{eq:piv-dis} and \eqref{eq:piv-dis2}. This time all these pivotal measures are supported on $\D$.

Finally, we define quantities related to crossing events. 
Fix $0<\ell'<\ell$ and $0<r'<r$ such that $\ell+r=1$ as in Section~\ref{sec:crossing}. Let $A_1,A_2,A_3,A_4$ be four points on $\bdy \D$ clockwise aligned such that $A_1=1$ and the $\nub_{\gffd}$-length of the clockwise and counterclockwise arc $A_1A_3$ is $\ell$ and $r$ respectively, while the $\nub_{\gffd}$-length of the clockwise arc $A_1A_2$ (respectively, $A_4A_1$) is $\ell'$ (respectively, $r'$). Recall the event $\Eb_{\op b}(v)$ and $\Eb_{\op w}(v)$ for $v\in \ol\D$ defined in \eqref{eq:crossing} and \eqref{eq:crossing3}. \nina{Recall the discrete counterparts $E_\ob(v)$ and $E_\ow(v)$ of these events from Section \ref{sec:crossing-discrete}.}

For $(M,\si)=(M_n,\si_n)$ as in the above coupling, let $a_1,a_{2},a_{3},a_{4}$ be the 1st, $\lceil\ell'h_n\rceil$-th, $\lceil\ell h_n\rceil$-th, and $\lceil(1-r')h_n\rceil$-th outer edges in clockwise order around $M_n$ respectively (starting from the root-edge $\frk e_n=a_1$). Let $A_{2,n}$ (respectively, $A_{4,n}$) be the endpoint of the edge $a_2$ (respectively, $a_4$) in the clockwise arc from $a_1$ to $a_2$ (respectively,  from $a_4$ to $a_1$). 
For a vertex $v\in M_n$, we denote the events $E_{\op b}(v)$ and $E_{\op w}(v)$ by  $E^n_{\op b}(v)$ and $E^n_{\op w}(v)$, respectively, to indicate the dependence on $n$.

\begin{theorem}	\label{thm:finite}
	In the coupling of $(M_n,\sigma_n,\frk v_n)_{n\in \N_+}$ and $(\gffd,\etad, \ub)$ described above, the following quantities converge jointly in probability as $n\rta\infty$.
	\begin{compactitem}
		\item[(i)] Area and boundary length measure: the measures $\mu_n$ and $\nu_n$ converge in the weak topology to the $\sqrt{8/3}$-LQG area and boundary length measures $\mub_{\gffd}$ and $\nub_{\gffd}$, respectively. Here we view all the measures (also the boundary measures) as measures on $\ol\D$ (rather than $\partial\D$).
		\item[(ii)] Space-filling percolation exploration: $\eta_n$ converges in the uniform topology to the space-filling SLE$_6$ $\etad$.
		
		\item[(iii)] Percolation cycles: the embedded percolation cycles $\ga_1^n,\ga_2^n,\dots$ converge to the CLE$_6$ loops $\gab_1,\gab_2,\dots$ as parametrized curves in $\D$.
		\item[(iv)] Pivotal measures: for any fixed $\ep>0$ and $i,j\in\N_+$, the pivotal measures $\nu^{\ep,1}_{j,n}$, $\nu^{\ep,2}_{j,n}$, $\nu_{i,j,n}$, and $\nu_\eps^n$ converge in the weak topology to $\nub^{\ep,1}_{j}$, $\nub^{\ep,2}_{j}$, $\nub_{i,j}$, and $\nub_\eps$. \nina{Furthermore, $\1_{\nu^{\ep,1}_{j,n}(\C)=0}$ converges to $\1_{\nub^{\ep,1}_{j}(\C)=0}$ and the analogous statement holds for the other three measures.}
	\end{compactitem}
	For any $k,k'\in\N_+$, let $e_1,\dots,e_{k}$ and $e'_1,\dots,e'_{k'}$ be $k+k'$ edges sampled independently at random from $E(M_n)$, such that the first $k$ (respectively, last $k'$) edges are sampled uniformly at random from the set of inner (respectively, outer) edges, respectively. 
	Then we may extend the coupling above so that the tuple $P_n := (\phi_n(e_1), \dots, \phi_n(e_{k}), \phi_n(e'_1), \dots, \phi_n(e'_{k'}))$ 
	converges in probability to a tuple $P$ which consists of $k+k'$ points on $\ol\D$ sampled independently ($k$ points from $\mub_{\gffd}$ and $k'$ points from $\nub_{\gffd}$).
	\begin{compactitem}
		\item[(v)] 
		DFS tree finite marginals: $\tau^*_n|_{P_n}$ converges to $\taud|_{P}$. Furthermore, for any fixed 
		$\ell \in (0,1)$, suppose that either $z=\etab(\ell \frk m)\in\D$, or $z\in \bdy\D$ is such that the $\nub_{\gffd}$-length of the counterclockwise arc from $1$ to $z$ is $\ell$. Then the curve $\wh\eta^z_n$ converges uniformly to $\etadh^z$. 
	\end{compactitem}
	\begin{compactitem}
		\item[(vi)] Crossing events: 
		The events $E^n_{\op b}(A_{4,n})$, $E^n_{\op b}(\frk v_n)$, $E^n_{\op w}(A_{2,n})$, $E^n_{\op w}(\frk v_n)$ converge to 
		$\Eb_{\op b}(A_4)$, $\Eb_{\op b}(0)$, $\Eb_{\op w}(A_{2})$, $\Eb_{\op w}(0)$ respectively.
	\end{compactitem}
\end{theorem}
Another natural random triangulation model is the critical Boltzmann disk.
\begin{definition}\label{def:bol}
	Given an integer $\ell\ge 2$, the law of the \emph{critical Boltzmann triangulation with boundary length $\ell$}
	is a probability measure on near-triangulations with a simple boundary of length $\ell$ such that the probability assigned to each such near-triangulation $\cM$ is proportional to $(2/27)^{\#V(\cM)}$.
\end{definition}
Since Theorem~\ref{thm:finite} holds for all $\frk m>0$, we have the following immediate corollary.
\begin{cor}\label{cor:bol}
	In Theorem~\ref{thm:finite}, if $M_n$ is a critical Boltzmann triangulation with boundary length $h_n+2$, and $\gffd,\etad,\Zd$ are as in Theorem~\ref{thm:mot-disk}, then all the convergence statements in Theorem~\ref{thm:finite} still hold as convergence in law. 
\end{cor}
\begin{remark}\label{rmk:CLE}
	As mentioned at the end of Section~\ref{sec:crossing}, one can define other crossing events by rotating $A_1,A_2,A_3,A_4$. However, since the encoding of these other events are more complicated in terms of the random walk/Brownian motion, it is more challenging to establish convergence for these alternatives. For example, if we re-root $M_n$ at $a_{2,n}$ and color the boundary properly, one would obtain a new random walk $\wt Z^n$ with the same law as $Z^n$. However, the joint convergence of $\wt Z^n$ and $Z^n$ to their continuum counterparts is not easy.
	However, in the  work \cite{ghs-metric-peano} by the second and third author and E. Gwynne, it is proved that the percolation cycles on the uniform triangulation with simple boundary converges as curve-decorated metric measure spaces. Moreover the convergence is jointly with the one in Theorem~\ref{thm:finite}. By the re-rooting invariance of $\CLE_6$, the joint convergence of $\wt Z^n$ and $Z^n$ follows. 
	This would give the joint convergence of the crossing events in Theorem~\ref{thm:finite} and their variants after re-rooting.
\end{remark}
We also remark that the convergence in law of crossing events similar to those in Theorem~\ref{thm:finite} were established in \cite{angel-scaling-limit}, where the limiting probabilities were expressed in terms of L\'evy processes instead of CLE$_6$ on $\sqrt{8/3}$-LQG. 
But the convergence of other observables in Theorem~\ref{thm:finite} characterizing the full scaling limit of critical planar percolation was not considered there.

\bigskip

Finally we state the sphere version of our convergence results. In the continuum, recall the setup in Theorem~\ref{thm:mot-sphere},
where $(\C \cup \{\infty\},\gffs)$ is a particular representative of the unit area $\sqrt{8/3}$-LQG sphere; 
$\etas$ is a space-filling $\SLE_6$ parametrized in such a way that $\etas(0)=\etas(1)=\infty$ and $\mub_{\gffs} (\etas([0,t])) =t$ for all $0\le t\le 1$; 
$\Zs$ is the boundary length process of $(\etas,\etas)$ which has the law of a Brownian excursion in the first quadrant with duration 1. 

In the discrete, for $n\in\N_+$ let $M_n$ be a uniformly chosen triangulation of the sphere with $n$ vertices and a (directed) root-edge $\frk e_n$. Let $\sigma_n$ be a coloring of $V(M_n)$ in black and white, where we require that root-edge is oriented from a white to a black vertex, while all the other vertices are colored black or white independently with equal probability.
Recall the bijection $\Phi_0$ in Corollary~\ref{cor:bij-path-returning-0}. Let $Z=\Phi^{-1}_0(M_n,\sigma_n)\in\cK^{(0,0)}$ and $Z^n$ be renormalized as in \eqref{eq:Zn}. Then $Z_n$ converges to $\Zs$ in law. In the rest of this section we consider a coupling of $(M_n,\sigma_n)_{n\in\N_+}$ and $(\gffs,\etas)$ such that $\lim_{n\to\infty}Z^n=\Zs$ in uniform topology almost surely.

Observe that the sphere bijection $\Phi_0$ is a special case of the disk bijection $\bPhi$ (after doing a local modification of the walk and the map, respectively) when we require the boundary length of the disk to be exactly 3. We may therefore define the percolation observables exactly as in the disk case. The main difference is that we choose to embed the limiting $\sqrt{8/3}$-LQG surface into $\C$ instead of $\D$. We also need to enumerate the percolation cycles differently, since there is no well-defined largest CLE$_6$ loop in the scaling limit. Any percolation cycle $\gam$ of $M_n$ divides the vertex set $V(M_n)$ into two disjoint sets $V^{\gam}_1$ and $V^{\gam}_2$ with union $V(M_n)$, and we order the percolation cycles so that $\min\{ \#V^{\gam}_1,\#V^{\gam}_2 \}$ is non-increasing with ties broken in an arbitrary way. 

In the continuum, recall the branching $\SLE_6$ $\taus$ and the $\CLE_6$ $\Gabs$ associated with $(\gffs,\etas)$. For a $k$-tuple $P$, the subtree $\taus|_P$ of $\taus$ is defined in the way as in the case of $\taub^*|_P$ in Section~\ref{sec:cone}.
We enumerate the $\CLE_6$ loops so that $\min\{\area(\gab),1-\area(\gab) \}$ is decreasing. For $j\neq i\in \N_+$, the pivotal measures $\nub^{\eps, 1}_{j},\nub^{\eps, 2}_{j}$ and $\nub_{j,i}$ are defined in the same way as in the cone and disk cases. 
\begin{theorem}\label{thm:sphere}
	In the coupling of $(M_n,\sigma_n)_{n\in \N_+}$ and $(\gffs,\etas)$ described above, the following quantities converge jointly in probability as $n\rta\infty$.
	\begin{compactitem}
		\item[(i)] Area measure: the vertex counting measure $\mu_n$ converges in the weak topology to the $\sqrt{8/3}$-LQG area measure $\mub_{\gffs}$.
		
		\item[(ii)] Space-filling percolation exploration: $\eta_n$ converges in the uniform topology to the space-filling SLE$_6$ $\etas$.
		\item[(iii)] Percolation cycles: the embedded percolation cycles $\ga_1^n,\ga_2^n,\dots$ converge  to the CLE$_6$ loops $\gab_1,\gab_2,\dots$ as parametrized curves in $(\BB S^2,d_{\BB S^2})$.
		\item[(iv)] Pivotal measures: for any fixed $\ep>0$ and $i,j\in\N_+$, the pivotal measures $\nu^{\ep,1}_{j,n}$, $\nu^{\ep,2}_{j,n}$, $\nu_{i,j,n}$, and $\nu_\eps^n$ converge in the weak topology to $\nub^{\ep,1}_{j}$, $\nub^{\ep,2}_{j}$, $\nub_{i,j}$, and $\nub_\eps$. \nina{Furthermore, $\1_{\nu^{\ep,1}_{j,n}(\C)=0}$ converges to $\1_{\nub^{\ep,1}_{j}(\C)=0}$, and the analogous statement holds for the other three measures.}
	\end{compactitem}
	For any $k\in\N_+$, let $e_1,\dots,e_{k}$ be $k$ independently sampled edges in $E(M_n)$.
	Then we may extend the coupling above so that the $k$-tuple $P_n :=(\phi_n(e_1), \dots, \phi_n(e_{k}))$ converges in probability to a tuple $P$ which consists of $k$ points sampled independently from $\mub_{\gffs}$.
	\begin{compactitem}
		\item[(v)] DFS tree finite marginals: $\tau^*_n|_{P_n}$ converges to $\taus|_{P}$. Furthermore, for any fixed $t\in(0,1)$ and $z=\etas(t)$ the curve $\wh\eta^z_n$ converges uniformly to $\wh\etab_{\op S}^z$.
	\end{compactitem}
	
\end{theorem}

\subsection{Towards dynamical percolation and Cardy embedding}\label{sec:flip}

Although our Theorems~\ref{thm1},~\ref{thm:finite}, and~\ref{thm:sphere} imply the convergence of many observables of percolated triangulations to their counterparts in $\sqrt{8/3}$-LQG, a major drawback is that the embedding $\phi_n$ is implicit and depending on more information than the map $M_n$ itself. In light of the conformal invariance of the scaling limit, a more natural embedding for $M_n$ would be discrete approximations to the Riemann mapping, such as circle packing or Tutte embedding with proper boundary conditions. In \cite{hs-quenched}, the second and third named author will introduce such an embedding called \emph{Cardy's embedding}, based on the probability of crossing events discussed in Section~\ref{sec:crossing}. More precisely, in the setting of Theorem~\ref{thm:finite}, for each $v\in M_n$, let $\cardy^n_{\op x}(v)=\P [E^n_{\op b}(v) ]$, where we set $\ell'=1/3$ and $\ell=2/3$ when defining $E^n_{\op b}(v)$.  
Moreover, let $\cardy^n_{\op y}(v)$ and $\cardy^n_{\op z}(v)$ be defined in the same way with $(a_{1},a_{2},a_{3})$ permuted to $(a_{2},a_{3},a_{1})$ and $(a_{3},a_{1},a_{2})$, respectively. Our \emph{Cardy embedding} is defined to be 
\[
\cardy^n(v)=(\cardy^n_{\op x}(v), \cardy^n_{\op y}(v),\cardy^n_{\op z}(v))\in [0,1]^3
\qquad \forall\; v\in V(M_n). 
\]
In \cite{smirnov-cardy}, it is proved that conditioning on $M_n$ being the triangular lattice restricted to a simply connected domain $D$ whose boundary is a continuous curve, $\cardy^n$ converges to the Riemann mapping from $D$ to the simplex $\Delta=\{(x,y,z): x+y+z=1\;\textrm{and }x,y,z\ge 0 \}$.
It will be proved in \cite{hs-quenched} that as $n\to\infty$, all the results in Theorem~\ref{thm:finite} hold if the embedding $\phi_n$ is replaced by $\cardy^n$. The only difference is that now the unit length $\sqrt{8/3}$-LQG disk is embedded in $\Delta$ instead of $\D$. Moreover, the graph distance of $M_n$ under the embedding $\cardy^n$ also has a scaling limit, which is the metric defined by Miller and Sheffield in \cite{lqg-tbm1,lqg-tbm2,lqg-tbm3} and is isometric to the Gromov Hausdorff limit of $M_n$ called the Brownian disk \cite{bet-mier-disk,aasw-type2}. 

The approach taken in \cite{hs-quenched} is the so-called \emph{dynamical percolation}, where one starts from a percolation but each vertex on $V(M_n)$ updates its color independently according to an exponential clock. Some elementary ergodic theory considerations imply that the convergence of the 
Cardy embedded triangulation to the $\sqrt{8/3}$-LQG disk (embedded in $\Delta$)
would follow from the mixing property of the dynamical percolation at the correct time scale.\footnote{One needs to be careful about the topology under which the dynamics lives in when invoking the ergodic theory. The 
	one based on the embedding $\phi_n$ in Theorem~\ref{thm:finite} would not be sufficient given the strong dependence of $\phi_n$ on $\sigma_n$. 
	By incorporating the metric structure of $M_n$, the work \cite{ghs-metric-peano} will allow us to work under a more amenable topology called the \emph{Gromov-Hausdorff-Prokhorov-Uniform} topology.} 
Reviewing the full landscape of this program would be too ambitious here. In this section, we make an important step (namely, we obtain the joint convergence of the percolation loop ensembles before and after the color-flip of a significant pivotal point, and simultaneously the joint convergence of the pivotal measures before and after the color-flip) and explain how it fits into the general program.

Recall from Section~\ref{sec:cont-piv} that, given an instance $\Gab$ of CLE$_6$ on $\D$ and a pivotal point $z$, we obtain a new loop configuration $\Gab_z$ by flipping the color of $z$ (the definition in Section~\ref{sec:cont-piv} were for the CLE$_6$ on $\C$, but the same definition works for CLE$_6$ on $\D$). 
Define the convergence of a loop configuration as in assertion (iii) of Theorem \ref{thm:finite},
that is, as the uniform convergence of the loops viewed as parametrized curves in $\C$.  
Define pivotal points and $\eps$-significance for $\Gab_z$ in the same way as for $\Gab$ (Section~\ref{sec:cont-piv}). Then the set of pivotal points for $\Gab$ and $\Gab_z$ is the same, but the set of $\eps$-significant pivotal points is different. By the first of these observations, for each $\eps>0$ we may define a measure $\wh\nub_{\op D,\eps}$ supported on the set of $\eps$-pivotals for $\Gab_z$, such that $\wh\nub_{\op D,\eps}$ agrees with the measure $\nub_{\op D,\eps'}$ for arbitrary $\eps'>0$ on the intersection of their supports.

We adopt similar definitions in the discrete. Let $(M_n,\sigma_n)\in\smTP$ be as in Theorem~\ref{thm:finite}. Let $\Ga^n$ be the collection of percolation cycles of $(M_n,\si_n)\in\smTP$, 
let $\Ga^n_{z_n}$ be the collection after flipping an $\eps n$-significant pivotal point $z_n$, and
let $\wh\nu_{\op D,\eps}^n$ denote counting measure on the $\eps n$-pivotal points of $\Ga^n_{z_n}$ such that each pivotal has mass $n^{-1/4}$.

\begin{proposition}\label{prop23}
	Recall the setting of Theorem~\ref{thm:finite}. For $\eps>0$, let $z_n\in\D$ be sampled uniformly at random from the set of $\eps n$-pivotal points of $(M_n,\sigma_n)$. Let $z\in\D$ be sampled from $\nub_{\op D,\eps}$ renormalized to be a probability measure. 
	We may extend the coupling of Theorem~\ref{thm:finite} so that $z_n\rta z$ in probability. Letting $\Ga^n_{z_n},\Gab_z,\wh\nu^n_{\op D,\eps},\wh\nub_{\op D,\eps}$ be as above, it holds that  $\Ga^n_{z_n}$ (resp.\ $\wh\nu^n_{\op D,\eps}$) converges in probability to $\Gab_z$ (resp.\ $\wh\nub_{\op D,\eps}$).
\end{proposition}

\begin{remark}\label{rmk:LDP}  
	It is well known that the scaling limit of the dynamical percolation on the triangular lattice is governed by the color updates on significant pivotal points \cite{gps-near-crit}. Proposition~\ref{prop23} can be used to show that the variant of dynamical percolation on planar triangulations where only $\eps n$-significant pivotal points are allowed to update has a scaling limit. In \cite{hs-quenched}, it will be proved that as $\eps\to 0$, this $\eps$-variant of the scaling limit converges to the so-called \emph{Liouville dynamical percolation} introduced in \cite{ghss18}, whose mixing property is also established in the same article by considering quad crossing events. The mixing property of the process is then used to deduce convergence of the triangulation to $\sqrt{8/3}$-LQG under the Cardy embedding. Theorem \ref{thm:finite}(vi) (plus the variants of this result with the vertices permuted) guarantees that we have convergence of the observables used to define the Cardy embedding, which is also essential to guarantee convergence of the embedded triangulation.
\end{remark}

\section{Proofs of the scaling limit results}\label{app:conv}

This section contains the proofs of the scaling results stated in Section \ref{sec:conv}. In Section \ref{sec:convinf} we prove Theorem \ref{thm1}. In Section \ref{sec:convfin} we deduce Theorems \ref{thm:finite} and \ref{thm:sphere} from this result. In Section \ref{subsec:flip} we prove Proposition \ref{prop23}.
We use the following notation throughout this section.
\begin{notation}
	If $A$ and $B$ are two quantities which depend on a parameter $s$ we write $A=o_s(B)$ if $A/B\rta 0$ as (depending on the context) $s\rta 0$ or $s\rta\infty$.
\end{notation}
\begin{notation}
	If $(X_n)_{n\in\N}$ are random variables (in a topological space) we write $X_n\Rightarrow X$ to indicate that $(X_n)_{n\in\N}$ converges in law to $X$ as $n\rta \infty$.
\end{notation}

\subsection{Infinite volume case}
\label{sec:convinf}

Throughout this section we work in the setting of Theorem~\ref{thm1}. In particular, $Z^n$ is the renormalized random walk in \eqref{eq:Zn}, coupled with $(\gff,\etab,\Zb)$ as in Theorem~\ref{thm:mot} so that $Z^n$ almost surely converges uniformly to $\Zb$ on compact sets. This section is devoted to proving Theorem~\ref{thm1}.

\subsubsection{Random walk}\label{sec:walk}
We start by introducing the discrete analog of the L\'evy process relative to time $u$ for $u\in \R$ (see Section~\ref{sec:dictionary-fl}). Let $\beta>0$ be as in Remark \ref{rmk:variance}.
Recall from \eqref{eq:Zn} that 
$$Z\equiv (Z_k)_{k\in \Z}=(2\sqrt{n/\beta}\,Z^n_{k/(3n)})_{k\in \Z} $$
is the unscaled walk which is associated with the word $w\in\{a,b,c\}^\ZZ$, and that the steps of $w$ are chosen uniformly and independently at random. 
Here the dependence of $Z$ and $w$ on $n$ is dropped for simplicity. 
For each $u\in \R$, we write $w=w^{-,u} w^{+,u}$ where the first step of 
$w^{+,u}$ is $w_{\lfloor3nu\rfloor}$. Recall Definition~\ref{def:pi-w-inf}. 
We denote by $\dots,T^u(-2),T^u(-1),T^u(0)\in\ZZ^{<0}$ the times associated with the spine steps of $w^{-,u}$ so that $\pi(w^{-,u})=\dots w_{T^u(-2)}w_{T^u(-1)}w_{T^u(0)}$. 
For $m\in\ZZ$, let $\ell^u_m:=-\#\{k\in\Z^{\leq 0}~|~T^u(k)\geq m\}$. Finally, let $\wh Z^u_m=Z_{T^u(m)}$. 

\begin{figure}
	\includegraphics[scale=1]{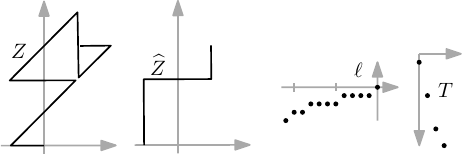}
	\centering
	\caption{\nina{Illustration of $Z$, $\wh Z$, $\ell$, and $T$ when $w^{-,u}=\dots acbbccaacca$ and $\pi(w^{-,u})=\dots abaa$.}}
	\label{fig:Zhat-ell-T}
\end{figure}

For all $t\in\R$ and $s\leq 0$, let
\begin{equation}
	\label{eq:Levy-walk}
	\wh Z^{u,n}_s:=(\wh L^{u,n}_s, \wh R^{u,n}_s)=\frac 12 \sqrt{\beta/n} \,\wh Z^u_{\lfloor sn^{3/4}\rfloor},
	\quad 
	T^{u,n}_s=(3n)^{-1}T^u(\lfloor sn^{3/4}\rfloor),
	\quad 
	\ell^{u,n}_t=n^{-3/4}\ell^u_{\lfloor 3nt\rfloor }.
\end{equation}
Then $\ell^{u,n}$, $T^{u,n}$, and $\wh Z^{u,n}$ are the discrete analogs of $\ellb^u,\Tb^u$, and $\wh \Zb^u$, respectively, defined in Sections~\ref{sec:dictionary-branch} and~\ref{sec:dictionary-fl}. \nina{See Figure \ref{fig:Zhat-ell-T} for an illustration.}
The main result in this subsection is the joint convergence of the triple $(\ell^{u,n},T^{u,n},\wh Z^{u,n})$ to their continuum counterparts.
\begin{lemma}\label{prop5} 
	In the setting of Theorem~\ref{thm1} with the above defined notions, for any fixed $u\in \R$, 
	$(\wh Z^{u,n}, T^{u,n} ,\ell^{u,n})$ converges in probability to $(\wh\Zb^u,\Tb^u,\ellb^u)$ on any fixed compact set, where the first two coordinates are in the Skorohod topology and the third coordinate is in the uniform topology.
\end{lemma}
To prove Lemma~\ref{prop5}, we consider the time-reversal of the walk $Z$. 
Define $\Wrev=(\cLrev_k,\cRrev_k)_{k\in\N}$ by $\Wrev_k=Z_{-k}$ for $k\in\N$. Then $\Wrev$ is a two-dimensional random walk with steps independently and uniformly distributed on $\{(0,-1),(-1,0),(1,1)\}$. Let $\Trev(0)=0$, and for $m\geq 1$ define $\Trev(m)=-T(-m+1)$ with $T$ as in Definition \ref{def:pi-w-inf}. Using the definition of $T$, we get 
\eqbn 
\Trev(m) = \inf\big\{ t> \Trev(m-1)\,|\, 
\cLrev_t \leq \cLrev_{\Trev(m-1)} 
\text{\,\,or\,\,}
\cRrev_t \leq \cRrev_{\Trev(m-1)} 
\big\}.
\eqen
Define $\Whrev=(\Whrev_m)_{m\in\N}$ by $\Whrev_m=\Wrev_{\Trev(m)}$ for $m\in\N$. 
\begin{lemma}\label{prop4} 
	The process $(\Trev(m))_{m\in\N}$ is a random walk with independent steps, whose step distribution is supported on $\N_+$ and satisfies the following asymptotics:
	\begin{equation}
		P(\Trev(1)=k)= \frac{1}{\sqrt{2}\,\Gamma(-3/4)}k^{-7/4}\,(1+O(1/\sqrt{k})).\label{eq:T-asympt}
	\end{equation}
	
	The process $(\Whrev_m)_{m\in\N}$ is a random walk with independent steps in $\bigcup_{k\in{-1}\cup \N_+}\{(k,0),(0,k)\}$. 
	The distribution of the steps are given by $\P[\Whrev_1=(k,0)]=\P[\Whrev_1=(0,k)]=\frac{1}{2}\P(J=k)$, where $J$ is a random variable satisfying $\P[J=-1] =2/3$ and for all $k>0$,
	\begin{align}
		P(J=k)=\frac{1}{4^{k}(2k-1)(k+1)}{2k\choose k}~=~ \frac{1}{2\sqrt{\pi}}k^{-5/2}\,(1+O(1/k))\label{eq:J-asympt}.
	\end{align}
\end{lemma} 
\begin{proof}
	We only prove \eqref{eq:T-asympt} and \eqref{eq:J-asympt} as other statements are easy observations. Let $\mK_0$ (respectively, $\mK_1$) be the set of finite walks starting at the origin, staying (respectively, staying strictly) in the first quadrant, and ending on the $y=0$ line. Observe that there is a bijection between $\mK_0$ and $\mK_1$: Given a walk in $\mK_0$ we obtain the associated walk in $\mK_1$ by adding $(1,1)$ in the beginning and $(0,-1)$ at the end. Note that the size of the excursion increases by two, and the height of the excursion increases by one when we apply the bijection. Also observe that the excursion we consider in the lemma is contained in $\mK_1$. In our analysis it will be more convenient to study walks in $\mK_0$. Therefore we define $J_0=J-1$ and $T_0=\Trev(1)-2$.
	Let 
	$$K(x,z)=\sum_{P\in \mK_0}x^{i(P)}z^{|P|},$$ 
	where $|P|$ is the number of steps and $(i(P),0)$ is the ending point of $P$.
	The series $K=K(x,z)$ is known to be algebraic \cite{Kreweras:walks} over $\QQ(x,z)$. The following equation is readily obtainable from \cite[Theorem 1]{MBM:Kreweras}:
	\begin{eqnarray*}
		&&16\,{K}^{6}{x}^{6}{z}^{10}-48\,{x}^{4}{z}^{8} ( x-2\,z ) {K}^{5}
		+8\,{x}^{2}{z}^{6} ( 6\,{x}^{3}{z}^{2}+7\,{x}^{2}-24\,xz+24\,{z}^{2} ) {K}^{4}\\
		&&-32\,{z}^{4} ( x-2\,z ) ( 3\,{x}^{3}{z}^{2}+{x}^{2}-2\,xz+2\,{z}^{2} ) {K}^{3}\\
		&&+{z}^{2} ( 48\,{x}^{4}{z}^{4}+64\,{x}^{3}{z}^{2}-264\,{x}^{2}{z}^{3}+192\,x{z}^{4}+9\,{x}^{2}-32\,xz+32\,{z}^{2} ) {K}^{2}\\
		&&- ( x-2\,z ) ( 48\,{x}^{2}{z}^{4}+16\,x{z}^{2}-72\,{z}^{3}+1 ) K+16\,{x}^{3}{z}^{4}+8\,{x}^{2}{z}^{2}-72\,x{z}^{3}+108\,{z}^{4}+x-2\,z=0
	\end{eqnarray*}
	
	Next we will prove that
	\eqb
	P(J_0=k)=\frac{2}{9}[x^k]K(x,1/3),\qquad\text{and}\qquad P(T_0=k)=\frac{2}{9}[z^k]K(1,z/3).
	\label{eq82}
	\eqe 
	There are constants $c_1,c_2>0$ such that
	\eqb
	K(x,1/3) = c_1 \sum_{k=0}^\infty \P[J_0=k]x^k,\qquad
	K(1,z/3) = c_2 \sum_{k=0}^\infty \P[T_0=k]z^k.
	\label{eq45}
	\eqe 
	In order to prove \eqref{eq82} it is sufficient to show that $c_1=c_2=9/2$. The probabilities in each sum in \eqref{eq45} sum to $1/3$ (since with probability 2/3 there is no excursion, due to the walk starting with $(0,-1)$ or $(-1,0)$), and therefore we have $K(1,1/3)=c_1/3=c_2/3$. We have $\P[T_0=0]=\P[T=2]=2/9$, since $T=2$ if the first two steps are $(1,1)$,$(0,-1)$ or $(1,1)$,$(-1,0)$. The number of walks in $\mK_1$ of duration 2 is equal to 1 (since this happens if the walk starts with $(1,1)$,$(0,-1)$), which implies $1=c_2\P[T=2]$. Combining the above we get $1=c_2\P[T_0=2]=2c_2/9$, which gives $c_1=c_2=9/2$ as desired, and \eqref{eq82} follows.
	
	Let 
	$$J_0\equiv J_0(x)=\sum_{k=0}^\infty P(J_0=k)x^k:=\frac{2}{9}K(x,1/3)$$ 
	and 
	$$T_0\equiv T_0(z)=\sum_{k=0}^\infty P(T_0=k)z^n:=\frac{2}{9}K(1,z/3).$$
	By specializing the above equation for $K(x,z)$, we get
	\eqbn
	\frac{1}{324}(9J_0^2x^2-18J_0x+12J_0+4x+24)(9J_0^2x^2-18J_0x+12J_0+4x-3)^2=0,
	\eqen
	and 
	\begin{eqnarray}\label{eq:T}
		T_0^{6}{z}^{10}+2\,{z}^{8} \left( 2\,z-3 \right) T_0^{5}+2/3\,{z}^{6}
		\left( 10\,{z}^{2}-24\,z+21 \right) T_0^{4}\nonumber\\
		+{\frac {16\,{z}^{4}
				\left( 2\,z-3 \right) \left( 5\,{z}^{2}-6\,z+9 \right) T_0^{3}}{27}}
		+1/27\,{z}^{2} \left( 80\,{z}^{4}-264\,{z}^{3}+288\,{z}^{2}-288\,z+243
		\right) T_0^{2}\\
		+{\frac { \left( 4\,z-6 \right) \left( 16\,{z}^{4}-72
				\,{z}^{3}+48\,{z}^{2}+27 \right) T_0}{81}}-{\frac {32\,{z}^{3}}{27}}+4/9
		+{\frac {496\,{z}^{4}}{729}}+{\frac {32\,{z}^{2}}{81}}-{\frac {8\,z}{
				27}}=0.\nonumber
	\end{eqnarray}
	The coefficients of $J_0(x)$ must be positive. Therefore the former equation gives
	\begin{equation}\label{eq:J}
		9	{x}^{2}{J_0}^{2}+ 6 \left(2 -3\,x\right) J_0+4\,x-3=0,
	\end{equation}
	since we can exclude the first factor as it would lead to $J_0= -2-(10/3)x-8x^2+....$, which has negative coefficients. From \eqref{eq:J}, we get
	$$J_0(x)={\frac { 2\left( 1-x \right) ^{3/2}}{3{x}^{2}}}-\frac{2}{3x^2}+\frac{1}{x}.$$ 
	Applying the binomial theorem and Stirling's formula gives \eqref{eq:J-asympt}.
	In order to obtain \eqref{eq:T-asympt} from \eqref{eq:T} we can apply the techniques of \cite[Chapter VII.7]{Flajolet:analytic}. More precisely, we first check (by computing the discriminant of \eqref{eq:T}) that $z=1$ is the unique dominant singularity, and then obtain the asymptotic behavior 
	$$T_0(z)=_{z\to 1} \frac{1}{3}-\frac{2\sqrt{2}}{3}(1-z)^{3/4}+\frac{4}{3}\,(1-z)+O((1-z)^{5/4}),$$
	from which \eqref{eq:T-asympt} follows.
\end{proof}

\begin{remark}\label{rmk:constant}
	Due to the explicit asymptotics of the tail of $T^{\op{rev}}$ in \eqref{eq:T-asympt} and the fact that $\Tb^u$ is a stable subordinator, 
	$T^{u,n}$ converges in law to a constant $c$ times $\Tb^u$ for each fixed $u\in \R$. Recall that $\ellb^u$ is only defined up to a multiplicative constant. We fix that multiplicative constant by requiring that $c=1$.
\end{remark}

\begin{proof}[Proof of Lemma~\ref{prop5}]
	Without loss of generality we assume $u=0$. To simplify the notation, we drop the dependence on $u=0$ and write $\wh Z^{0,n}$, $T^{0,n}$, $\ell^{0,n} $, $\wh\Zb^0$, $\Tb^0$, $\ellb^0$ as $\wh Z^{n}$, $T^{n}$, $\ell^{n} $, $\wh\Zb$, $\Tb$, $\ellb$, respectively.
	
	We first show that $(Z_t^n)_{t\in\R}$ and $(T_s^n)_{s\leq 0}$ converge jointly in law to $(\Zb_t)_{t\in\R}$ and $(\Tb_s)_{s\leq 0}$. 
	The convergence of the two marginal laws is immediate.
	We need to show that there is a unique joint subsequential limit which is as desired. By Skorokhod embedding, we may assume that $(Z^n, T^n)$ converges almost surely along a subsequence to a limit $(\Zb,\Tb')$ with $\Tb'\eqD \Tb$. To conclude we need to show that $(\Tb_s)_{s\leq 0}=(\Tb_s')_{s\leq 0}$ almost surely, where $\Tb$ is determined from $\Zb$ as in Section \ref{sec:cont-cle}. 
	Let 
	\[
	a^n_{-}=\sup\{s<-1:s\;\textrm{is in the image of}\; T^n \} \quad \textrm{and}\quad a^n_+=\inf\{s>-1:s \;\textrm{is in the image of}\; T^n \}.
	\]
	Let $\ab'_-$ and $\ab'_+$ be the limit of $a^n_-$ and $a^n_+$, respectively, and let 
	\[
	\ab_-=\sup\{s<-1:s\in\ans(0)\} \qquad \textrm{and}\qquad \ab_+=\inf\{s>-1:s\in\ans(0)\}.
	\]	
	In the discrete we have $L^n_t\geq L^n_{a^n_+}$ and $R^n_t\geq R^n_{a^n_+}$ for all $t\in(a^n_-,a^n_+)$. By convergence of $(L^n,R^n)$ to $(\Lb,\Rb)$, this gives 
	$\Lb_t\geq \Lb_{\ab'_+}$ and $\Rb_t\geq \Rb_{\ab'_+}$ for all $t\in(\ab_-,\ab_+)$, so $\ab'_+$ is an ancestor of all points in $(\ab'_-,\ab'_+)$. Furthermore, $\Zb|_{[ \ab'_-,\ab'_+ ]}$ is a cone excursion. 
	Therefore $( \ab'_-,\ab'_+ )\cap \ans(0)=\emptyset$, so $[ \ab'_-,\ab'_+ ]\subseteq [ \ab_-,\ab_+ ]$ and further $\ab'_+\le \ab_+$. Furthermore, since both $\Zb|_{[ \ab'_-,\ab'_+ ]}$ and $\Zb|_{[ \ab_-,\ab_+ ]}$ are cone excursions, we must have $\ab_-\le \ab'_-$. On the other hand, $\ab_-\eqD\ab'_-$ and $\ab'_+ \eqD\ab_+$, so we must have $\ab_-=\ab'_-$ and $\ab'_+= \ab_+$ almost surely. The same argument is true if we replace the time $-1$ in the definition of $a_-^n$, $a_+^n$, $\ab_-$, and $\ab_+$ by any other negative rational. This implies that the image of $\Tb'$ is equal to $\ans(0)$, which implies further (since $\Tb\eqD\Tb'$ and these two processes are stable subordinators) that $\Tb'=\Tb$.
	
	Observe that $( Z^n, T^n )\Rightarrow (\Zb,\Tb)$ implies $( Z^n,\wh Z^n, T^n )\Rightarrow (\Zb,\wh\Zb,\Tb)$, since $\wh Z^n$ is tight by Lemma \ref{prop4}, and since if $( Z^n, T^n )\rightarrow (\Zb,\Tb)$ almost surely, then almost surely for any fixed $u\geq 0$, 
	$$
	\wh Z^n_s
	=\frac 12 \sqrt{\beta/n} \wh Z_{\lfloor un^{3/4}\rfloor}
	=\frac 12 \sqrt{\beta/n} Z_{T(\lfloor un^{3/4}\rfloor)}
	\rta 
	\Zb_{\Tb_{u}}
	= \wh\Zb_s
	.
	$$
	
	Since $T^n\Rightarrow \Tb$, it is immediate from the definition of $T^n,\Tb,\ell^n,\ellb$ that $(\ell^n, T^n )\Rightarrow (\ellb,\Tb)$. 
	Recall that for sequences of random variable $(x_n)$, $(y_n)$, and $(z_n)$, if $(x_n,y_n)$ converges jointly to $(x,y)$, $(x_n,z_n)$ converges jointly to $(x,z)$, and $x_n$ determines $z_n$ and $x$ determines $z$, then $(x_n,y_n,z_n)$ converges jointly to $(x,y,z)$. Hence, combining the above results with the fact that $T^n$ (respectively, $\Tb$) determines $\ell^n$ (respectively, $\ellb$), we conclude the proof.
\end{proof}

\subsubsection{Envelope intervals}\label{subsub:env}
In this section, we extend the coupling in Theorem~\ref{thm1} and establish the convergence of envelope intervals. 
Before stating this result, we need to define a notion of convergence for intervals.
\begin{notation}
	When we discuss convergence of intervals in the remainder of this section we mean convergence for the following pseudometric $d$ on the set of intervals. For closed intervals $I=[a_1,a_2]$ and $J=[b_1,b_2]$ define $d_{\Io}(I,J):=|a_1-b_1|+|a_2-b_2|$. The metric is defined in the same way if one or both of the intervals is open or half-open. 
	\label{not1}
\end{notation}
For a bounded interval $I$, recall from Section \ref{sec:cont-cle} that $\li(I)$ is the smallest envelope interval containing $I$. For fixed $I\subset\R$ it is almost surely the case that $\li(I)$ is also the smallest cone interval containing $I$. 

In the discrete, we say that an interval $[s,t]$ is a \emph{cone interval} for $w$ if $ns,nt\in\ZZ$ and 
the steps $ns+1$ and $nt$ of the walk $w$ form a close-matching (in other words, the subwalk of $Z$ between time $ns$ and $nt$ is a cone excursion). For an interval $I$ we let $\ci_n(I)$ be the smallest cone interval of $w$ containing $I$.

\begin{lemma}\label{lem:env-int}
	The coupling in Theorem~\ref{thm1} can be extended in such a way that the convergence of $(\wh Z^{u,n}, T^{u,n} ,\ell^{u,n})$ in Lemma~\ref{prop5} holds almost surely for all $u\in \Q$. 
	Under such a coupling, for any fixed bounded interval $I\subset\R$, we have $\ci_n(I)\rta\li(I)$ almost surely.
\end{lemma}
\begin{proof} 
	The existence of a coupling such that the convergence of $(\wh Z^{u,n}, T^{u,n} ,\ell^{u,n})$ in Lemma~\ref{prop5} holds almost surely for all $u\in \Q$ is immediate by a diagonal argument. 
	
	We now consider a bounded interval $I$. 
	Observe that for $q\in\R$, $I$ and the image of $\Tb^q$ are disjoint if and only if there exists a cone excursion ending at a time $< q$ which contains $I$. By varying $q\in\Q$ 
	and recalling that $\li(I)$ is the smallest cone interval containing $I$, it follows that we have convergence of the right endpoint of $\ci_n(I)$ to the right endpoint of $\li(I)$. Letting $t_*$ denote the left endpoint of $\li(I)$ it holds that $-t_*$ is a stopping time for the time reversal $\Zb^{\op{rev}}=(\Zb_{-t})_{t\in\R}$ of $\Zb$.
	Therefore both coordinates of $(\Zb^{\op{rev}}_{t-t_*}-\Zb^{\op{rev}}_{-t_*})_{t\in\R}$ take both positive and negative values immediately after $t=0$. Since $Z^n\rta \Zb$ uniformly on compact intervals, this implies that the left endpoint of $\ci_n(I)$ converges to the left endpoint of $\li(I)$. This concludes the proof.
\end{proof}

We now briefly recall the description of discrete envelope excursions in terms of the walk $w$. Given $e\in E(M_n)$, let $k_2=\etae^{-1}(e)$ and let $w_{k_1}$ be the far-match of $w_{k_2}$. Recall the definition of $\treepar$ from Section~\ref{sec:cluster-tree-inf}. By Corollary \ref{cor9}, $e$ is an envelope edge and $w_{k_1}\dots w_{k_2}$ is an envelope excursion if and only if $\treepar(w_{k_2})$ has a different color than $w_{k_2}$. In this case we call $[(3n)^{-1}(k_1-1),(3n)^{-1}k_2]$ an \emph{envelope interval} of $w$.
By Claim \ref{claim:disjoint-intervals}, for any fixed interval $I$ there exist a \emph{smallest} envelope interval containing $I$, and we denote this interval by  $\li_n(I)$.
Note that in contrast to the continuum, $\ci_n(I)\neq \li_n(I)$.

The following lemma gives convergence of envelope intervals. Since the discrete description of envelope intervals in terms of the walk is \emph{not} a direct analog of the continuum description, we need to use probabilistic techniques to study the behavior of the walk near a \emph{typical} envelope interval. \nina{Recall that we did a rescaling by $0.5\sqrt{\beta/n}$ when defining $Z^n=(L^n,R^n)$ in \eqref{eq:Zn}, and therefore $L^n_{s_n^c}-L^n_{t_n^c}=0.5\sqrt{\beta/n}$ (respectively, $R^n_{s_n^c}-R^n_{t_n^c}=0.5\sqrt{\beta/n}$) corresponds to a delta of one for the unscaled process.}
\begin{figure}
	\centering
	\includegraphics[scale=1]{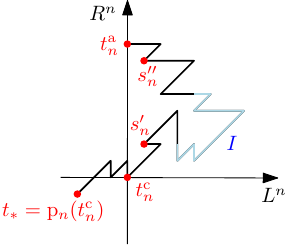}
	\caption{Illustration of the proof of Lemma \ref{prop2}. The label next to a point indicates the time at which the point is visited by $Z^n$.}
	\label{fig:env}
\end{figure}

\begin{lemma} 
	Let $I\subset\R$ be a fixed interval, and consider a coupling as in Lemma \ref{lem:env-int}.
	With probability converging to 1 as $n\rta\infty$, it holds that $\li_n(I)\subset \ci_n(I)$. Furthermore, $\li_n(I)$ converges to $\li(I)$ in probability. 
	
	Let $t_n^c$ (respectively, $s_n^c$) denote the right end-point of $\ci_n(I)$ (respectively, $\li_n(I)$). With probability converging to 1 as $n\rta\infty$, if $\ci_n(I)$ if a right (respectively, left) cone interval then $L^n_{s_n^c}-L^n_{t_n^c}=0.5\sqrt{\beta/n}$ (respectively, $R^n_{s_n^c}-R^n_{t_n^c}=0.5\sqrt{\beta/n}$).
	\label{prop2}
\end{lemma}
\begin{proof}
	We advise to look at Figure \ref{fig:env} while reading the proof. To simplify notation let $[t^a_n,t^c_n]=\ci_n(I)$, $[s_n,s^c_n]=\li_n(I)$, and $[t^a,t^c]=\li(I)$. Without loss of generality assume that $\li(I)$ is a right cone excursion, that is, $\Lb_{t^a}=\Lb_{t^c}$. By Lemma \ref{lem:env-int}, we may assume throughout the proof that $L^n_{t^a_n}=L^n_{t^c_n}$, since this holds with probability $1-o_n(1)$ by Lemma \ref{prop5}.
	
	First we introduce some notation. Assume $t\in (3n)^{-1}\Z$ is such that $w_{3tn}=c$. We say that $s\in (3n)^{-1}\Z$ is the close-match (respectively, far-match) of $t$ if $w_{3sn+1}$ is the close-match (respectively, far-match) of $w_{3tn}$. Observe that $t^{\op{a}}_n$ is the close-match of $t^{\op{c}}_n$, and that $s_n$ is the far-match of $s^{\op{c}}_n$. We say that $t\in (3n)^{-1}\Z$ is of \emph{$a$-type} (respectively, \emph{$b$-type}) if the far-match of $w_{3tn}$ is an $a$-step (respectively, $b$-step). In this case we say that the $a$-match of $t$ is equal to its far-match (respectively, near-match), while the $b$-match of $t$ is equal to its near-match (respectively, far-match). Define $\matchingpar_n(t)\in (3n)^{-1}\Z$ to be the smallest $t'>t$ such that $w_{3t'n}=c$, and such that the $a$-match (respectively, $b$-match) of $t'$ is smaller than the $b$-match (respectively, $a$-match) of $t$ on the event that $t$ is of type $a$ (respectively, $b$). Notice that $\matchingpar(w_{3tn})=w_{3t'n}$ with $\matchingpar$ as in Section \ref{sec:tree-cluster-walk} if and only if $\matchingpar_n(t)=t'$.
	
	We will prove that the following happens with probability $1-o_n(1)$: 
	\begin{compactitem}
		\item[(i)] $\matchingpar_n(t^c_n)$ is of the same type as $t^c_n$.
		\item[(ii)] Let $$s_n':=\inf\{ t\geq \sup(I)\,:\,L^n_t=L^n_{t_n^c}+\frac{1}{2}\sqrt{\beta/n} \}.$$
		Then $s_n'$ is the right endpoint of an envelope interval $J_n$ satisfying 
		$I\subset J_n\subset\ci_n(I)$, $J_n\rta\li(I)$ as $n\rta\infty$, and $\matchingpar_n(s'_n)=t_n^c$.
		\item[(iii)] The envelope interval $J_n$ in (ii) equals $\li_n(I)$.
	\end{compactitem}
	First we argue that (i-iii) imply the lemma. By (ii) and (iii), we see that $\li_n(I)=J_n\subset\ci_n(I)$. 
	Combining (ii) and (iii), we have that $\li_n(I)=J_n\rta\li(I)$. The last assertion of the lemma is satisfied by the definition of $s_n'$.
	
	Proof of (i): Conditioned on $\ci_n(I)$ and $Z^n|_{\ci_n(I)}$, the increments of $Z^n|_{\R\setminus \ci_n(I)}$ are uniform and independent. 
	Let $t_*:=\inf\{ t>t_n^c\,:\,R_t^n<R_{t_n^c}^n \}$.
	Then $\matchingpar_n(t_n^c)=t_*$ by definition of $\matchingpar_n$, since the $b$-match of $t_*$ is smaller than the $b$-match of $t^{\op{c}}_n$. With probability $1-o_n(1)$ we have 
	$t_*\in(t_n^c,t_n^c+n^{-0.99})$, $|L^n_{t_*}-L^n_{t_n^c}|<n^{-0.49}$, and $|R^n_{t_*}-R^n_{t_n^c}|=n^{-1/2}/2<n^{-0.49}$. By uniform convergence of $(L^n,R^n)$ to $(\Lb,\Rb)$ on compact sets, with probability $1-o_n(1)$, we have $\inf_{t\in(t_n^a-n^{-0.97},t_n^a)} (L_t^n-L^n_{t_n^c})<-n^{-0.49}$ and $\inf_{t\in(t_n^a-n^{-0.97},t_n^a)} (R_t^n-R^n_{t_n^c})>n^{-0.49}$, and the far-match of $t_n^c$ is smaller than $t_n^a-n^{-0.97}$. Since $\matchingpar_n(t_n^c)=t_*$ it follows that with probability $1-o_n(1)$, $\matchingpar_n(t_n^c)$ is of the same type as $t^c_n$.
	
	Proof of (ii): Let $s''_n=\sup\{ t<s'_n\,:\,L^n_t=L^n_{s_n'} \}$. By uniform convergence of $(L^n,R^n)$ to $(\Lb,\Rb)$ on compact sets, there exists a random constant $c>0$ independent of $n$, such that with probability $1-o_n(1)$, $L_t^n>L_{t_n^c}^n+c$ and $R_t^n>R_{t_n^c}^n+c$ for all $t\in I$. Therefore, by the intermediate value theorem applied to $L^n$, we have $s_n'\in(\sup(I),t_n^c)$ and $s''_n\in(t_n^a,\inf(I))$. Furthermore, we have $s_n''- t_n^a\to 0$ and $s'_n- t_n^c\to 0$, which gives that $J_n=[s''_n,s'_n]\rta\li(I)$. By definition of $\ci_n(I)$ as the minimal cone interval containing $I$, we see that $J_n$ is not a cone interval, $s''_n$ is the far-match of $s'_n$, and $s'_n$ is of a different type than $t^c_n$. Moreover, $\matchingpar_n(s'_n)=t^c_n$. Using (i), we conclude that $J_n$ is an envelope interval, which completes the proof of (ii).
	
	Proof of (iii): By (ii) we have $\li_n(I)\subseteq J_n\subset\ci_n(I)$ with probability $1-o_n(1)$. To conclude we need to show that $s^c_n=s'_n$.
	First we will show that 
	\begin{equation}\label{eq:match}
		\matchingpar_n(s^c_n)=t^c_n.
	\end{equation}
	On the one hand, we must have $\matchingpar_n(s^c_n)\geq t^c_n$, otherwise $\matchingpar_n(s^c_n)$ would be the end of a cone excursion contained in $\ci_n(I)$ and containing $I$, which is a contradiction to the definition of $\ci_n(I)$. On the other hand $\matchingpar_n(s^c_n)\leq t^c_n$, by the definition of $\matchingpar_n$, since both the near-match and the far-match of $t^c_n$ are smaller than $s_n$. 
	This gives \eqref{eq:match}. We have $\matchingpar_n(s^c_n)=\matchingpar_n(s_n')$ and $[s_n,s_n^c]\subseteq [s_n'',s_n']$, hence $s_n^c=s_n'$.
\end{proof}
\begin{assumption}
	In the remainder of the section we consider an extension of the coupling in Theorem~\ref{thm1} such that the convergence in Lemma~\ref{prop5} holds almost surely for all $u\in \Q$, and such that the convergence in Lemmas~\ref{lem:env-int} and~\ref{prop2} hold almost surely for all intervals $I$ whose endpoints are both rational.
	\label{a1}
\end{assumption}

In Theorem \ref{thm1} the enumeration of the percolation cycles is based on the number of vertices enclosed by the various cycles. To argue that the $j$th cycle in the discrete converges to the $j$th cycle in the continuum, the following lemma will be useful.
\begin{lemma} 
	Let $I\subset\R$ be a fixed interval, and let $\ga^n$ (respectively, $\gab$) be the percolation cycle with envelope interval $\li_n(I)$ (respectively, $\li(I)$). Then $\area_n(\ga^n)\rta\area(\gab)$ in probability as $n\rta\infty$.
	\label{prop10}
\end{lemma}
\begin{proof}
	To simplify notation, define $J_n=\li_n(I)$ and $J=\li(I)$, and let $u_n\in\R$ (respectively, $u\in\R$) denote the terminal endpoint of $J_n$ (respectively, $J$). Denote the complementary components of $\ans(u)\cap J$ by $J^{{\Lo}}_i$ and $J^{{\Ro}}_i$ for $i\in\N$, such that $\Zb|_{J^{{\Lo}}_i}$ (respectively, $\Zb|_{J^{{\Ro}}_i}$) is a left (respectively, right) cone excursion. By symmetry we can assume that $\gab$ is traced in counterclockwise direction. 
	In this case, the discussion above Definition~\ref{def:CLE} gives
	\[
	\area(\gab)=\sum_{i\in\N}|J^{{\Lo}}_i|,
	\] 
	where we use $|\cdot|$ to denote the length of an interval. 
	
	We will use a similar formula in the discrete setting to conclude. Since $\gab$ is traced in counterclockwise direction, $\Zb|_{J}$ is a right cone excursion. Hence by Lemma \ref{prop2}, with probability $1-o_n(1)$ the percolation cluster for which $\ga^n$ is the outside-cycle is white. We condition on this event in the remainder of the proof. Define $i$ and $k$ such that $\env_n(I)=[i/(3n),k/(3n)]$. By Lemma \ref{lem:envelope-matching-is-loop-building}(iv), $\frk L_\ell(w_i\dots w_{k-1})=\frk L(\ga)$. Furthermore, by Theorem \ref{thm:LR}, each step $x$ of $\hpi(w)$ of the form $\ba_k$ corresponds to both a bicolor triangle $t$ and to a unicolor face $f$ of $\LR(M,\sigma)$. The site-percolated near-triangulation $(M',\sigma')$ of $(M,\sigma)$ formed of $t$ and all the vertices and edges of $(M,\sigma)$ inside and on the boundary of $f$ is equal to $\bPhi(w')$, where $w'$ is the cone excursion of $w$ corresponding to $x$. By Remark \ref{rmk8}, the number of vertices in $(M',\sigma')$ is equal to 2 plus the number of $c$-steps in $w'$. All these vertices are enclosed by $\ga$, except one of the vertices on $t$. Let $J^{{\Lo}}_{i;n}$ denote intervals on the form $[i'/(3n),k'/(3n)]$ for $w'=w_{i'}\dots w_{k'}$ such that $[i',k']\subset[i,k]$, $Z^n|_{J^{{\Lo}}_{i;n}}$ is a left cone excursion, and $[i',k']$ is not contained in any other interval satisfying these properties. Assuming $J^{{\Lo}}_{i;n}$ and $J^{{\Lo}}_{i}$ are ranked by decreasing length, it follows from the scaling limit result for $T^n$ that $d_{\op{I}}(J^{{\Lo}}_{i;n},J^{{\Lo}}_{i} )\rta 0$ in probability for each $i$. Furthermore, since $\#c=(|w'|+h)/3$ where $h=R_{i'-1}-R_{k'}=o(|w'|)$, the number of $c$-steps in the word $w'$ associated with some interval $J^{{\Lo}}_{i;n}$, divided by $n$, converges to $|J^{{\Lo}}_i|$ in probability. Combining the above results, we get that $\area_n(\ga^n)$ converges in probability to $\area(\gab)=\sum_{i\in\N}|J^{{\Lo}}_i|$.
\end{proof}

\subsubsection{Definition of the continuum pivotal measure}\label{sec:piv-def}
Recall that Definitions~\ref{def:pivotal} and~\ref{def:piv-eps} are contingent on fixing unspecified normalizing constants, and on the covering Lemma \ref{lem:cover}. Before heading to technical discussions on pivotal measures, we first address these two issues and complete the definition of pivotal measures for $\CLE_6$ loops. 

Recall that we are still using the notation of Theorem~\ref{thm1} and Assumption \ref{a1}. For $u\in\R$, recall the quantities defined in \eqref{eq:Levy-walk}. In this section as well as in Sections~\ref{subsec:cycles} and~\ref{subsub:pivot}, when we consider $\wh Z^{u,n}$, $\wh L^{u,n}$, $\wh R^{u,n}$, $T^{u,n}$, and $\ell^{u,n}$ relative to a single time $u$, which is clear from the context, we simply write $\wh Z^n$, $\wh L^n$, $\wh R^n$, $T^n$, and $\ell^n$. Similarly, in the continuum we adopt the convention of Section~\ref{sec:dictionary-fl} and use 
$\wh \Zb$, $\wh \Lb$, $\wh \Rb$, $\Tb,\ellb$ to represent quantities relative to $u$. (In fact, the only $u$'s which will be considered are deterministic times and envelope closing times.)

Let us recall and complete the definition of pivotal measure introduced in Section~\ref{sec:cont-piv}. Let $(\wh\Lb_t,\wh\Rb_t)_{t\leq 0}$ denote the L\'evy processes relative to time $u$ (see Section \ref{sec:dictionary-fl}). As in \eqref{eq:foward}, let 
\[
\Ab_{\Lo}(s,u)=\{ t\in(s,0]\,|\,\inf_{t'\in[s,t]} \wh \Lb_{t'}=\wh \Lb_t \}
\]
be the set of forward running infima of $\wh\Lb$ relative to time $s$. As explained in Section~\ref{sec:cont-piv} below \eqref{eq:foward}, one can define a measure $\pb_{\Lo}(s,u)$ supported on $[s,0]$ via the local time at $\Ab_{\Lo}(s,u)$.
Define $\Ab_{\Ro}(s,u)$ and $\pb_{\Ro}(s,u)$ similarly and set 
\begin{equation}\label{eq:cont-piv-measure}
	\Ab(s,u)=\Ab_{\Lo}(s,u) \cup \Ab_{\Ro}(s,u) \qquad \textrm{and}\qquad \pb(s,u)=\pb_{\Lo}(s,u)+\pb_{\Ro}(s,u).
\end{equation} 
Recall from Section \ref{sec:dictionary} that the local time on the set of running infima of a $3/2$-stable L\'evy process with only negative jumps is only defined up to a multiplicative constant. Lemma \ref{lem:psu} will fix the convention for this constant in our paper in order to match the definitions from the discrete.

Let us now define the corresponding quantities in the discrete.
For $s<0$, let $A^n_{\Lo}(s,u)$ be the set of $t\in [s,0]$ such that $\wh L^n$ has a forward running infimum at time $t$ relative to time $s$, that is,
$$A^n_{\Lo}(s,u)=\{t\in(s,0] ~|~\forall t'\in[s,t), L^n_{t'}>L^n_{t}\}.$$
Roughly speaking, the set $A^n_{\Lo}(s,u)$ should be thought of as the set of times that correspond to cut-points of the white forested line $\frk L_\ell(w^-)$ that separate the vertex of the percolation path of $(M_n^-,\si_n^-)$ corresponding to time $s$ from the infinite line of $\frk L_\ell(w^-)$ (see Figure \ref{fig:LR-inf}).
We define $A^n_{\Ro}(s,u)$ similarly from  $\wh R^n$, and let $A^n(s,u)=A^n_{\Lo}(s,u)\cup A^n_{\Ro}(s,u)$.  

Then let $p_{\Lo}^n(s,u)$ (respectively, $p_{\Ro}^n(s,u)$) denote the measure on $(-\infty,0]$ given by counting measure on $A^n_{\Lo}(s,u)$ (respectively, $A^n_{\Ro}(s,u)$), where each point is given mass $n^{-1/4}$. In other words, an interval $I\subset(-\infty,0]$ is assigned mass the following mass for $p_{\Lo}^n(s,u)$
\begin{equation}\label{eq:running-inf} 
	n^{-1/4}\cdot\# (A^n_{\Lo}(s,u)\cap I),
\end{equation}
and a similar relation holds for $p_{\Ro}^n(s,u)$. 
Note that $p_{\Lo}^n(s,u)$ and $p_{\Ro}^n(s,u)$ are supported on $[s,0]$. 
Let 
\begin{equation}\label{eq:dis-piv-mea}
	p^n(s,u)=p^n_{\Lo}(s,u)+p^n_{\Ro}(s,u).
\end{equation}

\begin{lemma} \label{lem:psu} 
	In the above setting, it is possible to choose a normalizing constant when defining the local time on the set of strict running infima of a $3/2$-stable L\'evy process with only negative jumps such that
	\eqb
	\E[ \pb_{\op L}(-1,0) ] = \lim_{n\rta\infty} \E[ p^n_{\op L}(-1,0) ]. 
	\label{eq83}
	\eqe 
	Using this normalization, for any fixed deterministic $u\in\R$ and $s<0$, it holds that $p_{{\Lo}}^n(s,u)$ (respectively, $p_{{\Ro}}^n(s,u)$) converges to $\pb_{{\Lo}}(s,u)$ (respectively, $\pb_{{\Ro}}(s,u)$) as $n\rta\infty$ in probability.
\end{lemma}
\begin{proof}
	This is immediate by Lemma~\ref{prop4} and \cite[Theorem 2]{cd10}. 
\end{proof}

\begin{proof}[Proof of Lemma~\ref{lem:cover}]
	Let $\{B_{i}\}_{i\in\N_+}$ be an enumeration of all the bubbles in $\frk L_{\Xb}$.
	Let $t_i$ be the jump time of $\Xb$ corresponding to $B_i$.
	Let $N$ be a random integer such that 
	$\sum_{i=N+1}^{\infty}\mub_{\gff}(B_i)<\eps$. Consider an $\eps$-pivotal point $p\in \dbl_{\Xb,\eps}$. 
	Since $p$ is pivotal, there exists $\frk s_p <\frk t_p$ such that $p=\pi_{\Xb}(\frk s_p)=\pi_{\Xb}(\frk t_p)$,
	where  $\pi_{\Xb}$ is the quotient map as in Section~\ref{sec:dictionary-fl}. 
	Since $p$ is $\eps$-pivotal, there exists $i\le N$ such that $\frk s_p<t_i\le \frk t_p$.
	Moreover, if $t_i\neq \frk t_p$ then $p\in \dbl_{\Xb,_\eps}(t_i)$, while if $t_i=t_p$ then $p=\pi_{\Xb}(t_i)$.
	For all $i\in [N]$, one can choose $q_i\in \Q$ such that $\pi_{\Xb}(q_i)\in B_i$ and $\Xb_{q_i}> \Xb_{t_i}$. With this choice we get \(\dbl_{\Xb,\eps} \subseteq \cup^{N}_{i=1} \dbl_{\Xb} (q_i).\)
	This concludes the proof. \qedhere
\end{proof}

\subsubsection{Percolation cycles}
\label{subsec:cycles}

Combining Lemma \ref{prop:piv} and~\ref{lem:psu} give some notion of convergence for pivotal measure. However, in our proof of Theorem \ref{thm1} we need to obtain a version of Lemma~\ref{lem:psu} where $u$ is replaced by an envelope closing time. Controlling the behavior of the walk near an envelope closing time is challenging since the walk is conditioned on a complicated random walk event. We resolve this by studying critical percolation on critical Boltzmann triangulations (see Definition~\ref{def:bol}). 

The following lemma shows that critical Boltzmann triangulations arise naturally in the site-percolated UIPT. The continuum analog is that each connected component of the future wedge is an independent $\sqrt{8/3}$-LQG disk (see Section~\ref{sec:mot}).
\begin{lemma}
	Let $w\in\imK$ be chosen according to the uniform distribution. Write $w=w^-w^+$ and $w^+=w^+_1cw^+_2cw^+_3c\dots$ as in Section \ref{sec:bij-inf}. The percolated maps $\ol\Phi(w^+_1),\ol\Phi(w^+_2),\dots$ in $\smTP$ (which are used to define the future percolated near-triangulation) are independent and identically distributed. Conditioned on their boundary lengths, they have the law of critical Boltzmann triangulations decorated with a uniformly sampled percolation satisfying the root-interface condition (see Section \ref{subsec:finite} for the definition).
	\label{prop33}
\end{lemma}
\begin{proof}
	Since the $c$-steps in the decomposition $w^+=w^+_1cw^+_2cw^+_3c\dots$ correspond to simultaneous running infima for both coordinates, these $c$-steps are stopping times. Using this and that the steps of $w$ are independent and identically distributed, we get that $w^+_1,w^+_2,\dots$ are independent and identically distributed.
	This implies that the associated maps $\ol\Phi(w^+_1),\ol\Phi(w^+_2),\dots$ are also independent and identically distributed. Let $(M_1,\sigma_1)=\bPhi(w^+_1)$. Let $(\cM_0,s_0)$ be an arbitrary fixed percolated map in $\smTP$, let $v_0=\bPhi^{-1}(\cM_0,s_0)\in \smK$, and let $n$ (respectively, $m$) be the number of inner (respectively, outer) vertices of $\cM_0$. By Corollary \ref{cor:chordal-case}, 
	\eqb
	\P[ (M_1,\sigma_1)=(\cM_0,s_0) ] 
	= \P[ w_1\dots w_{|v_0|+1}=v_0c ]
	= 3^{-|v_0|-1} 
	= 3^{-(3n+2m-4)-1}. 
	\label{eq41}
	\eqe
	For fixed $m$ there are $m-1$ ways to choose $s_0$ restricted to the boundary of $\cM_0$, and there are $2^{-n}$ ways to choose $s_0$ restricted to the interior of $\cM_0$. Therefore, for fixed $m$, $\P[ (M_1,\sigma_1)=(\cM_0,s_0) ]$ is proportional to $2^{-n}\P[M_1=\cM_0]$. Using \eqref{eq41} it follows that, for fixed $m$, $\P[M_1=\cM_0]$ is proportional to $(2/27)^{n}$, which implies that $M_1$ has the law of a critical Boltzmann triangulation. By \eqref{eq41}, all the allowable percolation configurations occur with the same probability.
\end{proof}

Our next lemma will rule out certain pathological behaviors of the percolation cycle near the envelope closing time. Fix $\xi>0$.
Let $M_n'$ be a critical Boltzmann disk with boundary length $\lceil\xi n^{1/2}\rceil\geq 2$, and let $\sigma_n'$ be a uniformly sampled percolation satisfying the root-interface condition as in Lemma~\ref{prop33}. 
Denote the root-edge of $M_n'$ by $e_1$ and the edge at the other end of the percolation path of $\sigma_n'$ by $e_2$. Let $f_1\in F(M'_n)$ (respectively, $f_2\in F(M'_n)$) be the unique inner face of $M'_n$ which is adjacent to $e_1$ (respectively, $e_2$). Note that $(M'_n,\sigma'_n)\in\smTP$. By Corollary~\ref{cor:chordal-case}, we let $w^n=\ol\Phi^{-1}(M'_n,\sigma'_n)\in \smK$ and $\tau^*_n=\Delta_{M'_n}(\sigma'_n)$. Consider the branch in the DFS tree $\tau^*_n$ from $f_1^*$ to $f_2^*$, let $0=T(0)<\dots<T(k')\in\N$ be the set of times such that this branch consist of the faces (equivalently, vertices of the dual map $(M'_n)^*$) $\etavf(T(j))$, and let $t_1=(3n)^{-1}T(k')$. Define $\ell_m:=\min\{k\geq 0: T(k)\ge m \}$ for all $m\in\Z$. With $Z$ denoting the walk on $\Z^2$ associated with $w^n$ and
$\wh Z_u=(\wh L_u,\wh R_u)=Z_{T(u)}$ for $u\in\{0,\dots,k' \}$, define
\eqbn
Z^n_u = \frac 12 \sqrt{\beta/n} Z_{\lfloor 3nu\rfloor },\quad
\wh Z^n_u=\frac 12 \sqrt{\beta/n} \wh Z_{\lfloor un^{3/4}\rfloor},
\quad 
\ell^n_t=n^{-3/4}\ell_{\lfloor 3nt\rfloor }, 
\quad 
T^n_u=\frac{1}{3}n^{-1}T(\lfloor un^{3/4}\rfloor).
\eqen
Then define the measures $p^n_{{\Lo}}$, and $p^n_{{\Ro}}$ on $\R$ as in \eqref{eq:running-inf} and let $p^n=p^n_{{\Lo}}+p^n_{{\Ro}}$.
In other words, $p^n_{{\Lo}}$ (respectively, $p^n_{{\Ro}}$ ) counts the renormalized number of strict running infima for $\wh L^n$ (respectively, $\wh R^n$) in any given interval, such that each strict running infimum has mass $n^{-1/4}$.
\begin{lemma}
	In the setting above, for each $\eps>0$ there exists $\delta>0$ depending only on $\eps$ and $\xi$ but not on $n$ 
	such that
	\eqb
	\P\Big[\ell^n_{\delta}<\eps;\,
	p^n([0,\ell^n_{\delta}])<\eps;\,
	|\ell^n_{t_1}-\ell^n_{t_1-\delta}|<\eps;\,
	p^n([\ell^n_{t_1}-\ell^n_{t_1-\delta},\ell^n_{t_1}])<\eps
	\Big]>1-\eps.
	\label{eq33}
	\eqe
	\label{prop22}
\end{lemma}
We prove Lemma~\ref{prop22} by comparing the Boltzmann disk with the \emph{uniform infinite half-planar triangulation} (UIHPT), which is the weak local limit of a critical Boltzmann triangulation rooted at an outer edge \cite{angel-scaling-limit}. 
\begin{proof}[Proof of Lemma~\ref{prop22}]
	For $\delta>0$ define the following events $A_0=A_0(\delta)$ and $A'_0=A'_0(\delta)$.
	$$
	A_0(\delta)=\{ \ell^n_{\delta}<\eps\} \cap
	\{p^n_{{\Lo}}([0,\ell^n_{\delta}])<\eps \},
	\quad
	A'_0(\delta)=\{ |\ell^n_{t_1}-\ell^n_{t_1-\delta}|<\eps\} \cap \{p^n_{{\Lo}}([\ell^n_{t_1}-\ell^n_{t_1-\delta},\ell^n_{t_1}])<\eps \}.
	$$
	Assume that the numbers of left edges and the number of right edges of $M_n'$ are both bigger than $\lfloor 2\delta n^{1/2}/5\rfloor$. Note that this assumption is satisfied with probability $1-o_\delta(1)$.
	
	For $\delta,\delta',\delta''>0$, let $A(\delta,\delta',\delta'')$ denote the following event
	\eqb
	A(\delta,\delta',\delta'') = \{ T^n_{\delta'}>\delta \} \cap 
	\{ p^n_{{\Lo}}([0,\delta''])<\eps \} \cap 
	\Big\{ \inf_{t\in[0,\delta'']} \wh L^n_t < \inf_{t\in[0,\delta']} \wh L^n_t \Big\} \cap
	\Big\{ \inf_{t\in[0,\delta'']} \wh R^n_t < \inf_{t\in[0,\delta']} \wh R^n_t \Big\}.
	\label{eq101}
	\eqe
	First we will argue that for any $\delta',\delta''\in(0,\eps)$,
	\eqb
	\P[ A_0(\delta) ] \wedge \P[ A_0'(\delta) ] \geq \P[ A(\delta,\delta',\delta'') ].
	\label{eq96}
	\eqe
	It is immediate that $\P[ A_0(\delta) ] \geq \P[ A(\delta,\delta',\delta'') ]$ since $A(\delta,\delta',\delta'')\subset A_0(\delta)$. To prove that $\P[ A'_0(\delta) ] \geq \P[ A(\delta,\delta',\delta'') ]$ recall the space-filling exploration of the edges of $(M'_n,\sigma'_n)$ described in Section \ref{sec:exploration-tree-from-walk}. Consider two paths from $f_2$ to $f_1$: the time reversal of the space-filling path from $e_1$ to $e_2$, and the space-filling path from $e_2$ to $e_1$. The latter path is defined in exactly the same way as the former path, but for the map rooted at $e_2$ instead of $e_1$. The percolation interface between black and white is the same for the two paths, but for the former (respectively, latter) path the complementary components of the percolation path are visited after (respectively, before) their boundary is traced by the percolation interface. Due to the symmetry between $e_1$ and $e_2$, the latter process is equal in law to the forward exploration. Let $\check T^n$ and $\check\ell^n$ be defined just as $T^n$ and $\ell^n$, respectively, but for the percolated map rerooted at $e_2$. Observe that if $\check T^n_{\delta'}>\delta$ then $\{\check\ell^n_\delta<\eps \}$ and further $\{ |\ell^n_{t_1}-\ell^n_{t_1-\delta}|<\eps\}$. Recall the function $\wh\eta_{\op{v}}$ defined in Section \ref{subsec:pivot}, and assume that in our setting this function is parametrized by the integers between 1 and the length of the percolation interface. Observe that on the event $A(\delta,\delta',\delta'')$, by the requirement on $\wh L^n$ and $\wh R^n$ in \eqref{eq101}, and by the interpretation of $\wh L^n$ and $\wh R^n$ as boundary length processes, all pivotal points which were visited for the first time by the percolation interface before time $\delta'$ (in the sense that the pivotal point is given by $\wh\eta_{\op{v}}(i)$ for some $i<n^{3/4}\delta'$) are visited for the second time by the percolation interface before time $\delta''$ (in the sense that the pivotal point is given by $\wh\eta_{\op{v}}(i)$ for some $i<n^{3/4}\delta''$). If the pivotal point is visited for a second time at time $n^{-3/4}i$, then the measure $p^n$ will have a point mass at $n^{-3/4}i$. Therefore the event $A'_0(\delta)$ occurs if the event $A(\delta,\delta',\delta'')$ occurs for the percolated map rerooted at $e_2$. We conclude that \eqref{eq96} holds.
	
	By \eqref{eq96} and symmetry between $p_{\Lo}$ or $p_{\Ro}$, in order to conclude the proof it is sufficient to show that $\P[A(\delta,\delta',\delta'')]>1-\eps/4$.

	Consider an instance of the uniform infinite half-planar triangulation (UIHPT) $\wt M$ with root-edge $\wt e_1$ on the boundary. Then consider the following \emph{peeling process} of $\wt M$ (see also Section \ref{sec:LR-decomposition-walk}). The triangles of $\wt M$ are explored one by one, starting with the triangle which is adjacent to $\wt e_1$. The triangle we peel in step $u+1$ shares an edge with the triangle peeled in step $u$, and is connected to infinity by a path of unexplored triangles, where adjacent triangles on the path are required to share an edge. If there are two possible such triangles, we choose one uniformly at random. Observe that the sequence of triangles revealed in this peeling process have the same law as the triangles on the chordal percolation interface, assuming the left (respectively, right) frontier is white (respectively, black). Let $(\check L_t)_{t\in\N}$ (respectively, $(\check R_t)_{t\in\N}$) be the process describing how the length of the left (respectively, right) boundary evolves for this peeling process, and let $\check L^n,\check R^n$ denote the renormalized versions, that is, $\check L^n_t=\frac 12 \sqrt{\beta/n} \check L_{\lfloor tn^{3/4}\rfloor}$
	and
	$\check R^n_t=\frac 12\sqrt{\beta/n} \check R_{\lfloor tn^{3/4}\rfloor}$.
	
	The peeling process defines a percolation interface which separates triangulated disks with monocolored boundary from infinity. We will now relate the peeling process to the DFS defined in Section \ref{sec:exploration-tree-from-walk}, and explain how we can obtain a \emph{space-filling} exploration of $E(\wt M)$ by proceeding as in that section. 
	Recall from Section \ref{sec:exploration-tree-from-walk} that for a percolated map $(M,\sigma)\in\bmTP$ the ordering of the edges as defined by the word $\ol\Phi^{-1}(M,\sigma)\in\bmK$ via $\etae$ can be obtained by a DFS of the dual map $M^*$ for $M$. Since the UIHPT is the local limit of the critical Boltzmann triangulation, we may define a space-filling exploration of $(\wt M,\wt\sigma)$ by following the rules of Definition \ref{def:space-filling}.
	By Definition \ref{def:space-filling}(ii'), we see that this exploration is defined by following the percolation interface between white (left) and black (right), that is, by exploring the map via peeling as described above, except that every time we separate a finite submap $(\check M,\check{\sigma})$ from infinity by visiting a face $f$ for which all three incident vertices are on the boundary of the explored region, we do a space-filling exploration of $(\check M,\check{\sigma})$ before continuing to explore the unbounded unexplored part of $\wt M$. The space-filling exploration of $(\check M,\check{\sigma})$ is done by following the rules of Definition \ref{def:space-filling}.
	
	Define $\wt T$ and $\wt p$ similarly as in the exploration of the disk, that is, for $k\in\N$, $\wt T(k)$ is the number of edges explored by the space-filling path when we peel the $k$th face in the peeling process, $\wt p_L$ (respectively, $\wt p_L$) is counting measure on the running infima of $\check L$ (respectively, $\check R$), and $\wt p=\wt p_L+\wt p_R$. Let $\wt T^n$ and $\wt p^n$ denote the renormalized processes, and let $\wt A=\wt A(\delta,\delta',\delta'')$ denote the following event
	\eqbn
	\wt A=
	\wt A(\delta,\delta',\delta'')=
	\wt A_1(\delta,\delta') \cap \wt A_2(\delta'') \cap \wt A_3(\delta'') \cap \wt A_4(\delta',\delta''),
	\eqen
	where
	\eqbn
	\begin{split}
		\wt A_1(\delta,\delta') &= \bigg\{ \sum_{t\in[0,\delta']\,:\,\wt L^n_{t}\neq\wt L^n_{t^-} } (\wt T^n_{t}-\wt T^n_{t^-})>\delta \bigg\},\qquad
		\wt A_2(\delta'') = \{ \wt p^n_{{\Lo}}([0,\delta''])<\eps \},\\
		\wt A_3(\delta'') &= \Big\{ \sup_{t\in[0,\delta'']}|\wt L^n_t|\wedge|\wt R^n_t|<\eps^2/100 \Big\},\\
		\wt A_4(\delta',\delta'')&=
		\Big\{ \inf_{t\in[0,\delta'']} \check L^n_t < \inf_{t\in[0,\delta']} \check L^n_t \Big\} \cap
		\Big \{ \inf_{t\in[0,\delta'']} \check R^n_t < \inf_{t\in[0,\delta']} \check R^n_t \Big\}. 
	\end{split}
	\eqen
	
	Let $\wt e_3$ be the outer edge of $\wt M$ such that the left endpoint of $\wt e_3$ is $\lceil \xi n^{1/2}/5\rceil$ vertices to the left of the left endpoint of $\wt e_1$ along the boundary of $\wt M$. Consider the peeling process of $\wt M$ as described above starting from $\wt e_3$. Let $f$ be the first face we explore, and let $u$ be the vertex incident to $f$ which is not an endpoint of $\wt e_3$. Let $\wt B$ be the event that $u$ is an outer vertex of $\wt M$ which is exactly $\lceil\xi n^{1/2} \rceil-1$ vertices to the right of the right endpoint of $\wt e_3$, that is, when this face is peeled we enclose a disk of boundary length $\lceil\xi n^{1/2} \rceil$ with $\wt e_1$ on its boundary. On the event $\wt B$, the enclosed disk has the law of a critical Boltzmann triangulation \cite{angel-curien-uihpq-perc}. By the definition of $\wt A$ and the locality of percolation, this implies that $\P[A]\geq\P[\wt A|\wt B]$. Using this and Bayes' rule,
	\eqbn
	\P[A]\geq\P[\wt A|\wt B]=\frac{\P[\wt B|\wt A] \P[\wt A]}{ \P[\wt B] }.
	\eqen
	We will show that
	\[
	\textrm{(i) $\P[\wt A]>1-\eps/8\qquad$ and \qquad (ii) $\P[\wt B|\wt A]/\P[\wt B]>1-\eps/8$},
	\]
	which is sufficient to conclude the proof of the lemma.
	
	(i) By for example \cite[Section 3]{angel-peeling}, $(\wt L^n,\wt R^n)$ is a random walk with increments that are independent and identically distributed such that $\wt L^n\eqD \wt R^n$ and $\P[\wt L^n_1>a]\sim ca^{-5/2}$ for $a>1$ and some constant $c>0$. By \cite{cd10}, counting measure on the running infima of a random walk converging to a L\'evy process converges in law to the local time on the running infima of the limiting L\'evy process (if 0 is regular for the limiting process, which is the case for us). This result implies that $\P[ \wt A_2(\delta'') ]>1-\eps/50$ for sufficiently small $\delta''$. Using that $\wt L^n$ and $\wt R^n$ converge to L\'evy processes in the scaling limit, we also get that
	$\P[ \wt A_4(\delta',\delta'') ]>1-\eps/50$, 
	$\P[ \wt A_1(\delta,\delta') ]>1-\eps/50$, and
	$\P[ \wt A_3(\delta'') ]>1-\eps/50$ by first choosing $\delta''$ sufficiently small, then choosing $\delta'$ sufficiently small compared to $\delta''$, and finally choosing $\delta$ sufficiently small compared to $\delta'$. We conclude by a union bound that $\P[A(\delta,\delta',\delta'')]>1-\eps/8$.
	
	(ii) Define $k':=\lceil \xi n^{1/2}\rceil+\wt L^n_{\delta'}+\wt R^n_{\delta'}$, and observe that on the event $\wt A(\delta,\delta')$ we have $|k'-\lceil \xi n^{1/2}\rceil|<\eps^2/50$. 
	By for example \cite[Section 3]{angel-peeling}, $\P[\wt B]=c\lceil \xi n^{1/2}\rceil^{-5/2}(1+o_n(1))$ for some constant $c>0$. By the Markov property of the peeling process, $\P[ \wt B\,|\,\sigma(\wt A,k')]=c(k')^{-5/2}(1+o_{k'}(1))$ on the event that $\wt A=\wt A(\delta,\delta')$ occurs. It follows that on the event $\wt A$,
	\eqbn
	\frac{\P[\wt B\,|\,\wt A,k']}{ \P[\wt B] }
	=
	\frac{ \lceil \xi n^{1/2}\rceil^{-5/2}(1+o_n(1)) }{ (k')^{-5/2}(1+o_{k'}(1))}.
	\eqen
	On the event $\wt A$ and for large $n$, the ratio on the right side is larger than $1-\eps/8$.
\end{proof}

Recall the looptree perspective on CLE$_6$ described in Section~\ref{sec:cont-cle}. There are two looptrees associated with a $\CLE_6$ loop $\gab$: the \emph{left} and \emph{right} looptrees. The looptree $\frk L(\gab)$ is mutually absolutely continuous with respect to $\frk L_{\Xb}$ for $\Xb$ a $3/2$-stable L\'evy excursion with only negative jumps.
The other looptree is mutually absolutely continuous with respect to a part of $\frk L_{\Xb}$ (more precisely, its encoding process may be described as $\Xb$ restricted to the latter part of its domain of definition). 
In particular, we have the following properties:
\begin{itemize}
	\item[(1)]	
	almost surely $\pb(s,u)$ does not have an atom at 0;
	\item[(2)] 
	almost surely the quantum natural parametrization of $\gab$
	does not accumulate mass around $\etab(u)$ in the sense that if $\Tb$ is the L\'evy process relative to the envelope closing time of $\gab$ then $\Tb_{-\delta}\rta 0$ as $\delta\rta 0$.
\end{itemize}

\begin{figure}
	\centering
	\includegraphics[scale=0.9]{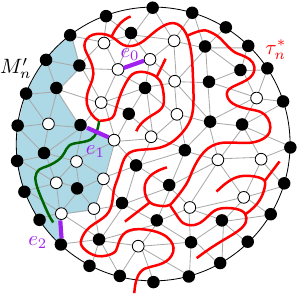}
	\caption{Illustration of the proof of Lemma \ref{prop21}. 
		The triangulation made of the blue faces is $M'_n$, and the red tree is the part of $\tau^*_n$ which is explored before the time $t_2^{n,0}$ at which $e_1$ is treated. Observe that the green chordal path in $M'_n$ is the last segment of the percolation cycle around the large white percolation cluster. The edge $e_0=\etae(0)$ is the root-edge. }
	\label{fig2}
\end{figure}

Using Lemmas \ref{prop33} and~\ref{prop22} along with properties (1) and (2) above, we obtain convergence of quantum natural time and pivotal measure when exploring the map toward an envelope closing edge.
\begin{lemma}\label{prop21}
	Fix $\eps>0$, and let $u$ (respectively, $u_n$) be an envelope closing time for a CLE$_6$ loop $\gab\in\Gab$ (respectively, percolation cycle $\ga^n$ of $(M_n,\si_n)$) such that the envelope interval of $\ga^n$ converges in probability to the envelop interval of $\gab$. Let $\ellb,\Tb$ (respectively, $\ell^n, T^n$) be defined as above, but with the Brownian motion $\Zb$ (respectively, the walk $Z^n$) recentered at $u$ (respectively, $u_n$). Then $\ell^n$ (respectively, $T^n$) converges in probability to $\ellb$ (respectively, $\Tb$) in the uniform topology (respectively, Skorokhod topology). Furthermore, for any fixed $s<0$, the measure $p_n(s,u_n)$ converges in probability to $\pb(s,u)$.	
\end{lemma}
\begin{proof} 
	Without loss of generality, we assume that $0\in \reg(\gab)$ and $\env(\gab)$ is a left cone interval. Define $t_2=u$ and let $t_1$ be such that $[t_1,t_2]=\env(\gab)$. Recall from Section~\ref{sec:cont-cle} the decomposition of $\gab$ into two segments $\gab_1$ and $\gab_2$ in the past and future wedge relative to $0$. See in particular Figure \ref{fig:cle-2parts}. Let $[t_1^0,t_2^0]$ be the maximal cone interval inside $\env(\gab)$ containing $0$. Then, $\etab([t_1^0,t_2^0])$ is the closure of the connected component of $\etab(\env(\gab))\setminus \gab$ containing $0$. Recall that $\gab_2$ is the segment of $\gab$ restricted to $\etab([t_2^0,t_2])$. Let $D$ be the interior of $\etab([t_2^0,t_2])$.
	Let $x=\etab(t_2^0)$ and $y=\etab(t_2)$. \nina{By Section~\ref{sec:cont-cle}, conditioning on $(D,x,y)$ the curve $\gab_2$ is a chordal $\SLE_6$ on $(D,x,y)$. Furthermore, since the space-filling SLE$_6$ can be obtained from the regular SLE$_6$ by filling in each complementary component of the curve by an independent space-filling SLE$_6$ loop, it also holds that $\etab|_{[t_2^0,t_2]}$ is a chordal space-filling $\SLE_6$ on the interior of $(D,x,y)$.} 
	
	Let $\{q_i\}_{i\in \N_+}$ be an increasing rational sequence greater than $t_2^0$ tending to $t_2=u$ with the following property. For each $i\in \N_+$, $\etab(q_i)$ is contained in a connected component $B_i$ of $D\setminus \gab_2$ whose boundary has nonempty intersection with $\bdy \etab(\env(\gab))$. Let $[\frk s_i,\frk t_i]$ be such that $\etab([\frk s_i,\frk t_i])$ is the closure of $B_i$.
	\nina{Observe that
		$\lim_{i\to\infty} \frk s_i=t_2$ since $q_i\in(\frk s_i,\frk t_i)$, $q_i\uparrow t_2$, and $(\frk s_i,\frk t_i)\in\frk I$ for $\frk I$ a countable collection of disjoint intervals $(\frk s,\frk t)\subset(t_2^0,t_2)$ for which $\frk t\neq t_2$. By continuity of $\ellb$, this gives further that $\ellb_{\frk s_i}\rta\ellb_{t_2}=0$.}
	Let $\ellb^i$ be the local time of $\ans(q_i)$. Then for each $t<\frk s_i$, we have that $\ellb_{\frk s_i}-\ellb_t= \ellb^i_{\frk s_i}-\ellb^i_t $. Hence, since $\ellb_{\frk s_i}\rta\ellb_{t_2}=0$,
	\begin{equation}\label{eq:ell}
		\textrm{$\ellb^i-\ellb^i_{\frk s_i}$ converges in probability to $\ellb$},
	\end{equation}
	uniformly on compact sets.
	Similarly, let $\pb^i$ be the Borel measure supported on $[s,\ellb_{\frk s_i}]$ defined as follows. Given any interval $I$, $\pb^i(I) $ equals the $\pb(s-\ellb_{\frk s_i} +\ellb^i_{\frk s_i} , q_i )$-mass of the shifted interval $I-\ellb_{\frk s_i}+\ellb^i_{\frk s_i}$.
	Then 
	\begin{equation}\label{eq:pb}
		\textrm{$\pb^i$ converge in probability to $\pb(s,u)$.}
	\end{equation} 
	
	In the discrete, for $q_1$ defined as above, let $t_2^n$ (respectively, $t_2^{n,0}$) be such that $3nt_2^n $ (respectively, $3nt_2^{n,0}$) is the index of the first (respectively, last) $c$-step after (respectively, before) $3nq_1$ in $w$ that does not have an $a$-match nor a $b$-match in $w^+$, so in particular $3nt_2^{n,0}$ and $3nt_2^n$ are two consecutive cut-times for $w^+$ (note that $t_2^{n}$ and $t_2^{n,0}$ exist with probability $1-o_n(1)$).  
	Let $t^n_1$ (respectively, $t_1^{n,0}$) be such that $3nt_1^{n}$ (respectively, $3nt_1^{n,0}$) is the index of the near match of $w_{3nt_2^{n}}$ (respectively, $w_{3nt_2^{n,0}}$). Equivalently, $[t^n_1-(3n)^{-1},t^n_2]$ 
	is the smallest cone interval which contains both 0 and $q_1$, while $[t^{n,0}_1-(3n)^{-1}, t^{n,0}_2]$  is the largest cone interval containing 0 but not $q_1$. By Lemma \ref{lem:env-int} and Assumption \ref{a1}, we have $[t_1^{n},t_2^n]\rta[t_1,t_2]$ and $[t_1^{n,0},t_2^{n,0}]\rta[t^0_1,t^0_2]$ in probability, which gives that $[t_2^{n,0},t_2^n]\rta [t_2^0, t_2]$ and $[t_1^n, t_1^{n,0}]\rta [t_1, t_1^0]$ in probability. Moreover, suppose $I$ is an interval with rational endpoints such that $\env(I)=[t_1,t_2]$. Then with probability $1-o_n(1)$ we have $\ci_n(I) = [t_1^n,t_2^n]$. Given $i\in \N_+$, let $\frk t_i^n$ be such that $3n\frk t_i^n$ is the index of the first $c$-step after $3nq_i$ without a $b$-match in $w|_{[t_2^{n,0},t_2^n ]}$. Let $\frk s_i^n$ be such that $3n\frk s_i^n$ is the index of the near match of $w_{3n\frk t_i^n}$. Then $[\frk s_i^n, \frk t_i^n]$ converges to $[\frk s_i,\frk t_i]$ in probability. 
	
	By Lemma~\ref{prop33}, $w|_{[t_2^{n,0},t_2^n-(3n)^{-1}]}$ gives a critical Boltzmann disk $M_n'$ decorated with a uniformly sampled percolation $\sigma'_n$ satisfying the root-interface condition (see Figure \ref{fig2}). Moreover, the total number of outer edges in $M_n'$ divided by $\sqrt n$ converges in probability as $n\to \infty$. Since $\lim_{n\to\infty} \frk s_i^n=\frk s_i$ in probability and $\lim_{i\to\infty} s_{i} =t_2$ for each $i\in \N$, combining Lemmas~\ref{prop5}, ~\ref{lem:psu}, and~\ref{prop22} and Assertions~\eqref{eq:ell} and~\eqref{eq:pb}, we see that 
	Lemma~\ref{prop21} holds with $u_n$ replaced by the right endpoint of $\ci_n(I)$, which is $t_2^n$. 
	
	We will now argue that the lemma still holds with $u_n$ instead of $t_2^n$. Since $[t_1^n-(3n)^{-1},t_2^n]=\li_n([0,q_1])$ and $u_n$ is the terminal endpoint of $\env_n([0,q_1])$, Lemma \ref{prop2} gives that $u_n<t_2^n$ with probability converging to 1 as $n\rta\infty$, and that $|u_n-t_2^n|$ converges to 0 as $n\rta\infty$. Consider the process $T^{n,u_n}=T^n$ (respectively, $T^{n,t_2^n}$) relative to $u_n$ (respectively, $t_2^n$), which defines a parametrization of the times $t$ at which $\eta_{\op{e}}(3nt)$ is on the path of the exploration tree toward $\eta_{\op{e}}(3n u_n)$ (respectively, $\eta_{\op{e}}(3n t_2^n)$). We first claim that the range of $T^{n,u_n}$ is contained in the range of $T^{n,t^n_2}$ with probability $1-o_n(1)$. The claim follows from the last assertion of Lemma \ref{prop2}, since this implies that with probability $1-o_n(1)$ there is no cone excursion ending between $u_n$ and $t_2^n$ and starting before $u_n$, so any letter which is enclosed by a near-matching before the $3nt_2^n$th letter is also enclosed by a near-matching before the $3nu_n$th letter. Since $|u_n-t_2^n|$ converges to 0 as $n\rta\infty$, by Lemma \ref{prop22}, and since we argued above that the lemma holds with $t_2^n$ instead of $u_n$, our claim implies that the lemma also holds with $u_n$.
\end{proof}

\begin{lemma} 
	Let $j\in \N_+$, and let $\li_j^n\subset\R$ (respectively, $\li_j\subset\R$) be the envelope interval corresponding to $\ga_j^n$ (respectively, $\gab_j$). Then $\li_j^n\rta\li_j$ in probability.
	\label{lem:env-cov}
\end{lemma}
\begin{proof} 
	Recall from \eqref{eq:value} that a percolation cycle $\ga$ may be associated with a value val$(\ga)$. For any fixed interval $I$ let val$_n(I)$ (respectively, val$(I)$) be the value of the percolation cycle $\ga^n$ (respectively, CLE$_6$ loop $\gab$) with envelope interval $\env_n(I)$ (respectively, $\env(I)$). Recall from Lemma \ref{prop10} that $\area_n(\ga^n)\rta\area(\gab)$ in probability as $n\rta\infty$, and recall from Lemma \ref{prop2} that $\env_n(I)\rta\env(I)$ in probability. 
	Let $T^n$ and $\Tb$ be defined by \eqref{eq:Tu} and \eqref{eq:T-cle}, respectively, relative to the envelope closing time of $\ga^n$ and $\gab$. Then $T^n\rta\Tb$ by Lemma \ref{prop21}. 
	
	By symmetry we can assume that $\Zb|_{\env(I)}$ is a left cone excursion. Furthermore, we may assume without loss of generality that $\ga^n$ is the outside-cycle of a black cluster, since this holds for all sufficiently large $n$. If $J\subset \env(I)$ is a complementary component of the range of $\Tb$, then $\Zb|_J$ is a cone excursion. The set $\etab(J)$ is in the interior of $\gab$ if and only if this cone excursion is a right cone excursion. A similar description holds in the discrete. Indeed, by Theorem \ref{thm:LR}, if $J^n$ is an interval ending before the envelope closing time of $\ga^n$ such that  $Z^n|_{J^n}$ is a cone excursion, and if $J^n$ is not contained inside any larger interval satisfying these two properties, then $Z^n|_{J^n}$ encodes a bubble $B$ of one of the looptrees of the spine-looptrees decomposition relative to the envelope closing time of $\ga^n$. Furthermore, $Z^n|_{J^n}$ is a right (respectively, left) excursion if and only if $B$ is part of the black (respectively, white) looptree. Therefore, $B$ is enclosed by $\ga^n$ if and only if $Z^n|_{J^n}$ is a right excursion and $J^n\subset\env_n(I)$. It follows from the above discussion that for any fixed intervals $I,I'\subset\R$ such that $I\subset I'$ and for all sufficiently large $n$, the percolation cycle with envelope interval $\env_n(I)$ is enclosed by the percolation cycle with envelope interval $\env_n(I')$ if and only if this holds for the associated CLE$_6$ loops.
	
	Let $\ga'$ and $\ga'_n$ be as in \eqref{eq:value0} and \eqref{eq:value}, respectively. Then the envelope interval of $\ga'_n$ converges to the envelope interval of $\ga'$ for the metric $d_{\op{I}}$ as $n\rta\infty$, by Lemma \ref{prop10}, and since the preceding paragraph implies that for any fixed $I\subset\R$ and all sufficiently large $n$, the percolation cycle with envelope interval $\env_n(I)$ surrounds the origin if and only if this holds for the CLE$_6$ loop with envelope interval $\env(I)$. By the preceding paragraph we also get that the envelope interval of $\op{anc}(\ga^n,\ga'^n)$ converges toward the envelope interval of $\op{anc}(\gab,\gab')$. Now applying Lemma \ref{prop10}, we get that for any fixed interval $I$,
	\eqb\label{eq32}
	\val_n(I)\rta \val(I)
	\eqe
	in probability as $n\rta\infty$.
	
	Let $\eps>0$.  Let $k_1\in\N$ be such that with probability at least $1-\eps$, 
	the event $E(k_1)= \{\val(\gab_j)<2^{k_1-3}\}$ occurs. 
	Let $\mcl I_k$ denote the set of closed intervals such that both endpoints are contained in $(2^{-k}\Z)\cap[-k,k]$. 
	Then choose an integer $k_2>k_1$ such that with probability at least $1-\eps$ the following event $E(k_2)$ occurs:
	\begin{itemize}
		\item[(i)] for all $j'\in[j]$, if $\env_{j'}=[s_1,s_2]$ is such that $\Zb|_{\env_{j'}}$ is a left (respectively, right) cone excursion and 
		$s_0=\inf\{t<s_1\,:\,\Rb_s>\Rb_{s_2}\,\forall s\in[t,s_1] \}$ (respectively, $s_0=\inf\{t<s_1\,:\,\Lb_s>\Lb_{s_2}\,\forall s\in[t,s_1] \}$), then $|s_0-s_1|>2^{-k_2+3}$, and
		\item[(ii)] for any loop $\gamma''$ containing the loop $\gamma'$ in \eqref{eq:value0} such that $\env(\gamma'')\not\subset[-k_2,k_2]$ we have $\area(\gamma'')>2^{k_1}$.
	\end{itemize} 
	Let $E(k_1,k_2)=E(k_1)\cap E(k_2)$.
	By a union bound, $\P[E(k_1,k_2)]\geq 1-2\eps$.
	
	We will now observe that on the event $E(k_1,k_2)$, for all $j'\in[j]$ there exists $I\in\mcl I_2$ such that $\env(\gamma_{j'})=\env(I)$. It is sufficient to show that if $\gamma=\gamma_{j'}$ satisfies $\env(\gamma)\neq\env(I)$ for all $I\in\mcl I_2$ then $\val(\gamma)>2^{k_1-3}$ since this contradicts $E(k_1)$. If $\gamma=\gamma_{j'}$ satisfies $\env(\gamma)\neq\env(I)$ for all $I\in\mcl I_2$, at least one of the following hold: 
	(a) $\env(\gamma)<2^{-k_2+1}$, 
	(b) there is a $I\in\mcl I_1$ such that $\env(I)\subsetneq \env(\gamma)$ and $\env(\gamma)\setminus\env(I)$ has Lebesgue measure smaller than $2^{-k_2+1}$, or 
	(c) $\env(\gamma)\not\subset[-k_2,k_2]$. In either case we have $\val(\gamma)>2^{k_1-3}$: 
	In case (a) we use that the value of a percolation cycle is always at least the inverse length of its envelope interval,
	in case (b) we use (i) in the definition of $E(k_2)$, which implies that (b) cannot occur, and 
	in case (c) we use (ii) in the definition of $E(k_2)$ and $\op{val}(\gab)>\area (\sma(\gab,\gab'))$.

	Recall that with probability $1-o_n(1)$, for all $I\in\mcl I_{k_2}$ we have $\val_n(I)\rta \val(I)$ and $\env_n(I)\rta\env(I)$. Furthermore, on $E(k_1,k_2)$, for all $j'\in[j]$ there exists $I\in\mcl I_2$ such that $\env(\gamma_{j'})=\env(I)$. Since the percolation cycles and CLE$_6$ loops are ordered by their value, this implies that with probability $1-o_n(1)$, among the percolation cycles whose envelope interval can be written on the form $\env_n(I)$ for $I\in\mcl I_2$, there are at least $j$ distinct percolation cycles with value smaller than $2^{k_1-2}$, so $\val(\ga_j^n)<2^{k_1-2}$.
	
	To conclude the proof it is sufficient to argue that with probability at least $1-3\eps$ for all sufficiently large $n$, for each $j'\in [j]$ we have $\env^n_{j'}=\env_n(I_{j'})$ for some fixed $I_{j'}\in\mcl I_{k_2}$. This is sufficient since we know that $\env_n(I)\rta\env(I)$ in probability for each $I\in\mcl I_{k_2}$.

	We will proceed by contradiction. We condition on the event $E(k_1,k_2)$, and assume we can find arbitrarily large $n$ for which there exists $j'\in[j]$ such that $\env^n_{j'}\neq\env_n(I)$ for all $I\in\mcl I_{k_2}$. Recall that $\val(\ga_j^n)<2^{k_1-2}$. 
	Let $j'\in[j]$ and $n$ be such that $\env^n_{j'}\neq\env_n(I)$ for all $I\in\mcl I_{k_2}$. Since the value of a percolation cycle is at least the inverse of the length of its envelope interval, this implies that $\env_{j'}^n$ has length at least $2^{-k_1+2}$ for all $j'\in[j]$. Similarly, using (ii) in the definition of $E(k_2)$ and that $\val(\ga^n_{j'})>\area(\op{anc}(\ga^n_{j'},\ga'_n))$ (with $\ga'_n$ as in \eqref{eq:value}), we get that $\env_{j'}^n\subset[-k_2,k_2]$ for all $j'\in[j]$. Since $\env_{j'}^n\subset[-k_2,k_2]$ has length at least $2^{-k_1+2}$ but is not equal to $\env_n(I)$ for any $I\in\mcl I_{k_2}$, there exists $I\in\cI_{k_2}$ such that $\env(I)\subset\env_{j'}^n$, and such that the distance between the left (respectively, right) endpoint of $\env_{j'}^n$ and $\env(I)$ is less than $2^{-k_2}$. 
	
	Let $[t_1,t_2]=\env_n(I)$. We may assume without loss of generality that $t_2$ is of $b$-type. By the definition of $E(k_2)$, for all sufficiently large $n$,
	\eqb
	\not\exists\,(t'_1,t'_2) \in [t_1-2^{-k_2},t_1]\times [t_2,t_2+2^{-k_2}]
	\quad\text{such\,\,that}\quad
	t'_1 \text{\,\,is\,\,the\,\,$a$-match\,\,of\,\,} t'_2.
	\label{eq94}
	\eqe
	By \eqref{eq94}, the right endpoint of $\env_{j'}^n$ must be of $b$-type, so $\ga^n_{j'}$ is the outside-cycle of a black cluster.
	Then $\matchingpar_n(t_2)$ is the end of a left cone excursion containing $\env_n(I)$. Since $\matchingpar_n(t_2)$ and $t_2$ are of different types by the definition of an envelope interval, $\matchingpar_n(t_2)$ is not equal to the envelope closing time of $\ga^n_{j'}$. Furthermore, $\matchingpar_n(t_2)$ must be strictly smaller than the envelope closing time of $\ga^n_{j'}$, by the definition of $\matchingpar_n(t_2)$ and since the two endpoints of $\env_{j'}^n$ define a $b$-match. It follows that the cone interval $J$ ending at $\matchingpar_n(t_2)$ satisfies $\env_n(I)\subset J\subset\env_{j'}^n$. Let $J'$ be the cone interval chosen as large as possible such that $J\subseteq J'\subset \env_{j'}^n$. By \eqref{eq94}, $J'$ is a left cone interval. By Theorem \ref{thm:LR}, the cone excursion $Z^n|_{J'}$ encodes a bubble $B$ of the left (i.e., white) looptree in the spine-looptrees decomposition relative to the envelope closing time of $\ga^n_{j'}$. In particular, the bubble $B$ is \emph{not} a bubble of $\frk L(\ga)$ (Definition \ref{def:lt-cluster}) since this looptree is black. Therefore the area of $\ga$ is bounded above by the length of its envelope interval minus the length of $J'$, so $\area(\ga)<2\cdot 2^{-k_2}$, which is a contradiction.
\end{proof}

\subsubsection{Pivotal points}\label{subsub:pivot}
In this section we prove the convergence of pivotal measures. We adopt the notation $\wh Z^n$, $\wh L^n$, $\wh R^n$, $T^n$, $\ell^n$ of Section~\ref{sec:piv-def} (relative to a given time $u$), and their continuum counterparts $\wh \Zb$, $\wh \Lb$, $\wh \Rb$, $\Tb$, $\ellb$. Recall also the measures $p_{\Lo}^n(s,u)$, $p_{\Ro}^n(s,u)$, and $p^n(s,u)$ defined in Section~\ref{sec:piv-def}, and their continuum counterparts $\pb_{\Lo}(s,u)$, $\pb_{\Ro}(s,u)$, and $\pb(s,u)$. We will use the function $\wh\eta_{\op{v}}:\ZZ^{\leq 0}\to V(M_n)$ defined in Section \ref{subsec:pivot} in order to identify the times sampled from the measures $p_{{\Lo}}^n(s,u)$  or $p_{{\Ro}}^n(s,u)$ with vertices of the triangulation $M_n$.

We first prove a lemma asserting that a pivotal point sampled from one of the measures $p_{{\Lo}}^n(s,u)$ or $p_{{\Ro}}^n(s,u)$ is $\eps n$-significant with probability $1-o_\eps(1)$, where the $o_\eps(1)$ is uniform in $n$. 
In the statement and proof of the lemma, if $\sigma$ is a measure of finite total mass then $x\sim\sigma$ means that $x$ is sampled from $\sigma$ renormalized to be a probability measure. 
\begin{lemma}\label{prop19} 
	\begin{compactitem}
		\item[(i)] For fixed $u\in\R$, $s<0$, and $\xi>0$ there is an $\eps>0$ depending on $s,\xi$ but not on $n$
		such that if $t^n_+\sim p_{\Lo}^n(s,u)$, then $\wh\eta_{\op{v}}(n^{3/4}t^n_+)$ is an $\eps n$-pivotal point with probability at least $1-\xi$. 
		\item[(ii)] For $\xi>0$, $j\in\N_+$, and (possibly random) $s<0$, there is an $\eps>0$ depending on $s,\xi,j$ but not on $n$
		such that if $u_n$ is the envelope closing time of the percolation cycle $\gab_j^n$ and $t_+^n\sim p_{\Lo}^n(s,u_n)$, then $\wh\eta_{\op{v}}(n^{3/4}t_+^n)$ is an $\eps n$-pivotal point with probability at least $1-\xi$. 
		\item[(iii)] For fixed $u\in\R$ and $s<0$ let $t^n_1\sim p_{\Lo}^n(s,u)$ and define 
		$t^n_2=\inf\{t>t^n_1\,:\, \wh L^n_t<\wh L^n_{t^n_1} \}$ (on the event of probability $o_n(1)$ that this is not well-defined, define $t^n_2$ arbitrarily).
		Then, for any fixed $\eps>0$, $\P[T_{t^n_2}^n-T_{(t^n_1)^-}^n<\eps]=1-o_n(1)$, where the notation $T^n_{a^-}$ stands for $\lim_{b\uparrow a}T^n_{b}$.
	\end{compactitem} 
	The assertions (i-iii) still hold with $\mathrm L$ and $\wh L^n$ replaced by $\mathrm R$ and $\wh R^n$, respectively.
\end{lemma}
The following notion will be convenient when dealing with LQG looptrees.
Given an LQG looptree $\frk L$, we call the total $\mub_{\gff}$-area of the $\sqrt{8/3}$-LQG disks associated with the bubbles of $\frk L$ the \emph{significance} of $\frk L$ and we denote it by $\sig(\frk L)$. Given $\eps>0$, we say that $\frk L$ is \emph{$\eps$-significant} if $\sig(\frk L)\ge \eps$.
\begin{proof} 
	We first prove	(i).  
	Recall the notation $\fll,\flr,\pi_{\Lb},\pi_{\Rb}$ of Section~\ref{sec:dictionary-fl} (relative to the time $u$). Let $t_+\sim \pb_{\Lo}(s,u)$. 
	By Lemma \ref{lem:psu} we may extend our coupling such that $t^n_+\rta t_+$ almost surely. 
	Let
	$$s_+=\inf\{t\in [t_+,0] : \wh \Lb_{t'} >\wh \Lb_t \; \forall t'\in (t,0] \}\quad \textrm{and}\quad s_-= \inf\{ s <s_+\,:\, \wh\Lb_{s'}>\wh\Lb_{s_+}\,\forall s'\in(s,s_+) \}.$$
	Then $\wh\Lb|_{[s_-,s_+]}$ is a $3/2$-stable L\'evy excursion encoding a looptree $\frk L$ on $\fll$ such that $p:=\pi_{\Lb}(t_+) \in \frk L$. Let 
	$t_-= \inf\{ t <t_+\,:\, \wh\Lb_{t'}>\wh\Lb_{t_+}\,\forall t'\in(t,t_+) \}$
	which implies $\pi_{\Lb}(t_-)=p$ almost surely. Let $\frk L'$ be the looptree on $\fll$ encoded by $\wh \Lb|_{[t_-,t_+]}$, and let $\frk L^\pm$ be the looptrees on $\flr$ containing $\pi_{\Rb} (t_\pm)$. 
	We will now argue that $z:=\wh\etab(t_+)$ is an $\eps$-pivotal point if the following four events occur:
	(1) $\sig(\frk L')\ge \eps$; 
	(2) $\sig(\frk L) -\sig(\frk L') \ge \eps$;
	(3) $\sig(\frk L^+)\ge \eps$; 
	(4) $\sig(\frk L^-)\ge \eps$.
	
	By Lemma \ref{lem:pivotal-continuous}, all pivotal points for an instance of CLE$_6$ are of exactly one type 1-4, with types defined as in Section \ref{subsec:pivot}. 
	If $z$ is of type 1 then there is a CLE$_6$ loop $\gab$ such that $\phi_{\Lo} (\frk L)\subset \gab$. In this case, (1) and (2) imply that $z$ is an $\eps$-pivotal. 
	If $z$ is of type 2 then there is a CLE$_6$ loop $\gab$ such that $\phi_{\Lo}(\frk L'),\phi_{\Ro}(\frk L^\pm)\subset\gab$, and we have that $z$ is an $\eps$-pivotal if (1) and at least one of (3) and (4) hold.
	If $z$ is of type 3 or type 4 then there are distinct CLE$_6$ loops $\gab_+$ and $\gab_-$ such that $\phi_{\Ro}(\frk L^\pm)\subset\gab_\pm$, and $\gab_+$ and $\gab_-$ are nested if and only if $z$ is of type 4. If $z$ is of type 3 we have that $z$ is an $\eps$-pivotal if both (3) and (4) hold, while if $z$ is of type 4 it is sufficient if either (1) or (2) hold in addition to (3). In either case we see that (1)-(4) implies that $z$ is an $\eps$-pivotal.
	
	The events $(1)$ and $(2)$ occur with probability $1-o_\eps(1)$ since $\sig(\frk L) > \sig(\frk L')>0$ almost surely. The events $(3)$ and $(4)$ occur with probability $1-o_\eps(1)$ since $\wh\Lb$, hence $t_-$ and $t_+$, are independent of $\wh\Rb$. 
	\begin{figure}
		\centering
		\includegraphics[scale=.75]{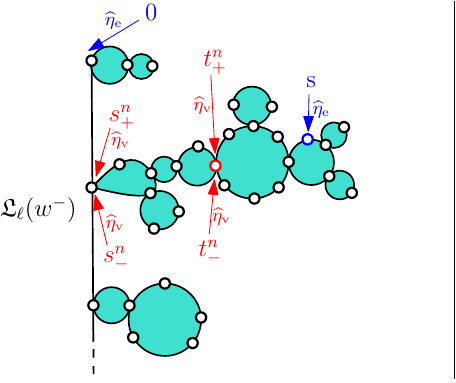}\quad\quad	
		\includegraphics[scale=1.3]{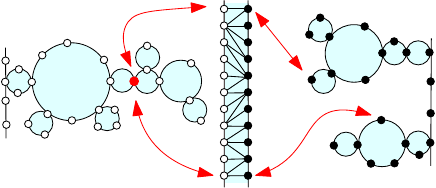}
		\caption{Left: The discrete forested line $\frk(w^-)$ (relative to time $u$) for $(M_n,\si_n)$, and the vertices corresponding to the times  $s, t_+^n,t_-^n,s_+^n,s_-^n$. An arrow labeled $\wh \etav$ indicates the vertex $\wh \etav(n^{3/4}t)$ associated to a time $t$. An arrow labeled $\wh \etae$ indicates the white endpoint of the edge $\wh \etae(n^{3/4}t)$ associated to a time $t$.
			Right: Illustration of the proof of Lemma \ref{prop19}(i). A uniformly sampled pivotal point from the measure $p_{{\Lo}}^n(s,u)$ (in red on the figure) is $\eps n$-significant with probability $1-o_\eps(1)$ uniformly in $n$. The red arrows indicate pairs of vertices which are identified.}
		\label{fig:pivotal-discrete}
	\end{figure}

	Define the discrete times $t_-^n,s_-^n,s_+^n$ in the same way as $t_-,s_-,s_+$ with $\wh L^n,\wh R^n$ in place of $\wh \Lb,\wh \Rb$. The meaning of these discrete times is represented in Figure \ref{fig:pivotal-discrete} (left).
	By convergence of $(\wh L^n,\wh R^n)$ to $(\wh \Lb,\wh \Rb)$, the discrete analog of the events (1-4) also hold with probability $1-o_\eps(1)$ for $\wh L^n, \wh R^n$, where the $o_\eps(1)$ in uniform in $n$. Note that when defining $\frk L^\pm$ there may be multiple vertices on the forested line encoded by $\wh R^n$ which are connected to the sampled pivotal point $\wh\eta_{\op{v}}(n^{3/4}t_+^n)$, but, by convergence of
	$(\wh L^n,\wh R^n)$ to $(\wh \Lb,\wh \Rb)$,
	with probability $1-o_n(1)$ these can be divided into two disjoint sets, such that all vertices in the same set are on the same looptree of the forested line encoded by $\wh R^n$. See the right part of Figure \ref{fig:pivotal-discrete} for an illustration.
	By the first assertion of Lemma \ref{prop:piv}, the vertex $v=\wh\eta_{\op{v}}(n^{3/4}t_+^n)$ is a pivotal point with probability $1-o_n(1)$.
	By considering separately the different types of pivotal points 1-4, we see that discrete analogs of (1-4) imply that $v$ is an $\eps n$-pivotal. This implies that the vertex $v$ is an $\eps n$-pivotal with probability $1-o_\eps(1)$. 
	For example, in the case of a pivotal point of type 1, and with notation as in the continuum case, the percolation cycle immediately surrounding $\frk L$ will be contained in $\cL_v\cap\Gamma$ (Definition \ref{def:pivot}), while we will have two macroscopic percolation cycles in $\cL_v\cap\Gamma_v$, one of which immediately surrounds $\frk L'$. 
	This concludes the proof of (i).
	
	Next, we prove (ii). By Lemma \ref{prop21}, $p_n(s,u_n)$ converges in probability to $\pb(s,u)$, which assigns mass 0 to 0. Therefore the mass assigned to $[-\delta,0]$ by $p_n(s,u_n)$ is $o_\delta(1)$ with probability $1-o_\delta(1)$, where the $o_\delta(1)$ is uniform in $n$. Recall that $\nu_n(s,u_n)$ is the pushforward of $p_n(s,u_n)$ under $\phi\circ \wh\eta_{\op{v}}$. 
	Uniformly over all $A\subset \C$, for all sufficiently large $n$ and for $v>0$, the difference between the $\nu_n(s,u_n)$-mass and the $\nu_n(s,u_n-v)$-mass of $A$ is bounded above by $o_{v}(1)$ with 
	probability $1-o_n(1)$. 
	Suppose a pivotal point is sampled from $\nu_n(s,q)$ for a fixed $q\in\R$. Then the argument in our proof of Assertion (i) shows that this pivotal point is an $\eps n$-pivotal with probability $1-o_\eps(1)$. By approximating $u_n$ by $q\in\Q$ the same holds for $\nu_n(s,u_n)$. 
	
	Lastly, Assertion (iii) is immediate by the convergence of $\wh L^n$, $T^n$, and $p_{\Lo}^n(s,u)$ established in Lemmas~\ref{prop5} and \ref{prop21}, since this gives that $t_2^n-t_1^n$ converge to zero in probability, and that $T^n_{t_1^n+\delta}-T^n_{t_1^n-\delta}=o_\delta(1)$ with probability at least $1-\delta$ for all $\delta\in(0,1)$ and $n\in\N$ sufficiently large.
\end{proof}
Recall that for $s<0$ and $u\in\R$ we defined the sets $\Ab_{\op L}(s,u), \Ab_{\op R}(s,u), \Ab (s,u)\subset(-\infty,0]$ in Section~\ref{sec:piv-def}. 
For $t\leq 0$ recall the time $\frk s_t$ defined by \eqref{eq:st}. If $\frk s_t\neq t$, then $\Xb=\wh\Lb|_{[\frk s_t,t]}$ is a L\'evy excursion which encodes a looptree $\frk L=\frk L_{\Xb}$. 
Let $\Ab^{\op{L}}_\eps(s,u)$ be the union of $\Ab_{\op{L}}(s,u)$ and the set of $t\leq 0$ such that $\frk s_t\neq t$ and the associated looptree $\frk L$ is $\eps$-significant. 
Define $\Ab^{\op{R}}_\eps(s,u)$ similarly, and set $\Ab_\eps(s,u)=\Ab^{\op{L}}_\eps(s,u)\cup \Ab^{\op{R}}_\eps(s,u)$.
By the looptree perspective on $\CLE_6$ loops described in Section~\ref{sec:cont-cle} and Lemma~\ref{lem:cover}, we have the following lemma.
\begin{lemma}
	Fix $\eps>0$. Then the following holds for all CLE$_6$ loops $\gab$ almost surely. Let $u$ be the envelope closing time of $\gab$, and define $s<0$ such that $|s|$ is the quantum natural length of $\gab$. Then we can almost surely find rationals $q_1,\dots,q_k<0$ such that 
	\eqb
	\Ab_\eps(s,u) \subset
	\bigcup_{i=1}^k \Ab(q_i,u) \cup \Ab(s,u).
	\label{eq46} 
	\eqe
	\label{prop30}
\end{lemma}
For $j\in \N_+$, let $P^{n}_{\eps,j}\subset V(M_n)$ be the set of $\eps n$-pivotal points associated with the percolation cycle $\gamma^n_j$. Recall the sets $S_{\ga_j^n,\ga'}$ defined in Lemma \ref{prop:piv}, each with cardinality at most three.
\begin{lemma} 
	For $j\in\N_+$ let $u$ (respectively, $u_n$) be the envelope closing time of $\gab_j$ (respectively, $\ga_j^n$), and let $-s>0$ (respectively, $-s_n>0$) be the quantum natural length (respectively, renormalized length) of $\gab_j$ (respectively, $\ga_j^n$). Let $q_1,\cdots,q_k\in\Q$ be as in Lemma~\ref{prop30} for $\gab=\gab_j$. With probability converging to 1 as $n\rta\infty$,
	\eqbn
	P^{n}_{\eps,j}\setminus \bigcup_{\ga'\in\Gamma\setminus \{\ga_j^n \}\,:\,\sig(\ga')\geq\ep} S_{\ga_j^n,\ga'} 
	\subset
	\left\{\wh\eta_{\op{v}}(n^{3/4}t)\,:\, t\in \left(\bigcup_{i=1}^k A^n(q_i,u_n)\right)\cup A^n(s_n,u_n)\right\}.
	\eqen
	Furthermore, for each $v\in P^{n}_{\eps,j}\setminus \bigcup_{\ga'} S_{\ga_j^n,\ga'}$ and $q\in\{q_1,\dots,q_k,s_n \}$ there is at most one $t\in A^n(q,u_n)$ such that $v=\wh\eta_{\op{v}}(n^{3/4}t)$.
	\label{prop16}
\end{lemma}

\begin{proof}
	For $\alpha>0$ we say that a looptree identified with a subset of vertices and edges of a map $M$ is an $\alpha$-looptree if the cycle immediately surrounding the looptree has area at least $\alpha$ (Definition \ref{def:pivot}). Let $A^{n,{\Lo}}_\eps(s,u_n)$ be defined in the exact same way as $\Ab^{{\Lo}}_\eps(s,u)$, that is, $t\in A^{n,{\Lo}}_\eps(s,u_n)$ if and only if $t$ is a strict running infimum for $L^n|_{[s,0]}$, or if for $t_-=\sup\{s<t:L^n_{s}=L^n_t \}$ the looptree with counterclockwise code $L^n|_{[t_-,t]}$ is an $\eps n$-looptree. Define $A^{n,{\Ro}}_\eps(s,u_n)$ similarly with $R$ instead of $L$, and set $A^{n}_\eps(s,u_n)= A^{n,{\Lo}}_\eps(s,u_n)\cup A^{n,{\Ro}}_\eps(s,u_n)$. To conclude the proof it is sufficient to show the following: 
	(i) \(A^{n}_{\eps}(s_n,u_n) \subset \left(\bigcup_{i=1}^k A^n(q_i,u_n)\right)\cup A^n(s_n,u_n),\) and 
	(ii) $A^{n}_{\eps,j}\setminus \bigcup_{\ga'} S_{\ga,\ga'} \subset \{ \wh\eta_{\op{v}}(n^{3/4}t)\,:\,t\in A^{n,{\Lo}}_\eps(s_n,u_n)\}$. 
	
	Assertion (i) follows from Lemma~\ref{prop30} and Assumption~\ref{a1}. In fact, if for arbitrarily large $n$ we can find $t_n\in A^{n,{\Lo}}_{\eps}(s_n,u_n)\setminus A^n(s_n,u_n)$ such that $t_n\not\in\bigcup_{i=1}^k A^n(q_i,u_n)$ with uniformly positive probability, then as $n\to\infty$, a compactness argument implies that with positive probability there exists $t\in \Ab^{{\Lo}}_{\eps}(s,u)\setminus \Ab(s,u)$ such that 
	$t\not\in\bigcup_{i=1}^k \Ab(q_i,u)$, which contradicts Lemma~\ref{prop30}. 
	
	Assertion (ii) follows from the last assertion of Lemma \ref{prop:piv}.
\end{proof}

\subsubsection{Proof of Theorem~\ref{thm1}} 
\label{sec:proof-main}

We are now well equipped to prove Theorem~\ref{thm1}. 

Recall that in \eqref{eq35}, for a given metric space $(B,d_B)$, we defined a metric $d_{\op p}$ for comparing parametrized curves on $B$. Also recall that $d_{\BB S^2}$ denotes the spherical metric on $\C$.	
\begin{fact} \label{fact: curve}
	For $n\in \N_+$, let $I_n\subset\R$ be intervals and let $f_n:I_n\to \R$ be such that $f_n$ converges in the $d_p$-metric if we equip $\R$ with the Euclidean metric.
	Let $J_n\subset \R$ be intervals such that $f(I_n)\subset J_n$ for all $n\in \R$. Given a metric space $(B,d_B)$, let 
	$g_n:J_n\to B$ be such that $g_n$ converges in the $d_{\op{p}}$-metric. 
	Then $g_n\circ f_n$ also converges in the $d_{\op{p}}$-metric. 
\end{fact}

Recall that in \eqref{def-phie} when defining the embedding $\phi^n|_{E(M_n)}$ of the edges of $M_n$ we identified an edge $e=\{u_1,u_2\}$ with the mid-point between $\phi^n(u_1)$ and $\phi^n(u_2)$. However, in our proof it is more convenient to consider an embedding where $e$ is mapped to $\etab((3n)^{-1}\etae^{-1}(e))\in\C$. We will use Lemma \ref{prop6} to argue that the two embeddings are asymptotically equivalent, in the sense that uniformly over all choices of $e$, the distance between the two embeddings of $e$ is $o_n(1)$. The proof of Lemma \ref{prop6} uses Lemma \ref{prop35}, which we state and prove first.

Recall the future/past decomposition described in Section \ref{subsec:future-past}. For $i\in\N_+$ let $\xi^\Lo_+(i)$ (respectively, $\xi^\Ro_+(i)$) denote the $i$th vertex along the left (respectively, right) boundary of the future near-triangulation, such that the end-points of the root-edge are $\xi^\Lo_+(1)$ and $\xi^\Ro_+(1)$. Similarly, let $\xi^\Lo_-(i)$ (respectively, $\xi^\Ro_-(i)$) denote the $i$th vertex along the left (respectively, right) boundary of the \emph{past} near-triangulation, such that the end-points of the top-edge are $\xi^\Lo_-(1)$ and $\xi^\Ro_-(1)$. For a bounded interval $J\subset\R$ let $V^{\Lo}_+(J)\subset V(M_n)$ (respectively, $V^{\Ro}_+(J)\subset V(M_n)$) be the set of vertices $v$ on the left (respectively, right) boundary of the future near-triangulation such that $(3n)^{-1}\eta_{\op{vf}}^{-1}(v)\in J$.
Let $V^{\Lo}_-(J)\subset V(M_n)$ (respectively, $V^{\Ro}_-(J)\subset V(M_n)$) be the set of vertices $v$ on the left (respectively, right) boundary of the past near-triangulation such that for some edge $e$ with end-point $v$ in the interior of the past near-triangulation, we have $(3n)^{-1}\eta_{\op{e}}^{-1}(e)\in J$. 
\begin{lemma}
	For any bounded interval $J\subset(-\infty,0]$ or $J\subset[0,-\infty)$ the following convergence results hold in probability for the Hausdorff distance
	\eqbn
	\begin{split}
		\{ 0.5\sqrt{\beta/n}\cdot (\xi^\Lo_-)^{-1}(v)\,:\,v\in V^{\Lo}_-(J) \} 
		&\rta
		\ol{\{ \Lb_t\,:\,t\in J, \Lb_t=\inf_{t'\in[t\wedge 0,t\vee 0]} \Lb_{t'} \}}\qquad\text{for\,\,}J\subset(-\infty,0],\\
		\{ 0.5\sqrt{\beta/n}\cdot (\xi^\Lo_+)^{-1}(v)\,:\,v\in V^{\Lo}_+(J) \} 			&\rta			\ol{\{ \Lb_t\,:\,t\in J, \Lb_t=\inf_{t'\in[t,0]} \Lb_{t'} \}}\qquad\text{for\,\,}J\subset[0,\infty).
	\end{split}
	\eqen
	By symmetry, the same result holds with $\op{L}$ and $\Lb$ replaced by $\op{R}$ and $\Rb$, respectively.
	\label{prop35}
\end{lemma}
\begin{proof}
	By uniform convergence of $Z^n$ to $\Zb$ on compact sets it is sufficient to establish the following for $J\subset[0,\infty)$
	\eqb
	\{ 0.5\sqrt{\beta/n}\cdot (\xi^\Lo_+)^{-1}(v)\,:\,v\in V^{\Lo}_+(J) \}
	= \{ L^n_t\,:\,t\in J, L^n_t<\inf_{t'\in[0,t)} L^n_{t'} \}
	\label{eq99}
	\eqe
	and the following for $J\subset(-\infty,0]$, where $(L_i)_{i\in\Z^{<0}}$ is the unscaled lattice walk
	\eqb
	\begin{split}
		\text{With $v\in V^{\Lo}_-(J)$ and $e=\eta_{\op{e}}(i)$ as in the definition of $V^{\Lo}_-(J)$ we have}\\
		\text{$i=3tn$, $L_i\in\{ (\xi_-^{\Lo})^{-1}(v)+1,(\xi_-^{\Lo})^{-1}(v)+2 \}$, and $\min_{i'\in\{i,\dots,-1 \}} L_{i'}\geq L_i-1$.}
	\end{split}
	\label{eq100}
	\eqe
	These identities follow from the interpretation of $Z^n$ as a boundary length process and the bijection between vertices/faces (respectively, edges) and steps of the walk (see Remark \ref{rk:phi-well-def} and Definitions \ref{def:eta} and \ref{def:eta-extend}). 
	In particular, for any vertex $v$ on (say) the left boundary of the future percolated near-triangulation, if $v=\eta_{\op{vf}}(i)$ then in step $i$ of the bijection we apply the mapping $\ol\phi_c$ (see Section \ref{subsec:bij-chordal} for the definition) for the case where there is no active left edge, and the edge between $(\xi^\Lo_-)^{-1}(v)$ and $(\xi^\Lo_-)^{-1}(v)+1$ is the top-edge immediately before this step. This gives \eqref{eq99}.
	For any vertex $v$ on the left boundary of the past site-percolated triangulation, if $e$ is an interior edge of the past site-percolated triangulation with end-point and $\eta_{\op{e}}(i)=e$, then we have $L_{i}=-(\xi^\Lo_-)^{-1}(v)+1$ (respectively, $L_{i}=-(\xi^\Lo_-)^{-1}(v)+2$) if $w_i\in\{b,c \}$ (respectively, $w_i=a$). This gives \eqref{eq100}. 
\end{proof}


\begin{lemma} 
	Given $e\in E(M_n)$ with endpoints $u_1,u_2\in V(M_n)$, define
	\eqbn
	t^0_e:=\frac{1}{3n}\etae^{-1}(e),\quad 
	t^1_e:=\frac{1}{3n}\etavf^{-1}(u_1),\quad
	t^2_e:=\frac{1}{3n}\etavf^{-1}(u_2).
	\eqen
	This is well-defined if we require that $t^1_{e}<t^2_{e}$. Then it holds almost surely that
	\eqb
	\sup_{e\in E(M_n) }
	d_{\BB S^2}(\etab(t^0_e),\etab(t^1_e))+
	d_{\BB S^2}(\etab(t^0_e),\etab(t^2_e))
	\rta 0 
	\qquad
	\text{as}\,\, 
	n\rta\infty.
	\nonumber
	\eqe 
	\label{prop6}
\end{lemma}
\begin{proof}
	Let $\ep>0$. It is sufficient to show that with probability at least $1-\ep$ for sufficiently large $n$,
	\eqb
	\sup_{e\in E(M_n) }
	d_{\BB S^2}(\etab(t^0_e),\etab(t^1_e))<\ep,
	\label{eq34}
	\eqe
	since the case of $t^2_e$ (instead of $t^1_e$) can be treated in the exact same way. For all $\wt\eps>0$ there exists $K>0$ such that for all $t\notin[-K,K]$, $\etab(t)$ is in the $\wt\eps$-ball around $\infty$ for the metric $d_{\BB S^2}$. Therefore it is sufficient to prove that the following holds with probability at least $1-\ep$ for an arbitrary fixed $K$ and sufficiently large $n$, 
	\eqb
	\sup_{e\in E(M_n),\{t_0^1,t^1_e\}\cap [-K,K]\neq\emptyset }
	d_{\BB S^2}(\etab(t^0_e),\etab(t^1_e))<\ep.
	\label{eq92}
	\eqe	
	We start by proving the weaker statement that the following holds with probability at least $1-\ep$ for all sufficiently large $n$ 
	\eqb
	\sup_{e\in E(M_n),\{t_0^1,t^1_e\}\subset [-K,K] }
	d_{\BB S^2}(\etab(t^0_e),\etab(t^1_e))<\ep.
	\label{eq93}
	\eqe
	Notice that we always have $t_e^0\leq t_e^1$ (with equality if and only if $w_{3t_e^1}=c$). Since $\etab|_{[-K,K]}$ admits a modulus of continuity we may assume that for some $a>0$ depending only on $\eps$ we have $t_e^1-t_e^0>a$. Therefore, by invariance under recentering at any fixed time $q\in\R$, when proving \eqref{eq93} we may assume that $t_e^0<0<t_e^1$. In other words, it is sufficient to show that the following holds with probability at least $1-\ep$ for all sufficiently large $n$
	\eqb
	\sup_{e\in E(M_n),~t_e^0\in [-K,0),~t^1_e\in (0,K] }
	d_{\BB S^2}(\etab(t^0_e),\etab(t^1_e))<\ep.
	\label{eq95}
	\eqe

	For $\delta>0$ let $E_1=E_1(K,\eps,\delta)$ be the event that for all $t\in[-K-\delta,K]$ the set $\etab([t,t+\delta])\subset\C$ has diameter smaller than $\ep/2$ for the spherical metric. Then $\lim_{\delta\rta 0}\P[E_1]=1$. Fix $\delta>0$ such that $\P[E_1]>1-\ep/2$. 
	
	Let $I_1,\dots,I_S$ for $S\in\N$ be an enumeration of the intervals of the form $[(k-1)\delta,k\delta]$ for $k\in\Z$ which intersect $[-K,K]$. Assume the intervals are ordered such that the right end-point of $I_j$ is equal to the left end-point of $I_{j+1}$. We say that two intervals $I_j\subset(-\infty,0]$ and $I_{j'}\subset[0,\infty)$ are \emph{adjacent} if $\etab(I_j)$ and $\etab(I_{j'})$ share a non-trivial boundary arc (i.e., a connected set with more than one point). Define
	\eqbn
	U^\Xo_{j}=\{ \Xb_t\,:\,t\in J, \Xb_t=\inf_{t'\in[t\wedge 0,t\vee 0]} \Xb_{t'} \},\qquad\text{for\,\,}(\Xo,\Xb)=(\Lo,\Lb),(\Ro,\Rb). 
	\eqen
	By the mating-of-trees construction in Section~\ref{sec:mot} (see also \cite[equation (1.3)]{ghs-dist-exponent}), if $I_j\subset(-\infty,0]$ and $I_{j'}\subset[0,\infty)$ are not adjacent then $U^\Lo_{j}$ and $U^{\Lo}_{j'}$ (respectively, $U^\Ro_{j}$ and $U^{\Ro}_{j'}$) have positive distance almost surely in the sense that 
	\eqbn
	\inf\{|x-x'|\,:\,x\in U^{\Xo}_j,x'\in U^{\Xo}_{j'} \}>0,
	\qquad\text{for\,\,}(\Xo,\Xb)=(\Lo,\Lb),(\Ro,\Rb). 
	\eqen
	For $(\Xo,\Xb)=(\Lo,\Lb),(\Ro,\Rb)$ define
	\eqbn
	\begin{split}
		U_{j}^{\Xo,n} = \{ 0.5\sqrt{\beta/n}\cdot (\xi^\Xo_-)^{-1}(v)\,:\,v\in V^{\Xo}_-(I_j) \}\qquad\text{for\,\,}&I_j\subset(-\infty,0],\\
		U_{j'}^{\Xo,n} = \{ 0.5\sqrt{\beta/n}\cdot (\xi^\Xo_+)^{-1}(v)\,:\,v\in V^{\Xo}_+(I_{j'}) \}\qquad\text{for\,\,}&I_{j'}\subset[0,\infty).
	\end{split}
	\eqen
	For $\wt\delta>0$ let $E_2=E_2(K,\ep,\delta,\wt\delta,n)$ be the event that for all $I_j\subset(-\infty,0]$ and $I_{j'}\subset[0,\infty)$ which are not adjacent, $U_{j}^{\Xo,n}$ and $U_{j'}^{\Xo,n}$ have distance greater than $\wt\delta$ for $\Xo=\Lo,\Ro$, that is,
	\eqb
	\inf\{|x-x'|\,:\,x\in U^{\Xo,n}_j,x'\in U^{\Xo,n}_{j'} \}>\wt\delta.
	\label{eq98}
	\eqe
	Then $\lim_{\wt\delta\rta 0}\lim_{n\rta\infty} \P[E_2]=1$. Fix $\wt\delta>0$ and $n_0\in\N$ such that $\P[E_2]>1-\ep/2$ for all $n\geq n_0$. 
	
	Assume $E_1$ and $E_2$ both occur. Let $e\in E(M_n)$ and $u_1\in V(M_n)$ be as in the statement of the lemma such that $t_e^0\in[-K,0)$ and $t_e^1\in(0,K]$. Let $j,j'$ be such that $t_e^0\in I_j$ and $t_e^1\in I_{j'}$. Then $u_1$ is on the boundary of the future wedge, while $e$ is an interior edge of the past wedge. Without loss of generality assume that $u_1$ is on the left boundary of the future wedge. Since the past and the future near-triangulations are ``glued together'' according to boundary length, $0.5\sqrt{\beta/n}\cdot(\xi_+^\Lo)^{-1}(u_1)\in U_j^{\Lo,n}\cap U_{j'}^{\Lo,n}$. In particular, the condition \eqref{eq98} is not satisfied, so by occurrence of $E_2$ the intervals $I_j$ and $I_{j'}$ are adjacent. By occurrence of $E_1$ this implies that $d_{\BB S^2}(\etab(t_e^0),\etab(t_e^1))< 2\cdot\ep/2$. Since $e$ and $u_1$ were arbitrary this implies that \eqref{eq95} is satisfied, so
	\eqbn
	\P\bigg[ \sup_{e\in E(M_n),t_e^0\in [-K,0),t^1_e\in (0,K] }
	d_{\BB S^2}(\etab(t^0_e),\etab(t^1_e))\geq \ep\bigg]
	\leq \P[ E_1^c\cup E_2^c ] \leq \ep.
	\eqen
	We have established \eqref{eq95}, which implies \eqref{eq93}.
	
	In order to prove \eqref{eq92}, first observe that for any $\eps>0$ and fixed $K>0$ we can find $K'>K$ such that with probability at least $1-\eps$, the following event $E$ occurs
	\eqbn
	\begin{split}
		E =&\, \bigg\{ 
		\Big(\inf_{t\in[-K',-K]} \Lb_t\Big)
		\vee
		\Big(\inf_{t\in[K,K']} \Lb_t\Big)
		< 
		\inf_{t\in[-K,K]} \Lb_t
		\bigg\} \\
		&\qquad\qquad\cap
		\bigg\{ 
		\Big(\inf_{t\in[-K',-K]} \Rb_t\Big)
		\vee
		\Big(\inf_{t\in[K,K']} \Rb_t\Big)
		< 
		\inf_{t\in[-K,K]} \Rb_t
		\bigg\} 
		.
	\end{split}
	\eqen 
	For $t_e^0$ and $t_e^1$ as above assume that $w_{3t_e^0}\in\{ a,b\}$ (as remarked above, if $w_{3t_e^0}=c$ then $t_e^0=t_e^1$, so this case is immediate when proving \eqref{eq92}). By the interpretation of $Z^n$ as a boundary length process and assuming without loss of generality that $u_1$ is on the left frontier of the map at time $t_e^1$, we have 
	$t_e^1=\inf\{ t>t_e^0\,:\,L^n_t<L^n_{t_e^0} \}$.
	By this identity, on the event $E$ and with probability converging to 1 as $n\rta\infty$, if $\{t_e^0,t_e^1 \}\cap[-K,K]\neq\emptyset$ then $\{t_e^0,t_e^1 \}\subset[-K',K']$. Therefore \eqref{eq92} follows from \eqref{eq93}. 
\end{proof}

\begin{proof}[Proof of Theorem \ref{thm1}(i-iii) and (v)]
	We verify the four assertions separately.
	
	(i) \nina{For any ball $D\subset\C$ the measure $\mu_n|_D$ renormalized to be a probability measure it tight since $D$ is bounded. Furthermore, the random variable $\mu_n(D)$ is tight by \eqref{eq:mu-def} and since $D\subset\etab([-T,T])$ for some random $T>0$. Therefore we know that $\mu_n$ converges subsequentially for the weak topology, and to conclude it is sufficient to establish uniqueness of subsequential limits.} 
	Since the set of space-filling SLE$_6$ segments $\etab([q_1,q_2])$, $q_1,q_2\in\Q$, is a $\Pi$-system that generates the Borel $\sigma$-algebra (since any open subset of $\C$ can be written as a countable union of SLE$_6$ segments), it is sufficient to show that \nina{for the coupling considered in the rest of this section,} for any fixed $q_1,q_2\in\Q$, we have $\mu_n(D)\rta\mub(D)$ almost surely for $D:=\etab([q_1,q_2])$. Defining 
	\eqbn
	X_n:=\#\Big\{ v\in V(M_n)\,:\,\frac{1}{3n}\etavf^{-1}(v)\in[q_1,q_2] \Big\},\qquad
	Y_n:=\#\Big\{ e\in E(M_n)\,:\,\frac{1}{3n}\etae^{-1}(e)\in[q_1,q_2] \Big\},
	\eqen
	in order to conclude it is sufficient to show that almost surely,
	\eqb\label{eq31}
	\lim_{n\rta\infty}\frac{X_n}{Y_n}\rta \frac 13,\qquad \textrm{and}\qquad
	\lim_{n\rta\infty} \frac{Y_n}{3n}\rta q_2-q_1.
	\eqe
	The second part of \eqref{eq31} is immediate by the definition of $\phi_n$. For the first part, note that $X_n$ (respectively, $Y_n$) represents the number of $c$-steps (respectively, steps) in a word of length $\Theta(n)$. By concentration for the sum of Bernoulli$(1/3)$ random variables, for any fixed $\eps>0$, 
	\begin{equation}\label{eq:con}
		\P[|X_n/Y_n-1/3|\ge \eps ] \textrm{ decays exponentially in $n$, with a rate depending on $\eps$.}
	\end{equation}

	(ii) Recall that we defined $\eta_n(t)=\phi_n(\etae(3tn))$ for $tn\in\Z$, where $\phi_n(e)$ for $e=\{u,v\}\in E(M_n)$ is the midpoint between $\phi_n(u)=\etab( (3n)^{-1}\etavf^{-1}(u) )$ and $\phi_n(v)=\etab( (3n)^{-1}\etavf^{-1}(v) )$. If we instead had defined $\eta_n(t)=\etab((3n)^{-1}\etae(tn))$ then the result would have been immediate by continuity of $\etab$. We deduce the result for our choice of $\phi_n$ by using Lemma \ref{prop6}. 
	
	(iii) Recall that in Section~\ref{sec:cont-cle} we showed that it is possible to define the local time $\ellb^u$, and its inverse $\Tb^u$ (defined as in \eqref{eq:Tu}) relative to a time $u$ defined as the envelope closing time of a $\CLE_6$ loop. For $j\in \N_{+}$, $\Tb^{u_j}$ is defined relative to the envelope closing time $u_j$ of $\gab_j$. Recall the definition of $T^{n,j}$ given in \eqref{eq:T-cle}. By Assumption \ref{a1} and Lemma \ref{prop21}, we get that $T^{n,j}\rta \Tb^{u_j}$ for the metric $d_{\op{p}}$ given by \eqref{eq35}. By an application of Fact~\ref{fact: curve}, we get that $\ga_1^n,\ga_2^n,\dots$ converge to $\gab_1,\gab_2,\dots$ as parametrized curves. By Lemmas \ref{prop10} and \ref{lem:env-cov}, $\area_n(\ga_j)\rta\area(\gab_j)$. 
	
	It remains to prove that $\1_{\ga_i^n\subset\reg(\ga_j^n)}$ converges to  $\1_{\gab_i\subset\reg(\gab_j)}$. Without loss of generality, assume that $\gab_j$ is oriented counterclockwise and let $s$ denote the quantum natural length of $\gab_j$. Consider the forested line relative to the envelope closing time $u$ of $\gab_j$ and let $\frk L$ be the looptree on this forested line for which $\gab_j$ is the outside cycle. We have $\gab_i\subset\reg(\gab_j)$ if and only if $\gab_i$ is contained inside one of the bubbles of $\frk L$; equivalently, if and only if $\env(\gab_i)\subset J$ for $J$ a complementary component of $\op{range}(\Tb^u|_{[-s,0]})$ (so $\etab(J)$ is a bubble of $\frk L$).
	With probability $1-o_n(1)$ the description of the discrete event $\ga_i^n\subset\reg(\ga_j^n)$ in terms of $\env_n(\ga^n_i)$ and $T^{u,n}$ is the same, which implies that $\1_{\ga_i^n\subset\reg(\ga_j^n)}\rta\1_{\gab_i\subset\reg(\gab_j)}$.
	
	(v) The existence of tuples $(e^n_1,\ldots,e^n_k)$ and $P=(z_1,\dots,z_k)$ with the correct marginal laws such that $P_n=(\phi_n(e^n_1),\ldots,\phi_n(e^n_k))$ converges to $P$ is immediate by (i) and (iii). Fix $j\in[k]$, and let $u_{n,j}=(3n)^{-1} \etavf^{-1}(e^n_j)$. Then $u_{n,j}$ converge almost surely to $u_j:=\sup\{ t\in\R\,:\, \etab(t)=z_j \}$.
	Let $T^{n}$ and $\Tb$ be as in \eqref{eq:Levy-walk} and Section~\ref{sec:dictionary-branch}, respectively, where we have recentered at the time $u_{n,j}$ and $u_j$, and we do not indicate the $j$ dependence to simplify notation. Then we have convergence of $T^n$ to $\Tb$ in probability by Lemma~\ref{prop5}. Combining this result with Fact~\ref{fact: curve}, \eqref{eq:eta-def}, and the 
	convergence of $\eta_n$, gives convergence of the branch from $\infty$ to each point of $P_n\cup\{0 \}$.

	To conclude we need to argue that the path of $\tau^*_n$ between any pair of points $\phi_n(e_i^n),\phi_n(e_j^n)\in P_n$ converges. This does not follow from the preceding paragraph, since we do not know that the point where the two branches from infinity merge, also converges. For $z\in\C$, let $t_{z}=\sup\{t\in\R\,:\,\etab(t)=z \}$, and assume $t_{z_i}<t_{z_j}$. In the continuum there are disjoint time sets $U,U_i,U_j\subset\R$ such that $U\cup U_i=\ans(t_{z_i})$ and $U\cup U_j=\ans(t_{z_j})$, and such that $\sup(U)\leq\inf(U_i)<\sup(U_i)\leq \inf(U_j)$. By Theorems \ref{thm:spinelt-inf} and \ref{thm:dfs-inf} (see also \eqref{eq:branching-times}) there are time sets $U^n,U^n_i,U^n_j\subset\R$ with the same properties in the discrete setting. In particular, $t\in U^n\cup U^{n}_i$ (respectively, $t\in U^n\cup U^{n}_j$) if and only if the dual of $\eta_{\op{e}}(3nt)$ is on the branch of $\tau^*_n$ from $\infty$ to $\phi_n(e_i^n)$ (respectively, $\phi_n(e_j^n)$). By convergence of the process $T^n$ defined relative to $u_{n,i}$ and $u_{n,j}$, respectively, counting measure on $U^n,U^n_i$, and $U^n_j$, converges to the local time of $U,U_i$, and $U_j$, respectively. Using this we get convergence of the branch in $\tau^*_n$ connecting $\phi_n(e_i^n)$ and $\phi_n(e_j^n)$, since this branch is the concatenation of the branch from $\phi_n(e_i^n)$ to $\eta_n(\sup (U^n))$ of length $n^{-3/4}|U^n_i|$, and the branch from $\eta_n(\sup (U^n))$ to $\phi_n(e_j^n)$ of length $n^{-3/4}|U^n_j|$.
\end{proof}

Theorem \ref{thm1}(iv), which has not yet been proved, is immediate by the following lemma. 
Recall the notation $\Gab_z$ introduced in Section \ref{sec:cont-piv} for a pivotal point $z$. We adopt the analogous notation $\Ga^n_{z_n}$, where  $z_n\in \C$ is the image by $\phi_n$ of a pivotal point of $(M_n,\si_n)$.
Convergence of 	 $\Ga^n_{z_n}=(\wh{\ga}^n_j)_j$ to $\Gab_{z}=(\wh{\gab}_j)_j$ means uniform convergence of each percolation cycle $\wh{\ga}^n_j$ to $\wh{\gab}_j$, where percolation cycles are viewed as parametrized curves in $\C$ and the percolation cycles of $\Ga^n_{z_n}$ and $\Gab_{z}$ are ordered by their value.
\begin{lemma}\label{prop34}
	Assertion (iv) of Theorem \ref{thm1} holds, that is, for any fixed $\ep>0$ and $i,j\in\N_+$, the pivotal measures $\nu^{\ep,1}_{j,n}$, $\nu^{\ep,2}_{j,n}$, $\nu_{i,j,n}$, and $\nu^{\ep}_{n}$ converge \nina{in probability} for the weak topology to $\nub^{\ep,1}_{j}$, $\nub^{\ep,2}_{j}$, $\nub_{i,j}$, and $\nub^{\ep}$, respectively. \nina{Furthermore, $\1_{\nu^{\ep,1}_{j,n}(\C)=0}$ converges in probability to $\1_{\nub^{\ep,1}_{j}(\C)=0}$, and the analogous statement holds for the three other measures.}
	
	For $j\in\N_+$ and $\eps>0$, let $z_n\in\C$ (respectively, $z\in\C$) have the law of a uniformly sampled $\eps n$-pivotal point associated with $\ga^n_j$ (respectively, $\gab_j$). There exists a coupling of $z_n$ and $z$ such that $z_n\rta z$ almost surely.
	In such a coupling let $\Ga^n_{z_n}$ (respectively, $\Gab_z$) be the collection of percolation cycles after flipping the color of $z_n$ (respectively, $z$). Then $\Ga^n_{z_n}$ converges in probability to $\Gab_z$.
\end{lemma}
\begin{proof}
	Throughout the proof we identify points $\phi_n(v)\in\C$ with the associated vertex $v$ of $M_n$. Let $u$ (respectively, $u_n$) denote the envelope closing time of $\gab_j$ (respectively, $\ga_j^n$), and let $q\in(-\infty,0)\cap\Q$ be such that $|q|$ is at most the quantum natural length of $\gab_j$. Let $\nub(q,u)$ be the pushforward of the measure $\pb(q,u)$ defined in \eqref{eq:cont-piv-measure}, and recall that $\nu_n(s,u_n)$ is the pushforward of $p_n(s,u_n)$ under $\phi_n\circ \wh\eta_{\op{v}}(n^{3/4} \, \cdot)$. By Lemmas \ref{prop21} and \ref{lem:env-cov}, we can find a coupling where $y_n\in \C$ (respectively, $y\in \C$) is sampled from $\nu_n(q,u_n)$ (respectively, $\nub(q,u)$) renormalized to be a probability measure, such that $y_n$ converges to $y$ almost surely. By Lemma \ref{prop19}(i), the pivotal point $y_n$ is an $\eps n$-pivotal point with probability $1-o_\eps(1)$. Let $\frk L=\frk L(\ga^n_j)$ be the looptree associated with $\ga^n_j$ (Definition \ref{def:lt-cluster}). Upon taking a subsequence we may assume that all the pivotal points $y_n$ are of the same type, and that the cluster for which $\ga_j^n$ is the outside-cycle is always of the same color. In the remainder of this paragraph (except at the very end) all convergence statements concern convergence along this subsequence. Assume the cluster for which $\ga^n_j$ is the outside-cycle is white, and that $y_n$ is a vertex on this looptree (equivalently, $y_n$ is of type 1); the other cases will be discussed in the next paragraph. Under these assumptions, there is a finite collection of white looptrees $\frk L'_1,\frk L'_2,\dots$ 
	such that if $i_1\neq i_2$ then the only vertex on both $\frk L'_{i_1}$ and $\frk L'_{i_2}$ is $y_n$, 
	the union of the vertices of $\frk L'_1,\frk L'_2,\dots$ is the set of vertices of $\frk L$, and
	for each $\frk L'_i$ there is a unique bubble $B_i$ which has $y_n$ on its boundary. By \cite[Theorem 2]{camia-newman-full} (which we recalled in Lemma \ref{lem:pivotal-continuous}), the limiting CLE$_6$ has no triple points, so except for at most two values of $i$ the outside-cycle of $\frk L'_i$ has diameter $o_n(1)$ when embedded into $\C$. Since the $\sqrt{8/3}$-LQG area measure $\mub$ has no atoms, if the outside-cycle of a looptree has diameter $o_n(1)$ when embedded into $\C$ then the looptree encloses less than $\eps n$ vertices for all sufficiently large $n$. Therefore, with probability $1-o_n(1)-o_\eps(1)$ (with $o_\eps(1)$ uniform in $n$) exactly two of the looptrees $\frk L'_1,\frk L'_2,\dots$ (say, $\frk L'_1$ and $\frk L'_2$) enclose at least $\eps n$ vertices and have diameter at least $\ep$ when embedded into $\C$, while each of the other looptrees enclose $o_n(n)$ vertices and have diameter $o_n(1)$. 
	The bubbles $B_1$ and $B_2$ enclose area $o_n(1)$ by Lemma \ref{prop19}(iii). It follows that with probability $1-o_n(1)-o_\eps(1)$, $\cL_{y_n}\cap\Gam^n_{y_n}$ contains exactly two percolation cycles $\ga'_n$ and $\ga''_n$ enclosing at least $\eps n$ vertices (since they enclose at least $\frk L'_i\setminus B_i$). Since $\etab$ is a continuous curve when parametrized by LQG area measure, and since the edges in $B_1$ (resp.\ $B_2$) are visited consecutively in the space-filling exploration of $M_n$, we see see that the diameter of $B_1$ and $B_2$ is $o_n(1)$ when embedded into $\C$. Combining the above, the two percolation cycles $\ga'_n$ and $\ga''_n$ converge to non-trivial percolation cycles $\gab'$ and $\gab''$ in the scaling limit, such that $\gab$ is the concatenation of $\gab'$ and $\gab''$ and $y$ is the common point of $\gab'$ and $\gab''$. 
	In particular, $\Gam^n_{y_n}$ converges to $\Gab_{y}$ in the scaling limit, at least along the considered subsequence. Since $\Gab_y$ is a deterministic function of $\Gab$ and $y$ this implies that we do not only have subsequential convergence of $\Gam^n_{y_n}$ to $\Gab_{y}$, but convergence along $n\in\N$. Furthermore, $y$ has the same type 1-4 as $y_n$ with probability $1-o_n(1)$.
	
	Now we will discuss the case of pivotal points of types 2-4. Again we may assume without loss of generality that $y_n$ is white, but now $y_n$ is \emph{not} a vertex on the cluster for which $\gamma_j^n$ is the outside cycle, so this cluster is black. Again, locally near $y_n$ we will see a picture similar to that illustrated in Figure \ref{fig:pivotal-discrete} (except for the case of pivotal points of type 2, where there will be one instead of two black looptrees on the black forested line which contain a vertex adjacent to $y_n$). Let $\frk L$ be the white looptree on the white forested line relative to time $u_n$ which contains $y_n$, and (as before) let $\frk L'_1,\frk L'_2,\dots$ be the white looptrees such that if $i_1\neq i_2$ then the only vertex which is on both $\frk L_{i_1}$ and $\frk L_{i_2}$ is $y_n$, and for each $\frk L'_i$ there is a unique bubble $B_i$ which has $y_n$ on its boundary. As before, \cite[Theorem 2]{camia-newman-full} implies that all except two of these looptrees have diameter $o_n(1)$, since there will be one or two percolation cycle(s) which trace the interface between the white and black looptrees nearby $y_n$ on the forested line; if $y_n$ is of type 2 we use that no CLE$_6$ loops have triple points, and if $y_n$ is of type 3 or 4 we use that there is no point which is both a double point of a CLE$_6$ loop and a point where two CLE$_6$ loops meet. We can now complete the argument similarly as in the paragraph above.

	By Theorem \ref{thm1}(i) and since the $\mub$-mass of the $\delta$-neighborhood of any CLE$_6$ loop converges to 0 as $\delta\rta 0$, we see that the if $y_n$ is $\eps n$-pivotal for arbitrarily large $n$ then $y$ is almost surely an $\eps$-pivotal, and vice versa. We have proved that if $E_n(y_n)$ (respectively, $E(y)$) denotes the event that $y_n$ (respectively, $y$) is an $\eps$-pivotal of type 1 associated with $\ga^n_j$ (respectively, $\gab_j$), then the pair $(y_n,\1_{E_n(y_n)})$ converge jointly to the pair $(y,\1_{E(y)})$ in the scaling limit. This implies further that if $\wh\nu_{j,n}^{\eps,1}(q,u_n)$ (respectively, $\wh\nub_{j}^{\eps,1}(q,u)$) denotes $\nu(q,u)$ (respectively, $\nub(q,u)$) restricted to the set of $\eps$-pivotal points of type 1, then $\wh\nu_{j,n}^{\eps,1}(q,u_n)$ converges to $\wh\nub_{j}^{\eps,1}(q,u)$ in the scaling limit. By Lemma \ref{prop30} we know that the $\eps$-pivotal points associated with $\gab_j$ can be covered by finitely many sets $\Ab(q_0,u),\dots,\Ab(q_k,u)$ (with $q_0$ denoting the quantum natural length of $\gab_j$ and the other $q'_i$s as above). By Lemma \ref{prop16} the same result holds in the discrete for the same choice of times $q_i$ (if we ignore the sets $S_{\ga^n_j,\ga'}$, which have negligible mass in the scaling limit). Defining
	\eqb
	\wh \nub^{\eps,1}_j = \sum \wh\nub_{j}^{\eps,1}(q_i,u),
	\label{eq97}
	\eqe
	the measure $\nub^{\eps,1}_j$ can be obtained from $\wh \nub^{\eps,1}_j$ by reweighing by a Radon-Nikodym derivative $\wh\etab^{\etab(u)}\circ f$ for $f:[-q_0,0]\to [0,\infty)$. The function $f$ has a particularly simple form, since it can be chosen to be piecewise constant with values strictly between $1/(k+1)$ and 1. More precisely, for each $i$ and $i'$ such that $q_i<q_{i'}$ the interval $[-q_i,0]$ can be written as a union of three closed intervals with disjoint interior $J_{i},J_{i'},J_{i,i'}$ (ordered by increasing left endpoint), such that $\wh\nub_{j}^{\eps,1}(q_i,u)$ is supported on the image of $J_i\cup J_{i,i'}$ under $\wh\etab^{\etab(u)}$, and $f(t)$ is equal to the inverse of the number of $i$ such that $t\in J_i\cup J_{i,i'}$. This description follows by using that the pullback of $\wh \nub^{\eps,1}_j$ by $\wh\etab^{\etab(u)}$ is supported on either the set of local running infima for $\wh \Lb$ or the set of local running infima for $\wh\Rb$. We can describe $\nub^{\eps,1}_{j,n}$ in the exact same way in the discrete setting by reweighing a measure $\wh\nub^{\eps,1}_{j,n}$ of the form \eqref{eq97} by a Radon-Nikodym derivative $\phi_n\circ f_n$. By using that $f_n\rta f$ and $\wh\nu_{j,n}^{\eps,1}(q_i,u_n)\rta \wh\nub_{j}^{\eps,1}(q_i,u)$ we get that $\nu^{\eps,1}_{j,n}\rta \nub^{\eps,1}_j$. 
	
	\nina{Next we will} argue that $\nu^{\ep}_{n}$ converges in the weak topology to $\nub^{\ep}$. Let 
	$\nu^\eps_{i,j,n}$ (respectively, $\nub^\eps_{i,j}$) denote 
	$\nu_{i,j,n}$ (respectively, $\nub_{i,j}$) restricted to the set of $\eps n$-pivotals (respectively $\eps$-pivotals). By Theorem \ref{thm1}(i,iii), since $\nu_{i,j,n}$ converges weakly to $\nub_{i,j}$, and since the $\mub$-mass of the $\delta$-neighborhood of any CLE$_6$ loop converges to 0 as $\delta\rta 0$, we see that $\nu^\eps_{i,j,n}$ converges weakly to $\nub^\eps_{i,j}$. We conclude by using that $\nu^{\ep}_{n}$ (respectively, $\nub^{\ep}$) is the sum of the measures 
	$\nu^{\eps,1}_{j,n},\nu^{\eps,2}_{j,n},\nu^{\eps}_{i,j,n}$ (respectively, 
	$\nub^{\eps,1}_{j},\nub^{\eps,2}_{j},\nub^\eps_{i,j}$) summed over all $i,j\in\N_+$.
	
	\nina{It remains to argue convergence of the event $\nu^{\ep,1}_{j,n}(\C)=0$ to the event $\nub^{\ep,1}_{j}(\C)=0$, and the analogous statement for the other two measures. We first consider the case of the measure $\nu^{\ep,1}_{j,n}$ and start by arguing that if $\nub^{\ep,1}_{j}(\C)\neq 0$ then $\nu^{\ep,1}_{j,n}(\C)\neq 0$ almost surely for all sufficiently large $n$. In the notation of the first three paragraphs of the proof, if $\nub^{\ep,1}_{j}(\C)\neq 0$ then there exists some $q\in\{q_1,\dots,q_k \}$ such that $E(y)$ occurs with positive probability for $y$ sampled from $\nub(q,u)$. Since we have convergence of $(y_n,\1_{E_n(y_n)})$ to $(y,\1_{E(y)})$, we see that $E_n(y_n)$ also happens with uniformly positive probability for all sufficiently small $n$, so $\nu^{\ep,1}_{j,n}(\C)\neq 0$ almost surely for all sufficiently large $n$. 
		
		To conclude we need to show that if there exists arbitrarily large $n$ such that $\nu^{\ep,1}_{j,n}(\C)\neq 0$ then we have $\nub^{\ep,1}_{j}(\C)\neq 0$ almost surely. Let $\cN\subset\N_+$ be an infinite set such that $\nu^{\ep,1}_{j,n}(\C)\neq 0$ for all $n\in\cN$, and for each $n\in\cN$ let $z_n$ be sampled from $\nu^{\ep,1}_{j,n}$. By compactness, $z_n$ converges along a subsequence $\cN'\subset\cN$ to some $z\in\C$. By convergence of $\gamma_j^n$ and $\mu_n$ to $\gab_j$ and $\mub_{\gff}$, respectively, we get that $z$ is almost surely an $\eps$-pivotal. In particular, we see that the set of $\eps$-pivotals of type 1 for $\gab_j$ is non-empty, and to conclude it is sufficient to show that this implies a.s.\ positivity of the total mass of $\nu^{\ep,1}_{j,n}$. 
		
		We first observe the following: If a $3/2$-stable L\'evy processes with only negative jumps has a running infimum in some open interval $I$, then the local time at its running infimum in $I$ is positive a.s. This observation follows for example by stopping the L\'evy process the first time $t$ it reaches a running infimum in $I$ and applying a 0-1 law to argue that the local time is a.s.\ positive in any interval of the form $[t,t+s]$ for $s>0$. 
		
		Recall from Section \ref{sec:cont-cle} that $\gab_j$ is associated with a looptree $\frk L(\gab_j)$.
		By the observation made about L\'evy processes in the previous paragraph, it is sufficient to show that if $(\Xb_s)_{s\in[0,T]}$ is the L\'evy excursion encoding $\frk L(\gab_j)$ then there is some open interval $I\subset[0,T]$ such that all running infima of $\Xb$ relative to the left end-point of $I$ correspond to $\eps$-pivotal points for $\gab_j$. 
		Note that if $s$ is the second time at which a pivotal point is visited (equivalently, $s$ is a local running infimum for $\Xb$), then $\frk L(\gab_j)$ is split into two looptrees if we cut at this pivotal point: one looptree $\frk L_s$ encoded by an excursion $\Xb|_{[s',s]}$ (where $s'=\sup\{ t\leq s\,:\,\Xb_t<\Xb_s \}$) and one looptree $\frk L'_s$ which contains the root of the original looptree $\frk L(\gab_j)$ such that $\frk L(\gab_j)$ can be obtained by concatenating $\frk L_s$ and $\frk L'_s$.
		Let $(\cF_s)_{s\in[0,T]}$ denote the filtration such that $\cF_s$ is the smallest $\sigma$-algebra containing information about $\Xb|_{[0,s]}$ along with the LQG disks associated with the jumps of $\Xb|_{[0,s]}$. Then $\Xb$ is a Markov process for the filtration $(\cF_s)$.
		The event that $\frk L_s$ has area at least $\eps$ (i.e., the event that the total LQG area of the LQG disks of $\frk L_s$ is at least $\eps$) is measurable with respect to $\cF_s$. For a fixed rational $q>0$ let $\tau$ be the stopping time given by the first time $s>q$ for which $\frk L_s$ has area at least equal to $\eps$. By the Markov property of $\Xb$, the process $\Xb$ has running infima relative to time $\tau$ in $[\tau,\tau+s']$ for any $s'>0$, and for such a running infimum $s''\in[\tau,\tau+s']$ the looptree $\frk L_{s''}$ (resp.\ $\frk L'_{s''}$) has an area which is $o_{s'}(1)$ larger (resp.\ smaller) than $\frk L_\tau$. The area of $\frk L'_\tau$ is exactly equal to $\eps$ with probability 0; this can e.g.\ be seen since the total area of an LQG disk has a continuous density. It follows that if $\tau$ corresponds to an $\eps$-pivotal point then the running infima in $[\tau,\tau+s']$ are also $\eps$-pivotal points for sufficiently small $s'$, so we can conclude the proof by setting $I=(\tau,\tau+s')$. By Lemma \ref{lem:cover} (alternatively, Lemma \ref{prop30}), by varying $q$ we can guarantee that if $\gab_j$ has an $\eps$-pivotal point of type 1 then this will be captured by the procedure we just described.
		
		We have concluded the proof for the case of the measure $\nu^{\ep,1}_{j,n}$, and will consider the measures $\nu^{\ep,2}_{j,n}$, $\nu_{i,j,n}$, and $\nu^{\ep}_{n}$ now. We will only point out the steps where the proof differs from the case of $\nu^{\ep,1}_{j,n}$. The case of $\nu^{\ep,2}_{j,n}$ is very similar to the case of $\nu^{\ep,1}_{j,n}$; the main difference is that instead of considering the looptree $\frk L(\gab_j)$ encoded by a L\'evy excursion $\Xb$ starting and ending at 0, we have a process $\Xb$ starting at a positive value and staying positive until it hits 0 (see Sections \ref{sec:cont-cle} and \ref{sec:cont-piv} and the right part of Figure \ref{fig:past-wedge}). 
		
		To treat the cases of $\nu_{i,j,n}$ and $\nu^{\ep}_{n}$ it is sufficient to treat the measure $\nu^{\ep}_{i,j,n}$ defined above, i.e., it is sufficient to show that $\1_{\nu^{\ep}_{i,j,n}(\C)}$ converges to $\1_{\nub^{\ep}_{i,j}(\C)}$. Notice that to determine whether a point of intersection between two loops is an $\eps$-pivotal or not it is sufficient to consider the area enclosed by each of the two loops and the difference between the two areas. In particular, either all or none of the points of intersection between the two loops are $\eps$-pivotals. By assumption, in our case there is at least one $\eps$-pivotal point in the intersection between $\gab_i$ and $\gab_j$, so all the points of intersection must be $\eps$-pivotal points. We also observe that there must be more than a single point of intersection between $\gab_i$ and $\gab_j$ since a single point of intersection corresponds to having a 6-arm event, which has exponent $35/12>2$ \cite{smirnov-werner-percolation} and therefore a.s.\ does not occur. Combining this with the fact that pivotal points of types 3 and 4 are (global) running infima of a $3/2$-stable L\'evy process $\Xb$ starting at a positive value and staying positive until it hits 0, we get the existence of an interval $I$ such that all running infima of $\Xb$ during this interval correspond to pivotal points between $\gab_i$ and $\gab_j$, which allows us to conclude the proof as in the case of pivotals of type 1.}
\end{proof}

\subsection{Finite volume case}
\label{sec:convfin} 
The proof of Theorem~\ref{thm:finite} is very similar to the proof of Theorem \ref{thm1} modulo two new inputs: Proposition~\ref{prop3} and Lemma~\ref{lem:chordal-conv}. Therefore, in this section we will only elaborate on the differences. We also omit the proof of Theorem \ref{thm:sphere} as it can be proved in the exact same way as Theorem~\ref{thm:finite}.

Recall the setting of Section~\ref{sec:finite} for Theorem~\ref{thm:finite}, where we study $(M_n,\sigma_n)$, $Z^n$, and $(\gffd,\etad,\Zd)$.
Let $Z=(L,R)$ be the unscaled lattice walk related to $Z^n$ by \eqref{eq:Zdn}. Then $Z$ is associated with a word $w\in \smK$. Here we still adopt the convention of hiding the dependence on $n$ for unscaled discrete quantities. 
Let $\frk v_n$ be a uniformly sampled inner vertex of $M_n$ and let $u_n$ be such that $\etavf(3nu_n)=\frk v_n$. 
Let $\ub$ be the almost surely unique time such that $\etad(\ub)=0$. 
The following lemma ensures the existence of the coupling in Theorem~\ref{thm:finite}.
\begin{lemma}\label{lem:uniform}
	The pair $(Z^n,u_n)$ converges to $(\Zd,\ub)$ in law.
\end{lemma}
\begin{proof}
	We first claim that the analog of \eqref{eq31} holds in the finite volume setting. The second part of \eqref{eq31} is verified exactly as before. To show the the first part, let $E_n=\{Z_{3n+2h_n}=(0,-h_n), Z_k\in [0,\infty)^2\;\forall k\in \{0,\dots,3n+2h_n\} \}$. Under the law of the infinite volume setting, where $Z$ has independent and identically distributed increments, $\P[E_n]$ decays polynomially in $n$. Therefore \eqref{eq:con} still holds in the finite volume setting since the two settings differ by a conditioning on $E_n$, while the probability in \eqref{eq:con} is exponentially small.
	
	Recall that $\ub$ is uniform in $(0,\frk m)$ and independent of $\Zd$. Let $\wh u_n$ be sampled independently from $(M_n,\sigma_n)$ such that $3n\wh u_n$ is uniform in $\{1,\dots,3n+2h_n \}$. Let $\wh e_n=\etae(3n\wh u_n)$. Note that $\wh e^n$ is an edge sampled uniformly at random from $M_n$ (except that the top-edge cannot be sampled). Conditioning on $M_n,\sigma_n, \wh e_n$, uniformly pick one of the two endpoints of $\wh e_n$ and denote it by $\wh v_n$. Then $\wh v_n$ is sampled from the uniform measure on the vertex set of $M_n$ weighed by the degree (where we do not count the top-edge when considering the degree of its two end-points). 
	
	By the finite volume analog of \eqref{eq31}, the Prokhorov distance between $(Z^n,u_n)$ and $(Z^n,\wh u^n)$ is $1-o_n(1)$.
	Since $(Z^n,\wh u_n)$ converge to $(\Zd,\ub)$ in law, we are done. 
\end{proof}

As in the UIPT case, the area convergence part of Assertion (i) and Assertion (ii) in Theorem~\ref{thm:finite} are almost immediate by the definition of $\phi_n$. To prove the former statement, we proceed as in the infinite volume case, and note that the analog of \eqref{eq31} holds as explained in the proof of Lemma \ref{lem:uniform}. To prove Assertion (ii), first observe that Lemma \ref{prop6} also holds in the finite volume setting. Indeed, the proof works in the same way once we use Theorem~\ref{thm:mot-disk} instead of Theorem~\ref{thm:mot}. Lemma~\ref{prop6} implies convergence of the space-filling percolation exploration as in the infinite volume case.

To prove the rest of Theorem~\ref{thm:finite}, we will use the following result, which first appeared in \cite[Proposition 5.19]{lsw-schnyder-wood}. 
The proposition allows to transfer properties of infinite volume maps to finite volume maps. In \cite{lsw-schnyder-wood} the result was stated for general centered random walks whose increments have some finite exponential moment, but for us it is sufficient to consider walks with increments $(1,0)$, $(0,1)$, and $(-1,-1)$. Recall that $\frk m>0$ denotes the $\sqrt{8/3}$-LQG area of the $\sqrt{8/3}$-LQG disk. Let $\cP^n$ be the law of $Z^n|_{[0,\frk m]}$ with $Z^n_0=(1,0)$ and independent and identically distributed steps, and let $\cP$ be the law of $\Zb|_{[0,\frk m]}$ for $\Zb$ a planar Brownian motion with correlation $1/2$ started from $\Zb_0=(1,0)$. For $\ep\geq 0$ let $E^{n,\ep}$ be the event that the walk stays in $(-\ep,\infty)^2$ for $t\in(0,1)$ and $|Z_{\frk m}^n|\leq \ep$, and define $\cP^{n,\ep}$ to be $\cP^n$ given $E^{n,\ep}$. Define $\cP^{\ep}$ similarly in the continuum, observing that $\cP^{0}$ may be defined as the limit of $\cP^{\ep}$ as $\eps\rta 0$. Then $\cP^{n,0}$ converges weakly to $\cP^0$.
\begin{prop}[\cite{lsw-schnyder-wood}]
	\nina{Fix $\xi\in(0,1/2)$. Suppose $Y^n$ (respectively, $\Yb$) is a random variable which is measurable with respect to $Z^n|_{[\xi,\frk m-\xi]}$ (respectively, $\Zb|_{[\xi,\frk m-\xi]}$). If the $\cP^n$-law of $(Z^n,Y^n)$ converges weakly to the $\cP$-law of $(\Zb,\Yb)$, then the $\cP^{n,0}$-law of $(Z^n,Y^n)$ converges weakly to the $\cP^0$-law of $(\Zb,\Yb)$.}
	\label{prop3}
\end{prop}
\begin{lemma}
	Assertions (iii) and (iv) in Theorem~\ref{thm:finite} hold.
\end{lemma}
\begin{proof}
	(iii): We will argue convergence in probability of the processes $T^{j,n}$. This is sufficient to conclude, since convergence of $T^{j,n}$ implies convergence of $\ga_j^n$ by the same argument as in the infinite volume case by Fact~\ref{fact: curve}. Fix $m\in\N_+$ and $\eps>0$. It is sufficient to prove the result only for $j=1,\dots,m$ and only on an event of probability $1-\eps$. 
	
	Let $\delta>0$ be such that with probability at least $1-\eps/10$ we have $\area(\gab_m)>2\delta$. If we sample a point uniformly at random from the $\sqrt{8/3}$-LQG area measure there will be a CLE$_6$ loop (in fact, infinitely many loops) surrounding this point almost surely. Therefore we can find $m'\in\N_+$ such that except on an event of probability $\eps/10$ the set of points which are \emph{not} enclosed by $\gab_1,\dots,\gab_{m'}$ have $\sqrt{8/3}$-LQG area at most $\delta$. Then let $\xi>0$ be such that except on an event of probability $\eps/10$ the envelope intervals of the percolation cycles $\gab_1,\dots,\gab_{m'}$ are contained in $[2\xi,\frk m-2\xi]$. 
	
	For $j=1,\dots,m'$ we can find intervals $I_j$ with rational endpoints such that $\env(I_j)$ is the envelope interval of $\gab_j$. Envelope intervals are encoded by the random walk and the Brownian excursion in a local way, in the sense that for two fixed intervals $I\subset J$ the event
	$\env(I)\subset J$ is measurable with respect to $\Zb|_J$, and the event $\env_n(I)\subset J$ is measurable with respect to $Z^n|_J$ with probability $1-o_n(1)$. By using this, Lemma \ref{prop2}, and Proposition \ref{prop3}, the intervals $\env_n(I_j)$ converge to the intervals $\env(I_j)$ in probability. This result, Lemma \ref{prop21}, and Proposition \ref{prop3} imply we can find $k^n_1,\dots,k^n_{m'}\in\N_+$ such that $T^{k_j,n}$ converges to $\Tb^{j}$ in probability for $j=1,\dots,m'$. To conclude the proof it is sufficient to show that with probability to converging to 1 we have $k_j=j$ for $j=1,\dots,m$. By convergence of $\mu_n$ to $\mub$ and convergence of $\area_n(\ga_{k_j})$ to $\area(\gab_j)$ for $j=1,\dots,m'$, with probability at least $1-3\cdot\eps/10$ for sufficiently large $n$, the number of vertices which are \emph{not} enclosed by any of the percolation cycles $\ga_{k_j}$ is at most $1.1\delta n$, and the percolation cycles $\ga_{k_j}$ for $j=1,\dots,m$ all enclose at least $1.9\delta n$ vertices. Since the percolation cycles are ordered by enclosed area, this implies that $k_1,\dots,k_m\in\{1,\dots,m' \}$. By Proposition \ref{prop3}, $\area_n(\ga_{k_j})$ converges to $\area(\gab_j)$ for $j=1,\dots,m'$. Therefore we have $k_j=j$ for $j=1,\dots,m$.
	
	(iv): We proceed similarly as in (iii). First we use Proposition \ref{prop3} to obtain the finite volume version of Lemmas~\ref{prop21} and~\ref{prop16}, which can be viewed as statements about random walk. This further implies the convergence of the pivotal measures $\nu^{\eps,1}_{j,n}$, $\nu^{\eps,2}_{j,n}$, and $\nu_{i,j,n}$. 
\end{proof}

By Corollary~\ref{cor:chordal-case} and Remark~\ref{rmk8}, monocolored outer edges on the percolated disk $(M_n,\sigma_n)$ are in one-to-one correspondence to $c$-steps in $w$ that only have an $a$-match. For $k\in \{1,\dots,h_n \}$, let $\lambda(k)$ be the index of the $c$-step corresponding to the $k$-th monocolored outer edge when enumerating the edges in clockwise order around $\partial M_n$. Then by Fact \ref{fact:c},
\[
\lambda(k) = \inf\{m\ge 0 : R_m < k\}.
\]
For $s\in [0,1]$, let 
\begin{equation}\label{eq:boundary}
	\lambda^n(s)=(3n)^{-1}\lambda(\lfloor s h_n\rfloor)\qquad \textrm{and}\qquad \bm{\lambda}(s) = \inf\{t\ge 0: \Rd_t\le s\}.
\end{equation}
Since $Z^n$ converges to $\Zd$, by Donsker's invariance principle, if we uniformly sample $U\in (0,1)$, then $\lambda^n(U)$ converges to $\bm{\lambda}(U)$ in probability. 
This implies that the boundary measure $\nu_n$ on $M_n$ converges to $\nub_{\gffd}$ in probability. 

In order to prove Assertions (v) and (vi) in Theorem~\ref{thm:finite}, in addition to Proposition~\ref{prop3}, we need to deal with the branches of the DFS trees whose terminal points are at the boundary of $M_n$. In particular, we need to rule out pathological behaviors at the endpoints. This is essentially done in Lemma~\ref{prop22}. Let us explain it in more detail.

Fix constants $\ell,r>0$ such that $\ell+r=1$. Let $\wt e_n$ be the $\lceil\ell h_n\rceil$-th outer edge when clockwise tracing the outer edges from the root edge $\frk e_n$. (The edge $\wt e_n$ is the same as the edge $a_{3}$ defined in Section~\ref{sec:finite},
but we use $\wt e_n$ since our discussion is meant to be for a generic branch of $\tau^*_n$ terminating at the boundary.) Let $u$ be such that $3n u$ is the index of the $\lceil\ell h_n\rceil$-th $c$-step in $w$ which does not have a $b$-match. 

Let us consider the future/past decomposition of $w$ relative to $u$. Write $w=w^-w^+$ where the first step in $w^+$ is $w_{3nu}$. We denote by $\dots,T(-2),T(-1),T(0)\in\ZZ^{<3nu}$ the times associated with the spine steps of $w^{-}$ so that $\pi(w^{-})=\dots w_{T(-2)}w_{T(-1)}w_{T(0)}$. For $m\in\ZZ^{\leq 0}$, let $\ell_m:=-\#\{k\in\Z^{\leq 0}~|~T(k)\geq m\}$. Finally, let $\wh Z_m=Z_{T(m)}+Z_{T(0)}$ where the shift by $Z_{T(0)}$ is to be consistent with \eqref{eq:Zdn}. 
We also define the rescaled versions of $T,\ell,\wh Z$ as in \eqref{eq:Levy-walk} and denote them by $T^n,\ell^n,\wh Z^n$. 

In the continuum, let $\wt x\in \bdy \D$ be such that the $\nub_{\gffd}$-length of the clockwise arc from $1$ to $\wt x$ equals $\ell$. 
Let $\wt u$ be the almost surely unique time such that $\etad(\wt u) =\wt x$.
Let $\ellb$, $\Tb,$ and $\wh \Zb$ be the local time, inverse local time, and L\'evy process, respectively, relative to $\wt u$. (Again $\wt x$ and $\wt u$ are just $A_3$ and $t_3$ but for the same reason as above we relabel them.) 
Now Assertions (v) and (vi) of Theorem~\ref{thm:finite} are easy consequences of the following.
\begin{lemma} \label{lem:chordal-conv}
	\begin{enumerate}
		\item [(i)]
		In the setting of Theorem~\ref{thm:finite} and the two paragraphs above, 
		the triple $(\wh Z^{n}, T^{n} ,\ell^{n})$ converges in probability to $(\wh\Zb,\Tb,\ellb)$, where the first two coordinates are equipped with the Skorokhod topology and the third coordinate is equipped with the uniform topology. 
		\item [(ii)] Given a constant $\ell'$ such that $0<\ell'<\ell$, let $\bm {s}$ and $\bm{t}$ be defined as $t_2$ and $t_2'$, respectively, in \eqref{eq:t3'}, with $\wt u$ in place of $t_3$. Let $s$ (respectively, $t$) be the index of the last\footnote{We remark that the $s$ here corresponds to the time $t_2$ defined in Section \ref{sec:crossing-discrete}.} (respectively, first) $c$-step in $w^+$ before (respectively, after) $\lambda ( \lceil(1-\ell') h_n \rceil )$ that has neither an $a$-match nor a $b$-match within $w^+$. Let $s^n=(3n)^{-1}s$ and $t^n=(3n)^{-1} t$. Then $\lim_{n\to\infty} s^n=\bm{s}$ and $\lim_{n\to\infty} t^n=\bm{t}$ in probability.
	\end{enumerate}
\end{lemma}
\begin{proof}
	To prove Assertion (i), we first consider the setting of Corollary~\ref{cor:bol} and not the setting of Theorem~\ref{thm:finite}. In the corollary, $M_n$ is a critical Boltzmann triangulation with boundary length $h_n+2$, and the definitions of $(\wh Z^{n}, T^{n} ,\ell^{n})$ and $(\wh\Zb,\Tb,\ellb)$ are changed accordingly. In this case, 
	on the one hand, away from the two endpoints of the interval of definition of 	 $(\wh Z^{n}, T^{n} ,\ell^{n})$, Assertion (i) follows from Lemma~\ref{prop5} and Proposition~\ref{prop3}. More precisely, we use a variant of Proposition~\ref{prop3} for walks of random duration, but this variant is immediate from the fixed duration variant since we can condition on the duration. On the other hand, Lemma~\ref{prop22} yields that the boundary contribution is negligible. This proves Assertion (i) in the setting of Corollary~\ref{cor:bol}.
	
	Back to the original setting of Theorem~\ref{thm:finite},
	we observe that the law of random walk $\wh Z^n$ in this theorem restricted to its initial interval ending at $u$ is absolutely continuous with respect to the law in the random area setting considered in the above paragraph. This allows us to transfer Assertion (i) from the Boltzmann case to the fixed size case. 
	
	Assertion (ii) in Lemma~\ref{lem:chordal-conv} is a direct consequence of the convergence of $Z^n$ to $\Zd$.
\end{proof}

\begin{lemma}
	Assertions (v) and (vi) in Theorem~\ref{thm:finite} hold.
\end{lemma}
\begin{proof}
	(v): The convergence of the branches away from $\bdy \D$ follows from Proposition~\ref{prop3} and Assertion (iii) of Theorem~\ref{thm1}. The fact that the boundary effect is negligible follows from Assertion (i) of Lemma~\ref{lem:chordal-conv} and its proof. 
	
	(vi): Given the mating-of-trees description of the crossing events given in Section~\ref{sec:crossing}, and their discrete counterparts given in Section \ref{sec:crossing-discrete},
	Assertion (vi) follows from Lemmas~\ref{lem:uniform} and~\ref{lem:chordal-conv}.
\end{proof}

\subsection{Proof of Proposition \ref{prop23}}\label{subsec:flip}
\begin{figure}
	\centering
	\includegraphics[scale=1]{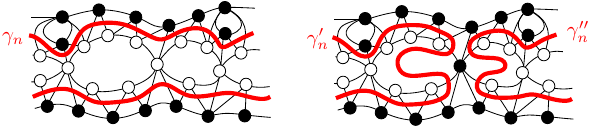}
	\caption{The figure shows the percolation interface between two (or more) clusters before (left) and after (right) the color of a pivotal point has been flipped. It follows from Lemma \ref{prop19}(iii) that the red ``excursions'' inside the two white bubbles on the right figure are microscopic, which we use in the proof of Proposition \ref{prop23}.
	}
	\label{fig:flip}
\end{figure}

In this section, we finish the proof of Proposition~\ref{prop23}. We retain the notions defined in Section~\ref{sec:flip}. The technical bulk of this section is the convergence of pivotal measures after flipping the color of a vertex. We start by proving an infinite volume version of Proposition \ref{prop23}. The loops in $\Gab_z$ are ordered by the same rule as the loops in $\Gab$, that is, by their value.
\begin{lemma}
	Consider the setting of Lemma \ref{prop34}.
	For $j\in\N_+$ and $\eps>0$ let $z_n\in\C$ (respectively, $z\in\C$) have the law of a uniformly sampled $\eps n$-pivotal (respectively, $\eps$-pivotal) point associated with $\ga^n_j$ (respectively, $\gab_j$), such that $z_n\rta z$ almost surely. Let $\Ga^n_{z_n}$ (respectively, $\Gab_z$) be the collection of percolation cycles after flipping the color of $z_n$ (respectively, $z$). Then $\Ga^n_{z_n}$ converges in probability to $\Gab_z$. Furthermore, for any $j\in\N_+$, the counting measure on the $\eps n$-pivotal points associated with the $j$th percolation cycle of $\Ga^n_{z_n}$ converge in probability to the $\eps$-pivotal measure associated with the $j$th loop of $\Gab_{z}$. 
	\label{prop25}
\end{lemma}
\begin{proof}
	It was proved in Lemma \ref{prop34} that $\Ga^n_{z_n}$ converges in probability to $\Gab_z$. 
	Now we will prove convergence of the pivotal measure of $\Gamma^n_{z_n}$. Fix $\eps'>0$. Let $\Delta_n$ be the set of vertices which are $\eps n$-significant for $\Ga^n_{z_n}$ but which are not $\eps' n$-significant for $\Ga^n$.
	We claim that the desired convergence follows from 
	\eqb
	n^{-1/4}\cdot \# \Delta_n
	= o_{\eps'}(1),
	\label{eq85}
	\eqe
	where $o_{\eps'}(1)$ may depend on $\eps$ and $j$ but not on $n$. Before proving \eqref{eq85}, we will show how we use it to conclude the proof of the lemma.

	Recall that the percolation cycles in $\Ga^n_{z_n}$ and $\Gab_z$ are ordered by their value. However, it is sufficient to prove the lemma for an arbitrary other ordering, and in the proof we will reorder the cycles in $\Ga^n_{z_n}$ and $\Gab_z$ so the ordering is closer to the ordering of the cycles in $\Ga^n$ and $\Gab$. Fix $k\in\Z^{\geq 3}$. Let $D^n(\eps,k)\subset\C$ (respectively, $D^n_{z_n}(\eps,k)$) be the set of $\eps n$-pivotal points associated with the first $k$ percolation cycles of $\Gamma^n$ (respectively, $\Gam^n_{z_n}$). Here the percolation cycles in $\Gam^n$ are ordered using their value as defined by \eqref{eq:value}, while the percolation cycles in $\Gam^n_{z_n}$ are ordered as follows. We may assume that the set $\cL_{z_n}$ contains exactly three percolation cycles of area at least $\eps n$, since this holds with probability $1-o_n(1)$ by \cite[Theorem 2]{camia-newman-full}; see the first paragraph in the proof of Lemma \ref{prop34} for an argument. If there is exactly one percolation cycle in $\Gamma^n_{z_n}\cap \cL_{z_n}$ of area at least $\eps n$ then we let this percolation cycle be the first percolation cycle of $\Gamma^n_{z_n}$, and we let the remaining percolation cycles have the same relative ordering as for $\Gam^n$. 
	If there are two percolation cycles in $\Gamma^n_{z_n}\cap \cL_{z_n}$ of area at least $\eps n$ then we let the percolation cycles in $\Gamma^n_{z_n}\cap \cL_{z_n}$ be the first and second, respectively, percolation cycle of $\Gamma^n_{z_n}$, and we let the remaining percolation cycles have the same relative ordering as for $\Gam^n$. The percolation cycles of $\Gamma^n_{z_n}\cap \cL_{z_n}$ of area less than $\eps$ will be given the largest possible rank such that all percolation cycles with smaller rank have a smaller value. The percolation cycles of $\Gab_z$ are ranked in a similar way in the continuum.
	
	Let $k>3$. Let $z'_n\in\C$ be sampled uniformly at random from $D^n_{z_n}(\eps,k-2)$. Let $z'$  be a point sampled from the measure on the $\eps$-pivotal points associated with the first $k-2$ CLE$_6$ loops of $\Gab_{z}$. We want to prove that $z_n'$ converges in law to $z'$.
	With the new ordering of $\Gam_{z_n}^n$, for any $\ga\in \Gamma_{z_n}^n\cap \Gam^n$, the rank of $\ga$ in $\Gam_{z_n}^n$ and the rank of $\ga$ in $\Gam^n$ differ by at most 2. Therefore $D^n_{z_n}(\eps,k-2)\setminus D^n(\ep',k)\subset \Delta_n$. 
	By \eqref{eq85}, for $\eps'\in(0,\eps)$ and with probability $1-o_n(1)$,
	\eqb
	z'_n\in D^n(\eps',k)\quad \text{with probability}\quad 1-o_{\eps'}(1),
	\label{eq84}
	\eqe 
	where $o_{\eps'}(1)$ may depend on $\eps$ and $k$ but not on $n$.
	Let $z''_n$ be sampled uniformly at random from $D^n(\eps',k)\cap D^n_{z_n}(\eps,k-2)$, and let $z''$ be sampled uniformly at random from the $\eps$-pivotal measure supported on the continuum counterpart of this set. 
	Let $z'''_n$ be sampled uniformly at random from $D^n(\eps',k)$, and let $E_n$ be the event that $z'''_n\in D^n(\eps',k)\cap D^n_{z_n}(\eps,k-2)$. Define $z'''$ and $E$ similarly in the continuum. Then $z'''_n$ and $\1_{E_n}$ converge jointly to $z'''$ and $\1_E$ by Lemma \ref{prop34}. Since $z''_n$ (respectively, $z''$) has the law of $z'''_n$ (respectively, $z'''$) conditioned on the event $E_n$ (respectively, $E$), the convergence of $z'''_n$ to $z'''$ implies that $z''_n$ converges in law to $z''$. 
	By \eqref{eq84} and its continuum counterpart, the convergence of $z''_n$ to $z''$ implies that $z'_n$ converges in law to $z'$. This shows that the \emph{probability measure} associated to $D^n_{z_n}(\eps,k-2)$ converges to its continuum counterpart. Lastly, since the $\eps n$-pivotal counting measure associated with $\Gamma_{z_n}^n$ is identical to the $\eps n$-pivotal counting measure associated with $\Gamma^n$ on the intersection of their supports, this shows that the (non-normalized) measure associated to $D^n_{z_n}(\eps,k-2)$ converges to its continuum counterpart. This proves Lemma~\ref{prop25}.
	\begin{figure}
		\centering
		\includegraphics[scale=1]{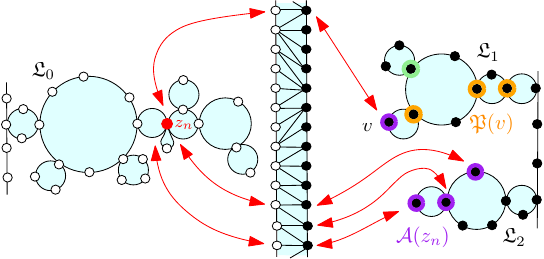}
		\caption{Illustration of objects defined in the proof of Lemma \ref{prop25}. The pivotal points $\frk P(v)$ for $v\in\mcl A(z_n)$ are of type (vi). The pivotal point marked in green is of type (v). Most new $\eps n$-pivotal points which appear after flipping the color of $z_n$ are of one of these kinds.}
		\label{fig:pivlocal}
	\end{figure}	
	It remains to prove \eqref{eq85}. We advise to study Figure \ref{fig:pivlocal} and the left part of Figure \ref{fig:flip} while reading the proof. Without loss of generality, assume $z_n$ is white. The pivotal point $z_n$ is \emph{associated} with some percolation cycle $\ga$ (recall the definition of association from Section \ref{subsec:pivot}). Consider the pair of forested lines relative to the envelope closing time of $\ga$. Let $\frk L_0$ on the left forested line containing $z_n$ (Definition~\ref{def:lt-cluster}) as in the left part of Figure \ref{fig:pivlocal}. We prove \eqref{eq85} by considering separately different classes of vertices which are candidates for vertices contained in $\Delta_n$. Before defining these classes we will describe the typical picture we see on the map locally near $z_n$, and we need to introduce some notation.
	
	We say that a vertex $v$ is \emph{on} a looptree $\frk L$ (respectively, a bubble $B$ of a looptree) if it is contained in the vertex set of $\frk L$ (respectively, $B$). We say that a vertex $v$ is \emph{inside} a bubble $B$ of a looptree if it is not on $B$ but is separated from $\infty$ by $B$. We say that a vertex $v$ on a looptree $\frk L$ is \emph{pivotal for $\frk L$} if we can find two bubbles $B_1$ and $B_2$ of $\frk L$ such that $v$ is on both $B_1$ and $B_2$.
	
	We will now argue that the percolated map $(M_n,\sigma_n)$ in a neighborhood around $z_n$ is rather similar to the neighborhood around a pivotal point sampled from a measure $p^n(s,u)$ for $u\in\R$ and $s<0$ deterministic. More precisely, the following hold with probability $1-o_n(1)$. Let $\wh L^n$ and $\wh R^n$ be the renormalized L\'evy walks relative to the envelope closing time of $\ga$. Then $\frk L_0$ is a discrete looptree on the discrete forested line encoded by $\wh L^n$. If we flip the color of $z_n$ then $\frk L_0$ is split into exactly two looptrees which are macroscopic (i.e., looptrees which enclose at least $\eps'n$ vertices with probability $1-o_{\eps'}(1)$, uniformly in $n$), plus possibly some looptrees which enclose area $o_n(1)$; this follows by properties of the scaling limit $\wh\Lb$ of $\wh L^n$. Let $\mcl A(z_n)$ denote the set of black vertices that are adjacent to $z_n$ and \emph{not} inside any bubble of $\frk L_0$. \nina{In the scaling limit, $z=\wh\etab^{w}(t_1)=\wh\etab^{w}(t_2)$ for $t_1<t_2<0$ such that
		$\wh\Lb(t_1)=\wh\Lb(t_2)$ and
		$(\wh\Lb-\wh\Lb(t_1))|_{[t_1,t_2]}$ is a positive excursion from 0. Furthermore, $z$ is on exactly two of the looptrees encoded by $\wh\Rb$
		if and only if the time-reversal $t\mapsto\wh\Rb(-t)$ of $\wh\Rb$ has a running infimum in $[-t_2,-t_1]$; otherwise $z$ is on exactly one of the looptrees encoded by $\wh\Rb$.
		Using this property of the limiting process,} 
	we get that with probability $1-o_n(1)$ there are either one or two looptrees such that all vertices of $\mcl A(z_n)$ are on one of these looptrees. We assume in the remainder of the proof that there are two (not one) such looptrees, but the case of one looptree can be treated in a similar way. Denote the two looptrees by $\frk L_1$ and $\frk L_2$. Let $\mcl A(z_n)=\mcl A_1(z_n)\cup \mcl A_2(z_n)$, where vertices in $\mcl A_1(z_1)$ (respectively, $\mcl A_2(z_2)$) are on $\frk L_1$ (respectively, $\frk L_2$).

	Let $v$ be a vertex on $\frk L_1$, and let $v_0$ denote the root-vertex of $\frk L_1$. We can find bubbles $B_1,\dots,B_k$ of $\frk L_1$ such that $B_i$ and $B_{i+1}$ share a vertex for all $i$, $v$ is on $B_1$, $v_0$ is on $B_k$, and $B_i\neq B_j$ for all $i\neq j$. Let $\frk P(v)\subset V(M_n)$ denote the set of vertices $v'$ on $\frk L_1$ for which we can find an $i$ such that $v'$ is on both $B_i$ and $B_{i+1}$. Define $\frk P(v)$ in the exact same way if $v$ is a vertex on $\frk L_2$.
	
	When studying vertices in $\Delta_n$ we consider the following classes of vertices $v_*\in V(M_n)$ separately:
	\begin{compactitem}
		\item[(i)] $v_*$ is inside a bubble $B$ of $\frk L_0$, $\frk L_1$, or $\frk L_2$, such that $z_n$ is \emph{not} on $B$,
		\item[(ii)] $v_*$ is inside a bubble $B$ of $\frk L_0$, such that $z_n$ is on $B$, 
		\item[(iii)] $v_*$ is on $\frk L_0$,
		\item[(iv)] $v_*$ is on $\frk L_1$ or $\frk L_2$ but is not a pivotal point for this looptree,
		\item[(v)] $v_*$ is a pivotal point for $\frk L_1$ or $\frk L_2$, but is not contained in $\frk P(v)$ for any $v\in\mcl A(z_n)$,
		\item[(vi)] $v_*\in\frk P(v)$ for some $v\in\mcl A(z_n)$, and
		\item[(vii)] $v_*$ is not covered by any of the cases above.
	\end{compactitem}
	
	By Definition \ref{def:pivot}, notice that the significance of a vertex $v_*$ of class (i), (iii), (iv), or (vii) for $\Gam_{z_n}^n$ is smaller than or equal to the significance of $v_*$ for $\Gam^n$. Therefore, and since $\eps'<\eps$, the set considered in \eqref{eq86} contains no vertices of these classes. 
	Moreover, with probability $1-o_n(1)$ vertices of class (ii) enclose area $o_n(1)$, which implies that such vertices are not $\eps' n$-significant pivotal points for $\Gam_{z_n}^n$ with probability $1-o_n(1)$. Hence it only remains to deal with vertices of class (v) and (vi). Let $\mcl V_{\op{(v)}}\subset V(M_n)$ (respectively, $\mcl V_{\op{(vi)}}\subset V(M_n)$) denote the set of vertices $v_*$ satisfying (v) (respectively, (vi)). 
	The above shows that  with probability $1-o_n(1)$,
	\eqb
	n^{-1/4}\cdot \# \big( \Delta_n \setminus (\mcl V_{\op{(v)}} \cup \mcl V_{\op{(vi)}} )\big) = o_n(1).
	\label{eq86}
	\eqe
	
	Next we deal with  vertices of class (vi). We will prove that with probability $1-o_n(1)$,
	\eqb
	n^{-1/4}\cdot \# (\Delta_n \cap \mcl V_{\op{(vi)}} ) = o_{\eps'}(1).
	\label{eq87}
	\eqe
	This is equivalent to showing that with probability $1-o_n(1)$,
	\eqb
	n^{-1/4} \cdot \# \bigcup_{v\in\mcl A(z_n) } \frk P(v) \cap \Delta_n = o_{\eps'}(1).
	\label{eq88}
	\eqe 
	First we argue tightness of $\#\mcl A(z_n)$. Let $t_-,t_+<0$ be chosen as small as possible such that $\wh L^n|_{[t_-,t_+]}$ encodes a looptree rooted at $z_n$ which encloses area at least $\eps$, and such that $\wh L^n_t>\wh L^n_{t_-}=\wh L^n_{t_+}$ for all $t\in(t_-,t_+)$. Observe that in any bounded neighborhood of $n^{3/4}t_+$ (respectively, $n^{3/4}t_-$), $(\wh L,\wh R)$ converges in law for the total variation distance to a bi-infinite walk with independent and identically distributed increments as for the time-reversal of the walk in Lemma \ref{prop4}, except that the walk is conditioned to have a strict running infimum at time $n^{3/4}t_+$ (respectively, a strict running infimum in backwards direction at time $n^{3/4}t_-$). Furthermore, note that by Lemma \ref{prop:piv} and with $\wh\eta_{\op{v}}$ defined as above this lemma,
	if $t'>t_+$ is such that $\inf_{s\in[n^{3/4}t_+,n^{3/4}t']} \wh L_s<\wh L_{n^{3/4}t_+}$ then there are no $t''>t'$ such that $\wh\eta_{\op{v}}(t'')=z_n$,
	if $t'<t_-$ is such that $\inf_{s\in[n^{3/4}t',n^{3/4}t_-]} \wh L_s<L_{nt_-}$ then there are no $t''<t'$ such that $\wh\eta_{\op{v}}(t'')=z_n$,
	and by the definition of $t_\pm$ there are no $t''\in(t_-,t_+)$ such that $\wh\eta_{\op{v}}(t'')=z_n$. Combining these two results we get tightness of $\#\mcl A(z_n)$. In fact, we get convergence in law of $\#\mcl A(z_n)$, but this stronger result is not needed.
	
	For each $v\in\mcl A(z_n)$, $\#\frk P(v)$ has magnitude of order $n^{1/4}$, since it follows from the definition of the mappings $\ccode$ and $\ccwcode$ (see also Lemma \ref{prop:piv}), that there is a bijection between $\#\frk P(v)$ and the set of strict running infima for $\wh R^n$ relative to $\wh\eta_{\op{v}}^{-1}(v)$ (if we exclude the root-vertex of the looptree $\frk L_1$ or $\frk L_2$ containing $v$). Furthermore, with probability $1-o_n(1)$,
	\eqbn
	n^{-1/4} \cdot \# \big( \frk P(v) \cap \Delta_n\big) = o_{\eps'}(1).
	\eqen
	Combining this with tightness of $\#\mcl A(z_n)$, we get \eqref{eq88}.

	It remains to consider vertices of class (v). We will prove that with probability $1-o_n(1)$,
	\eqb
	n^{-1/4}\cdot \# \big(\Delta_n \cap \mcl V_{\op{(v)}} \big) = o_{\eps'}(1).
	\label{eq89}
	\eqe
	Vertices $v_*\in \Delta_n$ of class (v) must satisfy the following, where we assume without loss of generality that $v_*$ is on $\frk L_1$. Since $v_*$ is a pivotal point for $\frk L_1$, we can find a finite collection of looptrees $\frk L'_1,\frk L'_2,\dots$ such for any $i\neq j$ the only vertex on both $\frk L'_i$ and $\frk L'_j$ is $v_*$, and such that the union of the vertices on $\frk L'_1,\frk L'_2,\dots$ equals the set of vertices on $\frk L_1$. Since the significance of $v_*$ increases when the color of $z_n$ is flipped, and by the definition of points of class (v), one of these looptrees (say, $\frk L'_1$) encloses area at least $\eps$, and another looptree (say, $\frk L'_2$) contains the root of $\frk L_1$ and all vertices in $\mcl A(z_0)$, and encloses area less than $\eps'$. Uniformly over all $v\in\mcl A(z_n)$, the bubbles $B_1,\dots,B_k$ on the path in the definition of $\frk P(v)$ enclose at least area $\eps'$ with probability $1-o_{\eps'}(1)$. Therefore, with probability $1-o_{\eps'}(1)$ a looptree $\frk L'_2$ as just described does not exist, which implies \eqref{eq89}.
	
	Combining \eqref{eq86}, \eqref{eq87}, and \eqref{eq89} we obtain \eqref{eq85}, which concludes the proof of the lemma. 
\end{proof}

\begin{proof}[Proof of Proposition \ref{prop23}]
	This is immediate by an application of Lemma \ref{prop25} and Proposition~\ref{prop3}. 
\end{proof}

\appendix

\section{\xin{Recovering the DFS tree and the walk from the percolated UIPT}} \label{sec:recover}
\xin{In this appendix  we prove the last sentence of Theorem~\ref{thm:UIPT} and Theorem~\ref{thm:dfs-inf}(i), which is summarized in the following lemma.}
\begin{lemma}\label{prop9}
	Let $(M,\si)$ be an instance of the percolated UIPT. Then there is a.s.\ a unique $w\in\imK$ and a unique $\tau^*\in\DFS_{M^*}$ such that $\Phi^\infty(w)=(M,\si)$ and $\si=\Lambda_{M}(\tau^*)$. 
\end{lemma}
Note that the combination of Lemmas~\ref{prop17} and~\ref{prop9} imply that almost surely each of the following objects determine the other two: $w$, $(M,\sigma,\frk e)$, and $(M,\tau^*,\frk e)$.

\begin{figure}[h]
	\centering
	\includegraphics[scale=1]{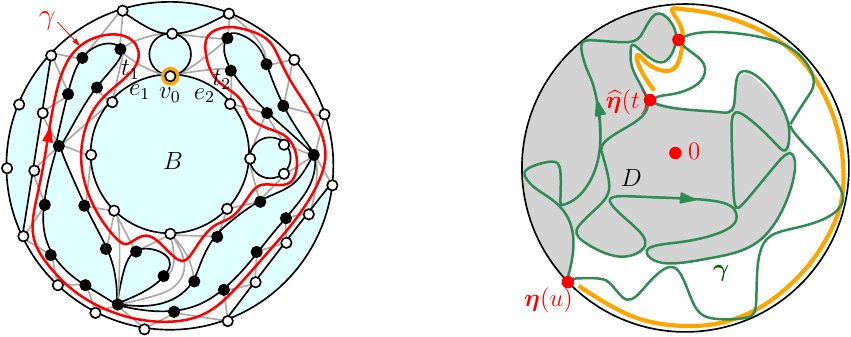}
	\caption{Illustration of the proof of Lemma~\ref{prop9}. Left: Illustration of the notation for the event $A(k)$. 
		Right: Illustration of the event described in part (ii) of the proof. The CLE$_6$ loop $\gab$ is traced in clockwise direction starting and ending at $\etab(u)$. The gray region is the region covered by $\etab$ between time $u$ and the time at which the complementary component $D$ of $\gab$ has just been filled in. The yellow curve indicates the left boundary of $\etab$ at this time. The two red points in the upper part of the figure are pivotal points.}
	\label{fig:meas} 
\end{figure}

In our proof of Lemma \ref{prop9}, we will use the following result, which follows directly from the description in Section \ref{sec:cont-cle}.
\begin{lemma}\label{prop83}
	Recall the definitions of Section \ref{sec:cont-cle} about the conformal loop ensemble $\CLE_6$ $\Gab$ on $\C$.
	Let $\gab\in\Gab$ be a CLE$_6$ loop oriented clockwise, let $s$ be the quantum natural length of $\gab$, and let $\wh\Zb=(\wh\Lb,\wh\Rb)$ be the L\'evy process relative to the envelope closing time $u$ of $\gab$. Then the set of pivotal points of type 2 associated with $\gab$ is given by
	\eqbn
	\big\{ \wh\etab^{\etab(u)}(t)\,:\, t\in(-s,0), \exists \delta>0 
	\text{\,\,such\,\,that\,\,}
	\wh\Lb_t
	=\inf_{t'\in[t-\delta,t]} \wh\Lb_{t'}
	>\inf_{t'\in[-s,t]} \wh\Lb_{t'}
	\big\}.
	\eqen
	If $\gab$ is oriented counterclockwise, then the same result holds with $\wh\Rb$ instead of $\wh\Lb$.
\end{lemma}

\begin{proof}[Proof of Lemma \ref{prop9}]
	Let $w\in\imK$ be chosen according to the uniform distribution. Let $(M,\si)=\Phi^\infty(w)$, and let $\etae$ be the corresponding bijection $\ZZ\to E(M)$.
	Let $\tau^*$ be the set of edges $e^*$ of $M^*$ such that $w_{\etae^{-1}(e)}$ is an $a$-step or a $b$-step. From Lemma~\ref{prop17} we know that $(M,\si)$ has the law of the percolated UIPT, and that $\tau^*$ is in $\DFS_{M^*}$ and satisfies $\Lambda_{M}(\tau^*)=\si$.
	We want to show that a.s.\ for $w$, if $\tw\in\imK$ and $\ttau^*\in \DFS_{M^*}$ are such that $\Phi^\infty(\tw)=(M,\si)$ and $\Lambda_{M}(\tau^*)=\si$, then $\tw=w$ and $\ttau^*=\tau^*$. Let $B(k)$ be the event (depending on $w$) that for all $i\in \{-k,\dots,k \}$ we have $\tw_i=w_i$, and $\etae(i)$ is in $\ttau^*$ if and only if $\etae(i)$ is in $\tau^*$. It suffices to show that for all $k\in \NN$, $B(k)$ holds almost surely.
	
	Let us fix $k\in \NN$. We now define an event $A(k)$ (depending on $w$) such that $A(k)$ implies $B(k)$ and we will later show that $A(k)$ holds almost surely. Figure~\ref{fig:meas} (left) indicates our notation. 
	Recall that each cluster $C$ of $(M,\si)$ has an associated outside-cycle $\ga_C$ and an associated looptree $\frk L_C=\frk L(\ga_C)$ (Definition~\ref{def:lt-cluster}). By \emph{bubble} of $C$ we mean the submap of $M$ made of the vertices and edges which are on or inside one of the bubbles of the looptree $\frk L_C$. 
	We define $A(k)$ to be the event that there is a cluster $C$, having a bubble $B$ and a vertex $v_0$ (pivotal of type 2) of $B$ such that the following holds.
	First, the interior of the bubble $B$ contains all the edges $\etae(i)$ for $i\in\{-k,\dots,k\}$. 
	Second, flipping the color of $v_0$ splits the percolation cycle $\ga_C$ into several cycles $\ga_1\ldots \ga_i$, with $\ga_1$ in the inside-region of $\ga_2$, and  $B$ in the inside-region of $\ga_1$.
	
	Let us now show that $A(k)$ implies $B(k)$. We break this proof into a series of facts. For concreteness we suppose throughout that $C$ is a white cluster.
	\begin{fact} \label{fact1}
		The DFS tree $\ttau^*$ can be obtained from a DFS of $M^*$ where, each time the algorithm is in Case (a) of Definition~\ref{def:DFS-algo-inf} and several edges are possible to move the chip along, the choice is made according to the rule (ii) of Definition~\ref{def:Delta}.
	\end{fact}
	Since $\ttau^*\in\DFS_{M^*}$, we know that $\ttau^*$ is associated to a DFS of $M^*$. Now suppose by contradiction that $\ttau^*$ is not associated to a DFS of $M^*$ satisfying rule (ii) of Definition~\ref{def:Delta}. This implies that there is a DFS $Y$ of $M^*$ associated to $\ttau^*$, and a face $f$ of $M^*$ such that when the DFS $Y$ reaches the first vertex $v_f\in V(M^*)$ incident to $f$ the rule (ii) of Definition~\ref{def:Delta} is not respected. We can indeed assume that the rule is not respected at the first vertex $v_f$ incident to $f$, because if the rule is not respected at a vertex $v\neq v_f$ incident to $f$, the choice made by the DFS at $v$ does not affect the resulting DFS tree $\ttau^*$ (so we could change the choice at $v$ to ensure that rule (ii) is satisfied there). Let $e_L$ and $e_R$ be the left forward and right forward edges at $v_f$, and let $v_L$ and $v_R$ be the incident vertices.
	Let us assume for concreteness that $f$ is black. Since the DFS $Y$ violates rule (ii) at $v_f$, the chip moves through the edge $e_R$ in $Y$, so $v_R$ is the child of $v_f$ in  $\ttau^*$.
	Since $v_f$ is the first vertex incident to $f$ to be reached during $Y$, the vertex $v_R$ will be an ancestor of $v_L$ in $\ttau^*$. 
	Hence, it is easy to see that the color of the vertex $f^*$ is white in $\Lambda_M(\ttau^*)$ (since $e_L^*$ will be the parent-edge of $f^*$ in $\ttau$). We reach a contradiction, so Fact~\ref{fact1} holds.

	We set some additional notation indicated in Figure~\ref{fig:meas}. 
	Let $T$ be the set of triangles of $M$ outside of $B$ but incident to an edge of $B$.
	Let $e_1$ (resp. $e_2$) be the edges preceding (resp. following $v_0$) in clockwise order along the boundary of $B$, and let $t_1\in T$ (resp. $t_2\in T$) be the triangle incident to $e_1$ (resp. $e_2$) lying outside $B$. Let $\ga_1'$ be the part of $\ga$ going from $t_1^*$ to $t_2^*$ in counterclockwise around~$B$. 
	
	\begin{fact} \label{fact3} In the tree $\ttau^*$, $t_1^*$ is an ancestor of $t_2^*$, moreover $\ga_1'\subset \ttau^*$. In particular $t_2^*$ is the descendant of all the vertices of $M^*$ dual to triangles in $T$. Moreover, denoting $s_2$ the triangle of $M$ incident to $e_2$ inside $B$, $t_2^*$ is the parent of $s_2^*$, and  $s_2^*$ is the ancestor of all the vertices of $M^*$ inside $B$.
	\end{fact}
	Consider a DFS $X$ associated to $\ttau^*$ as specified by Fact~\ref{fact1}.
	Let $u_0\in V(M^*)$ be the first vertex on $\ga$ visited during this DFS. We claim that $u_0$ does not belong to $\ga_1'$. Indeed $\ga_1'$ is a part of $\ga_1$ which is nested inside $\ga_2$. This implies that the the vertices of $M^*$ on $\ga\setminus \ga_1'$ together with the interior of the face $v_0^*\in F(M^*)$ separates $\ga_1'$ from $\infty$. Thus $u_0$ does not belong to $\ga_1'$. By definition of the DFS $X$, once the DFS $X$ reaches $u_0$, it will follow the edges of $\ga$, with the white cluster $C$ on its left, until it has gone through all the edges of $\ga$ except one. Since $u_0$ does not belongs to $\ga_1'$, the DFS $X$ will reach $t_1^*$ before $t_2^*$ and will follow the edges of $\ga_1'$ from $t_1^*$ to $t_2^*$. Hence $\ga_1'\subset \ttau^*$, and $t_1^*$ is an ancestor of $t_2^*$. Moreover, at the time the DFS $X$ reaches $t_2^*$, the vertices of $M^*$ inside $B$ are unvisited and reachable from $t_2^*$. Moreover, any path of unvisited vertices starting at $t_2^*$ and ending inside $B$ goes through $s_2^*$. This shows that $t_2^*$ is the parent of $s_2^*$, and  $s_2^*$ is the ancestor of all the vertices of $M^*$ inside $B$.

	\begin{fact} \label{fact4}
		The restrictions of the trees $\ttau^*$ and $\tau^*$ to the submap $B$ coincide.
	\end{fact}
	By Fact~\ref{fact3}, when the DFS of $\ttau^*$ enters $B$ it is by following the edge $e_2^*$ from $t_2^*$ to $s_2^*$, and at that time all the vertices of $M^*$ inside $B$ are unvisited, while all the adjacent vertices of $M^*$ not in $B$ (that is the vertices dual to the triangles in $T$) are visited. Hence the DFS of $\tau^*$ will visit all the vertices inside $B$ before backtracking from $s_2^*$ to $t_2^*$. 
	Since $\ttau^*$ is arbitrary, the same holds for the DFS associated to $\tau^*$. Since the DFS associated to $\ttau^*$ and $\tau^*$ follow the same rule (ii) of Definition~\ref{def:Delta}, they will perform exactly identically while inside $B$, which proves Fact~\ref{fact4}.

	\begin{fact} \label{fact5}
		For all $i\in \{-k,\dots,k \}$, $\tw_i=w_i$.
	\end{fact}
	Let $\wt \etae$ (resp. $\wt \etavf$) be the bijection giving the order of creation of the edges (resp. vertices and faces) of $M$ during the application of $\Phi^\infty(\tw)$.
	Recall that the order of creation of the triangles of $M$ given by $\wt \etavf$ corresponds to the order of visit of the vertices of $M^*$ for a DFS of $M^*$ associated with $\tau^*$ (although it is not equal to the DFS defined by Fact~\ref{fact1}). Hence, by Fact~\ref{fact4} all the triangles in $T$ are created after all the triangles inside $B$ during  $\Phi^\infty(\tw)$. This shows that the preimage $\{\wt \etae^{-1}(e),~e\in E(B)\}$ of the edges of $B$ forms an interval of integers $I=\{a,a+1,\ldots,b\}$.
	Moreover we claim that the word $\tw_{I}=\tw_a\ldots \tw_b$ is in $\bmK$, and $\bPhi(\tw_{I})=(B',\si')$, where $(B',\si')$ is obtained from the restriction of $(M,\si)$ to $B$ by flipping the color of $v_0$ (and taking $e_2$ as the root-edge). Indeed, we know from the definition of $\bPhi$ that the color of $v_0$ in $\bPhi(\tw_{I})$ is black. Moreover the color of the vertices in $\bPhi(\tw_{I})$ and in $(M,\si)$ have to be the same except for those on the boundary of the past-triangulation corresponding to the subword of $\tw^-=\ldots \tw_{b-1}\tw_b$, and the only such vertex is $v_0$. Since $\tw$ is arbitrary the same holds for $w$, so there is an interval of integers $J$ such that $\bPhi(w_{J})=(B',\si')$. By injectivity of $\bPhi$ we get $w_{J}=\tw_{I}$. Moreover, by definition, $\etae(0)=\wt \etae(0)$ is the root-vertex of $(M,\si)$, which is inside $B$. So $I=J\supseteq \{-k,\ldots,k\}$. This proves Fact~\ref{fact5}, and completes the proof of the fact that $A(k)$ implies $B(k)$.\\

	It only remains to prove that $A(k)$ holds almost surely. Consider words $w^n\in\imK$ coupled in such a way that the convergence in Theorem~\ref{thm1} is almost sure. Precisely, each word $w^n\in\imK$ has uniform distribution, and they are coupled in such a way that the associated renormalized walks $Z^n$ converge to a limiting Brownian motion $\Zb$ as $n\rta\infty$. In order to prove that $A(k)$ holds almost surely for a uniformly random word $w\in\imK$ it suffices to prove that for all $\eps>0$, there exists $n=n(\eps)$ such that $w^n$ satisfies $A(k)$ with probability at least $1-\eps$. To prove this fact we will use Theorem~\ref{thm1}, and adopt the notation of that theorem. 
	
	We claim that in the limiting CLE$_6$ $\Gab$, with probability 1 we can find some loop $\gab\in\Gab$ and a bounded complementary component $D$ of $\gab$, such that $0\in D$ and $D$ is not in the inside-region of $\gab$. The latter property means that if $\gab$ is traced in clockwise (resp.\ counterclockwise) direction then $\partial D$ is traced in counterclockwise (resp.\ clockwise) direction. The claim follows since for each annulus of the form $\{ k<|z|<2k \}$ for $k\in\Z^{>0}$ it holds with positive probability  that we can find an appropriate loop $\gab$ contained in the annulus, and the event that this occurs is independent for each annulus.
	
	For concreteness, we assume that $\gab$ is traced in clockwise direction. Let $u$ be the envelope closing time of $\gab$, and let $z=\etab(u)$. The situation is represented in Figure~\ref{fig:meas} (right). Define 
	$\Tb=(\Tb_t)_{t\leq 0}$, 
	$\wh\etab^{z}=(\wh\etab^{z}(t))_{t\leq 0}$, and
	$\wh\Zb=(\wh\Lb_t,\wh\Rb_t)_{t\leq 0}$ 
	as in Section~\ref{sec:cont-cle} when we recenter $\Zb$ at the envelope closing time $u$, and let $s$ be the quantum natural length of $\gab$. Let $t<0$ be the time for $\wh\etab^{z}$ at which $\wh\etab^{z}$ encloses $D$. If we flip the color of $\etab(t)$, then $\gab$ is split into two nested loops, and $D$ is a complementary component of the inner loop. In particular, $\etab(t)$ is a pivotal point of type~2 for $\gab$. By Lemma \ref{prop83} this means that $\wh\Lb$ has a local running infimum at time $\ellb_t$ which is \emph{not} a global running infimum relative to time $-s$. Furthermore, by properties of L\'evy excursions, there exists $\delta>0$ such that the mass assigned to such running infima of $\wh\Lb$ (equivalently, to pivotal points of type 2) in the interval $[\ellb_t,\ellb_t+\delta]$ is positive. If we flip the color of such a pivotal point then $\gab$ is split into two nested loops, and $D$ is a complementary component of the inner of these loops. By the convergence of the pivotal measure and of the CLE$_6$ loops (Theorem~\ref{thm1}), this implies the existence of $C,B,v_0$ as prescribed in the event $A(k)$ for the percolated UIPT $(M_n,\si_n)=\Phi^\infty(w^n)$ with probability converging to $1$ as $n\rta\infty$. More precisely, the preceding convergence result implies the existence of $C,B,v_0$ if we remove the requirement that the pivotal point $v_0$ belongs to the bubble $B$. But it is easy to see that if there is a such a vertex $v_0$ (possibly not in $B$), then there is another vertex $v_0'$ belonging to $B$ such that $C,B,v_0'$ are as prescribed in the event $A(k)$.
\end{proof}

\bigskip

\noindent {\bf Acknowledgements: } We thank Ewain Gwynne and Scott Sheffield for helpful discussions and thank Cyril Banderier and Ellen Powell for comments on the draft. \nina{We also thank the anonymous referee for careful reading of the paper and numerous helpful comments.} Olivier Bernardi was partially supported by NSF grants DMS-1400859 and  DMS-1800681. Nina Holden was partially supported by a doctoral research fellowship from the Norwegian Research Council and partially supported by Dr.\ Max R\"ossler, the Walter Haefner Foundation, and the ETH Z\"urich Foundation. 
Xin Sun was supported by Simons Foundation as a Junior Fellow at Simons Society of Fellows  and by NSF grants DMS-1811092 and by Minerva fund at Department of Mathematics at Columbia University.

\bibliographystyle{hmralphaabbrv}
\bibliography{biblio,biblong,bibshort,bib-perc,biblio-perco}

\newcommand{\etalchar}[1]{$^{#1}$}
\def\cprime{$'$} \def\cprime{$'$}
\begin{thebibliography}{CDCH{\etalchar{+}}14}

\bibitem[AAB17]{Albenque:limit-simple-triang}
M.~Albenque and L.~Addario-Berry.
\newblock The scaling limit of random simple triangulations and random simple
  quadrangulations.
\newblock {\em Ann. Probab.}, 45(5):2767--2825, 2017.

\bibitem[AC15]{angel-curien-uihpq-perc}
O.~Angel and N.~Curien.
\newblock Percolations on random maps {I}: {H}alf-plane models.
\newblock {\em Ann. Inst. Henri Poincar\'e Probab. Stat.}, 51(2):405--431,
  2015, \arxiv{1301.5311}. \MR{3335009}

\bibitem[AHS17]{ahs-sphere}
J.~Aru, Y.~Huang, and X.~Sun.
\newblock Two perspectives of the 2{D} unit area quantum sphere and their
  equivalence.
\newblock {\em Comm. Math. Phys.}, 356(1):261--283, 2017, \arxiv{1512.06190}.
  \MR{3694028}

\bibitem[AHS20]{aasw-type2}
M.~Albenque, N.~Holden, and X.~Sun.
\newblock Scaling limit of triangulations of polygons.
\newblock {\em Electron. J. Probab.}, 25:Paper No. 135, 43, 2020. \MR{4171388}

\bibitem[Ald91a]{aldous-crt1}
D.~Aldous.
\newblock The continuum random tree. {I}.
\newblock {\em Ann. Probab.}, 19(1):1--28, 1991. \MR{1085326 (91i:60024)}

\bibitem[Ald91b]{aldous-crt2}
D.~Aldous.
\newblock The continuum random tree. {II}. {A}n overview.
\newblock In {\em Stochastic analysis ({D}urham, 1990)}, volume 167 of {\em
  London Math. Soc. Lecture Note Ser.}, pages 23--70. Cambridge Univ. Press,
  Cambridge, 1991. \MR{1166406 (93f:60010)}

\bibitem[Ald93]{aldous-crt3}
D.~Aldous.
\newblock The continuum random tree. {III}.
\newblock {\em Ann. Probab.}, 21(1):248--289, 1993. \MR{1207226 (94c:60015)}

\bibitem[Ang03]{angel-peeling}
O.~Angel.
\newblock Growth and percolation on the uniform infinite planar triangulation.
\newblock {\em Geom. Funct. Anal.}, 13(5):935--974, 2003, \arxiv{0208123}.
  \MR{2024412}

\bibitem[{Ang}05]{angel-scaling-limit}
O.~{Angel}.
\newblock {Scaling of Percolation on Infinite Planar Maps, I}.
\newblock {\em ArXiv Mathematics e-prints}, December 2005,
  \arxiv{math/0501006}.

\bibitem[{Aru}17]{aru-survey}
J.~{Aru}.
\newblock {Gaussian multiplicative chaos through the lens of the 2D Gaussian
  free field}.
\newblock {\em ArXiv e-prints}, September 2017, 1709.04355.

\bibitem[AS03]{angel-schramm-uipt}
O.~Angel and O.~Schramm.
\newblock Uniform infinite planar triangulations.
\newblock {\em Comm. Math. Phys.}, 241(2-3):191--213, 2003. \MR{2013797
  (2005b:60021)}

\bibitem[BBMR15]{BeBoRa-arxiv}
O.~Bernardi, M.~Bousquet-M\'elou, and K.~Raschel.
\newblock Counting quadrant walks via {T}utte's invariant method.
\newblock \href{http://arxiv.org/abs/1708.08215}{Arxiv:1708.08215}, 2017.

\bibitem[BDG04]{BDFG:mobiles}
J.~Bouttier, P.~{Di Francesco}, and E.~Guitter.
\newblock Planar maps as labeled mobiles.
\newblock {\em Electron. J. Combin.}, 11(1):{R}69, 2004.

\bibitem[{Ben}17]{benoist-lqg-chaos}
S.~{Benoist}.
\newblock {Natural parametrization of SLE: the Gaussian free field point of
  view}.
\newblock {\em Electron. J. Probab.}, 23, 2018.

\bibitem[Ber96]{bertoin-book}
J.~Bertoin.
\newblock {\em L\'evy processes}, volume 121 of {\em Cambridge Tracts in
  Mathematics}.
\newblock Cambridge University Press, Cambridge, 1996. \MR{1406564 (98e:60117)}

\bibitem[Ber99]{bertoin-sub}
J.~Bertoin.
\newblock Subordinators: examples and applications.
\newblock In {\em Lectures on probability theory and statistics
  ({S}aint-{F}lour, 1997)}, volume 1717 of {\em Lecture Notes in Math.}, pages
  1--91. Springer, Berlin, 1999. \MR{1746300 (2002a:60001)}

\bibitem[Ber07]{bernardi-dfs-bijection}
O.~Bernardi.
\newblock Bijective counting of {K}reweras walks and loopless triangulations.
\newblock {\em J. Combin. Theory Ser. A}, 114(5):931--956, 2007.

\bibitem[Ber17]{berestycki-gmt-elementary}
N.~Berestycki.
\newblock An elementary approach to {G}aussian multiplicative chaos.
\newblock {\em Electron. Commun. Probab.}, 22:Paper No. 27, 12, 2017,
  \arxiv{1506.09113}. \MR{3652040}

\bibitem[Bet10]{bet-torus1}
J.~Bettinelli.
\newblock Scaling limits for random quadrangulations of positive genus.
\newblock {\em Electron. J. Probab.}, 15:no. 52, 1594--1644, 2010. \MR{2735376
  (2011h:60070)}

\bibitem[Bet12]{bet-torus2}
J.~Bettinelli.
\newblock The topology of scaling limits of positive genus random
  quadrangulations.
\newblock {\em Ann. Probab.}, 40(5):1897--1944, 2012, \arxiv{1012.3726}.
  \MR{3025705}

\bibitem[BF12]{OB-EF:girth}
O.~Bernardi and {\'E}.~Fusy.
\newblock Unified bijections for maps with prescribed degrees and girth.
\newblock {\em J. Combin. Theory Ser. A}, 119:1351--1387, 2012.

\bibitem[BKR17]{bkr-gessel-walk}
A.~Bostan, I.~Kurkova, and K.~Raschel.
\newblock A human proof of {G}essel's lattice path conjecture.
\newblock {\em Trans. Amer. Math. Soc.}, 369(2):1365--1393, 2017,
  \arxiv{1309.1023}. \MR{3572277}

\bibitem[BM05]{MBM:Kreweras}
M.~Bousquet-M\'elou.
\newblock Walks in the quarter plane: {K}reweras' algebraic model.
\newblock {\em Ann. Appl. Probab.}, 15(2):1451--1491, 2005.

\bibitem[BM16]{mbm-gessel-walk}
M.~Bousquet-M\'elou.
\newblock An elementary solution of {G}essel's walks in the quadrant.
\newblock {\em Adv. Math.}, 303:1171--1189, 2016, \arxiv{1503.08573}.
  \MR{3552547}

\bibitem[BM17]{bet-mier-disk}
J.~Bettinelli and G.~Miermont.
\newblock Compact {B}rownian surfaces {I}: {B}rownian disks.
\newblock {\em Probab. Theory Related Fields}, 167(3-4):555--614, 2017,
  \arxiv{1507.08776}. \MR{3627425}

\bibitem[BP93]{bp93}
R.~Burton and R.~Pemantle.
\newblock Local characteristics, entropy and limit theorems for spanning trees
  and domino tilings via transfer-impedances.
\newblock {\em Ann. Probab.}, 21(3):1329--1371, 1993. \MR{1235419}

\bibitem[BSS14]{bss-lqg-gff}
N.~{Berestycki}, S.~{Sheffield}, and X.~{Sun}.
\newblock {Equivalence of Liouville measure and Gaussian free field}.
\newblock {\em ArXiv e-prints}, October 2014, \arxiv{1410.5407}.

\bibitem[Car92]{cardy-formula}
J.~L. Cardy.
\newblock Critical percolation in finite geometries.
\newblock {\em J. Phys. A}, 25(4):201--206, 1992.

\bibitem[CD10]{cd10}
L.~Chaumont and R.~A. Doney.
\newblock Invariance principles for local times at the maximum of random walks
  and {L}\'evy processes.
\newblock {\em Ann. Probab.}, 38(4):1368--1389, 2010. \MR{2663630}

\bibitem[CDCH{\etalchar{+}}14]{cdhks14}
D.~Chelkak, H.~Duminil-Copin, C.~Hongler, A.~Kemppainen, and S.~Smirnov.
\newblock Convergence of {I}sing interfaces to {S}chramm's {SLE} curves.
\newblock {\em C. R. Math. Acad. Sci. Paris}, 352(2):157--161, 2014.
  \MR{3151886}

\bibitem[CK14]{curien-kortchemski-looptree-def}
N.~Curien and I.~Kortchemski.
\newblock Random stable looptrees.
\newblock {\em Electron. J. Probab.}, 19:no. 108, 35, 2014, \arxiv{1304.1044}.
  \MR{3286462}

\bibitem[CK15]{curien-kortchemski-looptree-perc}
N.~Curien and I.~Kortchemski.
\newblock Percolation on random triangulations and stable looptrees.
\newblock {\em Probab. Theory Related Fields}, 163(1-2):303--337, 2015,
  \arxiv{1307.6818}. \MR{3405619}

\bibitem[CLR90]{Cormen:introduction-algorithms}
T.~Cormen, C.~Leiserson, and R.~Rivest.
\newblock {\em Introduction to Algorithms}.
\newblock MIT press, first edition, 1990.

\bibitem[CN06]{camia-newman-full}
F.~Camia and C.~M. Newman.
\newblock Two-dimensional critical percolation: the full scaling limit.
\newblock {\em Comm. Math. Phys.}, 268(1):1--38, 2006. \MR{2249794}

\bibitem[CS12]{chelkak-smirnov-ising}
D.~Chelkak and S.~Smirnov.
\newblock Universality in the 2{D} {I}sing model and conformal invariance of
  fermionic observables.
\newblock {\em Invent. Math.}, 189(3):515--580, 2012. \MR{2957303}

\bibitem[DKRV16]{dkrv-lqg-sphere}
F.~David, A.~Kupiainen, R.~Rhodes, and V.~Vargas.
\newblock Liouville quantum gravity on the {R}iemann sphere.
\newblock {\em Comm. Math. Phys.}, 342(3):869--907, 2016, \arxiv{1410.7318}.
  \MR{3465434}

\bibitem[DMS14]{wedges}
B.~{Duplantier}, J.~{Miller}, and S.~{Sheffield}.
\newblock {Liouville quantum gravity as a mating of trees}.
\newblock {\em ArXiv e-prints}, September 2014, \arxiv{1409.7055}.

\bibitem[DRV16]{drv-torus}
F.~David, R.~Rhodes, and V.~Vargas.
\newblock Liouville quantum gravity on complex tori.
\newblock {\em J. Math. Phys.}, 57(2):022302, 25, 2016, \arxiv{1504.00625}.
  \MR{3450564}

\bibitem[DS11]{shef-kpz}
B.~Duplantier and S.~Sheffield.
\newblock Liouville quantum gravity and {KPZ}.
\newblock {\em Invent. Math.}, 185(2):333--393, 2011, \arxiv{1206.0212}.
  \MR{2819163 (2012f:81251)}

\bibitem[Dub09a]{dubedat-duality}
J.~Dub{\'e}dat.
\newblock Duality of {S}chramm-{L}oewner evolutions.
\newblock {\em Ann. Sci. \'Ec. Norm. Sup\'er. (4)}, 42(5):697--724, 2009,
  \arxiv{0711.1884}. \MR{2571956 (2011g:60151)}

\bibitem[Dub09b]{dubedat-coupling}
J.~Dub{\'e}dat.
\newblock S{LE} and the free field: partition functions and couplings.
\newblock {\em J. Amer. Math. Soc.}, 22(4):995--1054, 2009, \arxiv{0712.3018}.
  \MR{2525778 (2011d:60242)}

\bibitem[Duq03]{Duquesne03}
T.~Duquesne.
\newblock A limit theorem for the contour process of conditioned
  {G}alton-{W}atson trees.
\newblock {\em Ann. Probab.}, 31(2):996--1027, 2003. \MR{1964956}

\bibitem[DW15]{dw-cones}
D.~Denisov and V.~Wachtel.
\newblock Random walks in cones.
\newblock {\em Ann. Probab.}, 43(3):992--1044, 2015, \arxiv{1110.1254}.
  \MR{3342657}

\bibitem[FIM99]{Fayolle:walks-quarter-plane}
G.~Fayolle, R.~Iasnogorodski, and V.~Malyshev.
\newblock {\em Random walks in the quarter-plane}, volume~40 of {\em
  Applications of Mathematics (New York)}.
\newblock Springer-Verlag, Berlin, 1999. \MR{1691900 (2000g:60002)}

\bibitem[FS09]{Flajolet:analytic}
P.~Flajolet and R.~Sedgewick.
\newblock {\em Analytic combinatorics}.
\newblock Cambridge University Press, 2009.

\bibitem[Ges86]{Gessel:Kreweras}
I.~Gessel.
\newblock A probabilistic method for lattice path enumeration.
\newblock {\em J. Statist. Plann. Inference}, 14:49--58, 1986.

\bibitem[GH20]{gh-displacement}
E.~Gwynne and T.~Hutchcroft.
\newblock Anomalous diffusion of random walk on random planar maps.
\newblock {\em Probab. Theory Related Fields}, 178(1-2):567--611, 2020.
  \MR{4146545}

\bibitem[GHMS17]{kappa8-cov}
E.~Gwynne, N.~Holden, J.~Miller, and X.~Sun.
\newblock Brownian motion correlation in the peanosphere for {$\kappa>8$}.
\newblock {\em Ann. Inst. Henri Poincar\'e Probab. Stat.}, 53(4):1866--1889,
  2017, \arxiv{1510.04687}. \MR{3729638}

\bibitem[GHS16]{ghs-bipolar}
E.~{Gwynne}, N.~{Holden}, and X.~{Sun}.
\newblock {Joint scaling limit of a bipolar-oriented triangulation and its dual
  in the peanosphere sense}.
\newblock {\em ArXiv e-prints}, March 2016, \arxiv{1603.01194}.

\bibitem[GHS19a]{ghs-dist-exponent}
E.~Gwynne, N.~Holden, and X.~Sun.
\newblock A distance exponent for {L}iouville quantum gravity.
\newblock {\em Probab. Theory Related Fields}, 173(3-4):931--997, 2019.
  \MR{3936149}

\bibitem[GHS19b]{ghs-metric-peano}
E.~{Gwynne}, N.~{Holden}, and X.~{Sun}.
\newblock {Joint scaling limit of site percolation on random triangulations in
  the metric and peanosphere sense}.
\newblock {\em arXiv e-prints}, May 2019, \arxiv{1905.06757}.

\bibitem[GHS19c]{peano-survey}
E.~{Gwynne}, N.~{Holden}, and X.~{Sun}.
\newblock {Mating of trees for random planar maps and Liouville quantum
  gravity: a survey}.
\newblock {\em arXiv e-prints}, page arXiv:1910.04713, October 2019,
  1910.04713.

\bibitem[GHS20]{ghs-map-dist}
E.~Gwynne, N.~Holden, and X.~Sun.
\newblock A mating-of-trees approach for graph distances in random planar maps.
\newblock {\em Probab. Theory Related Fields}, 177(3-4):1043--1102, 2020.
  \MR{4126936}

\bibitem[GHSS19]{ghss18}
C.~{Garban}, N.~{Holden}, A.~{Sep{\'u}lveda}, and X.~{Sun}.
\newblock {Liouville dynamical percolation}.
\newblock {\em ArXiv e-prints}, May 2019, \arxiv{1905.06940}.

\bibitem[GKMW18]{gkmw-burger}
E.~Gwynne, A.~Kassel, J.~Miller, and D.~B. Wilson.
\newblock Active {S}panning {T}rees with {B}ending {E}nergy on {P}lanar {M}aps
  and {SLE}-{D}ecorated {L}iouville {Q}uantum {G}ravity for {$\kappa > 8$}.
\newblock {\em Comm. Math. Phys.}, 358(3):1065--1115, 2018, \arxiv{1603.09722}.
  \MR{3778352}

\bibitem[GM17a]{gwynne-miller-char}
E.~{Gwynne} and J.~{Miller}.
\newblock {Characterizations of SLE$_{\kappa}$ for $\kappa \in (4,8)$ on
  Liouville quantum gravity}.
\newblock {\em ArXiv e-prints}, January 2017, \arxiv{1701.05174}.

\bibitem[GM17b]{gm-spec-dim}
E.~{Gwynne} and J.~{Miller}.
\newblock {Random walk on random planar maps: spectral dimension, resistance,
  and displacement}.
\newblock {\em ArXiv e-prints}, November 2017, \arxiv{1711.00836}.

\bibitem[GM18]{gwynne-miller-sle6}
E.~Gwynne and J.~Miller.
\newblock Chordal {SLE$_6$} explorations of a quantum disk.
\newblock {\em Electron. J. Probab.}, 23:1--24, 2018, \arxiv{1701.05172}.

\bibitem[GMS17]{gms-tutte}
E.~{Gwynne}, J.~{Miller}, and S.~{Sheffield}.
\newblock {The Tutte embedding of the mated-CRT map converges to Liouville quantum gravity}.
\newblock {\em Ann. Probab.}, 2021.

\bibitem[GMS19]{gms-burger-cone}
E.~Gwynne, C.~Mao, and X.~Sun.
\newblock Scaling limits for the critical {F}ortuin--{K}asteleyn model on a
  random planar map {I}: {C}one times.
\newblock {\em Ann. Inst. Henri Poincar\'{e} Probab. Stat.}, 55(1):1--60, 2019,
  \arxiv{1502.00546}. \MR{3901640}

\bibitem[GPS13]{gps-pivotal}
C.~Garban, G.~Pete, and O.~Schramm.
\newblock Pivotal, cluster, and interface measures for critical planar
  percolation.
\newblock {\em J. Amer. Math. Soc.}, 26(4):939--1024, 2013, \arxiv{1008.1378}.
  \MR{3073882}

\bibitem[GPS18]{gps-near-crit}
C.~Garban, G.~Pete, and O.~Schramm.
\newblock The scaling limits of near-critical and dynamical percolation.
\newblock {\em J. Eur. Math. Soc. (JEMS)}, 20(5):1195--1268, 2018. \MR{3790067}

\bibitem[GS15]{gms-burger-finite}
E.~{Gwynne} and X.~{Sun}.
\newblock {Scaling limits for the critical Fortuin-Kasteleyn model on a random
  planar map III: finite volume case}.
\newblock {\em ArXiv e-prints}, October 2015, \arxiv{1510.06346}.

\bibitem[GS17]{gms-burger-local}
E.~Gwynne and X.~Sun.
\newblock Scaling limits for the critical {F}ortuin-{K}asteleyn model on a
  random planar map {II}: local estimates and empty reduced word exponent.
\newblock {\em Electron. J. Probab.}, 22:Paper No. 45, 56, 2017,
  \arxiv{1505.03375}. \MR{3661659}

\bibitem[HLS18]{natural}
N.~{Holden}, X.~{Li}, and X.~{Sun}.
\newblock {Natural parametrization of percolation interface and pivotal
  points}.
\newblock {\em ArXiv e-prints}, April 2018, 1804.07286.

\bibitem[HS18]{hs-Euclidean}
N.~Holden and X.~Sun.
\newblock S{LE} as a mating of trees in {E}uclidean geometry.
\newblock {\em Comm. Math. Phys.}, 364(1):171--201, 2018, \arxiv{1610.05272}.
  \MR{3861296}

\bibitem[HS19]{hs-quenched}
N.~{Holden} and X.~{Sun}.
\newblock {Convergence of uniform triangulations under the Cardy embedding}.
\newblock {\em ArXiv e-prints}, May 2019, \arxiv{1905.13207}.

\bibitem[KKZ09]{kkz-gessel-walk}
M.~{Kauers}, C.~{Koutschan}, and D.~{Zeilberger}.
\newblock {Proof of Ira Gessel's lattice path conjecture}.
\newblock {\em Proceedings of the National Academy of Science},
  106:11502--11505, July 2009, \arxiv{0806.4300}.

\bibitem[KMSW19]{kmsw-bipolar}
R.~Kenyon, J.~Miller, S.~Sheffield, and D.~B. Wilson.
\newblock Bipolar orientations on planar maps and {${\rm SLE}_{12}$}.
\newblock {\em Ann. Probab.}, 47(3):1240--1269, 2019, \arxiv{1511.04068}.
  \MR{3945746}

\bibitem[Kor13]{k13}
I.~Kortchemski.
\newblock A simple proof of {D}uquesne's theorem on contour processes of
  conditioned {G}alton-{W}atson trees.
\newblock In {\em S\'eminaire de {P}robabilit\'es {XLV}}, volume 2078 of {\em
  Lecture Notes in Math.}, pages 537--558. Springer, Cham, 2013. \MR{3185928}

\bibitem[KPZ88]{kpz-scaling}
V.~Knizhnik, A.~Polyakov, and A.~Zamolodchikov.
\newblock {Fractal structure of 2D-quantum gravity}.
\newblock {\em {Modern Phys. Lett A}}, 3(8):819--826, 1988.

\bibitem[KR12]{KR-12}
I.~Kurkova and K.~Raschel.
\newblock On the functions counting walks with small steps in the quarter
  plane.
\newblock {\em Publ. Math. Inst. Hautes \'Etudes Sci.}, 116:69--114, 2012.
\newblock \href{http://arxiv.org/abs/1107.2340}{ArXiv:1107.2340}. \MR{3090255}

\bibitem[Kre65]{Kreweras:walks}
G.~Kreweras.
\newblock Sur une classe de probl\`emes li\'es au treillis des partitions
  d'entiers.
\newblock {\em Cahier du B.U.R.O.}, 6:5--105, 1965.

\bibitem[KS16]{ks-ising}
A.~{Kemppainen} and S.~{Smirnov}.
\newblock {Conformal invariance in random cluster models. II. Full scaling
  limit as a branching SLE}.
\newblock {\em ArXiv e-prints}, September 2016, 1609.08527.

\bibitem[KZ08]{kauers07v}
M.~Kauers and D.~Zeilberger.
\newblock The quasi-holonomic {A}nsatz and restricted lattice walks.
\newblock {\em J. Difference Equ. Appl.}, 14:1119--1126, 2008.
\newblock \href{http://arxiv.org/abs/0806.4318}{ArXiv:0806.4318}.

\bibitem[Law05]{lawler-book}
G.~F. Lawler.
\newblock {\em Conformally invariant processes in the plane}, volume 114 of
  {\em Mathematical Surveys and Monographs}.
\newblock American Mathematical Society, Providence, RI, 2005. \MR{2129588
  (2006i:60003)}

\bibitem[LG07]{LeGall:limitmaps}
J.-F. Le~Gall.
\newblock The topological structure of scaling limits of large planar maps.
\newblock {\em Inventiones Mathematica}, 169:621--670, 2007.

\bibitem[LG13]{legall-uniqueness}
J.-F. Le~Gall.
\newblock Uniqueness and universality of the {B}rownian map.
\newblock {\em Ann. Probab.}, 41(4):2880--2960, 2013, \arxiv{1105.4842}.
  \MR{3112934}

\bibitem[LSW04]{lsw-lerw-ust}
G.~F. Lawler, O.~Schramm, and W.~Werner.
\newblock Conformal invariance of planar loop-erased random walks and uniform
  spanning trees.
\newblock {\em Ann. Probab.}, 32(1B):939--995, 2004, \arxiv{math/0112234}.
  \MR{2044671 (2005f:82043)}

\bibitem[LSW17]{lsw-schnyder-wood}
Y.~{Li}, X.~{Sun}, and S.~S. {Watson}.
\newblock {Schnyder woods, SLE(16), and Liouville quantum gravity}.
\newblock {\em ArXiv e-prints}, May 2017, \arxiv{1705.03573}.

\bibitem[LV16]{lawler-viklund-16}
G.~F. {Lawler} and F.~{Viklund}.
\newblock {Convergence of loop-erased random walk in the natural
  parametrization}.
\newblock {\em Duke Math. J.}, 2021.

\bibitem[Mie13]{miermont-brownian-map}
G.~Miermont.
\newblock The {B}rownian map is the scaling limit of uniform random plane
  quadrangulations.
\newblock {\em Acta Math.}, 210(2):319--401, 2013, \arxiv{1104.1606}.
  \MR{3070569}

\bibitem[MM06]{marc-mokk-tbm}
J.-F. Marckert and A.~Mokkadem.
\newblock Limit of normalized quadrangulations: the {B}rownian map.
\newblock {\em Ann. Probab.}, 34(6):2144--2202, 2006, \arxiv{math/0403398}.
  \MR{2294979 (2007m:60092)}

\bibitem[MS16a]{lqg-tbm2}
J.~{Miller} and S.~{Sheffield}.
\newblock {Liouville quantum gravity and the Brownian map II: geodesics and
  continuity of the embedding}.
\newblock {\em ArXiv e-prints}, May 2016, \arxiv{1605.03563}.

\bibitem[MS16b]{lqg-tbm3}
J.~{Miller} and S.~{Sheffield}.
\newblock {Liouville quantum gravity and the Brownian map III: the conformal
  structure is determined}.
\newblock {\em Probab. Theory Related Fields}, 179(3):1183--1211, 2021.

\bibitem[MS16c]{ig1}
J.~Miller and S.~Sheffield.
\newblock Imaginary geometry {I}: interacting {SLE}s.
\newblock {\em Probab. Theory Related Fields}, 164(3-4):553--705, 2016,
  \arxiv{1201.1496}. \MR{3477777}

\bibitem[MS16d]{ig2}
J.~Miller and S.~Sheffield.
\newblock Imaginary geometry {II}: {R}eversibility of
  {$\operatorname{SLE}_\kappa(\rho_1;\rho_2)$} for {$\kappa\in(0,4)$}.
\newblock {\em Ann. Probab.}, 44(3):1647--1722, 2016, \arxiv{1201.1497}.
  \MR{3502592}

\bibitem[MS16e]{ig3}
J.~Miller and S.~Sheffield.
\newblock Imaginary geometry {III}: reversibility of {$\mathrm{SLE}_\kappa$}
  for {$\kappa\in(4,8)$}.
\newblock {\em Ann. of Math. (2)}, 184(2):455--486, 2016, \arxiv{1201.1498}.
  \MR{3548530}

\bibitem[MS17]{ig4}
J.~Miller and S.~Sheffield.
\newblock Imaginary geometry {IV}: interior rays, whole-plane reversibility,
  and space-filling trees.
\newblock {\em Probab. Theory Related Fields}, 169(3-4):729--869, 2017,
  \arxiv{1302.4738}. \MR{3719057}

\bibitem[MS19]{sphere-constructions}
J.~Miller and S.~Sheffield.
\newblock Liouville quantum gravity spheres as matings of finite-diameter
  trees.
\newblock {\em Ann. Inst. Henri Poincar\'{e} Probab. Stat.}, 55(3):1712--1750,
  2019, \arxiv{1506.03804}. \MR{4010949}

\bibitem[MS20]{lqg-tbm1}
J.~Miller and S.~Sheffield.
\newblock Liouville quantum gravity and the {B}rownian map {I}: the {${\rm
  QLE}(8/3,0)$} metric.
\newblock {\em Invent. Math.}, 219(1):75--152, 2020. \MR{4050102}

\bibitem[Mul67]{mullin-maps}
R.~C. Mullin.
\newblock On the enumeration of tree-rooted maps.
\newblock {\em Canad. J. Math.}, 19:174--183, 1967. \MR{0205882 (34 \#5708)}

\bibitem[MW17]{miller-wu-dim}
J.~Miller and H.~Wu.
\newblock Intersections of {SLE} {P}aths: the double and cut point dimension of
  {SLE}.
\newblock {\em Probab. Theory Related Fields}, 167(1-2):45--105, 2017,
  \arxiv{1303.4725}. \MR{3602842}

\bibitem[Pol81a]{polyakov-qg1}
A.~M. Polyakov.
\newblock Quantum geometry of bosonic strings.
\newblock {\em Phys. Lett. B}, 103(3):207--210, 1981. \MR{623209 (84h:81093a)}

\bibitem[Pol81b]{polyakov-qg2}
A.~M. Polyakov.
\newblock Quantum geometry of fermionic strings.
\newblock {\em Phys. Lett. B}, 103(3):211--213, 1981. \MR{623210 (84h:81093b)}

\bibitem[Pol90]{poly-2dqg}
A.~Polyakov.
\newblock Two-dimensional quantum gravity. {S}uperconductivity at high {$T_C$}.
\newblock In {\em Champs, cordes et ph\'enom\`enes critiques ({L}es {H}ouches,
  1988)}, pages 305--368. North-Holland, Amsterdam, 1990. \MR{1052937}

\bibitem[PS03]{Schaeffer:triangulation}
D.~Poulalhon and G.~Schaeffer.
\newblock A bijection for triangulations of a polygon with interior points and
  multiple edges.
\newblock {\em Theoret. Comput. Sci.}, 307(2):385--401, 2003.

\bibitem[Ric18]{Richier:looptrees}
L.~Richier.
\newblock The incipient infinite cluster of the uniform infinite half-planar
  triangulation.
\newblock {\em Electron. J. Probab.}, 23:Art. no. 89, 2018.

\bibitem[RS05]{schramm-sle}
S.~Rohde and O.~Schramm.
\newblock Basic properties of {SLE}.
\newblock {\em Ann. of Math. (2)}, 161(2):883--924, 2005, \arxiv{math/0106036}.
  \MR{2153402 (2006f:60093)}

\bibitem[RV14]{rhodes-vargas-review}
R.~Rhodes and V.~Vargas.
\newblock Gaussian multiplicative chaos and applications: {A} review.
\newblock {\em Probab. Surv.}, 11:315--392, 2014, \arxiv{1305.6221}.
  \MR{3274356}

\bibitem[Sch98]{Schaeffer:these}
G.~Schaeffer.
\newblock {\em Conjugaison d'arbres et cartes combinatoires al{\'e}atoires}.
\newblock PhD thesis, Universit{\'e} Bordeaux I, 1998.

\bibitem[Sch00]{schramm0}
O.~Schramm.
\newblock Scaling limits of loop-erased random walks and uniform spanning
  trees.
\newblock {\em Israel J. Math.}, 118:221--288, 2000, \arxiv{math/9904022}.
  \MR{1776084 (2001m:60227)}

\bibitem[She07]{shef-gff}
S.~Sheffield.
\newblock Gaussian free fields for mathematicians.
\newblock {\em Probab. Theory Related Fields}, 139(3-4):521--541, 2007.
  \MR{2322706 (2008d:60120)}

\bibitem[She09]{shef-cle}
S.~Sheffield.
\newblock Exploration trees and conformal loop ensembles.
\newblock {\em Duke Math. J.}, 147(1):79--129, 2009, \arxiv{math/0609167}.
  \MR{2494457 (2010g:60184)}

\bibitem[She16a]{shef-zipper}
S.~Sheffield.
\newblock Conformal weldings of random surfaces: {SLE} and the quantum gravity
  zipper.
\newblock {\em Ann. Probab.}, 44(5):3474--3545, 2016, \arxiv{1012.4797}.
  \MR{3551203}

\bibitem[She16b]{shef-burger}
S.~Sheffield.
\newblock Quantum gravity and inventory accumulation.
\newblock {\em Ann. Probab.}, 44(6):3804--3848, 2016, \arxiv{1108.2241}.
  \MR{3572324}

\bibitem[Shi85]{shimura-cone}
M.~Shimura.
\newblock Excursions in a cone for two-dimensional {B}rownian motion.
\newblock {\em J. Math. Kyoto Univ.}, 25(3):433--443, 1985. \MR{807490
  (87a:60095)}

\bibitem[Smi01]{smirnov-cardy}
S.~Smirnov.
\newblock Critical percolation in the plane: conformal invariance, {C}ardy's
  formula, scaling limits.
\newblock {\em C. R. Acad. Sci. Paris S\'er. I Math.}, 333(3):239--244, 2001,
  \arxiv{0909.4499}. \MR{1851632 (2002f:60193)}

\bibitem[SS09]{ss-dgff}
O.~Schramm and S.~Sheffield.
\newblock Contour lines of the two-dimensional discrete {G}aussian free field.
\newblock {\em Acta Math.}, 202(1):21--137, 2009, \arxiv{math/0605337}.
  \MR{2486487 (2010f:60238)}

\bibitem[SS11]{ss-planar-perc}
O.~Schramm and S.~Smirnov.
\newblock On the scaling limits of planar percolation.
\newblock In {\em Selected works of {O}ded {S}chramm. {V}olume 1, 2}, Sel.
  Works Probab. Stat., pages 1193--1247. Springer, New York, 2011,
  \arxiv{1101.5820}.
\newblock With an appendix by Christophe Garban. \MR{2883400}

\bibitem[SW01]{smirnov-werner-percolation}
S.~Smirnov and W.~Werner.
\newblock Critical exponents for two-dimensional percolation.
\newblock {\em Math. Res. Lett.}, 8(5-6):729--744, 2001, \arxiv{math/0109120}.
  \MR{1879816 (2003i:60173)}

\bibitem[Wer04]{werner-notes}
W.~Werner.
\newblock Random planar curves and {S}chramm-{L}oewner evolutions.
\newblock In {\em Lectures on probability theory and statistics}, volume 1840
  of {\em Lecture Notes in Math.}, pages 107--195. Springer, Berlin, 2004,
  \arxiv{math/030335}. \MR{2079672 (2005m:60020)}

\end{thebibliography}

\end{document}